\pdfoutput=1
\ifdefined\pdfinfoomitdate\pdfinfoomitdate=1\fi
\ifdefined\pdfsuppressptexinfo\pdfsuppressptexinfo=-1\fi
\ifdefined\pdftrailerid\pdftrailerid{}\fi
\documentclass[11pt,reqno,a4paper]{amsart}
\usepackage[T1]{fontenc}
\usepackage[utf8]{inputenc}
\usepackage{lmodern}
\usepackage[a4paper,left=22mm,right=22mm,top=25mm,bottom=25mm,headsep=6mm,footskip=10mm]{geometry}
\usepackage[protrusion=true,expansion=false]{microtype}
\usepackage{mathtools,amssymb,amsfonts,amsthm,mathrsfs}
\usepackage{enumitem,booktabs,array,tabularx,graphicx,placeins,xcolor,needspace,longtable}
\usepackage{pdflscape}
\usepackage{tikz}
\usetikzlibrary{arrows.meta,positioning,calc}
\usepackage{aliascnt}
\usepackage[hidelinks]{hyperref}
\usepackage[nameinlink,capitalise,noabbrev]{cleveref}
\allowdisplaybreaks
\numberwithin{equation}{section}
\newtheorem{theorem}{Theorem}[section]
\newaliascnt{proposition}{theorem}\newtheorem{proposition}[proposition]{Proposition}\aliascntresetthe{proposition}
\newaliascnt{lemma}{theorem}\newtheorem{lemma}[lemma]{Lemma}\aliascntresetthe{lemma}
\newaliascnt{corollary}{theorem}\newtheorem{corollary}[corollary]{Corollary}\aliascntresetthe{corollary}
\theoremstyle{definition}
\newaliascnt{definition}{theorem}\newtheorem{definition}[definition]{Definition}\aliascntresetthe{definition}
\newaliascnt{conjecture}{theorem}\newtheorem{conjecture}[conjecture]{Conjecture}\aliascntresetthe{conjecture}
\newaliascnt{question}{theorem}\newtheorem{question}[question]{Question}\aliascntresetthe{question}
\theoremstyle{remark}
\newaliascnt{remark}{theorem}\newtheorem{remark}[remark]{Remark}\aliascntresetthe{remark}
\crefname{proposition}{Proposition}{Propositions}
\crefname{lemma}{Lemma}{Lemmas}\crefname{corollary}{Corollary}{Corollaries}
\crefname{definition}{Definition}{Definitions}\crefname{remark}{Remark}{Remarks}
\crefname{question}{Question}{Questions}\crefname{conjecture}{Conjecture}{Conjectures}
\crefname{equation}{Equation}{Equations}\crefname{section}{Section}{Sections}
\crefname{appendix}{Appendix}{Appendices}
\setlist[enumerate]{label=\textup{(\roman*)},leftmargin=*,itemsep=2pt}
\newcommand{\C}{\mathbb C}\newcommand{\R}{\mathbb R}\newcommand{\Z}{\mathbb Z}\newcommand{\D}{\mathbb D}

\newcommand{\interior}{\operatorname{int}}
\newcommand{\intr}{\operatorname{int}_X}
\newcommand{\re}{\operatorname{Re}}\newcommand{\im}{\operatorname{Im}}

\newcommand{\file}[1]{{\small\nolinkurl{#1}}}

\newcommand{\Fh}{\mathcal F_h}
\newcommand{\certmark}{}
\newenvironment{certificate}{\par\smallskip\noindent\textit{Verification.}\ }{\par}

\newcommand{\suite}[1]{{\small\nolinkurl{#1}}}
\newcommand{\companioncheck}[1]{\par\noindent\textit{Companion:} record~\ref{#1} of~\cite{Companion26}.}
\def\CompanionDOI{10.5281/zenodo.24117729}
\graphicspath{{figures/}}
\hypersetup{pdftitle={The closure of the nonreal roots of bounded-height polynomials is connected exactly at pronic heights},pdfauthor={Bernat Espigulé},pdfsubject={Connectedness and components of bounded-height polynomial root closures},pdfkeywords={bounded-height polynomials, pronic numbers, restricted power series, connected components, double zeros, connectedness loci, Pisot numbers}}
\theoremstyle{plain}
\newtheorem{maintheorem}{Theorem}

\title[Nonreal roots of bounded-height polynomials]{The closure of the nonreal roots of bounded-height polynomials is connected exactly at pronic heights}
\author[B. Espigul\'e]{Bernat Espigul\'e}
\address{Department of Computer Science, Applied Mathematics and Statistics, Universitat de Girona, Catalonia, Spain}
\email{bernat@espigule.com}
\date{10 November 2026}
\subjclass[2020]{Primary 30C15; Secondary 11C08, 11R06, 30B10, 28A80, 54F15}
\keywords{bounded-height polynomials, pronic numbers, restricted power series, connected components, double zeros, connectedness loci, Pisot numbers}
\begin{document}
\begin{abstract}
Let $\mathcal C_h$ be the closure of the nonreal roots of integer polynomials
with coefficients of modulus at most $h$. With the nonzero real roots
included, the closure is connected at every height; without them,
connectedness becomes arithmetic. We prove that $\mathcal C_h$ is connected
exactly when $h=q(q-1)$ is a pronic number $2,6,12,20,\dots$, and has at
least five components otherwise. At height two this confirms a
twenty-year-old numerical observation of Shmerkin and Solomyak on double
zeros of power series: apart from $1/2$ they form the interval
$[1/\rho,1)=\mathcal C_2\cap(0,1)$, with $\rho$ the tribonacci constant.
Disconnection is almost always infinite. For $(p-1)^2<h\le p^2-3$,
$h\ne p(p-1)$, infinitely many components with nonempty interior accumulate
at $\pm p$ and $\pm1/p$, and when $h<p(p-1)$ these four points are
themselves components. Next to each square $q^2\ge16$, at least two of the
three heights $q^2-2,q^2-1,q^2$ have infinitely many components accumulating
at a common point of $(1/(q+1),1/q)$. Explicit series give infinitely
many components at every nonpronic height below $167$ other than
$1,3,4,9,25,63$ and $81$, and in general at most $\sqrt X+1$ nonpronic
heights $h\le X$ can have finitely many components, although
$\mathcal C_h$ is a Peano compactum at every height. The proofs compare values and
derivatives of all admissible power series near a real double zero; the
comparison sets are open and closed in the whole core, so components
separated near the axis cannot reconnect elsewhere. On the real axis, the
double zeros at heights $1,3,4,7$ form sets with infinitely many
components; at height one this proves a conjecture of Shmerkin and
Solomyak and one of Solomyak on the connectedness locus for pairs of
affine maps.
\end{abstract}

\maketitle
\section{Introduction}\label{sec:introduction}

For each positive integer $h$, the closure of the nonzero roots of integer
polynomials with coefficients in $[-h,h]$ is connected. What happens if we
discard the real roots \emph{before} taking the closure? Real points may
still be limits of nonreal roots, but the real intervals that connect the
full locus need not survive. We prove that the resulting set is connected
exactly at the pronic heights $2,6,12,20,\ldots$, the products of
consecutive positive integers. At a set of heights of natural density one
it has infinitely many connected components.

The first pronic height was visible on the real axis twenty years ago.
For power series with constant term $1$ and coefficients in
$\{-2,\dots,2\}$, Shmerkin and Solomyak found that the double zero $1/2$ is
isolated and observed numerically that the remaining double zeros in
$(0,1)$ appeared connected~\cite[Proposition~6.3 and Remark~6.4]{SS06}.
Proposition~\ref{prop:pronic-continuation} confirms this: they form the
interval $[1/\rho,1)$, where $\rho=1.839\ldots$ is the tribonacci constant,
the real root of $\rho^3=\rho^2+\rho+1$, and this interval is
exactly the set of points of $(0,1)$ that are limits of nonreal roots of
integer polynomials of height at most two. It is the real trace of a
connected planar set, and Theorem~\ref{root:thm:main} below shows that this
connectedness occurs exactly at the pronic heights. At height one the same
authors conjectured that the double zeros form a set with infinitely many
components; Theorem~\ref{thm:real-traces} proves this.

The height parameters below are integers. For $h\ge1$, put
\begin{equation}\label{root:eq:polynomial-family}
 \mathcal P_h=\left\{\sum_{j=0}^{d}a_jz^j:
 d\ge1,\ a_d\ne0,\ a_j\in\mathbb Z,\ |a_j|\le h\right\},
 \qquad
 \mathcal R_h=\{z:p(z)=0\text{ for some }p\in\mathcal P_h\},
\end{equation}
and define
\begin{equation}\label{root:eq:closures}
 \mathcal Z_h=\overline{\mathcal R_h},\qquad
 \mathcal C_h=\overline{\mathcal R_h\setminus\mathbb R}.
\end{equation}
All closures are in $\mathbb C$. Here height refers to the coefficients
of a defining polynomial, not necessarily its minimal polynomial.
We call $\mathcal C_h$ the \emph{planar core}.

For the full root closure the answer is known and uniform: the nonzero
part $\mathcal Z_h\setminus\{0\}$ is connected and locally connected at
every height, by its interpretation as a connectedness locus for affine
iterated function systems (Section~\ref{sec:foundations}). The planar
core behaves differently. Our first result says that its connectedness
is an arithmetic property of the height.

\begin{maintheorem}[Connectedness]\label{root:thm:main}\label{thm:eventual}
For every integer $h\ge1$,
\[
 \mathcal C_h\text{ is connected}
 \quad\Longleftrightarrow\quad
 h=q(q-1)\quad\text{for an integer }q\ge2.
\]
At these heights $\mathcal C_h$ is also locally connected; at every
other height it has at least five connected components. The full
nonzero root closure $\mathcal Z_h\setminus\{0\}$ is connected and
locally connected at every height.
\end{maintheorem}

\begin{figure}[tbp]
\centering
\includegraphics[width=\textwidth]{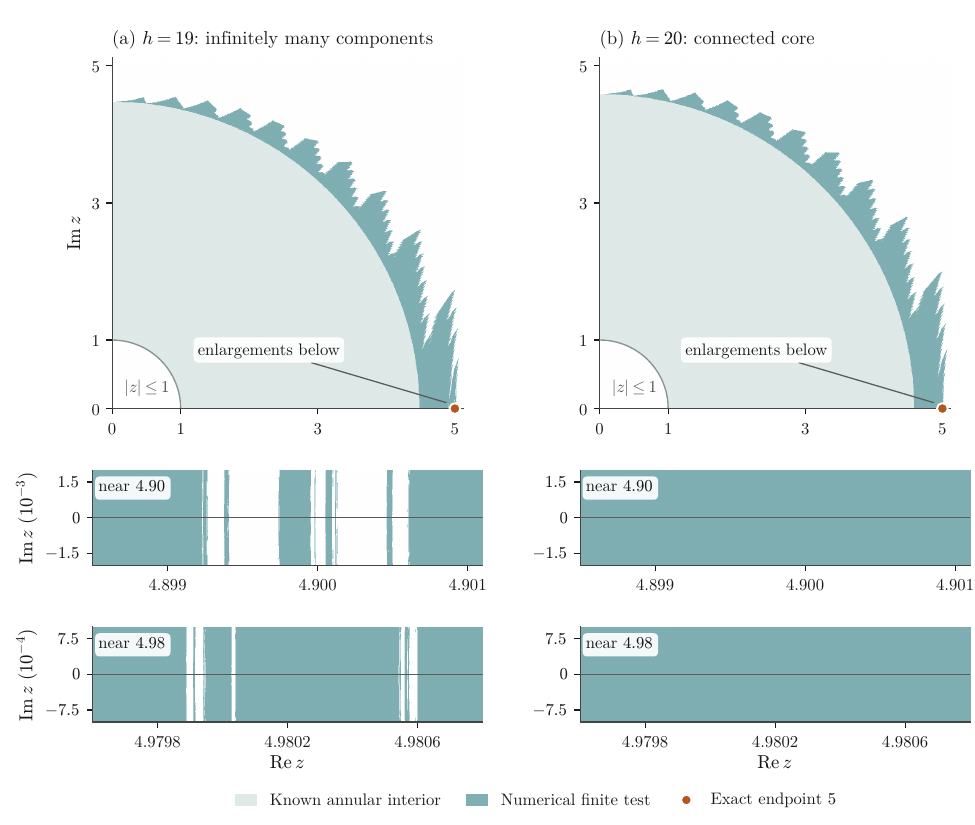}
\caption{Reconnection at consecutive heights $19$ and $20$.
The exterior parts of the cores ($|z|>1$) are shown in the first
quadrant, with matched enlargements near their common real endpoint $5$.
The enlargements straddle the real axis and scale the axes independently.
At height $19$, infinitely many components accumulate at the singleton
component $\{5\}$; at height $20$, the core is connected.
These assertions follow from Theorems~\ref{root:thm:main}
and~\ref{root:thm:infinite}. Shading illustrates finite-depth numerical
tests and does not certify component counts.}
\label{fig:reconnection}
\end{figure}

\begin{figure}[tbp]
\centering
\includegraphics[width=\textwidth]{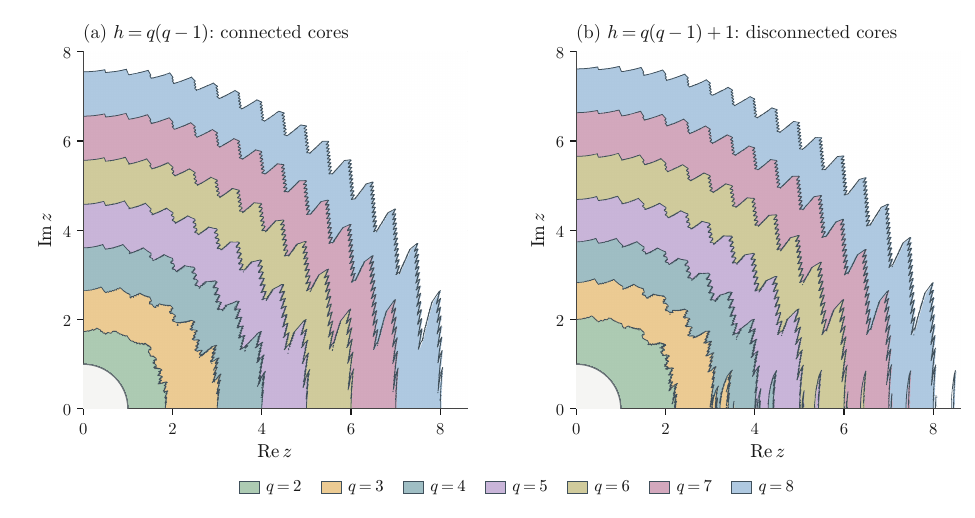}
\caption{Connectedness is not monotone in the height.
The exterior first-quadrant views show the nested cores at
$h=q(q-1)$ (left) and $h=q(q-1)+1$ (right), for $2\le q\le8$.
The same colour identifies $q$ in both panels; smaller-height sets
are drawn over larger ones. Every core on the left is connected,
whereas every core on the right is disconnected by
Theorem~\ref{root:thm:main}. The coloured regions and outlines are
finite-depth numerical illustrations.}
\label{fig:nested-families}
\end{figure}

The connected heights $2,6,12,20,30,42,\ldots$ are the \emph{pronic}
heights. The cores are nested as $h$ increases, but connectedness
is not monotone. Real intervals account for the difference between the
two sets in the theorem: they can join distinct components of the core.
At a pronic height, the remaining real intervals are attached at single
endpoints and contract onto the connected core.
Figure~\ref{fig:reconnection} shows this change from height $19$ to $20$;
Figure~\ref{fig:nested-families} places it in the sequence of pronic
heights and their successors.

Our second result gives infinitely many components throughout both
sides of each pronic height, apart from a short interval around the
next square.

\Needspace{14\baselineskip}
\begin{maintheorem}[Infinitely many components]\label{root:thm:infinite}
Let $p\ge3$ and $h$ be integers such that
\[
 (p-1)^2+1\le h\le p^2-3,\qquad h\ne p(p-1).
\]
Then $\mathcal C_h$ has infinitely many connected components, including
infinitely many with nonempty planar interior. Every neighbourhood of
each of the four points
\[
 p,\quad -p,\quad p^{-1},\quad -p^{-1}
\]
contains infinitely many such components. If $h<p(p-1)$, these four
points are singleton components. If $h>p(p-1)$, they belong to the
component containing the unit circle. In both cases the core is not
locally connected at these points.
\end{maintheorem}

Theorem~\ref{root:thm:infinite} is most striking just below a pronic
height. At height $19$ the point $5$ is a connected component of the
core by itself, and infinitely many components with interior accumulate
on it; at height $20$ everything has joined into a single continuum
(Figure~\ref{fig:reconnection}). The same happens at $\pm p$ and
$\pm1/p$ at every height strictly between $(p-1)^2$ and $p(p-1)$.

Near a square the rational anchors of Theorem~\ref{root:thm:infinite} are
no longer available. We call $q^2-2$, $q^2-1$, $q^2$ the
\emph{square-adjacent} heights. For at least two of the three we find the
same accumulation point in the open interval $(1/(q+1),1/q)$ between
the two rational anchors.

\begin{maintheorem}[Square-adjacent heights]\label{root:thm:gap}
Let $q\ge4$. At least two of the three cores $\mathcal C_{q^2-2}$,
$\mathcal C_{q^2-1}$, $\mathcal C_{q^2}$ have infinitely many connected
components, including infinitely many with nonempty planar interior.
They accumulate at a common point $t_*\in(1/(q+1),1/q)$ and at its
negative and reciprocals, where these cores are not locally connected.
All three cores have this property for $q\in\{4,6,12\}$. At every
square-adjacent height with $4\le q\le12$ other than $25$, $63$ and
$81$, the core has infinitely many components with nonempty planar
interior, accumulating at a point of $(1/(q+1),1/q)$; the same holds at
the heights $7$ and $8$, with $q=3$, and at $194$, $195$, $223$, $224$,
$225$, $254$ and $255$. At every square-adjacent height $h$ with $q\ge4$, consider
the series $1-2qz+\sum_{j\ge2}a_jz^j$ with integer coefficients
$h-2q+2\le a_j\le h$ and a double zero in $(1/(q+1),1/q)$. If one of
them is not eventually constant with coefficient $h-2q+2$, $h-q+1$ or
$h$, then $\mathcal C_h$ has infinitely many components; this is the
case unless one of them is eventually constant with coefficient $h-q+1$.
\end{maintheorem}

\begin{corollary}[Density one]\label{cor:density-one}
The set of heights at which $\mathcal C_h$ has infinitely many connected
components has natural density one. Among nonpronic heights $h\le X$, at
most $\sqrt X+1$ can fail to have this property. Possible exceptions lie
in $\{1,3,4,9,25,63,81\}$ and in the triples $\{q^2-2,q^2-1,q^2\}$
with $q\ge13$, with at most one exception in each triple.
\end{corollary}
\begin{proof}
The nonpronic heights not covered by Theorem~\ref{root:thm:infinite}
are exactly
\begin{equation}\label{eq:uniform-remaining-heights}
 \{1\}\ \cup\ \{q^2-1,q^2:q\ge2\}
 \ \cup\ \{q^2-2:q\ge3\}.
\end{equation}
Theorem~\ref{root:thm:gap} settles the heights $7$ and $8$, all but at
most one height of each triple with $q\ge4$, and every height of the
triples with $4\le q\le12$ other than $25$, $63$ and $81$. This leaves
the seven listed heights and at most one height in each triple with
$q\ge13$. Arrange these possible exceptions in a sequence: the seven
listed heights, then one height from each triple with
$q=13,14,\ldots$. Its $j$-th term is at least $(j-1)^2$, by inspection
for $j\le7$ and because $q=j+5$ for $j\ge8$.
Hence at most $\sqrt X+1$ exceptions are at most $X$.
The pronic heights themselves contribute only $O(\sqrt X)$, so the
complement has natural density zero.
\end{proof}

\begin{figure}[tbp]
\centering
\includegraphics[width=\textwidth]{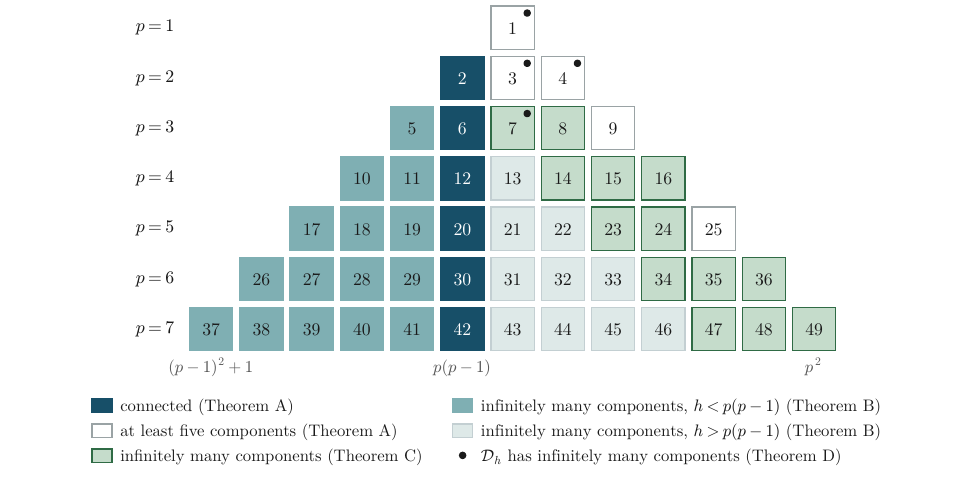}
\caption{What Theorems~\ref{root:thm:main}--\ref{thm:real-traces} prove
for $h\le49$. Row $p$ lists the heights $(p-1)^2<h\le p^2$, shifted so that
the pronic heights $p(p-1)$ form one column. For $h<p(p-1)$ the four points
$\pm p,\pm1/p$ are singleton components of $\mathcal C_h$; for $h>p(p-1)$
they lie in the component containing the unit circle. The green cells are
settled by the explicit witnesses of
Propositions~\ref{prop:explicit-witnesses}, \ref{prop:square-witnesses}
and~\ref{prop:heights-seven-eight}.
In the white cells, $h=1,3,4,9,25$, infinitude is open
(Conjecture~\ref{conj:universal-infinitude}); a dot marks the heights at
which the real double-zero set $\mathcal D_h$ has infinitely many components
(Theorem~\ref{thm:real-traces}), which does not decide the planar question.}
\label{fig:height-status}
\end{figure}

The count bounds actual exceptions, not the heights still undecided
individually: at most one of $167,168,169$ can be an exception, but the
argument does not decide which.
Figure~\ref{fig:height-status} shows what the theorems prove up to
height $49$; Section~\ref{sec:remaining-scope} discusses the remaining
heights and reduces the square-adjacent ones to a question about double
zeros of explicit integer polynomials.

\subsection{The arithmetic mechanism}
Taking reciprocals relates the root closure to the restricted series
\[
 \mathcal F_h=\left\{1+\sum_{j\ge1}a_jz^j:
                   a_j\in[-h,h]\cap\mathbb Z\right\}.
\]
Let $\mathcal D_h\subset(0,1)$ be their real multiple-zero set, and
let $\mathcal A_h=\mathcal C_h\cap(0,1)$ be the \emph{active} real
set. A point is active when nonreal roots approach it. A real collision
of conjugate roots gives $\mathcal A_h\subseteq\mathcal D_h$, but
the inclusion can be strict.

For $q\ge3$, the extremal series
\[
 F_q(z)=\frac{(1-(q+1)z)^2}{1-z}
       =1-(2q+1)z+q^2\sum_{j\ge2}z^j
\]
explains three distinct events. The point $1/(q+1)$ first becomes a
double zero at height $q^2$, when its positive tail is saturated and
the zero is inactive. At height $q^2+1$, perturbations $F_q+z^N$
produce nearby nonreal zeros, so the point becomes active. The core
reconnects at height $q(q+1)$, when integer digits allow no real double
zero closer to $0$ (Lemma~\ref{lem:pronic-double-zero-bound}) and every
point of $[1/(q+1),1)$ is active (Proposition~\ref{prop:pronic-core});
at the heights from $q^2+1$ to $q(q+1)-1$, the point $1/(q+1)$ is a
singleton component (Theorem~\ref{root:thm:infinite}). Throughout these
heights, $q+1$ remains the positive outer endpoint of the core.

\begin{figure}[tbp]
\centering
\includegraphics[width=\textwidth]{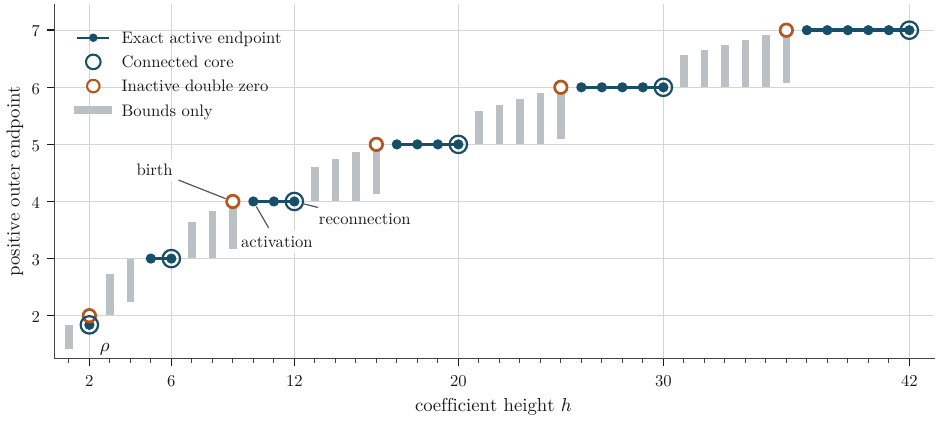}
\caption{The positive outer active endpoint for $1\le h\le42$.
Solid points give exact endpoints; circles around them mark connected
cores. Open circles mark inactive double zeros. Grey bars give bounds
only. For $q\ge3$, each constant segment begins at $q^2+1$ and ends
at the pronic height $q(q+1)$.}
\label{fig:real-axis-staircase}
\end{figure}

Perturbing a double zero can produce nearby nonreal roots without
producing distinct components. Admissible series may share a zero,
and connections can leave the local disk. We fix a reference series
with a real double zero, the \emph{anchor}, and compare every
representation in a finite coefficient block with every admissible
representation outside it. The first departure from the reference has
the form $z^NC(z)$. Shared conjugate pairs approaching the anchor would
force a coefficient limit of normalized differences to have both value
and derivative zero there, contradicting the comparison inequalities. Compactness and a quadratic
root bound then make the local root sets clopen in $\mathcal C_h$
(Lemmas~\ref{lem:finite-block-transfer}
and~\ref{lem:compact-negative-separation}). Each such set contains whole
global components.

What this argument needs from a single series is local: the comparison
inequalities at one double zero, and infinitely many levels at which
one coefficient can be lowered and the later ones raised to compensate
(Proposition~\ref{prop:local-criterion}). For the series
$(1-pz)^2/(1-z)$ this compensation is possible exactly when $h>p(p-1)$,
which is how the pronic heights enter above the anchor $1/p$; the
components then accumulate at a point of the component containing the
unit circle. At height $p(p-1)-1$, we use instead reference series with
double zeros converging to $1/p$. Two coefficient markers in each
reference isolate a finite block whose root set contains $1/p$, and the
comparison is uniform over this convergent family. The local root sets
are clopen and their diameters tend to zero, so $\{1/p\}$ is a
component.

The comparison also applies when most coefficient choices are removed
(Corollary~\ref{cor:sparse-alphabets}). For polynomials with constant
coefficient $1$ and all other coefficients in $\{-5,4,5\}$, the
nonreal-root closure has a nonisolated singleton component at $1/3$;
with the alphabet $\{-7,8,9,13\}$, infinitely many components
accumulate at $1/4$. The separating sets were constructed against every
representation in the full alphabet, so restricting the alphabet
preserves separation, while the explicit perturbations preserve the
accumulating roots.

At the three square-adjacent heights next to $q^2$, neither rational
anchor is available: at $1/q$ the first comparison inequality fails for
the digit range that the construction needs, while $1/(q+1)$ is not a
double zero at heights $q^2-2$ and $q^2-1$ (Lemma~\ref{prop:double-zero})
and is the saturated, inactive double zero of $F_q$ at height $q^2$
(Proposition~\ref{prop:endpoint-activation}). We work instead with the
series $1-2qz+\sum_{j\ge2}a_jz^j$ with digits $m\le a_j\le h$, $m=h-2q+2$,
whose double zeros $t=1/\beta$ lie in the anchor interval where
$2(\beta-q)(\beta-1)^2\in[m,h]$, inside $(1/(q+1),1/q)$. At width $2q-2$,
the comparison inequalities hold throughout this interval. The tail sums
$S_N=\sum_{i\ge1}(a_{N+i}-m)t^i$ measure the room for changing the
coefficient $a_N$ while keeping a double zero nearby. When this room
exists at infinitely many levels $N$, one can change a digit arbitrarily
far out and separate components near $t$. A witness without this
property is eventually constant with one of three digits, $m$, $m+q-1$
or $h$, and a witness common to all three heights can be exceptional in
this way at no more than one of them. This gives the two-of-three
conclusion. Of the three exceptional types, only the middle one can make
the set of anchors finite, which turns the remaining cases into a
question about double zeros of integer polynomials
(Proposition~\ref{prop:diophantine-reduction}). When $1/t$ is a Pisot
number, a finite graph decides existence and liveness, and every witness
is eventually periodic (Proposition~\ref{prop:pisot-anchors});
Table~\ref{tab:explicit-witnesses} and
Proposition~\ref{prop:square-witnesses} give explicit examples, all of
the form $g(z)^2/(1-z^r)$ with $g\in\mathbb Z[z]$. At the heights
$7$ and $8$ the comparison inequalities fail on part of the anchor
interval, but they hold at the double zeros of two such series, to
which the local criterion applies
(Proposition~\ref{prop:heights-seven-eight}).

\subsection{Real double zeros}
By Lemma~\ref{low:lem:countable}, $\mathcal D_h\setminus\mathcal A_h$
consists of isolated points of $\mathcal D_h$. Shmerkin and Solomyak proved
that $\mathcal D_1$ is disconnected~\cite[Theorem~2.6]{SS06} and conjectured
that it has infinitely many components~\cite[Section~1]{SS06}; they point out that such results
bear on infinite Bernoulli convolutions, where double zeros are the
obstruction to the transversality method of~\cite{Sol95,Sol98}. At four
small heights the real multiple-zero set breaks into infinitely many
pieces. Let $\Delta_h$ be the set of $\lambda\in(-1,1)$ such that
$(\lambda,\lambda)$ is a limit of pairs $(\gamma,\lambda')\in(-1,1)^2$,
$\gamma\ne\lambda'$, of real zeros of a common series in $\mathcal F_h$.
Coefficient compactness and Rolle's theorem give
$\Delta_h\cap(0,1)\subseteq\mathcal D_h$.

\begin{maintheorem}[Real traces near the minimum]\label{thm:real-traces}
Let $h\in\{1,3,4,7\}$. Then $\mathcal D_h$, its complement
$(0,1)\setminus\mathcal D_h$ and $\Delta_h$ have infinitely many connected
components. The gaps of $\mathcal D_h$ that separate them accumulate at
a point of $\mathcal A_h$ within $10^{-87}$ of $\min\mathcal D_h$.
\end{maintheorem}

\begin{corollary}[Solomyak's diagonal trace]\label{cor:solomyak}
Let $\mathcal N=\{(\gamma,\lambda)\in(-1,1)^2:f(\gamma)=f(\lambda)=0
\text{ for some }f\in\mathcal F_1\}$ be the connectedness locus for pairs
of affine maps studied in~\cite{Sol15}. Then
$\overline{\mathcal N\setminus\operatorname{Diag}(\R)}\cap\operatorname{Diag}(\R)$
has infinitely many connected components.
\end{corollary}
\begin{proof}
Inside $(-1,1)^2$ this set is $\{(\lambda,\lambda):\lambda\in\Delta_1\}$.
By Theorem~\ref{thm:real-traces}, infinitely many components of $\Delta_1$
are separated by gaps of $\mathcal D_1$ near $\min\mathcal D_1$; the closure
can add only the points $\pm(1,1)$, which do not join them.
\end{proof}

For $h=1$, Theorem~\ref{thm:real-traces} proves the conjecture of Shmerkin
and Solomyak, and Corollary~\ref{cor:solomyak} proves the conjecture,
recorded by Solomyak~\cite[Section~1]{Sol15}, that the diagonal trace has
infinitely many components. The proof in Section~\ref{sec:real-traces} is
computer-assisted. Near the smallest double zero of series with a
prescribed initial block, the other double zeros of such series form
scaled slices of a self-affine set, one for each later position at which a
corner digit can be lowered. At explicit digit patterns one of these
slices lies alone in a short window, which forces a gap; after the last
such pattern the remaining digits lie in a set of small dimension, and a
finite branch and bound shows that the next block again contains a
pattern. The finite statements are checked in bounded integer arithmetic
(Section~\ref{app:certified} of the supplement). Rigorous computation of
this kind has a precedent for Bernoulli convolutions: Kleptsyn, Pollicott
and Vytnova proved, with a computer-assisted algorithm, that their
Hausdorff dimension is at least $0.96399$ for every parameter in
$(1/2,1)$~\cite{KPV22}. Theorem~\ref{thm:real-traces} concerns real traces.
Real gaps do not separate planar components, so it does not decide
Conjecture~\ref{conj:universal-infinitude} at these heights, nor whether
$\mathcal A_h$ itself has infinitely many components.

\subsection{Context}
Restricted polynomial root sets and their power-series limits were studied
by Odlyzko and Poonen~\cite{OP93}, who proved path-connectedness for the
closure of the zeros of polynomials with coefficients $0$ and $1$ and
constant term $1$, by Solomyak~\cite{Sol94}, and by Shmerkin and
Solomyak~\cite{SS06}. For the alphabet $\{-1,0,1\}$, the nonzero full locus
is related to the connectedness locus for pairs of similarities introduced
by Barnsley and Harrington~\cite{BH85}; Bousch proved it connected and
locally connected~\cite{Bousch92,Bousch93}, and for the consecutive integer
alphabets considered here these properties hold at every
height~\cite[Theorem~1]{EJS24}, by a theorem of Nakajima~\cite{Nak24}. Its
interior and holes were studied by Bandt~\cite{Bandt02}, Solomyak and
Xu~\cite{SX03} and Calegari, Koch and Walker~\cite{CKW17}, who proved Bandt's
conjecture that interior points are dense off the real axis and that the
locus has infinitely many holes; \cite{Esp26holes} extends the latter to
the collinear loci of every arity, and interior regions and finite-capture
methods for collinear systems were developed in~\cite{IJK95,EJS24,EJ26}. A
related planar set arises from the Galois conjugates of the growth rates of
postcritically finite quadratic maps~\cite{Thurston14}; its
three-dimensional extension, Thurston's Master Teapot, is
connected~\cite{BDLW21}. How closely roots of restricted polynomials
approach a given point was studied by Borwein and Pinner~\cite{BP97}, and
real multiple zeros of restricted series
in~\cite{BBBP98,BBBP98b,Pinner99,SS06,Sol15}; Kewarth~\cite{Kewarth26}
classifies connectedness of the full constant-term-one locus for symmetric
integer alphabets. All of these results concern the full locus, real
intervals included. Whether those intervals are what joins distinct planar
parts is the question answered here: deciding it requires the active real
trace $\mathcal A_h$, which Section~\ref{sec:connected-heights} computes at
the pronic heights and whose gaps Section~\ref{sec:component-cascades} shows
to separate planar components elsewhere.

\subsection{Organization of the proofs}
Section~\ref{sec:foundations} supplies the root-set dictionary and the
analytic bounds used throughout, and proves that the core is a Peano
compactum. The remaining dependencies are as follows.
\begin{center}
\small
\begin{tabularx}{\textwidth}{@{}>{\raggedright\arraybackslash}p{.24\textwidth}XX@{}}
\toprule
Conclusion & Main argument & Technical input\\
\midrule
Connected pronic cores & Section~\ref{sec:connected-heights}: arithmetic
endpoint and retraction & Trace tests and continuations in
Sections~\ref{sec:interfaces} and~\ref{sec:complete-traces} of the
supplement\\[3pt]
Infinitely many components & Section~\ref{sec:component-cascades}:
comparison and global separation & None: the same section contains the
fixed and moving reference constructions\\[3pt]
Sparse alphabets & Corollary~\ref{cor:sparse-alphabets}: restrict the
separating sets & Explicit perturbations retain roots in each set\\[3pt]
Remaining nonpronic heights & Section~\ref{sec:classification}:
planar separators & Section~\ref{sec:sep:framework} of the supplement\\[3pt]
Square-adjacent heights & Section~\ref{sec:gap-heights}: square-adjacent
witnesses, tail-sum dichotomy and transfer between heights & Exact
checks of finitely many inequalities and of the explicit witnesses
(Section~\ref{app:certified} of the supplement)\\[3pt]
Real traces, $h=1,3,4,7$ & Section~\ref{sec:real-traces}: structure at
relative tips, moves and landing & Exact certificates of the minimum, slice
bounds, Move and Landing Lemmas (Section~\ref{app:certified} of the
supplement)\\
\bottomrule
\end{tabularx}
\end{center}
The finite checks are replayed by the computational
companion~\cite{Companion26}. Neither the comparison argument of
Section~\ref{sec:component-cascades} nor the sparse-alphabet application
requires a finite certificate; density of planar interior is used only for
the interior-bearing components in Theorems~\ref{root:thm:infinite}
and~\ref{root:thm:gap}. Several ingredients do not depend on the
classification and may be useful elsewhere:
Lemma~\ref{lem:finite-block-transfer} separates global components of a root
closure for any coefficient alphabet, and
Proposition~\ref{prop:local-criterion} produces infinitely many of them
from a single series; Lemma~\ref{low:lem:countable} and
Proposition~\ref{prop:peano} describe the real multiple-zero set and the
local structure of the core at every height; and
Proposition~\ref{prop:pisot-anchors} turns the search for double zeros with
a Pisot reciprocal anchor into a finite computation. Section~\ref{sec:remaining-scope}
states the remaining question.

\section{Root closures and double zeros}\label{sec:foundations}

We first record the passage between polynomials, power series and the
connected full root locus. Retain $\mathcal F_h,\mathcal D_h,\mathcal A_h$
from the introduction. All closures are taken in $\C$ unless a relative
space is specified. Write $\D=\{|z|<1\}$, $\mathbb T=\{|z|=1\}$,
$X=\{|z|>1\}$, and
\[
 B_h=[-h,h]\cap\Z,\qquad
 \Sigma_h=\{c\in X:f(1/c)=0\text{ for some }f\in\mathcal F_h\}.
\]
A series $f\in\mathcal F_h$ is a \emph{witness} at $w\in\D$ if $f(w)=0$.
Coefficientwise convergence in $\mathcal F_h$ is locally uniform in $\D$,
together with every derivative. This compactness will be used throughout.

\subsection{Normalization and the full locus}

\begin{proposition}\label{prop:root-series}
$\Sigma_h$ is the relative closure in $X$ of the roots of monic
integer polynomials of height at most $h$. Every nonreal zero in $\D$
of a series in $\mathcal F_h$ belongs to $\mathcal C_h$.
\end{proposition}
\begin{proof}
Reversal puts each exterior monic root in $\Sigma_h$, and coefficient
compactness makes $\Sigma_h$ relatively closed. Conversely, polynomial
truncations of $f\in\mathcal F_h$ converge locally uniformly to $f$.
Rouch\'e's theorem on a small circle about a zero gives nearby polynomial
zeros. Reversing the truncations gives monic polynomials with the same
height bound. At a nonreal zero the circle can be chosen disjoint from
$\R$, so the approximating zeros remain nonreal.
\end{proof}

We use two established properties of this normalized locus:
\begin{equation}\label{eq:normalized-annulus}
 \{1<|c|<\sqrt{h+1}\}\subset\operatorname{int}\Sigma_h,
 \qquad \Sigma_h\text{ is connected and locally connected}.
\end{equation}
For clarity, the dictionary is as follows. If $n=h+1$ and $E(c,n)$ is
the attractor of the maps $w\mapsto t+c^{-1}w$, with
$t\in\{-n+1,-n+3,\ldots,n-1\}$, then its connectedness locus is
$\Sigma_h$. Consecutive-piece intersection is equivalent to
$2c\in E(c,n)-E(c,n)$; dividing its difference address by $2c$ gives
the series equation. In this normalization, the dictionary, annulus
and connectedness assertions are Proposition~2, Proposition~5 and
Theorem~1 of~\cite{EJS24}, respectively.
We retain the polynomial height $h$ as our parameter.

The following elementary estimate removes the possible nonunit leading
coefficient. It will also normalize the witnesses used in
Section~\ref{sec:component-cascades}.

\begin{lemma}[Conjugate-pair bound]\label{lem:conjugate-pair-bound}
If $f(z)=1+\sum_{j=1}^{d}b_jz^j$ has real coefficients $|b_j|\le A$
and a nonreal zero $\lambda$, then
\[
 |\lambda|^{-4}\le4(1+A^2).
\]
\end{lemma}
\begin{proof}
Put $r=1/\sqrt2$ and $\mu=\lambda/r$. Factor the conjugate pair from
$f(rz)$, then reverse this quadratic factor. Reversal preserves its
modulus on $\mathbb T$, so Parseval preserves the sum of squared
coefficients. The new constant coefficient is $|\mu|^{-2}$. Therefore
\[
 |\mu|^{-4}\le 1+\sum_{j=1}^d b_j^2r^{2j}\le1+A^2.
\]
Since $|\mu|^{-4}=|\lambda|^{-4}/4$, this is the desired estimate.
\end{proof}

\begin{proposition}\label{prop:global-monic-reduction}\label{prop:annulus}
Monicity does not change $\mathcal Z_h$ or $\mathcal C_h$.
Writing $\iota(z)=1/z$ and
$\mathcal C_h^{\rm ext}=\mathcal C_h\cap X$, one has
\begin{align}
 \mathcal Z_h\setminus\{0\}
   &=\Sigma_h\cup\mathbb T\cup\iota(\Sigma_h),\label{eq:full-gluing}\\
 \mathcal C_h
   &=\mathcal C_h^{\rm ext}\cup\mathbb T
       \cup\iota(\mathcal C_h^{\rm ext}).\label{eq:core-gluing}
\end{align}
Both sets are invariant under conjugation, negation and inversion.
Moreover,
\begin{gather}
 \mathcal Z_h\cap\R=\{0\}\cup[-h-1,-(h+1)^{-1}]
                          \cup[(h+1)^{-1},h+1],\label{eq:real-slice}\\
 \{(h+1)^{-1/2}<|z|<\sqrt{h+1}\}
       \subset\operatorname{int}\mathcal C_h.\label{eq:core-annulus}
\end{gather}
The nonzero full root closure is connected and locally connected.
\end{proposition}
\begin{proof}
For $1<x\le h+1$, the intervals $(a+J)/x$, $a\in B_h$, cover
$J=[-h/(x-1),h/(x-1)]$. Expanding $-1\in J$ recursively gives a
witness in $\mathcal F_h$ at $1/x$. The reverse radius bound follows
from $1\le h/(x-1)$. Thus
$\Sigma_h\cap\R=[-h-1,-1)\cup(1,h+1]$.

Let $c$ be an exterior root of a height-$h$ polynomial whose leading
coefficient has modulus $k\ge2$. Reversal and division by that coefficient
give a witness at $1/c$ with coefficient bound $h/k$.
If $c$ is nonreal, the preceding lemma gives
\[
 |c|^4\le4(1+(h/k)^2)\le h^2+4<(h+1)^2,
\]
so $c\in\Sigma_h$ by~\eqref{eq:normalized-annulus}. If $c$ is real,
the elementary bound $|c|\le1+h/k\le h+1$ gives the same conclusion.
The case $k=1$ is Proposition~\ref{prop:root-series}.
Hence $\mathcal Z_h\cap X=\Sigma_h$, including equality of the
nonreal-root closures there.

Inside $\D$, the monic polynomials $p(z)+z^N$, $N>\deg p$, converge
locally uniformly to any fixed $p\in\mathcal P_h$; Rouch\'e gives monic
approximations to each of its roots, preserving nonreality. Roots of
unity handle $\mathbb T$. This proves monic reduction.
After removing a maximal power of $z$, polynomial reversal preserves
height and inverts every nonzero root. The radius estimate applied in
both directions gives $(h+1)^{-1}\le|z|\le h+1$ for all nonzero roots.
Thus $0$ is isolated, and inversion passes to both closures. This proves
\eqref{eq:full-gluing}--\eqref{eq:real-slice}.

The exterior annulus consists of interior points of $\Sigma_h$ and
hence of its nonreal-root closure. Inversion gives the inner annulus;
their closures supply $\mathbb T$ and a neighbourhood of it. Finally,
\eqref{eq:full-gluing} joins the connected, locally connected exterior
locus to its reciprocal across this annular neighbourhood.
\end{proof}

Other restrictions on the defining polynomials are harmless as well:
Hokken~\cite{Hokken26} proves that the closure of the zeros of
polynomials with coefficients in a fixed integer set is unchanged when
the polynomials are required to have square discriminant.

\begin{proposition}[Density of interior]\label{cor:whole-plane-regularity}
For every $h\ge1$,
\[
 \mathcal C_h=\overline{\operatorname{int}\mathcal Z_h}
            =\overline{\operatorname{int}\mathcal C_h}.
\]
\end{proposition}
\begin{proof}
The nonreal interior-density results of
\cite[Theorem~1.1]{Kewarth26b} and~\cite{RC26} give the first
identity off $\mathbb T$ through the series dictionary and inversion;
the annular neighbourhood handles $\mathbb T$. Every real interior
point of $\mathcal Z_h$ is also approximated by nonreal points. Finally,
$\operatorname{int}\mathcal Z_h\subset\mathcal C_h\subset\mathcal Z_h$
gives equality of the two interiors.
\end{proof}

This last input is used to strengthen the separated components in
Section~\ref{sec:component-cascades} to components with planar interior.
The connectedness classification and the existence of infinitely many
components do not require this strengthening.

Although the core is not locally connected at most heights
(Theorems~\ref{root:thm:infinite} and~\ref{root:thm:gap}), its failure
of local connectedness is of the simplest kind. Following
Loridant, Luo and Yang~\cite{LLY19}, a compact metric space is a
\emph{Peano compactum} if every component is locally connected and, for
every $\varepsilon>0$, only finitely many components have diameter at
least $\varepsilon$.

Local connectedness at a point means that the point has arbitrarily
small connected neighbourhoods in the set; a compact metric space with
this property at every point is locally connected~\cite[Chapter~V]{Nadler92}.

\begin{proposition}[Peano structure of the core]\label{prop:peano}
For every $h\ge1$, $\mathcal C_h$ is a Peano compactum. The core is not
locally connected at a point $x$ if and only if a sequence of distinct
components of $\mathcal C_h$ accumulates at $x$, and every such $x$ is
real.
\end{proposition}
\begin{proof}
Let $N$ be the boundary, relative to $\R$, of the closed set
$\mathcal C_h\cap\R$; it is a closed nowhere dense subset of $\R$. We
first show that $\mathcal C_h$ is locally connected at every point of
$\mathcal C_h\setminus N$. If $x\in\mathcal C_h$ is nonreal, choose a
disk $B$ about $x$ disjoint from $\R$. A point of $\mathcal Z_h\cap B$
is a limit of roots lying in $B$, hence of nonreal roots, so
$\mathcal Z_h\cap B=\mathcal C_h\cap B$. If $x$ is real and not in
$N$, it is an interior point of $\mathcal C_h\cap\R$ relative to $\R$,
so for a small disk $B$ every real point of $\mathcal Z_h\cap B$ lies in
$\mathcal C_h$, and the nonreal ones do as before; again
$\mathcal Z_h\cap B=\mathcal C_h\cap B$. Since $0\notin\mathcal C_h$
by~\eqref{eq:real-slice}, Proposition~\ref{prop:annulus} gives local
connectedness of $\mathcal Z_h$ at $x$, hence of $\mathcal C_h$.

Next, let $K$ be a compact subset of the plane which is not locally
connected at a point $x$. We claim that either distinct components of
$K$ accumulate at $x$, or $x$ lies in a nondegenerate continuum at each
point of which $K$ fails to be locally connected. Choose $r>0$ such that
no connected subset of $K\cap\overline{B(x,r)}$ contains $x$ in its
interior relative to $K$, and let $Q$ be the component of $x$ in
$K\cap\overline{B(x,r)}$. Then there are points $y_n\to x$ of $K$
outside $Q$; let $Q_n$ be the component of $y_n$ in
$K\cap\overline{B(x,r)}$. A component containing infinitely many $y_n$
would contain $x$ and equal $Q$, so after passing to a subsequence the
$Q_n$ are pairwise distinct. If $Q_n$ does not meet the circle
$\partial B(x,r)$, then $Q_n$ is a component of $K$: a continuum in $K$
containing $Q_n$ and a point outside $\overline{B(x,r)}$ would, by the
boundary bumping theorem~\cite[Chapter~V]{Nadler92}, contain a
continuum joining $Q_n$ to the circle inside $\overline{B(x,r)}$, which
would then lie in $Q_n$. So either infinitely many $Q_n$ are components
of $K$ accumulating at $x$, or infinitely many $Q_n$ meet the circle;
in the latter case a subsequence converges in the Hausdorff metric to a
continuum $C$ containing $x$ and a point of the circle. Every point $y$
of $C$ in the open disk $B(x,r)$ is a point of non-local-connectedness
of $K$: a connected $V\subseteq K\cap B(x,r)$ containing $y$ in its
relative interior would lie in one component of
$K\cap\overline{B(x,r)}$ and meet $Q_n$ for all large $n$. By boundary
bumping applied to $C$ and $\overline{B(x,r/2)}$, the component of $x$
in $C\cap\overline{B(x,r/2)}$ is a nondegenerate continuum of such
points. This proves the claim.

Apply the claim to $K=\mathcal C_h$ at a point $x$ where
$\mathcal C_h$ is not locally connected. By the first paragraph, every
point of non-local-connectedness lies in $N$, and a nondegenerate
continuum in $N\subset\R$ would be an interval, contradicting nowhere
density. Hence distinct components of $\mathcal C_h$ accumulate at
$x$, and $x\in N$ is real. Conversely, if distinct components $K_n$
accumulate at $x$, then $\mathcal C_h$ is not locally connected at $x$:
a connected subset of $\mathcal C_h$ containing a neighbourhood of $x$
would lie in the component of $x$ and meet every $K_n$ with $n$ large.

Now let $K$ be a component of $\mathcal C_h$ and $x\in K$. If
$\mathcal C_h$ is locally connected at $x$, so is $K$, because a
connected subset of $\mathcal C_h$ containing a neighbourhood of $x$ lies
in $K$. Thus the points of non-local-connectedness of $K$ lie in $N$.
If $K$ is not locally connected at $x$, apply the claim to $K$ with
$r<\operatorname{diam}K/2$: the first alternative is impossible for a
continuum, and the second gives a nondegenerate continuum inside $N$.
So every component is locally connected. Finally, if $K_n$ were
distinct components of diameter at least $\varepsilon$, a subsequence
would converge to a continuum $C$ of diameter at least $\varepsilon$
each of whose points is an accumulation point of distinct components,
hence a point of $N$; this is again impossible.
\end{proof}

By \cite[Theorem~3 and Lemma~3.8]{LLY19}, the conclusion is equivalent
to the following finiteness property, visible in the figures: for every closed
strip bounded by two parallel lines, only finitely many components of
its intersection with $\mathcal C_h$ meet both lines. At the heights of
Theorems~\ref{root:thm:infinite} and~\ref{root:thm:gap}, infinitely
many components accumulate at real points, and their diameters tend to
zero.

\subsection{Real limits and coefficient saturation}

Recall that $\mathcal D_h$ is the real multiple-zero set in $(0,1)$
and $\mathcal A_h$ is its active subset. Inversion identifies
$\mathcal A_h$ with the positive exterior active trace. Both sets are
relatively closed, by coefficient compactness and the definition of the core.

\begin{lemma}\label{prop:double-zero}\label{prop:universal}
One has $\mathcal A_h\subseteq\mathcal D_h$. Every $t\in\mathcal D_h$
satisfies $t^{-1}\le1+\sqrt h$. Equality forces
\[
 f(z)=\frac{(1-t^{-1}z)^2}{1-z},\qquad
 a_1=-(1+2\sqrt h),\quad a_j=h\quad(j\ge2).
\]
\end{lemma}
\begin{proof}
At a real limit of nonreal roots, take a coefficientwise limit of their
normalized exterior witnesses. The conjugate secants tend to $f'(t)$,
so both $f(t)$ and $f'(t)$ vanish. Subtracting the value equation from
the weighted derivative equation gives
\[
 1=\sum_{j\ge2}(j-1)a_jt^j
       \le \frac{ht^2}{(1-t)^2}.
\]
All weights are positive, so equality forces $a_j=h$ for $j\ge2$;
the value equation then determines $a_1$.
\end{proof}

For real coefficients in $[-h,h]$, the bound $t^{-1}\le1+\sqrt h$ and
its extremal series are due to Beaucoup, Borwein, Boyd and
Pinner~\cite{BBBP98b}, who proved that for $h\ge3+2\sqrt2$ the smallest
disk containing two zeros of such a series with constant term $1$ has
radius $(1+\sqrt h)^{-1}$. The integer alphabet changes what happens at
the extremal point. The saturated series with $t^{-1}=q+1$ is admissible
at height $q^2$, but its double zero is inactive there
(Proposition~\ref{prop:endpoint-activation}); nonreal roots reach $q+1$
only from height $q^2+1$ on, and the core reconnects only at $q(q+1)$
(Theorem~\ref{root:thm:main} and Figure~\ref{fig:real-axis-staircase}).

For the binary alphabet, the inclusion $\mathcal A_1\subseteq\mathcal D_1$
is Lemma~2.5 of Solomyak~\cite{Sol98} (see
\cite[Lemma~2.9]{SS06}), and $\mathcal D_1$ is the set $\Omega_2$
of~\cite{SS06}; Calegari, Koch and Walker~\cite[Section~10]{CKW17} note that
it is not known whether some point of $\Omega_2$ lies outside the closure of
the interior of the binary locus. By Proposition~\ref{cor:whole-plane-regularity}
that question asks whether $\mathcal D_1\setminus\mathcal A_1$ is
nonempty. The next lemma shows that this difference consists of isolated
points of $\mathcal D_h$, at every height, so the binary question would
be settled by showing that $\mathcal D_1$ has no isolated points;
Proposition~\ref{prop:pronic-continuation} identifies the difference as
$\{1/2\}$ at height two, and the binary case remains open.

\begin{lemma}[Inactive double zeros are isolated]\label{low:lem:countable}
Every point of $\mathcal D_h\setminus\mathcal A_h$ is an isolated
point of $\mathcal D_h$, and $\mathcal D_h\subseteq\mathcal A_{h+1}$.
Consequently $\mathcal D_h\setminus\mathcal A_h$ is at most countable,
and every nondegenerate interval in $\mathcal D_h$, including its
endpoints in $(0,1)$, is contained in $\mathcal A_h$.
\end{lemma}
\begin{proof}
Let $f\in\mathcal F_h$ have a zero of order $\nu\ge2$ at $t$, and
write $f(t+w)=bw^\nu+O(w^{\nu+1})$, $b\ne0$.
Let $\epsilon\in\{-1,1\}$ be such that $bw^\nu+\epsilon$ has a nonreal
zero. With $r_N=t^{N/\nu}$,
\[
 t^{-N}\bigl(f(t+r_Nw)+\epsilon(t+r_Nw)^N\bigr)
        \longrightarrow bw^\nu+\epsilon
\]
locally uniformly, because $Nr_N\to0$.
If $f+\epsilon z^N$ is admissible for infinitely many $N$, Rouch\'e
gives nonreal zeros tending to $t$, so $t\in\mathcal A_h$. At height
$h+1$ every $f+\epsilon z^N$ is admissible, which proves
$\mathcal D_h\subseteq\mathcal A_{h+1}$.

Now let $t\in\mathcal D_h\setminus\mathcal A_h$. For every witness $f$
of $t$ as above, $a_N=\epsilon h$ for all large $N$. If $\nu\ge3$, both
signs $\epsilon$ give $bw^\nu+\epsilon$ a nonreal zero, so $a_N=h$ and
$a_N=-h$ eventually, which is absurd; hence $\nu=2$, $\epsilon$ is the
sign of $b=f''(t)/2$, and the coefficients of $f$ equal $\epsilon h$
from some index $N_0$ on. Suppose that distinct points
$t_n\in\mathcal D_h$ converge to $t$, with witnesses $f_n$ having zeros
of order at least two at $t_n$. Passing to a subsequence, $f_n$
converges coefficientwise to a series $f$ with $f(t)=f'(t)=0$, to which
the preceding paragraph applies; say $\epsilon=1$, so $f''(t)>0$ and
$a_N=h$ for $N\ge N_0$. For large $n$, the coefficients of $f_n$ of
index below $N_0$ agree with those of $f$, and all its coefficients are
at most $h$; hence $f_n(t)-f(t)=\sum_{N\ge N_0}(a_N(f_n)-h)t^N\le0$,
that is, $f_n(t)\le0$. On the other hand $f''>0$ on a neighbourhood
$U$ of $t$, so $f_n''>0$ on $U$ for large $n$; then $f_n$ is strictly
convex on $U$ with a zero of order at least two at $t_n\in U$, hence
nonnegative on $U$ with $t_n$ as its only zero there, and $t\ne t_n$
gives $f_n(t)>0$, a contradiction. For $\epsilon=-1$ all inequalities
reverse. Thus $t$ is isolated in $\mathcal D_h$.

A subset of $\R$ has at most countably many isolated points, and no
point of a nondegenerate interval in $\mathcal D_h$, nor an endpoint of
it in $(0,1)$, is isolated in $\mathcal D_h$.
\end{proof}

The lemma asserts approximation by nonreal zeros, not planar
interiority. This distinction is needed at height six, where a bounded
value--derivative orbit completes the active trace without a strict
planar-interior test at every real parameter.

\section{The connected heights}\label{sec:connected-heights}
Connectedness at $h=q(q-1)$ comes from a complete real interval to which
all nonreal pieces attach. There are two separate steps: an integer-digit
argument bounds the endpoint by $q$, and a coefficient continuation fills
the interval up to it. The coefficient construction is given in
Sections~\ref{sec:interfaces}--\ref{sec:complete-traces} of the supplement; the endpoint
and topological arguments are given here.

\subsection{The arithmetic endpoint}
The universal double-zero bound from Section~\ref{sec:foundations}
allows an endpoint as large as $1+\sqrt h$. At height $h=q(q-1)$,
integer coefficients give the smaller endpoint $q$.

\begin{lemma}[Pronic double-zero bound]\label{lem:pronic-double-zero-bound}
Let $q\ge2$ be an integer and $H=q(q-1)$. If
\[
 f(z)=1-\sum_{k\ge1}d_kz^k,\qquad d_k\in\mathbb Z,\quad |d_k|\le H,
\]
has a real double zero $t\in(-1,1)$, then $|t|\ge q^{-1}$.
\end{lemma}
\begin{proof}
Replace $f(z)$ by $f(-z)$ if necessary and suppose $x=t^{-1}>q$.
The double-zero equations give the tail states
\[
 P_n=\sum_{k\ge1}d_{n+k}x^{-k},\qquad
 V_n=\sum_{k\ge1}(1-k)d_{n+k}x^{-k},\qquad (P_0,V_0)=(1,1).
\]
Absolute convergence yields
\begin{equation}\label{eq:pronic-tail-recurrence}
 P_{n+1}=xP_n-d_{n+1},\quad V_{n+1}=xV_n+P_{n+1},\qquad
 |P_n|\le\frac H{x-1},\quad |V_n|\le\frac H{(x-1)^2}.
\end{equation}
Put $R=x-1$, $\delta=x-q$, and $A=(q-1)^2$.
The bound on $V_0$ gives $R^2\le H<q^2$, hence $0<\delta<1$.
Since $H<2R^2$, the bound on $V_1=2x-d_1$ implies
$d_1>2q-2$. On the other hand,
\[
 Rd_1\le Rx+H<2qR,\qquad 2qR-(Rx+H)=\delta(1-\delta)>0.
\]
Integrality forces $d_1=2q-1$. For $q=2$, this already exceeds $H$.

For $q\ge3$, define
\[
 U_n=RP_{n+1}+A,\qquad W_n=R^2V_{n+1}-A.
\]
Their bounds and initial values are
\[
 -(q-1)\le U_n\le H+A,\quad W_n\le q-1,\qquad
 U_0=\delta^2,\quad W_0=\delta(2R^2+x+q-2).
\]
In particular, $U_0\ge0$ and $x(q-1)U_0<\delta W_0$, since the
difference is $\delta^2R(2x-q)>0$.
Suppose a state satisfies
\begin{equation}\label{eq:pronic-forcing-cone}
 U\ge0,\qquad x(q-1)U<\delta W.
\end{equation}
Then $W>0$ and $xU<\delta$. Writing the next digit as $d=-A+j$,
$j\in\mathbb Z$, gives
\[
 U'=xU-Rj,\qquad W'=xW+xRU-R^2j.
\]
If $j\ge1$, then $U'<\delta-R=-(q-1)$; if $j\le-1$, then
$W'>R^2>q-1$. Both contradict the state bounds, so $j=0$.
Moreover,
\[
 \delta W'-x(q-1)U'
 =x\bigl(\delta W-x(q-1)U\bigr)+\delta xRU>0.
\]
Thus \eqref{eq:pronic-forcing-cone} persists, and induction gives
$U_n=x^n\delta^2$, contradicting $U_n\le H+A$.
\end{proof}

The forced digit is the essential arithmetic step. A continuously
variable coefficient could correct the expanding state; an integer
coefficient cannot do so without violating one of the two bounds.

\subsection{Completing the real trace}
The same value--derivative coordinates give a constructive way to prove
that a whole interval consists of double zeros. For $x>1$, start at
$(P_0,V_0)=(1,1)$ and choose digits $d_j\in B_h$ under
\begin{equation}\label{eq:bounded-continuation}
 P_{j+1}=xP_j-d_{j+1},\qquad V_{j+1}=xV_j+P_{j+1}.
\end{equation}
A bounded orbit yields a double zero at $1/x$. Indeed,
\[
 x^{-m}P_m=1-\sum_{j=1}^m d_jx^{-j},\qquad
 x^{-m}\bigl(V_m-(m+1)P_m\bigr)=\sum_{j=1}^m jd_jx^{-j}.
\]
Both expressions tend to zero. Thus
$1-\sum_{j\ge1}d_jz^j\in\mathcal F_h$ vanishes with its derivative
at $1/x$.

\begin{proposition}[Pronic continuation]\label{prop:pronic-continuation}
For $h=q(q-1)$ and $q\ge3$, a bounded orbit
of~\eqref{eq:bounded-continuation} exists for every $1<x<q$.
At height two, if $\rho$ is the real root of
$\rho^3=\rho^2+\rho+1$, then
\[
 \mathcal D_2=\{1/2\}\cup[1/\rho,1),\qquad
 \mathcal A_2=[1/\rho,1).
\]
\end{proposition}

The complete coefficient construction is given in Section~\ref{sec:complete-traces} of the supplement.
Here is its organization. Away from $q$, a finite prefix enters a strict
invariant rectangle for the two coordinates. Near $q$, long repetitions
of the digit $-(q-1)^2$ shadow its fixed state. Their return coordinates
are an exponential factor times an affine function of the repetition
length. The resulting parameter bands approach $q$ and overlap, so
there is no omitted interval between scales. Uniform estimates handle
large $q$ and long repetitions; exact polynomial inequalities handle
the finite complement. At height six, bounded continuation replaces
strict entry in the last interval. At height two, the continuation is
organized by the period-three witness
\[
 \frac{(1-z-z^2-z^3)^2}{1-z^3}
   =1-2z-z^2+z^3+\frac{z^4+z^5+2z^6}{1-z^3}\in\mathcal F_2.
\]
Sections~\ref{sec:complete-traces} and~\ref{app:certified} of the supplement include the
complete covers, their endpoints and the finite verification statements. In particular, the construction does
not infer interval coverage from sampled parameters.

\begin{proposition}[Connected pronic cores]\label{prop:pronic-core}
Let $h=q(q-1)$, $q\ge2$, and set $a=q$ for $q\ge3$, $a=\rho$ for $q=2$.
Then
\begin{equation}\label{eq:pronic-active-trace}
 \mathcal C_h\cap\R=[-a,-a^{-1}]\cup[a^{-1},a].
\end{equation}
The core is connected and locally connected. It is a strong deformation
retract of $\mathcal Z_h\setminus\{0\}$, and the complement of the core
in that set consists of the four intervals
\begin{equation}\label{eq:four-antennae}
 [-h-1,-a),\quad(-a^{-1},-(h+1)^{-1}],\quad
 [(h+1)^{-1},a^{-1}),\quad(a,h+1].
\end{equation}
\end{proposition}
\begin{proof}
For $q\ge3$, continuation and the sharp pronic bound give
$\mathcal D_h=[1/q,1)$. The endpoint is supplied by
$f_q=(1-qz)^2/(1-z)\in\mathcal F_h$.
Lemma~\ref{low:lem:countable} then gives
$\mathcal A_h=[1/q,1)$. For $q=2$, use
Proposition~\ref{prop:pronic-continuation}.
Negation and inversion prove~\eqref{eq:pronic-active-trace}.

Every nonreal point of $\mathcal Z_h$ already belongs to $\mathcal C_h$
by definition. Subtracting the active trace from~\eqref{eq:real-slice}
therefore gives~\eqref{eq:four-antennae}.
Fix the core pointwise and contract the closure of each displayed
interval linearly onto its endpoint in the core. The definitions agree
on their closed overlaps, giving a strong deformation retraction.
Connectedness follows from that of the full nonzero locus.
For local connectedness, given a neighbourhood $U$ of $y$ in the core,
choose a connected open neighbourhood $V$ of $y$ in the full locus
contained in the inverse image of $U$ under the retraction $r$.
Then $r(V)$ is a connected neighbourhood of $y$ contained in $U$.
\end{proof}

\subsection{Why the next endpoint appears before reconnection}
This short calculation explains the distinction between double-zero
existence, activation and reconnection in
Figure~\ref{fig:real-axis-staircase}. For $q\ge3$, set
\[
 F_q(z)=\frac{(1-(q+1)z)^2}{1-z}
       =1-(2q+1)z+q^2\sum_{j\ge2}z^j.
\]

\begin{proposition}[Endpoint activation]\label{prop:endpoint-activation}
At height $q^2$, the double zero $1/(q+1)$ is isolated in
$\mathcal D_{q^2}$ and is not active. For
$q^2+1\le h\le q(q+1)$, the greatest positive real core point is $q+1$.
\end{proposition}
\begin{proof}
Equality in the universal double-zero bound makes $F_q$ the unique
witness at $t=1/(q+1)$ in height $q^2$, with $F_q''(t)>0$.
Suppose qualifying zeros $z_j\ne t$ tend to $t$, where qualifying means
nonreal or real of multiplicity at least two. Compactness forces their
witnesses $g_j$ to converge to $F_q$. Factoring their nearby conjugate
pair, or a real double zero, gives
\[
 g_j(t)=|t-z_j|^2Q_j(t)>0,
 \qquad Q_j(t)\longrightarrow F_q''(t)/2>0.
\]
But the first coefficient eventually agrees with that of $F_q$, while
all later coefficients are at most $q^2$. Hence $g_j(t)\le F_q(t)=0$,
a contradiction. This proves isolation and inactivity.

At height $q^2+1$, every $F_q+z^N$, $N\ge2$, is admissible.
The quadratic scaling used in Lemma~\ref{low:lem:countable} gives
nonreal zeros tending to $t$, so $q+1$ is active. Nesting preserves this
up to height $q(q+1)$, and the pronic bound at that final height excludes
every larger positive real core point.
\end{proof}

Thus an endpoint can be present throughout a long range of heights while
its relation to the other components changes. The next section proves
that, before the pronic height, this endpoint is itself a singleton
component and is approached by infinitely many other components.

\section{Cascades of global components}
\label{sec:component-cascades}

Perturbing a double zero produces nearby nonreal roots, but does not
separate their connected components. A point may be a zero of several
admissible series, so disjoint families of coefficient sequences can have
overlapping root sets. Even a separation within a small disk could be
bypassed by a connection outside it. We address both difficulties by
constructing root sets that are open and closed in the entire core.
One comparison lemma supplies the required separation throughout both
height intervals in Theorem~\ref{root:thm:infinite}.

We use the normalized family $\mathcal F_H$ from
Section~\ref{sec:foundations}. A zero is called \emph{qualifying} if it is nonreal, or if it is real
and has multiplicity at least two.  This terminology records a necessary
condition at a real limit of nonreal roots; it does not assert that
every multiple real zero belongs to the core.

\subsection{From coefficient blocks to global components}

The normalization below ensures that the coefficient comparisons cover
every core point in the chosen neighbourhood, including real limits of
nonreal roots. Convergence of bounded coefficient sequences is locally
uniform in the unit disk, including all derivatives.

\begin{lemma}[Local witnesses and positive perturbations]
\label{lem:local-normalization}\label{lem:positive-unfolding}
Put $\sigma_H=(H^2+4)^{-1/4}$.
Every $z\in\mathcal C_H$ with $0<|z|<\sigma_H$ has a witness
in $\mathcal F_H$ for which it is a qualifying zero.
The qualifying-zero relation is closed on compact subdisks.

Suppose $F(t)=F'(t)=0$ and $b=F''(t)/2>0$, where $0<t<1$.  Let $C_N$
be power series whose coefficients are bounded independently of $N$,
with $C_N(t)\ge c>0$.  If $F+z^NC_N\in\mathcal F_h$ for arbitrarily
large $N$, then every sufficiently large such $N$ gives a nonreal
conjugate pair of zeros of $F+z^NC_N$ in $\mathcal C_h$, and these
pairs tend to $t$.
\end{lemma}
\begin{proof}
For a nonreal polynomial root, remove any factor $z^k$ and make the
constant coefficient positive.  If its value is an integer $d\ge2$,
dividing by $d$ and applying Lemma~\ref{lem:conjugate-pair-bound} gives
\[
 |z|^{-4}\le4\bigl(1+(H/d)^2\bigr)\le H^2+4.
\]
Thus roots in the stated strict disk have constant-term-one witnesses.
For a core point, take a convergent sequence of nonreal polynomial
roots and a coefficientwise convergent subsequence of these witnesses.
At a real limit, the conjugate roots coalesce: their divided difference
converges to the derivative, which therefore vanishes.  This also proves
closedness of the qualifying-zero relation, using derivative convergence
when the approximating zeros are already real and multiple.

For the last assertion put $\varepsilon_N=t^{N/2}$.  Uniformly for $w$
in compact sets,
\[
 t^{-N}(F+z^NC_N)(t+\varepsilon_Nw)-\bigl(bw^2+C_N(t)\bigr)
 \longrightarrow0,
\]
because $t^{-N}F(t+\varepsilon_Nw)\to bw^2$,
$(1+\varepsilon_Nw/t)^N\to1$ since $N\varepsilon_N\to0$, and
$C_N(t+\varepsilon_Nw)-C_N(t)=O(\varepsilon_N)$ by the bound on the
coefficients.  The zeros $\pm i\sqrt{C_N(t)/b}$ of the comparison
function stay in a compact subset of $\mathbb C\setminus\mathbb R$, and
on circles of a fixed small radius about them it is bounded away from
zero uniformly in $N$.  Rouch\'e's theorem gives the required pair.
Truncating the series and applying Rouch\'e again near either nonreal
zero places the pair in the polynomial-root core.
\end{proof}

\Needspace{5\baselineskip}
If $G\ne F$, its \emph{first departure} from $F$ is the expression
\[
 G-F=z^kC,\qquad C(0)\ne0,
\]
where $k$ is the first differing coefficient.  A finite departure block
is a finite union of coefficient cylinders specified by such a degree
and a set of allowed nonzero first differences. For example, the
negative departures with degrees in a given finite interval form one.

The next lemma assumes that no block witness shares a qualifying zero
with any outside witness in a fixed neighbourhood. Compactness and
rigidity then give relative openness of the block's nearby root set in
the core. A quadratic estimate at the reference double zero places this
set strictly inside the chosen disk for sufficiently late blocks.
It is therefore clopen in the entire core and contains whole global
components, although the set itself need not be connected.

\begin{lemma}[Global separation of finite blocks]
\label{lem:finite-block-transfer}
Let $F_r\in\mathcal F_H$ converge locally uniformly to $F_*$, and let
$t_r\to t_*\in(0,\sigma_H)$.  Suppose each $F_r$ has an exact double
zero at $t_r$, is the unique member of $\mathcal F_H$ with a multiple
zero there, and
\[
 F_*(t_*)=F_*'(t_*)=0,\qquad F_*''(t_*)>0.
\]
For each $r$ choose a finite departure block $\mathscr B_r$ relative
to $F_r$, whose least degree $n_r$ tends to infinity.  Assume that in
one fixed neighbourhood of $t_*$ no witness in $\mathscr B_r$ shares
a qualifying zero with any witness outside $\mathscr B_r$.

In a sufficiently small fixed closed disk $D$ about $t_*$, let $K_r$
be the qualifying-root projection of $\mathscr B_r$, intersected with
$\mathcal C_H$.  Then $K_r$ is compact, excludes $t_r$, and, for all
large $r$, is clopen in the entire core.  There are constants $c>0$
and $\rho<1$ such that
\begin{equation}\label{eq:block-shrinking}
 K_r\subseteq
 \left\{z:|z-t_r|^2\le\frac{2H\rho^{n_r}}{c(1-\rho)}\right\}.
\end{equation}
In particular their diameters tend to zero.  For any $h_*\le H$:
\begin{enumerate}[label=\textup{(\roman*)}]
\item If the $K_r$ are pairwise disjoint and infinitely many meet
$\mathcal C_{h_*}$, they contain distinct components with nonempty
planar interior, nested as the height increases from $h_*$ to $H$.
These components accumulate at $t_*$, where none of these cores is
locally connected.
\item If every $K_r$ contains $t_*$, and $t_*$ is nonisolated in
$\mathcal C_{h_*}$, then $\{t_*\}$ is a component at every height
$h_*\le h\le H$.  Every neighbourhood of it contains infinitely many
interior-bearing components, again chosen nested across these heights
and with their diameters at height $H$ tending to zero.  The cores are
not locally connected at $t_*$.
\end{enumerate}
\end{lemma}
\begin{proof}
On a small common disk,
\[
 F_r(z)=(z-t_r)^2Q_r(z),\qquad |Q_r(z)|\ge c>0,
\]
for all large $r$.  Indeed the quotients converge to
$F_*(z)/(z-t_*)^2$, which is nonzero near $t_*$.  Choose $D$ within
the normalization and separation neighbourhoods, and put
$\rho=\max_{z\in D}|z|<1$.
The coefficient block is compact, as is its qualifying-root projection
in the closed disk and closed core.  Rigidity excludes $t_r$.  At a
root of a departure of degree $k\ge n_r$,
\[
 c|z-t_r|^2\le |F_r(z)|
       \le\frac{2H\rho^k}{1-\rho},
\]
which proves \eqref{eq:block-shrinking}.

Fix $r$ and suppose points of $(\mathcal C_H\cap D)\setminus K_r$
converge to $z\in K_r$. Their normalized qualifying witnesses lie in
$\mathcal F_H\setminus\mathscr B_r$. A finite union of coefficient
cylinders is clopen in the compact coefficient space, so this complement
is compact. Pass to a convergent subsequence of the witnesses. The limit
still lies outside $\mathscr B_r$ and has a qualifying zero at $z$, by
closedness of that relation. A witness in $\mathscr B_r$ also has a
qualifying zero there, contradicting the separation hypothesis.
Thus $K_r$ is relatively open in $\mathcal C_H\cap D$.
The shrinking estimate eventually places it strictly inside $D$, so
it is open in the entire core.  Compactness makes it globally closed.

For (i), $K_r\cap\mathcal C_{h_*}$ is nonempty and relatively open.
By Proposition~\ref{cor:whole-plane-regularity}, it contains a planar
interior point.  Choose its component at height $h_*$ and, at each
larger height, the component containing it.  These components are
nested and remain in $K_r$, giving distinct components and the asserted
shrinking.  Their accumulation puts $t_*$ in every intermediate core.
A connected relatively open neighbourhood of $t_*$ could not meet a
clopen $K_r$ that excludes $t_*$, whereas every neighbourhood meets all
sufficiently late occupied sets.  At most one of the pairwise disjoint
$K_r$ can contain $t_*$, so discarding it proves failure of local
connectedness.

For (ii), the component of $t_*$ lies in every $K_r$ and is therefore
a singleton.  To obtain disjoint sets meeting the smallest core, form
successive finite intersections $V_j$ of the $K_r$.  Nonisolation lets
us choose $x_j\in V_j\cap\mathcal C_{h_*}$ different from $t_*$;
choose the next intersection to exclude $x_j$ and to have diameter
less than $1/j$.  The disjoint clopen differences
$V_j\setminus V_{j+1}$ meet $\mathcal C_{h_*}$ and shrink to $t_*$.
Apply the component selection just proved.  A nonisolated singleton
component has no connected relatively open neighbourhood, giving the
last assertion.
\end{proof}

\subsection{A common comparison for negative departures}

The following lemma verifies the no-shared-zero hypothesis by comparing
both the values and the derivatives of departure tails. A common nonreal
zero brings its conjugate with it; as the pair approaches the real anchor,
the limiting difference has a double zero. At a shared real multiple
zero, both the difference and its derivative already vanish. The
resulting relations contradict three strict coefficient inequalities.
Those inequalities also give uniqueness of the reference witness.
Their uniformity over a compact family will allow the reference itself
to move at the height $p(p-1)-1$ just below a pronic height.

\begin{lemma}[Uniform comparison for a compact family]
\label{lem:compact-negative-separation}
Fix integers $H\ge1$ and $0\le m\le H$.  Let $\mathscr F$ be a compact
family of pairs $(F,t)$, in coefficientwise and ordinary convergence,
where $F\in\mathcal F_H$,
\[
 F(t)=F'(t)=0,\quad F''(t)>0,\quad 0<t<1,
 \qquad m\le a_j(F)\le H\quad(j\ge2).
\]
Put $\ell=H-m$ and $\beta=1/t$.  Suppose throughout the family that
\begin{equation}\label{eq:compact-comparison-conditions}
 \frac{\ell}{\beta-1}<2,\qquad
 1-\frac{\ell}{(\beta-1)^2}>0,\qquad
 \Delta:=2\beta-\ell-\frac{2H}{(\beta-1)^2}>0.
\end{equation}
Each $F$ is the unique normalized height-$H$ witness with a multiple
zero at $t$.  In one neighbourhood of each anchor, of radius independent
of the reference, two witnesses cannot share a qualifying zero if their
first departures have opposite signs, or if they are both negative and
have degrees differing by at least two.  Adjacent negative departures
at degrees $j-1,j$ are also separated whenever $a_j(F)>m$.
Consequently, if $N<j$ and $a_N(F)>m$, $a_j(F)>m$, then in that
neighbourhood no negative departure from $F$ with degree in $[N,j-1]$
shares a qualifying zero with a witness that is not such a departure.
\end{lemma}
\begin{proof}
For a first departure $G-F=z^kC$, every nonconstant coefficient of $C$
lies in $[-2H,\ell]$.  Write $v=C(t)$ and $w=tC'(t)$.  The three cases
for its nonzero integer constant coefficient give
\begin{equation}\label{eq:three-departure-cones}
\begin{array}{ll}
 C(0)\ge1:&v-w\ge1-\ell/(\beta-1)^2>0,\\[2pt]
 C(0)\le-2:&v\le-2+\ell/(\beta-1)<0,\\[2pt]
 C(0)=-1:&w-2v\ge\Delta/\beta>0.
\end{array}
\end{equation}
For example, in the last line the coefficient of degree two cancels:
\[
 tC'(t)-2C(t)
 \ge2-\ell t-2H\sum_{j\ge3}(j-2)t^j
 =\Delta/\beta.
\]
A second multiple-zero witness would require $v=w=0$, which proves
rigidity.  All three margins are uniform by compactness.

Suppose the claimed uniform separation fails.  Take reference pairs
$(F_r,t_r)\to(F_*,t_*)$ and forbidden shared qualifying zeros tending
to $t_*$.  Every coefficient limit of their witnesses has a multiple
zero at $t_*$, by Lemma~\ref{lem:local-normalization}.  Rigidity makes
that limit $F_*$.  Discreteness of the coefficients then forces both
first departure degrees to tend to infinity.

The usual two-zero argument is uniform over this compact family.
In small real-symmetric disks the reference has just its double zero,
is strictly convex on the real diameter, and is positive at its ends.
Late departures converge with two derivatives; Rouch\'e gives exactly
two zeros, and strict convexity persists.  If such a departure has a
qualifying zero, its value at $t_r$ is positive: a nonreal conjugate
pair leaves no real zero in the disk, while a real double zero is its
unique real minimum; equality at $t_r$ would contradict rigidity.
Thus $C(t_r)>0$, and limiting departure values are nonnegative.
By \eqref{eq:three-departure-cones}, every negative departure has
constant coefficient $-1$.

Order the two departure degrees as $k_r\le l_r$.  First let the earlier
one be positive and the later one negative, ordering a mixed pair at
equal degrees this way.  If $l_r-k_r=s$ is fixed, coefficient limits
of their tails $C,D$ satisfy
\[
 v_C=t_*^sv_D,\qquad
 w_C=t_*^s(w_D+sv_D),\qquad v_C,v_D\ge0.
\]
These equations follow because the difference has coalescing conjugate
zeros, or a real multiple zero.  The two strict bounds give
\[
 v_C-w_C=t_*^s\bigl((1-s)v_D-w_D\bigr)
 \le-t_*^s\bigl((s+1)v_D+\Delta_*/\beta_*\bigr)<0,
\]
contrary to the positive-departure bound.  If instead the delays tend
to infinity, $C_r-z^{l_r-k_r}D_r$ converges locally uniformly to $C$.
Its coalescing zeros force $C(t_*)=C'(t_*)=0$, again impossible.

When the earlier departure is negative, compare the actual witnesses
$G_r$ and $\widetilde G_r$, which depart at degrees $k_r$ and $l_r$:
\[
 E_r(z)=\frac{G_r(z)-\widetilde G_r(z)}{z^{k_r}}.
\]
Its constant coefficient is at most $-1$, and all coefficients of degree
at least two lie in $[-2H,2H]$.  Its degree-one coefficient is at most
\[
\begin{array}{c|c}
 l_r-k_r\ge2&H-a_{k_r+1}(F_r)\le\ell\\
 l_r-k_r=1,\ \text{later positive}
   &H-a_{k_r+1}(F_r)-1\le\ell-1\\
 l_r-k_r=1,\ \text{later negative},\ a_{k_r+1}(F_r)>m
   &H-a_{k_r+1}(F_r)+1\le\ell.
\end{array}
\]
The final line uses integrality and the later constant coefficient
$-1$.  Passing to a coefficient limit $E$ gives
$E(t_*)=E'(t_*)=0$, for bounded or unbounded delays alike.  Yet
\[
 t_*E'(t_*)-2E(t_*)
 \ge2-\ell t_*-\frac{2Ht_*^3}{(1-t_*)^2}
 =\frac{\Delta_*}{\beta_*}>0.
\]
This contradiction proves every remaining asserted separation.  The
reference coefficients cancel once both witnesses have departed;
retaining this cancellation is what yields the third inequality in
\eqref{eq:compact-comparison-conditions}.

For the final assertion, a witness that is not such a departure is $F$
itself, a positive departure, or a negative departure of degree at most
$N-1$ or at least $j$.  In the disks above, $F$ has no zero other than
$t$, and rigidity excludes $t$ from every other witness.  The other
cases are the separations just proved, with the markers $a_N(F)>m$ and
$a_j(F)>m$ for the adjacent degrees.
\end{proof}

\subsection{A local criterion and the heights above a pronic number}

The next proposition isolates what the separation argument needs from
a single series: the comparison inequalities at one double zero, and
infinitely many levels at which one coefficient can be lowered and the
later ones raised to compensate.

\begin{proposition}[A local criterion]\label{prop:local-criterion}
Let $0\le m<h_*\le H$ be integers, put $\ell=H-m$, and let
$F\in\mathcal F_{h_*}$ have a double zero at $t=1/\beta\in(0,\sigma_H)$
and coefficients $a_j(F)\ge m$ for $j\ge2$. Suppose that the three
inequalities~\eqref{eq:compact-comparison-conditions} hold at $\beta$,
and that for some $\eta>0$ there are infinitely many indices $N$ with
\begin{equation}\label{eq:descending-condition}
 a_N(F)>m,\qquad \sum_{j\ge1}\bigl(h_*-a_{N+j}(F)\bigr)t^j\ge1+\eta .
\end{equation}
Then, at every height $h_*\le h\le H$, infinitely many distinct
connected components of $\mathcal C_h$ with nonempty planar interior
accumulate at each of the four points $\pm t$ and $\pm1/t$, at which
$\mathcal C_h$ is not locally connected. These components can be chosen
nested as $h$ increases, with their diameters at height $H$ tending to
zero. In particular $t\in\mathcal A_{h_*}$.
\end{proposition}
\begin{proof}
Negation and inversion transfer the statements at $t$ to the other
three points. Since $a_j(F)\ge m\ge0$ for $j\ge2$ and some $a_N(F)$ is
positive, $F''(t)>0$. Hence Lemma~\ref{lem:compact-negative-separation}
applies at height $H$ to the family consisting of the single pair
$(F,t)$: the series $F$ is the unique member of $\mathcal F_H$ with a
multiple zero at $t$, and $t$ has a separation neighbourhood.

Every index $N$ in~\eqref{eq:descending-condition} is followed by a
least index $j_N>N$ with $a_{j_N}(F)>m$, since there are infinitely
many such indices. Choose indices $N_1<N_2<\cdots$
satisfying~\eqref{eq:descending-condition} with $N_{i+1}\ge j_{N_i}$, and
let $\mathscr B_i\subset\mathcal F_H$ be the block of negative
departures from $F$ with degrees in $[N_i,j_{N_i}-1]$. These blocks are
pairwise disjoint, and by the last assertion of
Lemma~\ref{lem:compact-negative-separation} no witness in
$\mathscr B_i$ shares a qualifying zero in the separation neighbourhood
with a witness outside $\mathscr B_i$.
Lemma~\ref{lem:finite-block-transfer}, applied to the constant sequence
$(F,t)$, gives globally clopen subsets $K_i$ of $\mathcal C_H$. They are
pairwise disjoint, because a common point would be a shared qualifying
zero of two different blocks.

For occupancy, let $G_i=F+z^{N_i}C_i$ be the series that agrees with
$F$ below degree $N_i$, has the digit $a_{N_i}(F)-1\ge m$ in degree
$N_i$, and has the constant digit $h_*$ thereafter. Then
$G_i\in\mathscr B_i\cap\mathcal F_{h_*}$, the coefficients of $C_i$ are
bounded by $h_*$, and
\[
 C_i(t)=-1+\sum_{j\ge1}\bigl(h_*-a_{N_i+j}(F)\bigr)t^j\ge\eta>0 .
\]
Lemma~\ref{lem:positive-unfolding} provides nonreal zeros of $G_i$ in
$\mathcal C_{h_*}$ tending to $t$; they lie in $K_i$ for all large $i$.
Part~(i) of Lemma~\ref{lem:finite-block-transfer} gives the components,
their nesting and the failure of local connectedness, and
$t\in\mathcal A_{h_*}$ follows from the accumulation.
\end{proof}

\begin{corollary}[Components accumulating at a central anchor]
\label{cor:push-negative-lower-bulk}
For every integer $p\ge4$ and
\[
 p(p-1)+1\le h\le p^2-3,
\]
infinitely many distinct global components of $\mathcal C_h$ with
nonempty planar interior accumulate at each of $\pm p,\pm1/p$.
They can be chosen nested across this entire height interval, with
their diameters at its largest height tending to zero.  The four anchors
belong to the component containing the unit circle, and the cores are
not locally connected there.
\end{corollary}
\begin{proof}
Apply Proposition~\ref{prop:local-criterion} to
\[
 F(z)=\frac{(1-pz)^2}{1-z}
 =1-(2p-1)z+(p-1)^2\sum_{j\ge2}z^j,
 \qquad t=1/p,
\]
with $h_*=p(p-1)+1$, $H=p^2-3$ and $m=(p-1)^2-1$, so that
$\ell=2p-3$.  The three positive margins are
\[
 2-\frac\ell{p-1}=\frac1{p-1},\qquad
 1-\frac\ell{(p-1)^2}=\frac{(p-2)^2}{(p-1)^2},\qquad
 2p-\ell-\frac{2H}{(p-1)^2}=\frac{(p-3)^2}{(p-1)^2},
\]
and $t<\sigma_H$ since $p^4-H^2-4=6p^2-13>0$.  Every coefficient of $F$
of degree at least two exceeds $m$, and for every $N$ the sum
in~\eqref{eq:descending-condition} equals
$(h_*-(p-1)^2)/(p-1)=1+1/(p-1)$.  The connected core at height $p(p-1)$
contains the four anchors and the unit circle by
Proposition~\ref{prop:pronic-core}; inclusion in each larger core puts
them in its component containing the unit circle.
\end{proof}

With an arbitrary height $h$ in place of $h_*$, the sum
in~\eqref{eq:descending-condition} equals $(h-(p-1)^2)/(p-1)$, which
exceeds $1$ exactly when $h>p(p-1)$.  At the pronic height it equals
$1$ and the criterion just fails, as Theorem~\ref{root:thm:main}
requires.

\subsection{Singleton components below a pronic height}

At the height $p(p-1)-1$ we use nearby reference series.  A finite
negative block for each reference will contain the same limiting anchor.
This produces shrinking clopen neighbourhoods of that anchor without any
recurrence assumption on a single coefficient sequence.  The references
come from an intermediate-value argument, which
Section~\ref{sec:gap-heights} uses again.

For integers $0\le m'<H'$ and a finite integer sequence
$\pi=(a_0,\ldots,a_{N-1})$ with $N\ge2$, write
\[
 f^-_\pi(z)=\sum_{j<N}a_jz^j+\frac{m'z^N}{1-z},\qquad
 f^+_\pi(z)=\sum_{j<N}a_jz^j+\frac{H'z^N}{1-z}
\]
for the smallest and the largest of the \emph{completions}
$\sum_{j<N}a_jz^j+\sum_{j\ge N}a_jz^j$ of $\pi$ with integer digits
$a_j\in[m',H']$ for $j\ge N$.

\begin{lemma}[Completion with a prescribed minimum]\label{lem:ivt-completion}
Let $J=[a,b]\subset(0,1)$ satisfy $(H'-m')a\ge1-a$, and suppose that
\[
 f^-_\pi(a)\ge0,\qquad f^-_\pi(b)\ge0,\qquad
 \min_J f^-_\pi<0\le\min_J f^+_\pi .
\]
Then $\pi$ has a completion $F$ with $\min_JF=0$, and the minimum is
attained at an interior point $t$ of $J$. If moreover $a_j\ge0$ for
$2\le j<N$, then $F''>0$ on $(0,1)$, and $t$ is the only zero of $F$ in
$(0,1)$ and is exactly double.
\end{lemma}
\begin{proof}
Let $\pi'$ be an extension of $\pi$ by digits in $[m',H']$, with next
free degree $n\ge N$, and let $d<H'$ be a digit. On $J$,
\[
 f^+_{\pi'd}-f^-_{\pi'(d+1)}
 =z^n\Bigl((H'-m')\frac z{1-z}-1\Bigr)\ge0,
\]
by the hypothesis on $a$ and monotonicity of $z/(1-z)$. Therefore the
intervals $[\min_Jf^-_{\pi'd},\min_Jf^+_{\pi'd}]$, for
$d=m',\ldots,H'$, overlap consecutively, and their union is the interval
$[\min_Jf^-_{\pi'},\min_Jf^+_{\pi'}]$, since the extreme completions of
$\pi'$ coincide with those of its extreme children. Starting from
$\pi$, choose successively a child whose interval contains $0$. The
resulting completions converge uniformly on $J$ with all derivatives,
since the tails are bounded by $H'b^n/(1-b)$, and $\min_JF=0$ for the
limit $F$.

$F$ is not the completion $f^-_\pi$, whose minimum on $J$ is negative.
Hence some digit of $F$ beyond $\pi$ exceeds $m'$, and
$F-f^-_\pi=\sum_{j\ge N}(a_j(F)-m')z^j$ is strictly positive on
$(0,1)$. Thus $F(a)>0$ and $F(b)>0$, and the minimum $0$ of $F$ on $J$
is attained at an interior point $t$, where $F'(t)=0$. Under the last
hypothesis the coefficients of $F$ of degree at least two are
nonnegative and one of them is positive, so $F''>0$ on $(0,1)$; then
$t$ is the unique zero of $F$ there, and the zero is exactly double.
\end{proof}

\begin{theorem}[Below a pronic height]
\label{thm:moving-terminal-blocks}
For every integer $p\ge3$ and
\[
 (p-1)^2+1\le h\le p(p-1)-1,
\]
the four points $\pm p$ and $\pm1/p$ are singleton global components
of $\mathcal C_h$.  Every neighbourhood of each of them contains
infinitely many distinct global components with nonempty planar interior.
These components can be chosen nested as $h$ increases through the
stated interval, with their diameters at its largest height tending to
zero.  None of these cores is locally connected at the four anchors.
\end{theorem}
\begin{proof}
Put
\[
 \begin{gathered}
 H=p(p-1)-1,\quad A=(p-1)^2,\quad H_0=A+1,\quad m=(p-1)(p-2),\\
 t_0=1/p,\qquad F_0(z)=\frac{(1-pz)^2}{1-z}.
 \end{gathered}
\]
We construct references with a double zero $t_N>t_0$ tending to $t_0$,
all admissible already at height $H_0$; note that
$2p-1\le H_0\le H$ for $p\ge3$.
For sufficiently large $N$, let $r_N$ be the first root to the right
of $t_0$ of $1-pz+z^N$.  This root exists and tends to $t_0$, since
$1-pz+z^N$ converges with its derivative to $1-pz$ on a fixed small
interval about $t_0$, and its value at $t_0$ is positive.  Set
$I_N=[t_0,r_N]$ and consider the integer-tail completions
\begin{equation}\label{eq:moving-reference-family}
 F_N(z)=1-(2p-1)z+A\sum_{j=2}^{N-1}z^j
             +(A+1)z^N+\sum_{j>N}a_jz^j,
 \qquad m\le a_j\le A.
\end{equation}
The all-$m$ and all-$A$ completions are respectively
\[
 \frac{(1-pz)(1-pz+z^N)}{1-z},\qquad F_0(z)+z^N.
\]
The first vanishes at both endpoints of $I_N$ and is negative between
them, while the second is everywhere positive.  Since
$(A-m)t_0=1-t_0$, Lemma~\ref{lem:ivt-completion} with the digit range
$[m,A]$ gives a completion $F_N$ of the
form~\eqref{eq:moving-reference-family} whose minimum on $I_N$ is zero.
It is attained at an interior point $t_N$, which is the unique zero of
$F_N$ in $(0,1)$ and is exactly double, and $F_N''>0$ there.  We have
$t_N\to t_0$ and $F_N\to F_0$ locally uniformly.

The compact comparison lemma applies to this family after discarding
finitely many initial $N$, together with its limit $(F_0,t_0)$.
All coefficients of degree at least two belong to $[m,H]$, including
$A+1\le H$.  At the limit, $\ell=H-m=2p-3$, and the three margins are
\[
 2-\frac\ell{p-1}=\frac1{p-1},\qquad
 1-\frac\ell{(p-1)^2}=\frac{(p-2)^2}{(p-1)^2},\qquad
 2p-\ell-\frac{2H}{(p-1)^2}
      =\frac{(p-2)^2+1}{(p-1)^2}.
\]
They are strict for every $p\ge3$ and persist uniformly by continuity.
Also $p^2>H+1$ gives $p^4>H^2+4$, so a common normalization disk is
available.  We therefore have full-height rigidity and one common
separation neighbourhood for every reference, including $F_0$.

Let $j_N>N$ be the first later index with $a_{j_N}(F_N)>m$.  It exists
because $F_N$ is not the all-$m$ completion, and $a_N(F_N)=A+1>m$.  Let
$\mathscr B_N$ be the block of negative departures from $F_N$ with
degrees in $[N,j_N-1]$.  By the last assertion of
Lemma~\ref{lem:compact-negative-separation}, all hypotheses of
Lemma~\ref{lem:finite-block-transfer} are satisfied.

The resulting globally clopen set $K_N$ contains $t_0$: the witness
$F_0$ first departs negatively from $F_N$ at degree $N$, and $t_0$
belongs to $\mathcal C_{H_0}$ by the admissible perturbations
$F_0+z^k$.  These same perturbations and
Lemma~\ref{lem:positive-unfolding} show that $t_0$ is nonisolated
already at height $H_0$.  The diameters of $K_N$ tend to zero, regardless
of the possibly unbounded lengths of the selected finite blocks.
Part (ii) of the transfer lemma now gives the singleton, interior,
nesting and local connectedness conclusions at $t_0$.  Negation and
inversion give the other three anchors.
\end{proof}

\begin{proof}[Proof of Theorem~\ref{root:thm:infinite}]
Corollary~\ref{cor:push-negative-lower-bulk} and
Theorem~\ref{thm:moving-terminal-blocks} cover the two sides of the
stated interval; when $p=3$ its part above $p(p-1)$ is empty.
\end{proof}

Above the pronic height the components accumulate on a component
already present at that height; below it, at a singleton component that
remains separated until the next pronic height.  Both conclusions rest
on globally clopen sets, so no reconnection is possible outside the disk
used in the proof.

\subsection{Sparse coefficient alphabets}

Separation in the full height-$H$ core also controls root sets obtained
by deleting coefficient choices. The following consequence needs only
three or four nonconstant digits to produce the accumulating components.
For a finite integer alphabet $S$, define
\[
 \mathcal C(S)=\overline{\left\{z\in\mathbb C\setminus\mathbb R:
  1+\sum_{j=1}^{d}a_jz^j=0,\quad d\ge1,\quad
  a_j\in S,\quad a_d\ne0\right\}}.
\]
The constant coefficient is fixed at $1$; it need not belong to $S$.

\begin{corollary}[Persistence under changes of alphabet]
\label{cor:sparse-alphabets}
Put $A=(p-1)^2$ and $t=1/p$.
\begin{enumerate}[label=\textup{(\roman*)}]
\item Let $p\ge3$, $H=p(p-1)-1$, and
\[
 S^-_p=\{-(2p-1),A,A+1\}\subseteq S\subseteq[-H,H]\cap\mathbb Z.
\]
Then $\{t\}$ is a component of $\mathcal C(S)$.
\item Let $p\ge4$, $H=p^2-3$, and
\[
 S^+_p=\{-(2p-1),A-1,A,p(p-1)+1\}
       \subseteq S\subseteq[-H,H]\cap\mathbb Z.
\]
Then $\mathcal C(S)$ has infinitely many components accumulating at $t$.
\end{enumerate}
In both cases every neighbourhood of $t$ contains infinitely many
components containing nonreal points, and $\mathcal C(S)$ is not locally
connected at $t$. These components can be chosen to persist under every
allowed enlargement of $S$: if $S\subseteq T$, the component chosen for
$S$ lies in the corresponding component for $T$. Their diameters tend
to zero uniformly over all the allowed alphabets.
\end{corollary}
\begin{proof}
For either value of $H$, the set $\mathcal C(S)$ is a compact subset of
$\mathcal C_H$, and it increases with $S$. Thus the intersection of a
globally clopen subset of $\mathcal C_H$ with $\mathcal C(S)$ is clopen
in $\mathcal C(S)$.

For (ii), use the disjoint shrinking clopen sets constructed in the
proof of Proposition~\ref{prop:local-criterion} for the series of
Corollary~\ref{cor:push-negative-lower-bulk}; write $K_N$ for the set
whose block is the degree $N$. Their occupying witnesses are
\[
 \frac{(1-pz)^2}{1-z}
   +z^N\left(-1+\frac{pz}{1-z}\right).
\]
Their nonconstant coefficients are $-(2p-1)$ in degree one, $A$ in
degrees $2,\ldots,N-1$, $A-1$ in degree $N$, and $p(p-1)+1$
thereafter. Their nearby nonreal zeros therefore belong to
$\mathcal C(S^+_p)$: truncate the series and apply Rouch\'e's theorem
in a small disk disjoint from the real axis. Choose one such zero
$z_N\in K_N$. For every allowed $S$, the component of $z_N$ in
$\mathcal C(S)$ remains in $K_N$. These components are distinct,
shrink uniformly to $t$, and are nested when the alphabet increases.
Their limit $t$ belongs to $\mathcal C(S^+_p)$. Since $t\notin K_N$,
no connected relatively open neighbourhood of
$t$ can meet these sets, proving failure of local connectedness.

For (i), the proof of Theorem~\ref{thm:moving-terminal-blocks}
provides shrinking globally clopen neighbourhoods $K_N$ of $t$.
The witnesses $(1-pz)^2/(1-z)+z^k$ use only $S^-_p$, so their
nonreal zeros show that $t$ is nonisolated in $\mathcal C(S^-_p)$.
The intersection of the $K_N$ is $\{t\}$, proving the singleton
assertion for every allowed $S$. Successive finite intersections
$V_j$ of the $K_N$ can be chosen with diameters tending to zero so
that each difference $V_j\setminus V_{j+1}$ contains a nonreal point
of $\mathcal C(S^-_p)$. Select the component through that point at
each alphabet. The same clopen argument proves all the remaining
claims, including failure of local connectedness at the nonisolated singleton.
\end{proof}

For example, the alphabet $\{-5,4,5\}$ gives a singleton component at
$1/3$, and $\{-7,8,9,13\}$ gives infinitely many components
accumulating at $1/4$. No density of planar interior is needed for
these conclusions.

\section{Completing the connectedness classification}\label{sec:classification}

The preceding section proves more than disconnectedness at every height
in Theorem~\ref{root:thm:infinite}. The only heights still needed for
Theorem~\ref{root:thm:main} are the three heights next to each square,
together with a few small cases. We treat these by planar separators.
This avoids repeating the more general separation constructions on the
ranges already covered by infinitude. Section~\ref{sec:gap-heights}
proves infinitely many components at most of the square-adjacent
heights; the separators below cover all of them uniformly and remain
necessary at the heights not reached there.

\subsection{What the finite exclusions certify}
For $c=x+iy$, $x>1$, put $R=|c|^2$ and
\[
 T_h(c)=\left\{\sum_{j\ge0}d_jc^{-j}:d_j\in B_h\right\}.
\]
The normalized series criterion says $c\in\Sigma_h$ if and only if
$c\in T_h(c)$. Its inverse-digit orbit is $z\mapsto c(z-d)$.
For $y\ne0$, use the coordinates
\[
 V=\frac{\im z}{y},\qquad S=\frac{\im(cz)}{y}.
\]
Starting from $z_0=c$, they satisfy
\begin{equation}\label{eq:separator-recurrence}
 (V_0,S_0)=(1,2x),\qquad
 V_{j+1}=S_j-d_{j+1},\qquad S_{j+1}=2xV_{j+1}-RV_j.
\end{equation}
For a fixed prefix these are polynomials in $x,R$, so the formulas
continue to $y=0$. There they agree with
\eqref{eq:bounded-continuation} under $S_j=x(P_j+V_j)$.

Every point of $T_h(c)$ obeys the necessary bounds
\begin{equation}\label{eq:separator-bounds}
 |V|\le\frac h{(|c|-1)^2},\qquad
 |S|\le h+\frac h{(|c|-1)^2}.
\end{equation}
Indeed $|\im(c^{-k})|/|y|\le k|c|^{-k-1}$; summing over the
digits gives the first bound, and multiplying by $c$ gives the second.
Thus every possible next digit lies in a finite interval. A certificate
retains all of them and proves that every resulting branch eventually
violates~\eqref{eq:separator-bounds}.

Strict exclusion persists on a parameter neighbourhood. It therefore
excludes a point from $\mathcal C_h$, including when the point is real:
the real full locus is present there, but no nearby nonreal point can
belong to it. This is why these tests can certify a separator crossing
the real axis.

\begin{proposition}[Residual planar separation]\label{prop:residual-separators}
The exterior core $\mathcal C_h^{\rm ext}$ has at least three connected
components whenever
\[
 q^2-2\le h\le q^2\quad(q\ge3),\qquad
 \text{or}\qquad h\in\{1,3,4\}.
\]
\end{proposition}

The proofs and exact certificate statements are in
Sections~\ref{sec:sep:framework}--\ref{app:certified} of the supplement.
They consist of two ingredients: a curve disjoint from the core and
a core point in each of its separated regions. For a bounded Jordan
domain $U$ with $\partial U\cap\mathcal C_h=\varnothing$, the set
$U\cap\mathcal C_h$ is relatively clopen. A witness in $U$, its
negative in $-U$, and the central annulus then give three components.
The same reasoning applies to two reflected proper curves on the sphere.
No path-connectedness assertion about the core is needed.

For large $q$ the boundary inequalities are uniform in $q$.
The remaining parameter intervals are covered by rational subintervals.
After clearing positive denominators, an inequality is verified by
writing its polynomial in Bernstein form on the entire closed interval;
strictly positive coefficients certify positivity. The covers include
all joins, and the digit recursion keeps every admissible branch.
Interior-entry inequalities or a Rouch\'e enclosure certify the witnesses.
The exact-arithmetic checks of the companion~\cite{Companion26} verify
these assertions, rather than sampling the displayed curves.

\begin{figure}[htbp]
\centering
\includegraphics[width=.86\textwidth]{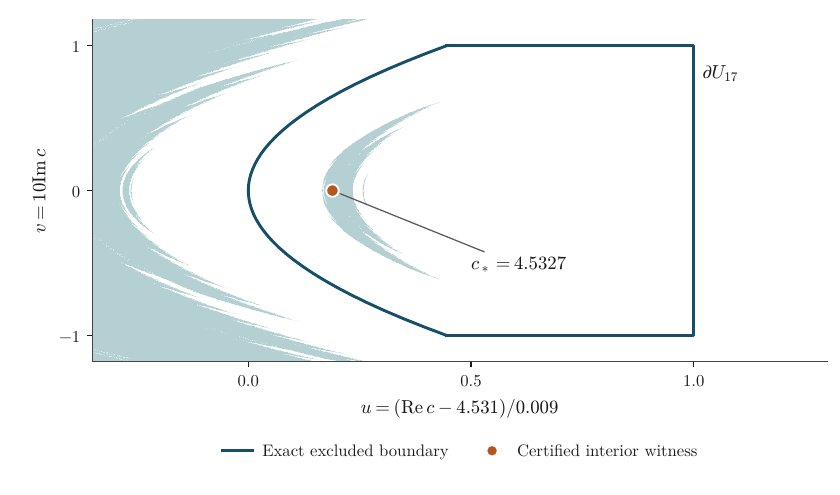}
\caption{The two parts of a planar separation certificate at height
$h=16$. In coordinates $u=(\re c-4.531)/0.009$, $v=10\im c$, the
domain is $|v|<1$, $(4/9)v^2<u<1$. Its blue boundary is disjoint from
the core, while the orange point $c_*=4.5327$ belongs to its interior
and to the planar interior of the root locus. Reflection supplies a
second separated region. All decimals specify exact rational numbers.
Pale shading is a depth-$48$ numerical illustration; the certificate,
not the shading, proves separation.}
\label{fig:bounded-separator}
\end{figure}

\subsection{Passing from the exterior to the whole plane}
Inversion doubles each component away from the central annulus.
The following elementary fact makes this precise.

\Needspace{14\baselineskip}
\begin{lemma}[Reciprocal gluing]\label{root:lem:reciprocal-components}
Let $K$ be bounded and relatively closed in $X$, invariant under
conjugation, and contain $\{1<|z|<b\}$ for some $b>1$.
Let $K_0$ be the component containing this annulus.
The components of $K\cup\mathbb T\cup\iota(K)$ are
\[
 K_0\cup\mathbb T\cup\iota(K_0),\qquad
 B,\ \iota(B)\quad(B\text{ a component of }K,\ B\ne K_0).
\]
In particular, connectedness is preserved; $m<\infty$ exterior
components give $2m-1$ components in the whole plane.
\end{lemma}
\begin{proof}
Fix $1<r_0<b$. The radial deformation
\[
 z\longmapsto\bigl((1-t)|z|+t\max\{|z|,r_0\}\bigr)\frac z{|z|}
 \qquad(0\le t\le1)
\]
retracts both $K$ and $K\cup\mathbb T$ onto
$K\cap\{|z|\ge r_0\}$. Its trajectories stay in the annulus or fix
their starting point, so adjoining $\mathbb T$ changes only $K_0$.
Every other component lies in $\{|z|\ge b\}$.

The continuous fold that fixes $|z|\ge1$ and sends
$z$ to $1/\bar z$ on $|z|\le1$ maps a connected subset into one
component of $K\cup\mathbb T$. The preimage of the central component
is connected through the annulus. The preimage of any other component
has an exterior part in $|z|\ge b$ and an interior part in
$|z|\le b^{-1}$; they are separated. Conjugation permutes the
components of $K$, giving the stated list.
\end{proof}

\begin{proof}[Proof of Theorem~\ref{root:thm:main}]
Proposition~\ref{prop:pronic-core} proves connectedness and local
connectedness at every $h=q(q-1)$. At all other heights covered by
Theorem~\ref{root:thm:infinite}, there are infinitely many components.
The remaining nonpronic heights are exactly
\[
 \{1,3,4\}\ \cup\ \bigcup_{q\ge3}\{q^2-2,q^2-1,q^2\}.
\]
Proposition~\ref{prop:residual-separators} gives at least three exterior
components at each of them. Apply the gluing lemma with
$b=\sqrt{h+1}$ and~\eqref{eq:core-gluing}. The whole-plane core then
has at least five components. The assertion about the full nonzero
closure was proved in Section~\ref{sec:foundations}.
\end{proof}

\section{The square-adjacent heights}\label{sec:gap-heights}

\begin{figure}[tbp]
\centering
\includegraphics[width=\textwidth]{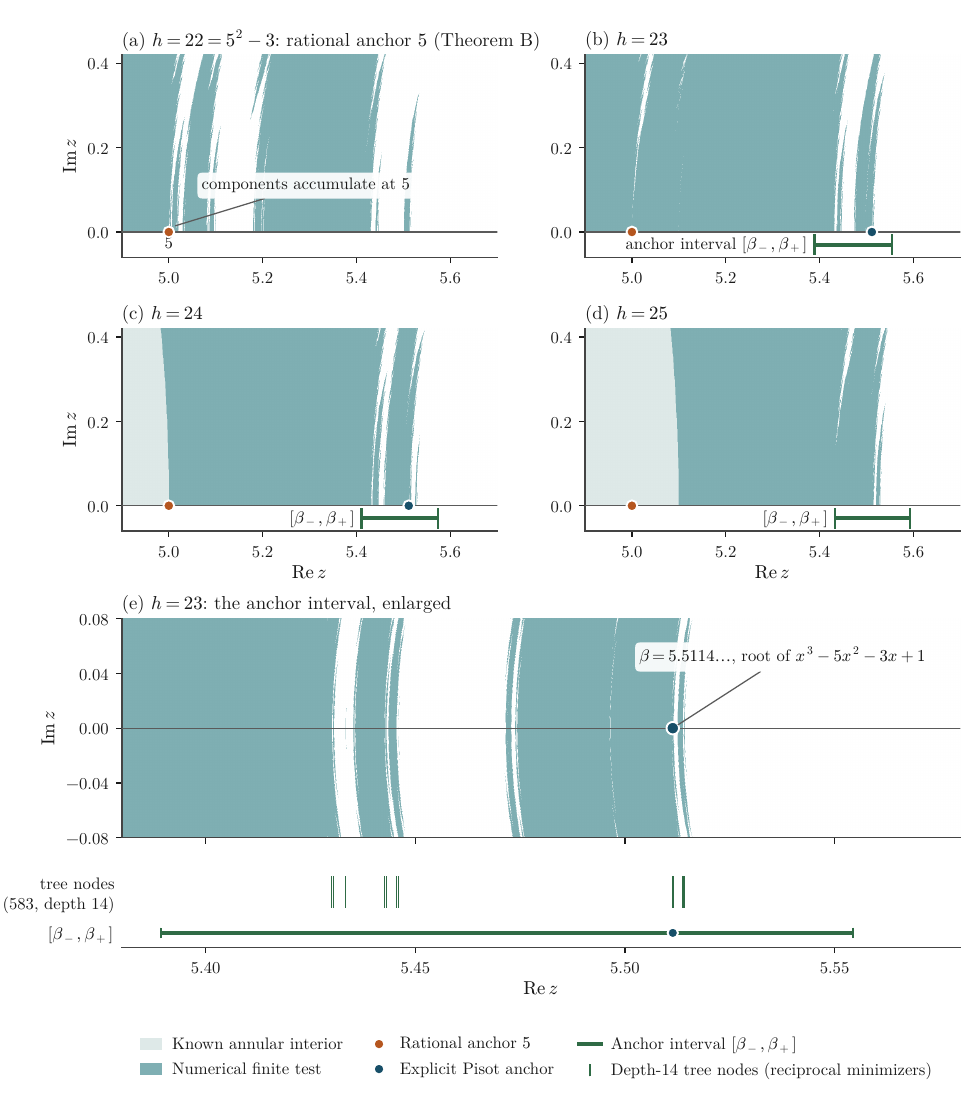}
\caption{The anchor interval next to the square $25$. Panels (a)--(d)
show the exterior core near the real axis at heights $22$ to $25$. At
$h=22=5^2-3$, components accumulate at the rational anchor $5$
(Theorem~\ref{root:thm:infinite}). At $h=25$ the point $5$ is
interior by~\eqref{eq:core-annulus}; at $h=23,24$ the plots suggest
the same. The bracket
marks the reciprocal anchor interval $[\beta_-,\beta_+]$ of
Lemma~\ref{lem:sa-anchors}. At $h=23,24$ the blue point is the
explicit Pisot reciprocal anchor of Proposition~\ref{prop:explicit-witnesses};
no such assertion is made here for $h=25$. Panel (e) enlarges the
interval at $h=23$. Its strip shows numerical reciprocal minimizers for
$583$ candidate prefixes at depth $14$. The shading and strip are
illustrative; the interval bounds and marked anchors have exact definitions.}
\label{fig:anchor-interval}
\end{figure}

Theorem~\ref{root:thm:infinite} leaves open the three heights next to
each square, $q^2-2$, $q^2-1$ and $q^2$. There the mechanism of
Section~\ref{sec:component-cascades} cannot be anchored at $1/q$: the
comparison margin $2-\ell/(\beta-1)$ of
Lemma~\ref{lem:compact-negative-separation} is not positive at $\beta=q$
for the digit range that the construction needs. Indeed $q$ itself soon
becomes an interior point of the core: for $q\ge4$ this first happens at
one of the three heights, since the core is not locally connected at
$q$ at height $q^2-3$ (Theorem~\ref{root:thm:infinite}) and contains a
neighbourhood of $q$ at height $q^2$ by~\eqref{eq:core-annulus}; which
of the three it is remains open
(Figure~\ref{fig:anchor-interval}). Nor can the mechanism be anchored
at $1/(q+1)$: the saturated series $(1-(q+1)z)^2/(1-z)$ is admissible
only from height $q^2$ on, and its double zero is inactive there
(Proposition~\ref{prop:endpoint-activation}).

We instead use series with a double zero strictly between $1/(q+1)$ and
$1/q$, the \emph{square-adjacent witnesses} defined below.
Theorem~\ref{thm:live-witness} shows that one witness outside three
exceptional types suffices for infinitely many components, and
Theorem~\ref{thm:two-of-three} finds such a witness common to at least
two of the three heights. In the explicit examples of
Propositions~\ref{prop:explicit-witnesses}
and~\ref{prop:square-witnesses} the reciprocal anchor is an algebraic
integer of degree at most nine.

Throughout the section, $q\ge4$ and $k\in\{0,1,2\}$ are fixed, and
\begin{equation}\label{eq:gap-parameters}
 h=q^2-2+k,\qquad m=q^2-2q+k,\qquad W=h-m=2q-2 .
\end{equation}
Thus $m=(q-1)^2-1+k$; the three values of $k$ give the three heights
next to the square $q^2$.

\subsection{The square-adjacent family and its anchors}

\begin{definition}\label{def:square-adjacent-family}
The \emph{square-adjacent family} at height $h$ is
\[
 \Phi_h=\Bigl\{\,1-2qz+\sum_{j\ge2}a_jz^j:\ a_j\in\mathbb Z,\ m\le a_j\le h\Bigr\}
 \subset\mathcal F_h .
\]
A \emph{square-adjacent witness} is a member of $\Phi_h$ with a real double zero
in $(0,1)$. Its double zero is called its \emph{anchor}. The set of all
anchors of square-adjacent witnesses is the \emph{anchor set}
$\mathcal T_h\subset(0,1)$.
\end{definition}

Every $F\in\Phi_h$ satisfies $F''(z)=\sum_{j\ge2}j(j-1)a_jz^{j-2}>0$
on $(0,1)$. Hence $F$ is strictly convex there, has at most one
double zero, and a square-adjacent witness is nonnegative on $(0,1)$ with its
anchor as the unique zero and the unique minimum point. The lowest
member of the family,
\begin{equation}\label{eq:lowest-member}
 M_m(z)=1-2qz+\frac{mz^2}{1-z}
 =\frac{1-(2q+1)z+(q^2+k)z^2}{1-z},
\end{equation}
is negative exactly on the interval $I_0=(t_-,t_+)$ with
\begin{equation}\label{eq:I0-endpoints}
 t_\pm=\frac{(2q+1)\pm\sqrt{4q+1-4k}}{2(q^2+k)} .
\end{equation}
Since every $F\in\Phi_h$ satisfies $F\ge M_m$ on $[0,1)$, all anchors
lie in $I_0$, and every member of $\Phi_h$ other than $M_m$ is strictly
positive at $t_-$ and at $t_+$.

Put
\begin{equation}\label{eq:phi-and-delta}
 \varphi(\beta)=2(\beta-q)(\beta-1)^2,\qquad
 \delta^*=\frac{\sqrt{q^2+4q-4}-q}{2},
\end{equation}
so that $\varphi$ increases on $[q,\infty)$ from $0$, and $\delta^*$ is
the positive root of $\delta^2+q\delta=q-1$. Let $\beta_-<\beta_+$ be
the solutions in $(q,\infty)$ of $\varphi(\beta_-)=m$ and
$\varphi(\beta_+)=h$.

\begin{lemma}[Anchors of square-adjacent witnesses]\label{lem:sa-anchors}
Every square-adjacent witness has anchor $t=1/\beta$ with
$\beta_-\le\beta\le\beta_+$, and
\[
 q+\tfrac14<\beta_-<\beta_+<q+\delta^*<q+1 .
\]
\end{lemma}
\begin{proof}
Let $F=1+\sum_{j\ge1}a_jz^j$ be a square-adjacent witness with anchor $t=1/\beta$.
In the notation of~\eqref{eq:pronic-tail-recurrence} with $d_j=-a_j$,
the double-zero equations give $P_0=V_0=1$,
$P_1=\beta-2q$ and $V_1=\beta V_0+P_1=2\beta-2q$, while the series
expression $V_1=\sum_{i\ge1}(i-1)a_{i+1}\beta^{-i}$ and the digit bounds
$m\le a_j\le h$ for $j\ge3$ give
\[
 \frac m{(\beta-1)^2}\le 2\beta-2q\le\frac h{(\beta-1)^2},
 \qquad\text{that is}\qquad m\le\varphi(\beta)\le h .
\]
Monotonicity of $\varphi$ gives $\beta_-\le\beta\le\beta_+$.
The lower bound follows from
$\varphi(q+\tfrac14)=\tfrac12(q-\tfrac34)^2<q^2-2q\le m$, valid for
$q\ge3$. For the upper bound we show $\varphi(q+\delta^*)>h$.
Since $\sqrt{q^2+4q-4}<q+2$, we have $\delta^*=2(q-1)/(\sqrt{q^2+4q-4}+q)
>(q-1)/(q+1)$, hence
\[
 \varphi(q+\delta^*)
 >\frac{2(q-1)}{q+1}\Bigl(\frac{(q-1)(q+2)}{q+1}\Bigr)^2
 =\frac{2(q-1)^3(q+2)^2}{(q+1)^3}
 \ge q^2\ge h\qquad(q\ge5),
\]
because $2(q-1)^3(q+2)^2-q^2(q+1)^3$ is positive at $q=5$ and increasing.
For $q=4$ one has $\delta^*=\sqrt7-2$ and
$\varphi(4+\delta^*)=2(\sqrt7-2)(\sqrt7+1)^2=8\sqrt7-4>16\ge h$.
Finally $\delta^*<1$ because $\delta^2+q\delta$ is increasing and equals
$q+1>q-1$ at $\delta=1$.
\end{proof}

The next lemma verifies the three inequalities of
Lemma~\ref{lem:compact-negative-separation} at every possible anchor,
so that the comparison argument applies to the entire square-adjacent
family.

\begin{lemma}[Comparison at the anchors]\label{lem:sa-cones}
Let $\mathscr F$ be a compact family of pairs $(F,t)$, where $F$ is a
square-adjacent witness at height $h$ with anchor $t$. Then
Lemma~\ref{lem:compact-negative-separation} applies to $\mathscr F$
with $H=h$, the lower digit bound $m$, and $\ell=W=2q-2$. Moreover
every anchor satisfies $t<\sigma_h=(h^2+4)^{-1/4}$, so
Lemma~\ref{lem:local-normalization} normalizes all witnesses near all
anchors.
\end{lemma}
\begin{proof}
Let $\beta=1/t\in[\beta_-,\beta_+]$. The first inequality
$\ell/(\beta-1)<2$ is $\beta>q$. The second, $(\beta-1)^2>2q-2$,
follows from $\beta>q$ and $(q-1)^2>2q-2$, which is $(q-1)(q-3)>0$.
For the third, write
$\Delta=2\beta-\ell-2h/(\beta-1)^2$ and observe that $\Delta>0$ is
equivalent to $(\beta-q+1)(\beta-1)^2>h$. Since
$(\beta-q)(\beta-1)^2=\varphi(\beta)/2\ge m/2$, it suffices that
$(\beta-1)^2>h-m/2=(q^2+2q-4+k)/2$. By Lemma~\ref{lem:sa-anchors},
$\beta-1>q-\tfrac34$, and
$(q-\tfrac34)^2\ge(q^2+2q-2)/2$ is equivalent to
$q^2-5q+\tfrac{25}8\ge0$, which holds for $q\ge5$. For $q=4$ one checks
directly that $\varphi(1+\sqrt{\smash[b]{h-m/2}})<m$ for $k=0,1,2$, which
places $1+\sqrt{\smash[b]{h-m/2}}$ below $\beta_-$: the three values are
$\varphi(1+\sqrt{10})=2(\sqrt{10}-3)\cdot10<8$,
$\varphi(1+\sqrt{21/2})=21(\sqrt{21/2}-3)<9$ and
$\varphi(1+\sqrt{11})=22(\sqrt{11}-3)<10$. All three inequalities hold
on the closed interval $[\beta_-,\beta_+]$, so the margins are uniform
over $\mathscr F$. For the normalization,
$\beta^4>(q+\tfrac14)^4>q^4+q^3>q^4+4\ge h^2+4$.
\end{proof}

By Lemma~\ref{lem:sa-cones} and the rigidity statement of
Lemma~\ref{lem:compact-negative-separation}, a square-adjacent witness is the
unique member of $\mathcal F_h$ with a multiple zero at its anchor. In
particular distinct square-adjacent witnesses have distinct anchors, and the map
from square-adjacent witnesses to $\mathcal T_h$ is a bijection. The set of square-adjacent
witnesses is closed in the coefficient topology, because a coefficient
limit of square-adjacent witnesses has a double zero at the limit of the anchors
by Lemma~\ref{lem:local-normalization}; hence $\mathcal T_h$ is compact.

\subsection{Witnesses from intermediate values}

Witnesses are produced by Lemma~\ref{lem:ivt-completion}. In this
section it is applied to prefixes $\pi=(1,-2q,a_2,\ldots,a_{N-1})$ with
digits $a_j\in[m,h]$ and to digit ranges $[m',H']\subseteq[m,h]$. The
completion $F$ it provides then belongs to $\Phi_h$, and the interior
minimum point is a double zero: $F$ is a square-adjacent witness with
anchor in $(a,b)$.

\begin{lemma}[Existence of square-adjacent witnesses]\label{lem:sa-existence}
The square-adjacent family $\Phi_h$ contains a square-adjacent witness. For $q\ge6$ the
family
\begin{equation}\label{eq:inner-family}
 \Psi_q=\Bigl\{\,1-2qz+\sum_{j\ge2}a_jz^j:\ q^2-2q+2\le a_j\le q^2-2\Bigr\}
 \subset\Phi_{q^2-2}\cap\Phi_{q^2-1}\cap\Phi_{q^2}
\end{equation}
contains a square-adjacent witness. The same holds for $q\in\{4,5\}$.
\end{lemma}
\begin{proof}
Apply Lemma~\ref{lem:ivt-completion} to the root prefix $\pi=(1,-2q)$,
the digit range $[m,h]$ and $J=[t_-,t_+]$. The smallest completion is
$M_m$, which vanishes at both endpoints and is negative inside. The
largest completion is $1-2qz+hz^2/(1-z)$, whose numerator
$1-(2q+1)z+(2q+h)z^2$ has discriminant $(2q-1)^2-4h<0$, so it is positive.
The contiguity condition $(2q-2)t_-\ge1-t_-$ is $t_-\ge1/(2q-1)$. Using
\eqref{eq:I0-endpoints}, it is equivalent to
\[
 2q^2-1-2k\ge(2q-1)\sqrt{4q+1-4k},
\]
and squaring gives
$4q^4-16q^3+(8+8k)q^2-16kq+4k^2+8k\ge0$. For $k=0$ the left side is
$4q^2(q^2-4q+2)$, for $k=1$ it is $4\bigl(q^2(q-2)^2-4q+3\bigr)$, and
for $k=2$ it is $4\bigl((q-1)^4-4q+7\bigr)$; each is positive for
$q\ge4$. This proves the first assertion.

For the second, put $m'=q^2-2q+2$, $H'=q^2-2$, and apply the lemma to
the root prefix with $J=[a,b]$, where $a=1/(2q-3)$ and $b=t'_+$ is the
larger zero of $M_{m'}$; the contiguity condition holds with equality
at $a$. The largest completion $1-2qz+(q^2-2)z^2/(1-z)$ is positive, as
its numerator has discriminant $9-4q<0$. The smallest completion is
$M_{m'}$, which vanishes at $b$ and is negative on $(t'_-,b)$. It is
nonnegative at $a$ precisely when $a\le t'_-$, that is,
\[
 2q^2-4q-7\ge(2q-3)\sqrt{4q-7}.
\]
For $q\ge7$ this follows from $\sqrt{4q-7}\le q-2$ and
$(2q-3)(q-2)\le2q^2-4q-7$, and for $q=6$ it reads $41\ge9\sqrt{17}$.
The interval $J$ contains points of $(t'_-,b)$, so
$\min_JM_{m'}<0$.

For $q\in\{4,5\}$ the root does not satisfy the endpoint condition,
and we apply the lemma to a deeper prefix. For $q=4$ take
$J=[1/5,1/3]$ and
\[
 \pi=(1,-8,12,13,14,12,13,10,13,13),
\]
and for $q=5$ take $J=[1/7,7/25]$ and
$\pi=(1,-10,20,22,18,21,22,21,17)$. In both cases the contiguity
condition holds with equality at the left endpoint, and the four
remaining conditions are inequalities between rational numbers, which
are verified in exact arithmetic (Proposition~\ref{cert:inner-prefix} of the supplement). For the lower bound on $f^+_\pi$ we use
convexity: rationals $z_1<z_2$ in $J$ with
$(f^+_\pi)'(z_1)\le0\le(f^+_\pi)'(z_2)$, $f^+_\pi(z_2)>0$ and
$f^+_\pi(z_1)+(f^+_\pi)'(z_1)(z_2-z_1)>0$ certify
$f^+_\pi>0$ on $J$.
\end{proof}

\begin{figure}[tbp]
\centering
\includegraphics[width=\textwidth]{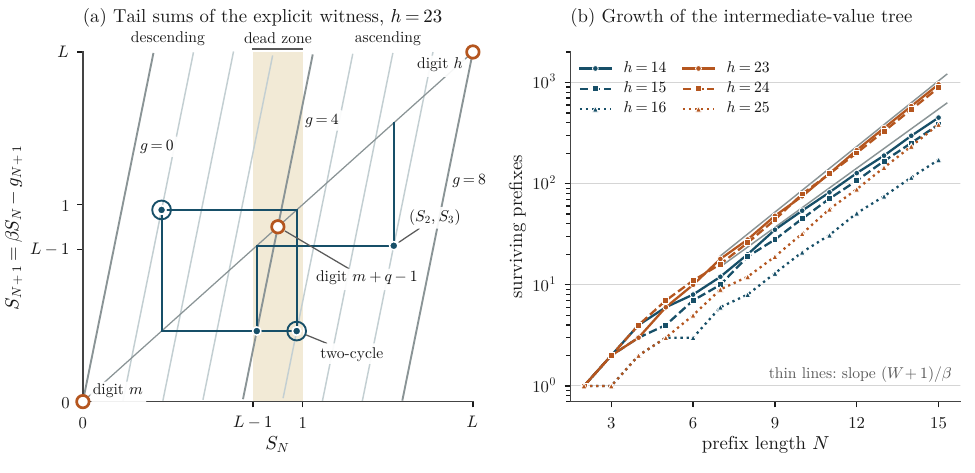}
\caption{(a) The tail-sum recursion~\eqref{eq:tail-recursion} at height
$23$ for the explicit witness of Proposition~\ref{prop:explicit-witnesses}:
the thin lines are the branches $S\mapsto\beta S-g$, $0\le g\le W$, with
$g=0$, $q-1$ and $W$ drawn darker; the shaded band is the dead zone
$[L-1,1)$; the hollow circles are the fixed points $0$, $(q-1)/(\beta-1)$
and $L$, which correspond to the exceptional tails with constant digit
$m$, $m+q-1$ and $h$. The cobweb follows the tail sums
$(S_N,S_{N+1})$ of the explicit witness from $(S_2,S_3)$ on: they end in
a two-cycle that alternates between the descending region and the dead
zone, so every other level is descending. (b) Numerical counts of the
surviving prefixes of the intermediate-value tree of the square-adjacent family,
by prefix length, at the six heights next to $16$ and $25$. The thin
lines illustrate growth by factors $(W+1)/4.47$ for $q=4$ and
$(W+1)/5.47$ for $q=5$ per additional digit, not a proved growth rate.
Only prefixes passing the exact tests of Section~\ref{app:certified} of the supplement
are used for the reported lower bounds.}
\label{fig:tail-dynamics}
\end{figure}

\subsection{Tail sums and the exceptional witnesses}

Let $F=1-2qz+\sum_{j\ge2}a_jz^j$ be a square-adjacent witness at height $h$
with anchor $t=1/\beta$. Write $g_j=a_j-m\in[0,W]$ for its
\emph{digit offsets}, and put
\begin{equation}\label{eq:tail-sums}
 L=\frac{Wt}{1-t}=\frac{W}{\beta-1},\qquad
 S_N=\sum_{i\ge1}g_{N+i}t^i\qquad(N\ge1).
\end{equation}
Then $0\le S_N\le L$ and
\begin{equation}\label{eq:tail-recursion}
 S_{N+1}=\beta S_N-g_{N+1}\qquad(N\ge1).
\end{equation}
By Lemma~\ref{lem:sa-anchors}, $1<L<2$. Call a level $N\ge2$
\emph{ascending} if $S_N\ge1$ and $g_N\le W-1$, and \emph{descending}
if $S_N<L-1$ and $g_N\ge1$; the \emph{margin} of a descending level is
$L-1-S_N$. The witness $F$ is \emph{live} if infinitely many of its
levels are ascending or descending, and \emph{exceptional} otherwise.

\begin{lemma}[Exceptional witnesses are eventually constant]\label{lem:exceptional}
An exceptional square-adjacent witness is eventually constant, with constant
digit $m$, $m+q-1$ or $h$.
\end{lemma}
\begin{proof}
Let $\beta$ be the reciprocal anchor. By Lemma~\ref{lem:sa-anchors},
$\beta<q+\delta^*$, that is $\beta^2-q\beta-(q-1)<0$. Expanding
$\beta(L-1)=\beta(2q-1-\beta)/(\beta-1)$ and
$(\beta-q+1)(\beta-1)=\beta^2-q\beta+q-1$, this single inequality is
equivalent to each of
\begin{equation}\label{eq:channel-inequalities}
 \beta(L-1)>q-1,\qquad \beta-(L-1)<q .
\end{equation}
Suppose no level $N\ge N_1$ is ascending or descending. For
$N\ge N_1$ this means: $S_N<1$ or $g_N=W$, and $S_N\ge L-1$ or $g_N=0$.

First let $N\ge N_1$ with $1\le g_N\le W-1$, so that $L-1\le S_N<1$.
Then \eqref{eq:tail-recursion} and $0\le S_{N+1}\le L$ give
\[
 g_{N+1}\ge\beta S_N-L\ge\beta(L-1)-L=W-\beta>q-3,\qquad
 g_{N+1}\le\beta S_N<\beta<q+1,
\]
so $q-2\le g_{N+1}\le q$, and in particular $1\le g_{N+1}\le W-1$.
Hence also $L-1\le S_{N+1}<1$, and then
\eqref{eq:tail-recursion} together with
\eqref{eq:channel-inequalities} gives
$q-2<\beta(L-1)-1<g_{N+1}<\beta-(L-1)<q$, so $g_{N+1}=q-1$. By
induction $g_j=q-1$ and $L-1\le S_j<1$ for all $j>N$. The affine
recursion $S_{j+1}-x^*=\beta(S_j-x^*)$ with $x^*=(q-1)/(\beta-1)$ then
forces $S_j=x^*$ for all $j\ge N$, and the digits are eventually
$m+q-1$.

Next let $g_N=0$, so $S_N<1$. Since $\beta<q+1\le W$, the next offset
satisfies $g_{N+1}\le\beta S_N<W$. If $g_{N+1}\ge1$, the first case
applies at level $N+1$. If $g_{N+1}=0$, then $S_{N+1}<1$ and
$S_N=S_{N+1}/\beta<\beta^{-1}$. Either some later offset is nonzero,
and the first nonzero offset after $N$ lies in $[1,W-1]$ by the same
bound, so the first case applies; or all later offsets vanish, which
means $S_N=0$ and the digits are eventually $m$.

Finally let $g_N=W$, so $S_N\ge L-1$. Then
$g_{N+1}\ge\beta S_N-L\ge W-\beta>q-3$, so $g_{N+1}\ge1$. If
$g_{N+1}\le W-1$, the first case applies. If $g_{N+1}=W$, then
$S_{N+1}=\beta S_N-W\ge L-1$. Either some later offset is smaller than
$W$, and the first such offset is at least $q-2$, so the first case
applies; or all later offsets equal $W$, and the recursion
$S_{j+1}-L=\beta(S_j-L)$ with $S_j\ge L-1$ forces $S_N=L$, so that the
digits are eventually $h$.
\end{proof}

The three exceptional types have different meanings. A witness with
constant digit $m$ or $h$ is the smallest or the largest completion of
its own prefix. A witness with constant digit $m+q-1$ is a polynomial
perturbation of the \emph{channel series}
\begin{equation}\label{eq:channel-series}
 1-2qz+(m+q-1)\frac{z^2}{1-z}
 =\frac{(1-qz)\bigl(1-(q+1)z\bigr)+(k-1)z^2}{1-z},
\end{equation}
and by~\eqref{eq:channel-inequalities} the offset $q-1$ is the only one
that can keep two successive tail sums in the \emph{dead zone}
$[L-1,1)$; the corresponding branch of~\eqref{eq:tail-recursion} has
the fixed point $(q-1)/(\beta-1)$.

\begin{lemma}[Margins]\label{lem:margins}
Put $y=W/(\beta_+-1)-1-(\beta_+-q)>0$ and
$\eta=\min\{\beta_--q,\ y/\beta_+\}$. A live square-adjacent witness has
infinitely many ascending levels, or infinitely many descending levels
of margin at least $\eta$.
\end{lemma}
\begin{proof}
Positivity of $y$ is the second inequality
in~\eqref{eq:channel-inequalities} at $\beta_+$. Let $F$ be live with
reciprocal anchor $\beta$ and write $\delta=\beta-q$ and
$y(\beta)=L-1-\delta$, which decreases in $\beta$, so $y(\beta)\ge y>0$;
also $y(\beta)<1$ since $L<2$. We show that a descending level $N$ of
margin $\varepsilon<\eta$ is followed by an ascending level or by a
descending level of margin larger than $\delta\ge\eta$; this proves the
lemma.

Write $S_N=L-1-\varepsilon$. Then
$\beta S_N=\beta(L-1)-\beta\varepsilon=(q-1)+y(\beta)-\beta\varepsilon$,
using $\beta(L-1)=W+L-\beta$. Since
$\beta\varepsilon<\beta_+\eta\le y\le y(\beta)$, the number
$\beta S_N$ lies in $(q-1,q)$ with fractional part
$y(\beta)-\beta\varepsilon$. By~\eqref{eq:tail-recursion} and
$0\le S_{N+1}\le L<2$, the offset $g_{N+1}$ lies in
$[\beta S_N-L,\beta S_N]$, so $g_{N+1}\in\{q-2,q-1\}$. If
$g_{N+1}=q-1$, then $S_{N+1}=y(\beta)-\beta\varepsilon<L-1-\delta$,
and level $N+1$ is descending with margin larger than $\delta$. If
$g_{N+1}=q-2$, then $S_{N+1}=1+y(\beta)-\beta\varepsilon>1$, and
level $N+1$ is ascending because $q-2\le W-1$.
\end{proof}

\subsection{From a live witness to infinitely many components}

\begin{theorem}[Live square-adjacent witnesses]\label{thm:live-witness}
Let $q\ge4$, let $h\in\{q^2-2,q^2-1,q^2\}$, and let $F_*$ be a live
square-adjacent witness at height $h$ with anchor $t_*$. Then $t_*\in\mathcal A_h$,
and infinitely many distinct connected components of $\mathcal C_h$ with
nonempty planar interior accumulate at each of the four points
$\pm t_*$ and $\pm1/t_*$, at which $\mathcal C_h$ is not locally
connected. If $F_*$ has infinitely many ascending levels, then these
four points are singleton components of $\mathcal C_h$.
\end{theorem}
\begin{proof}
Negation and inversion transfer the statements at $t_*$ to the other
three points, so we work at $t_*=1/\beta_*$. By
Lemma~\ref{lem:margins} there are two cases.

\emph{Ascending levels.}
Suppose infinitely many levels are ascending. For every ascending level
$N$ the series $F_*+z^N$ is admissible, because $g_N\le W-1$ means
$a_N(F_*)+1\le h$, and Lemma~\ref{lem:positive-unfolding} with $C_N=1$
gives nonreal zeros in $\mathcal C_h$ converging to $t_*$. Thus
$t_*\in\mathcal A_h$ and $t_*$ is nonisolated in $\mathcal C_h$.

For an ascending level $N$, let $\pi_N=(1,-2q,a_2,\ldots,a_{N-1})$ be
the prefix of $F_*$ and let $\pi'_N$ be the prefix $\pi_N$ followed by
the digit $a_N(F_*)+1\le h$. We apply Lemma~\ref{lem:ivt-completion}
with the digit range $[m,h]$ and $J=[t_-,t_+]$; the contiguity condition
was verified in Lemma~\ref{lem:sa-existence}. The smallest
completion of $\pi'_N$ is
\[
 f^-_{\pi'_N}=F_*+z^N-\sum_{j>N}(a_j(F_*)-m)z^j,\qquad
 f^-_{\pi'_N}(t_*)=t_*^N(1-S_N)\le0 .
\]
If $S_N>1$, its minimum on $J$ is negative. If $S_N=1$ and the minimum
were $0$, then $f^-_{\pi'_N}$ would be a square-adjacent witness with the anchor
$t_*$, contradicting the rigidity of $F_*$. The largest completion
$f^+_{\pi'_N}=F_*+z^N+\sum_{j>N}(h-a_j(F_*))z^j$ dominates $F_*\ge0$ on
$J$, and at the endpoints $f^-_{\pi'_N}$ exceeds $M_m$, which vanishes
there. The lemma gives a square-adjacent witness $F_N$ with prefix $\pi'_N$ and
anchor $t_N\in(t_-,t_+)$. The witnesses $F_N$ converge to $F_*$
coefficientwise, since they agree with $F_*$ below degree $N$, and
$t_N\to t_*$ by strict convexity.

Let $j_N$ be the least index above $N$ with $a_{j_N}(F_N)>m$; it exists
because the smallest completion of $\pi'_N$ has negative minimum. Let
$\mathscr B_N$ be the block of all negative departures from $F_N$ with
degrees in $[N,j_N-1]$. The witness $F_*$ belongs to $\mathscr B_N$:
its first departure from $F_N$ has degree $N$ and constant coefficient
$-1$. Lemma~\ref{lem:sa-cones} applies to the compact family
$\{(F_N,t_N)\}_N\cup\{(F_*,t_*)\}$ and yields one common separation
neighbourhood. Since $a_N(F_N)=a_N(F_*)+1>m$ and $a_{j_N}(F_N)>m$, the
last assertion of Lemma~\ref{lem:compact-negative-separation} gives the
hypotheses of Lemma~\ref{lem:finite-block-transfer} with $H=h_*=h$ and
$n_N=N\to\infty$. Every $K_N$ contains $t_*$, a qualifying zero of
$F_*\in\mathscr B_N$ lying in $\mathcal C_h$. Part~(ii) of that lemma
gives the singleton component $\{t_*\}$, the accumulating
interior-bearing components, and the failure of local connectedness.

\emph{Descending levels.}
Suppose instead that infinitely many levels are descending with margin
at least $\eta$. By Lemma~\ref{lem:sa-cones}, the three
inequalities~\eqref{eq:compact-comparison-conditions} hold at
$\beta_*$ with $H=h$ and $\ell=W$, and $t_*<\sigma_h$. A descending
level $N$ of margin at least $\eta$ has $a_N(F_*)>m$ and
\[
 \sum_{j\ge1}\bigl(h-a_{N+j}(F_*)\bigr)t_*^j=L-S_N\ge1+\eta .
\]
Proposition~\ref{prop:local-criterion} with $h_*=H=h$ gives the
accumulating interior-bearing components, the failure of local
connectedness, and $t_*\in\mathcal A_h$.
\end{proof}

\begin{lemma}[Exceptional witnesses are isolated]\label{lem:isolated-exceptional}
An exceptional square-adjacent witness is an isolated point of the set of square-adjacent
witnesses, and its anchor is an isolated point of $\mathcal T_h$.
Consequently $\mathcal T_h$ is infinite if and only if some square-adjacent
witness is live.
\end{lemma}
\begin{proof}
Let $F$ be exceptional with constant digit $c\in\{m,m+q-1,h\}$ from
degree $N_0$ on, and anchor $t=1/\beta$. Suppose square-adjacent witnesses
$G_n\ne F$ converge to $F$; their first departure degrees $M_n$ tend to
infinity, and their anchors $t_n$ tend to $t$. Write
$G_n=F+z^{M_n}C_n$ with $C_n(0)=c_n\in\mathbb Z\setminus\{0\}$; the
other coefficients of $C_n$ are differences of digits and lie in
$[-W,W]$, so $|C_n(z)-c_n|\le Wz/(1-z)$ for $0<z<1$. This bound
is $L<2$ at $t$ and tends to $L$ at $t_n$, hence is also $<2$ there for
large $n$. Since $G_n\ge0$ on $(0,1)$ and $F(t)=0$, we have $C_n(t)\ge0$,
hence $c_n\ge-1$; since $F\ge0$ and $G_n(t_n)=0$, we have
$C_n(t_n)\le0$, hence $c_n\le1$. Thus
$G_n$ has prefix $\pi$ followed by $c\pm1$, where $\pi$ is the prefix of
$F$ of length $M_n$, and $G_n$ is a completion of this prefix with
minimum $0$. We show that for large $M_n$ no such completion exists.

If $c=m+q-1$, the smallest completion of $(\pi,c+1)$ is
$F+z^{M_n}\bigl(1-(q-1)z/(1-z)\bigr)$. Its second factor is positive
for $z<1/q$, in particular near $t$, while away from $t$ the term
$F\ge0$ dominates the perturbation for large $M_n$; so this completion
is positive on $I_0$ and the prefix $(\pi,c+1)$ has no completion with
minimum $0$. The largest completion of $(\pi,c-1)$ is
$F-z^{M_n}\bigl(1-(q-1)z/(1-z)\bigr)$, which is negative at $t$. If
$c=m$, only the digit $c+1$ is available, and the smallest completion
$F+z^{M_n}$ is positive on $I_0$. If $c=h$, only $c-1$ is available, and
the largest completion $F-z^{M_n}$ is negative at $t$. In all cases
$G_n$ cannot exist for large $n$.

For the last assertion, let $F_*$ be live. Every ascending level $N$
produces a witness $F_N\ne F_*$ with prefix $\pi'_N$, as in the proof of
Theorem~\ref{thm:live-witness}. Every sufficiently late descending level
$N$ of margin at least $\eta$ produces a witness with prefix
$(\pi_N,a_N(F_*)-1)$ by Lemma~\ref{lem:ivt-completion}: the largest
completion of this prefix has value $t_*^N(-1+L-S_N)\ge\eta t_*^N$ at
$t_*$. Indeed, the largest completion is $F_*+z^N C_N(z)$, where
$C_N(t_*)\ge\eta$ and the coefficients of $C_N$ are uniformly bounded.
The derivative bound on a compact interval about $t_*$ gives
$C_N\ge\eta/2$ on a fixed smaller neighbourhood. Off that neighbourhood
in $[t_-,t_+]$, $F_*$ has a positive minimum and $z^N C_N$ tends
uniformly to zero. Thus the largest completion is positive throughout
$[t_-,t_+]$ for large $N$; the smallest is negative at $t_*$.
By Lemma~\ref{lem:margins}, infinitely many levels are of one of these
kinds, and distinct levels give distinct witnesses. Conversely, if every
square-adjacent witness is exceptional, the compact set of square-adjacent witnesses
consists of isolated points and is therefore finite.
\end{proof}

Of the three exceptional types, only the middle one can make the anchor
set finite: a witness with constant tail $m$ or $h$ always has a
neighbour in the tree of prefixes.

\begin{lemma}[Extreme tails]\label{lem:extreme-tails}
Let $F$ be a square-adjacent witness that is eventually constant with
digit $m$, respectively $h$, and let $N$ be the last level at which
its digit differs from that constant. Then there is a square-adjacent witness
$G$ that agrees with $F$ below degree $N$, has
$a_N(G)=a_N(F)-1$, respectively $a_N(G)=a_N(F)+1$, and whose digits
beyond level $N$ are not all equal.
\end{lemma}
\begin{proof}
The level $N\ge2$ exists, because the two constant series are $M_m$,
which is negative on $I_0$, and the largest member of $\Phi_h$, which
is positive. Let $t=1/\beta$ be the anchor of $F$. We apply
Lemma~\ref{lem:ivt-completion} with the digit range $[m,h]$ and
$J=[t_-,t_+]$, where $Wz\ge1-z$, as in Lemma~\ref{lem:sa-existence}.

Suppose first that the constant is $m$, and let $\pi$ be the prefix of
$F$ below degree $N$ followed by the digit $a_N(F)-1\ge m$. Then
\[
 f^-_\pi=F-z^N,\qquad
 f^+_\pi=F+z^N\Bigl(\frac{Wz}{1-z}-1\Bigr)\ge F\ge0\quad\text{on }J,
\]
while $f^-_\pi(t)=-t^N<0$ and $f^-_\pi\ge M_m$, which vanishes at
$t_\pm$. The lemma gives a witness $G$ with prefix $\pi$; let
$t'=1/\beta'$ be its anchor. If $a_j(G)=c$ for all $j>N$, then
\[
 G=F+z^N\Bigl(\frac{(c-m)z}{1-z}-1\Bigr).
\]
By Lemma~\ref{lem:sa-anchors}, $q-1<\beta-1<q$ and
$q-1<\beta'-1<q$. If $c-m\le q-1$, then
$G(t)=t^N\bigl((c-m)/(\beta-1)-1\bigr)<0$, which is impossible for a
witness. If $c-m\ge q$, then
$G(t')\ge t'^N\bigl((c-m)/(\beta'-1)-1\bigr)>0$, contradicting
$G(t')=0$.

If the constant is $h$, let $\pi$ be the prefix of $F$ below degree
$N$ followed by the digit $a_N(F)+1\le h$. Now $f^+_\pi=F+z^N>0$, and
$f^-_\pi=F+z^N(1-Wz/(1-z))$ has $f^-_\pi(t)=t^N(1-L)<0$ and
$f^-_\pi\ge M_m$. The lemma gives a witness $G$ with prefix $\pi$ and
anchor $t'=1/\beta'$. If $a_j(G)=c$ for all $j>N$, then
$G=F+z^N\bigl(1-(h-c)z/(1-z)\bigr)$, and either $h-c\ge q$ and
$G(t)<0$, or $h-c\le q-1$ and $G(t')>0$; both are impossible.
\end{proof}

\begin{corollary}[A sufficient criterion for infinitude]\label{cor:sa-criterion}
Let $q\ge4$ and $h\in\{q^2-2,q^2-1,q^2\}$. If the anchor set
$\mathcal T_h$ is infinite, or equivalently if some square-adjacent witness is not
eventually constant with digit $m$, $m+q-1$ or $h$, then $\mathcal C_h$
has infinitely many connected components with nonempty planar interior,
accumulating at a point of $(1/(q+1),1/q)$ at which the core is not
locally connected. In particular this holds if no square-adjacent witness
at height $h$ is eventually constant with digit $m+q-1$.
\end{corollary}
\begin{proof}
For the first assertion combine Lemmas~\ref{lem:exceptional}
and~\ref{lem:isolated-exceptional} with Theorem~\ref{thm:live-witness}.
For the second, suppose that no witness is eventually constant with
digit $m+q-1$, and let $F_0$ be a witness
(Lemma~\ref{lem:sa-existence}). If $F_i$ is exceptional, it is
eventually constant with digit $m$ or $h$ by
Lemma~\ref{lem:exceptional}, and Lemma~\ref{lem:extreme-tails} gives a
witness $F_{i+1}$ that is live, or exceptional with a larger last
nonconstant level than $F_i$. Thus some $F_i$ is live, or the $F_i$ are
pairwise distinct; in both cases $\mathcal T_h$ is infinite.
\end{proof}

\subsection{Transfer between neighbouring heights}

Shifting the height by one shifts the three exceptional digits by one,
so an exceptional witness at one height is usually a live witness at a
neighbouring height.

\begin{lemma}[Transfer]\label{lem:transfer}
Let $q\ge4$ and let $F$ be a square-adjacent witness at height $h=q^2-2+k$.
\begin{enumerate}
\item If $k\le1$, every tail digit $a_j$ ($j\ge2$) of $F$ is at least $m+1$, and $F$ is
exceptional at height $h$ with constant digit $m+q-1$ or $h$, then $F$
is a live square-adjacent witness at height $h+1$.
\item If $k\ge1$, every tail digit $a_j$ ($j\ge2$) of $F$ is at most $h-1$, and $F$ is
exceptional at height $h$ with constant digit $m$ or $m+q-1$, then $F$
is a live square-adjacent witness at height $h-1$.
\end{enumerate}
\end{lemma}
\begin{proof}
In case~(i), the digits of $F$ lie in $[m+1,h+1]$, which is the digit
range of the square-adjacent family at height $h+1$, so $F$ is a square-adjacent witness
there. Its constant digit is $m+q-1=(m+1)+q-2$ or $h=(h+1)-1$, neither
of which is one of the three exceptional digits $m+1$, $(m+1)+q-1$ and
$h+1$ of height $h+1$. By Lemma~\ref{lem:exceptional}, applied at
height $h+1$, $F$ is live there. Case~(ii) is symmetric: the digits lie
in $[m-1,h-1]$, and the constant digit $m=(m-1)+1$ or
$m+q-1=(m-1)+q$ is not among $m-1$, $(m-1)+q-1$ and $h-1$.
\end{proof}

\begin{theorem}[Two of three heights]\label{thm:two-of-three}
Let $q\ge4$. At least two of the three cores
$\mathcal C_{q^2-2}$, $\mathcal C_{q^2-1}$ and $\mathcal C_{q^2}$ have
infinitely many connected components, including infinitely many with
nonempty planar interior, accumulating at a common point
$t_*\in(1/(q+1),1/q)$ and at $-t_*$, $\pm1/t_*$, where these cores are
not locally connected. More precisely, if $F$ is a square-adjacent witness with
all digits in $[q^2-2q+2,q^2-2]$, then this holds at every height
$h\in\{q^2-2,q^2-1,q^2\}$ at which $F$ is not eventually constant with
digit $m_h$, $m_h+q-1$ or $h$, where $m_h=h-2q+2$, and these three
exceptional digits exclude at most one of the three heights.
\end{theorem}
\begin{proof}
Let $F$ be a square-adjacent witness in the family $\Psi_q$
of~\eqref{eq:inner-family}; it exists by Lemma~\ref{lem:sa-existence}.
Since $[q^2-2q+2,q^2-2]$ is contained in the digit range of the square-adjacent
family at each of the three heights, $F$ is a square-adjacent witness at each of
them. Within $[q^2-2q+2,q^2-2]$, the exceptional digits are
\[
 \{q^2-q-1,\ q^2-2\},\qquad \{q^2-q\},\qquad \{q^2-2q+2,\ q^2-q+1\}
\]
at the heights $q^2-2$, $q^2-1$ and $q^2$ respectively. These sets are
pairwise disjoint. If $F$ is eventually constant, its constant digit
belongs to at most one of them, and if it is not eventually constant it
belongs to none. By Lemma~\ref{lem:exceptional}, $F$ is live at every
height not excluded in this way, and Theorem~\ref{thm:live-witness}
gives the conclusions there, with the common anchor $t_*$ of $F$.
\end{proof}

\subsection{Pisot anchors and explicit witnesses}

The witnesses of Lemma~\ref{lem:sa-existence} need not be explicit. At
an anchor whose reciprocal is a Pisot number, the witness, if it exists,
is a rational function, and a finite computation decides whether it
exists and is live. This is the counterpart of the theorem of Bertrand
and Schmidt that, for a Pisot number $\beta$, the $\beta$-expansion of
every element of $\mathbb Q(\beta)\cap[0,1)$ is eventually
periodic~\cite{Bertrand77,Schmidt80}; here the value and the derivative
equations are imposed together, and rigidity replaces the greedy
algorithm.

\begin{proposition}[Witnesses at Pisot anchors]\label{prop:pisot-anchors}
Let $\beta\in[\beta_-,\beta_+]$ be a Pisot number, an algebraic integer
all of whose other conjugates lie in the open unit disk, and let $F$ be
a square-adjacent witness with anchor $1/\beta$. Then the digit sequence of $F$ is
eventually periodic, so $F$ is a rational function, and $F$ is live
unless its digits are eventually constant and equal to $m$, $m+q-1$ or
$h$. Moreover the digit sequence is the unique infinite path from the
initial state in a finite graph determined by $\beta$, $m$ and $h$; in
particular, whether a square-adjacent witness with anchor $1/\beta$ exists, and
whether it is live, is decided by a finite computation.
\end{proposition}
\begin{proof}
Write $t=1/\beta$ and $a_0=1$, $a_1=-2q$, and consider the partial sums
\[
 P_N=\sum_{j\le N}a_j\beta^{N-j},\qquad
 V_N=\sum_{j\le N}(N+1-j)a_j\beta^{N-j}\qquad(N\ge1),
\]
which satisfy $P_{N+1}=\beta P_N+a_{N+1}$ and $V_{N+1}=\beta V_N+P_{N+1}$,
and lie in $\mathbb Z[\beta]$. The equations $F(t)=0$ and $tF'(t)=0$ give
\begin{equation}\label{eq:pisot-tails}
 P_N=-\sum_{i\ge1}a_{N+i}t^i,\qquad V_N=\sum_{i\ge1}(i-1)a_{N+i}t^i,
\end{equation}
so $|P_N|\le h/(\beta-1)$ and $|V_N|\le h/(\beta-1)^2$. For an
embedding $\sigma$ of $\mathbb Q(\beta)$ with $|\sigma(\beta)|<1$, the
defining sums give $|\sigma(P_N)|\le h/(1-|\sigma(\beta)|)$ and
$|\sigma(V_N)|\le h/(1-|\sigma(\beta)|)^2$. Thus $P_N$ and $V_N$ are
algebraic integers of degree at most $[\mathbb Q(\beta):\mathbb Q]$ with
all conjugates bounded independently of $N$; the coefficients of their
minimal polynomials are bounded integers, so the pairs $(P_N,V_N)$ range
over a finite set, and some state repeats: $(P_M,V_M)=(P_N,V_N)$ with
$N<M$.

By~\eqref{eq:pisot-tails}, the state at a level determines the two tail
sums $\sum_{i\ge1}a_{N+i}t^i$ and $\sum_{i\ge1}(i-1)a_{N+i}t^i$, and
conversely the value and the derivative of a series at $t$ depend on its
digits beyond level $M$ only through these two sums. Hence the series $G$
with digits $a_1,\ldots,a_M$ followed by $a_{N+1},a_{N+2},\ldots$
satisfies $G(t)=F(t)=0$ and $G'(t)=F'(t)=0$. It belongs to $\Phi_h$, so
it is a square-adjacent witness with anchor $t$, and rigidity
(Lemma~\ref{lem:sa-cones}) gives $G=F$: the digits of $F$ are
periodic from level $N+1$ on, with period $M-N$. The claim about liveness
is Lemma~\ref{lem:exceptional}. For the last assertion, let the states
be the pairs $(P,V)$ in the finite set just described, with an edge
labelled $a\in[m,h]$ from $(P,V)$ to $(\beta P+a,\beta V+\beta P+a)$
whenever the target satisfies the bounds stated after~\eqref{eq:pisot-tails}.
The digit sequence of a witness is an infinite path from
$(P_1,V_1)=(\beta-2q,2\beta-2q)$. Conversely, along any infinite path
the states are bounded, so the series with those digits satisfies
$\sum_ja_jt^j=\lim_Nt^NP_N=0$ and
$\sum_jja_jt^j=\lim_N\bigl((N+1)\sum_{j\le N}a_jt^j-t^NV_N\bigr)=0$;
it has a double zero at $t$, so it is a square-adjacent witness with anchor $t$,
and rigidity makes it equal to $F$. Hence the path is unique, and a
witness with anchor $t$ exists if and only if some infinite path leaves
the initial state, which is decided by the finite graph.
\end{proof}

The searches in the companion~\cite{Companion26} enumerate Pisot
numbers of bounded degree in $[\beta_-,\beta_+]$ with numerical root
filters and follow the paths of the finite graph in floating-point
arithmetic. The witnesses of Table~\ref{tab:explicit-witnesses} are
then checked directly in exact arithmetic
(Proposition~\ref{cert:explicit-witnesses} of the
supplement), independently of the search.

\begin{proposition}[Explicit eventually periodic witnesses]\label{prop:explicit-witnesses}
For each row of Table~\ref{tab:explicit-witnesses}, the series
\[
 F(z)=1-2qz+\sum_{j=2}^{N-1}a_jz^j+\frac{z^N\sum_{i=0}^{r-1}c_iz^i}{1-z^r}
\]
is a live square-adjacent witness at each of the listed heights, with anchor
$t=1/\beta$, where $\beta$ is the Pisot number with the stated minimal
polynomial. Consequently the conclusions of Theorem~\ref{thm:live-witness}
hold at the heights
\[
 14,\ 15,\ 16,\ 23,\ 24,\ 34,\ 35,\ 36,\ 47,\ 48,\ 62,\ 80,\ 98,\ 119,\ 120,\ 194 .
\]
\end{proposition}

\begin{table}[htbp]
\centering\small
\caption{Explicit eventually periodic square-adjacent witnesses. The digits
$a_2,\ldots,a_{N-1}$ are followed by the period $c_0,\ldots,c_{r-1}$
repeated indefinitely. The reciprocal anchor $\beta$ is the largest
real root of the stated polynomial.}
\label{tab:explicit-witnesses}
\begin{tabular}{@{}cclllc@{}}
\toprule
$q$ & heights & digits $a_2,\ldots,a_{N-1}$ & period & minimal polynomial of $\beta$ & $\beta$\\
\midrule
4 & $14,15,16$ & $12,13,14,12$ & $13,10,13$ & $x^4-4x^3-2x^2-2x+1$ & $4.52842$\\
4 & $14,15$ & $12,14,9,14,13,14,11,10$ & $14,14,11,10$ & $x^5-4x^4-2x^3-x^2-2x+1$ & $4.51177$\\
5 & $23,24$ & $20,22,19$ & $16,20$ & $x^3-5x^2-3x+1$ & $5.51140$\\
5 & $24$ & $21,16,21,20,23$ & $20,21,16,24$ & $x^5-5x^4-2x^3-2x^2-2x+1$ & $5.44593$\\
6 & $34,35,36$ & $31,26,26,32,33$ & $30,34$ & $x^4-6x^3-3x^2+x-1$ & $6.44513$\\
7 & $47,48$ & $42,44,44$ & $36,45$ & $x^3-7x^2-4x+1$ & $7.51459$\\
8 & $62$ & $56,58,59,60,57$ & $60,56,48,61$ & $x^5-8x^4-4x^3-3x^2-3x+2$ & $8.51558$\\
8 & $62$ & $57,50,55,60,48,50,51$ & $48,52$ & $x^5-8x^4-4x^3+x^2-x+1$ & $8.46028$\\
9 & $80$ & $72,74,79$ & $64,80$ & $x^3-9x^2-5x+1$ & $9.51447$\\
10 & $98$ & $91,82,94,94,85,82,88$ & $80,89$ & $x^5-10x^4-5x^3+x^2-x+1$ & $10.46925$\\
11 & $119,120$ & $110,114,100,112,116$ & $108,117$ & $x^4-11x^3-6x^2+2x-1$ & $11.50697$\\
14 & $194$ & $183,172,174,172,192$ & $168,193$ & $x^4-14x^3-7x^2+2x-1$ & $14.47440$\\
\bottomrule
\end{tabular}
\end{table}

\begin{proof}
Write $P$ for the polynomial part and $T$ for the period polynomial,
so that $(1-z^r)F=R$ with $R=(1-z^r)P+z^NT\in\mathbb Z[z]$. A double
zero of $F$ in $(0,1)$ is the same as a double zero of $R$ there. For
each row, the greatest common divisor $g$ of $R$ and $R'$ is the
reversed form of the stated polynomial, it is irreducible, $g^2$ divides
$R$, and $g$ has exactly one zero in $(1/(q+1),1/q)$; these are finite
exact computations (Proposition~\ref{cert:explicit-witnesses} of the supplement), which also certify that $\beta$ is
a Pisot number, although the argument does not use this. All digits lie
in the digit range
of each listed height, and the period is not constant, so the witness
is not eventually constant. Lemma~\ref{lem:exceptional} shows that it is
live at each listed height, and Theorem~\ref{thm:live-witness} applies.
\end{proof}

The table also exhibits explicit irrational singleton components. For
the row $q=6$, with $\beta\approx6.44513$ and $m_h=h-10$, every late
digit $30$ is followed by $34$, so at that level $a_N=30\le h-1$ and
$S_N\ge(34-m_h)/\beta\ge8/\beta>1$: infinitely many levels are
ascending at each of the heights $34$, $35$, $36$, and by
Theorem~\ref{thm:live-witness} the four points $\pm\beta^{\pm1}$ are
singleton components of $\mathcal C_{34}$, $\mathcal C_{35}$ and
$\mathcal C_{36}$. Each row lists exactly the heights whose digit
range contains all of its digits.

\subsection{Witnesses of the form \texorpdfstring{$g(z)^2/(1-z^r)$}{g(z)\textasciicircum2/(1-z\textasciicircum r)}}\label{sec:heights-seven-eight}

Every witness of Table~\ref{tab:explicit-witnesses} is a quotient
$g(z)^2/(1-z^r)$ with $g\in\mathbb Z[z]$: in the proof of
Proposition~\ref{prop:explicit-witnesses}, $R=g^2$. For such a quotient
the double zero is visible, and the verification is a polynomial
identity: if $g^2=\sum_jc_jz^j$, the coefficients of the quotient are
$a_j=c_j+a_{j-r}$, with $a_j=c_j$ for $j<r$, and they repeat with period
$r$ beyond the degree of $g^2$. The next proposition lists thirteen
further series of this form. They come from an exhaustive search over
the polynomials $g$ of bounded degree, in integer arithmetic: the
condition $m\le a_j\le h$ confines each coefficient of $g$ in turn to at
most $q$ consecutive integers. The same search recovers every row of
Table~\ref{tab:explicit-witnesses}.

\begin{proposition}[Explicit witnesses that are squares]\label{prop:square-witnesses}
For each row of the following list, let $g(z)=z^dp(1/z)$, where $d$ is
the degree of $p$. The series $F=g^2/(1-z^r)$ equals
$1-2qz+\sum_{j\ge2}a_jz^j$, where the integers $a_j$ lie in the stated
interval and are eventually periodic with the stated period, and it has
a double zero at $1/\beta$, where $\beta$ is the only zero of $p$ in
$(1,\infty)$ and $q<\beta<q+1$.
\begin{center}\upshape\footnotesize\setlength{\tabcolsep}{4.5pt}
\begin{tabular}{@{}cclclll@{}}
\toprule
$q$ & heights & $p(x)$ & $r$ & digits in & period & $\beta$\\
\midrule
3 & $7$ & $x^6-3x^5-x^4-x^3-x^2-x+1$ & 5 & $[3,7]$ & $5,5,6,3,6$ & $3.40974$\\
3 & $8$ & $x^3-3x^2-2x+1$ & 2 & $[4,8]$ & $4,5$ & $3.49086$\\
7 & $49$ & $x^7-7x^6-3x^5-3x^4+x^3-x+1$ & 3 & $[37,49]$ & $39,39,43$ & $7.45408$\\
8 & $64$ & $x^7-8x^6-4x^5-3x^4+x^3-x+1$ & 3 & $[50,64]$ & $57,52,60$ & $8.50987$\\
9 & $79$ & $x^7-9x^6-4x^5-4x^4+2x^3-x+1$ & 3 & $[63,78]$ & $64,64,68$ & $9.46492$\\
10 & $98,99$ & $x^7-10x^6-5x^5-4x^4+2x^3-x+1$ & 3 & $[81,95]$ & $86,81,89$ & $10.51022$\\
10 & $100$ & $x^8-10x^7-5x^6+x^5-x^4+x^3-1$ & 2 & $[82,100]$ & $96,100$ & $10.46925$\\
11 & $121$ & $x^7-11x^6-5x^5-5x^4+2x^3-x+1$ & 3 & $[101,121]$ & $105,105,114$ & $11.47249$\\
12 & $142,143,144$ & $x^9-12x^8-6x^7-4x^6-4x^5+3x^4-x^2-x+1$ & 4 & $[124,140]$ & $136,132,124,137$ & $12.50722$\\
14 & $195$ & $x^7-14x^6-7x^5-6x^4+3x^3-x+1$ & 3 & $[169,195]$ & $176,169,184$ & $14.50995$\\
15 & $223,224$ & $x^7-15x^6-7x^5-6x^4+3x^3-x^2-x+1$ & 3 & $[196,221]$ & $212,196,217$ & $15.47656$\\
15 & $225$ & $x^6-15x^5-8x^4+2x^3-x^2+x-1$ & 2 & $[200,225]$ & $216,225$ & $15.50780$\\
16 & $254,255$ & $x^6-16x^5-8x^4+2x^3-x^2+x-1$ & 2 & $[225,250]$ & $234,250$ & $16.47833$\\
\bottomrule
\end{tabular}
\end{center}
In the rows with $q\ge7$, $F$ is a live square-adjacent witness at each
listed height. Consequently the conclusions of
Theorem~\ref{thm:live-witness} hold at the heights $49$, $64$, $79$,
$98$, $99$, $100$, $121$, $142$, $143$, $144$, $195$, $223$, $224$,
$225$, $254$ and $255$.
\end{proposition}
\begin{proof}
The digits follow from the recursion for the $a_j$. In each row
$p(q)<0<p(q+1)$, so $g$ has a zero $1/\beta$ with $q<\beta<q+1$, which
is a double zero of $F$. Since the $a_j$ are positive, $F''>0$ on
$(0,1)$, so this is the only zero of $F$, hence of $g$, in $(0,1)$. In
the rows with $q\ge7$ the stated interval is contained in $[h-2q+2,h]$
for each listed height $h$, so $F$ is a square-adjacent witness there.
Its period is not constant, so it is live by
Lemma~\ref{lem:exceptional}, and Theorem~\ref{thm:live-witness} applies.
\end{proof}

The series with $q=12$ is a witness common to its three heights, as in
Theorem~\ref{thm:two-of-three}.

For $q=3$ the inequalities~\eqref{eq:compact-comparison-conditions} hold
only on part of the interval of possible anchors, and the results of
this section do not apply to the whole family $\Phi_h$.
Proposition~\ref{prop:local-criterion} needs them at a single anchor,
and the first two series of Proposition~\ref{prop:square-witnesses}
satisfy all of its hypotheses.

\begin{proposition}[Heights seven and eight]\label{prop:heights-seven-eight}
For $h\in\{7,8\}$ let $\Gamma_h$ be the series of
Proposition~\ref{prop:square-witnesses} at height $h$,
\[
 \Gamma_7(z)=\frac{(1-3z-z^2-z^3-z^4-z^5+z^6)^2}{1-z^5},\qquad
 \Gamma_8(z)=\frac{(1-3z-2z^2+z^3)^2}{1-z^2},
\]
and $1/\beta_h$ its double zero, so that $\beta_7=3.40974\ldots$ and
$\beta_8=3.49086\ldots$. Then $\Gamma_h\in\mathcal F_h$, and infinitely many
connected components of $\mathcal C_h$ with nonempty planar interior
accumulate at each of the four points $\pm\beta_h^{\pm1}$, at which
$\mathcal C_h$ is not locally connected.
\end{proposition}
\begin{proof}
Expanding the squares and dividing gives
\begin{align*}
 \Gamma_7(z)&=1-6z+7z^2+4z^3+5z^4+7z^5+5z^6+5z^7
 +\frac{z^8(5+5z+6z^2+3z^3+6z^4)}{1-z^5},\\
 \Gamma_8(z)&=1-6z+6z^2+8z^3+4z^4+\frac{z^5(4+5z)}{1-z^2}.
\end{align*}
Thus $\Gamma_7\in\mathcal F_7$ with all coefficients of degree at least two
in $[3,7]$, and $\Gamma_8\in\mathcal F_8$ with all such coefficients in
$[4,8]$; in both cases $\ell=4$. The two polynomials $p$ change sign on
$(3.409,3.41)$ and on $(3.49,3.491)$ respectively, so
\[
 3.409<\beta_7<3.41,\qquad 3.49<\beta_8<3.491 .
\]
For $\beta$ in either interval and the corresponding height,
\[
 \frac{\ell}{\beta-1}<1.67,\qquad(\beta-1)^2>5.8>\ell,\qquad
 2\beta-\ell-\frac{2h}{(\beta-1)^2}>0.39,
\]
so the inequalities~\eqref{eq:compact-comparison-conditions} hold, and
$t_h=1/\beta_h<0.294<\sigma_h$, because $0.34^4<1/68$ and $0.37^4<1/53$.

It remains to verify~\eqref{eq:descending-condition}. For $\Gamma_7$ take
the indices $N\ge10$ divisible by $5$: then $a_N=6>3$, the next two
coefficients are $3$ and $6$, and the sum
in~\eqref{eq:descending-condition} is at least $4t_7+t_7^2>5/4$. For
$\Gamma_8$ take the even indices $N\ge6$: then $a_N=5>4$, the next two
coefficients are $4$ and $5$, and the sum is at least
$4t_8+3t_8^2>5/4$. Proposition~\ref{prop:local-criterion} applies with
$h_*=H=h$ and $\eta=1/4$.
\end{proof}

Some of these series belong to families. The series $\Gamma_8$ is the case
$q=3$ of $\bigl(1-qz-\tfrac{q+1}2z^2+z^3\bigr)^2/(1-z^2)$, $q$ odd,
whose coefficients lie in the digit range of a square-adjacent height
exactly when $3\le q\le9$; for $q=5,7,9$ it gives the first row with
$q=5$ and the rows with $q=7$ and $q=9$ of
Table~\ref{tab:explicit-witnesses}. Likewise
$\bigl(1-qz-\tfrac q2z^2+3z^3-2z^4+z^5-z^6+z^7\bigr)^2/(1-z^2)$ is a
live witness common to the three heights next to $q^2$ for
$q=22,24,26$, by the same verification. Such families serve only
finitely many $q$. Indeed $r\ge2$, because $a_1=c_1+1$ is odd when
$r=1$; and if $g=qu+v$, with fixed polynomials $u$, $v$ and a fixed
$r\ge2$, gave square-adjacent witnesses $g^2/(1-z^r)$ for infinitely many
$q$, comparing the coefficients of $q^2$ would give
$u(z)^2=z^2(1+z+\cdots+z^{r-1})$, which is impossible because the zeros
of $1+z+\cdots+z^{r-1}$ are simple.

\section{Real double zeros near the minimum}\label{sec:real-traces}

This section proves Theorem~\ref{thm:real-traces}. The real gaps it
produces lie inside a window of width at most $10^{-87}$ just right of
the minimum of $\mathcal D_h$. Near a point $v$ of $\mathcal D_h$
that is the smallest double zero of series with a prescribed initial
block $P$, the other double zeros of such series form \emph{groups},
one for each later position at which a corner digit can be lowered,
and each group is a scaled slice of a self-affine set
(Theorem~\ref{rt:thm:structure}). At certain digit patterns the group
that starts the next block lies alone in a window of relative width
$2\mu$; then $\mathcal D_h$ has a gap just left of the smallest double
zero $v'$ of the lowered block, and $v'$ has the same structure as $v$
(Corollary~\ref{rt:cor:move}). If the tail of $v$ has infinitely many
such patterns, the gaps accumulate at $v$. Otherwise we move to $v'$
from the last one. The digits that follow the last pattern lie in an
explicit set of small dimension, and a finite computation shows that
the tail of $v'$ again contains a pattern (Proposition~\ref{rt:prop:landing}).
The process therefore never stops. The local analysis holds at every
height; the heights $1,3,4,7$ enter only through finitely many
inequalities, certified in Section~\ref{app:certified} of the supplement
(Propositions~\ref{cert:rt-minimum}--\ref{cert:rt-landing}).

\subsection{Corner digits, slices and greedy tails}\label{rt:sec:digits}
Fix $h\ge1$. For $f=1+\sum_{j\ge1}a_jz^j\in\Fh$ write $b_j=h-a_j\in\{0,\ldots,2h\}$
for its \emph{corner digits} and $p_0(t)=1+ht/(1-t)$. Then
\begin{equation}\label{rt:eq:corner}
 f(t)=p_0(t)-\sum_{j\ge1}b_jt^j,\qquad
 p_k(t)=\sum_{j>k}b_jt^{j-k},\qquad q_k(t)=\sum_{j>k}(j-k)b_jt^{j-k}
\end{equation}
define the \emph{remainders} and \emph{moments} of $f$ after $k$ digits.
For $t\in(0,1)$ and $x\ge0$ let
\[
 S_t(x)=\Bigl\{\sum_{r\ge1}rc_rt^r:\ c_r\in\{0,\ldots,2h\},\ \sum_{r\ge1}c_rt^r=x\Bigr\}
\]
be the set of moments of the digit sequences of weight $x$; it is
compact, and $q_k(t)\in S_t(p_k(t))$. Splitting off the first digit
gives
\begin{equation}\label{rt:eq:slice-recursion}
 S_t(x)=x+t\bigcup_{c=0}^{2h}S_t(x/t-c),
\end{equation}
so $S_t(x)$ is the vertical slice over $x$ of the attractor of the
maps $(x,y)\mapsto(t(c+x),t(c+x+y))$, $0\le c\le2h$, and
$S_t(x)\subset[0,2ht/(1-t)^2]$. Put $\beta=1/t$ and $\bar d=\lfloor\beta\rfloor$.
The \emph{greedy expansion} of $x\in[0,2)$ is defined by $x_0=x$,
$g_k=\lfloor\beta x_{k-1}\rfloor$, $x_k=\beta x_{k-1}-g_k$. Its
remainders satisfy $0\le x_k<1$, so $g_k\le\bar d$ for $k\ge2$, while
$g_1\le\lfloor2\beta\rfloor\le2h$ at every height considered here; we
write $q^{\mathrm g}_t(x)=\sum_rrg_rt^r$ for its moment and
$d_\beta(1)$ for the greedy expansion of $1$. A finite word over
$\{0,\ldots,\bar d\}$ is \emph{weakly admissible} if each of its
suffixes is lexicographically at most the prefix of $d_\beta(1)$ of
the same length; every window of a greedy expansion of a number in
$[0,1)$ is weakly admissible~\cite{Parry60}.

\subsection{The minimum}\label{rt:sec:minimum}
For $0<t\le1/(h+1)$ every $f\in\Fh$ satisfies $f(t)\ge1-ht/(1-t)\ge0$, with
equality only for $t=1/(h+1)$ and $f=1-hz/(1-z)$, whose derivative does not
vanish there; hence $\mathcal D_h\subset(1/(h+1),1)$, and $\mathcal D_h$ is
relatively closed by coefficient compactness. For
$h\in\{1,3,4,7\}$ it is nonempty, and Proposition~\ref{cert:rt-minimum}
gives rationals $a_h<b_h$, with $b_h-a_h=2\cdot10^{-150}$ for $h=1$ and
$2\cdot10^{-90}$ otherwise, such that
\[
 t_0:=\min\mathcal D_h\in[a_h,b_h],
\]
together with digit words $P_0$ of lengths $n_0=49,24,19,15$ and longer
words of lengths $N_h=850,257,216,162$ with the following properties:
every series of $\Fh$ with a double zero in $W_h=[a_h-10^{-12},b_h+10^{-9}]$
begins with $P_0$, and every series with a double zero in $[a_h,b_h]$
begins with the long word. Table~\ref{rt:tab:data} lists the data.

For a word $P\supseteq P_0$ of length $n$ let $\Fh(P)$ be the series
of $\Fh$ that begin with $P$, let $\mathcal D_h(P)$ be the set of
double zeros in $W_h$ of members of $\Fh(P)$, and, when it is nonempty,
$m(P)=\min\mathcal D_h(P)$, a \emph{relative tip}. Thus $t_0=m(P_0)$.
Fix $v=m(P)$ with $|v-t_0|\le10^{-80}$ and a series $f_v\in\Fh(P)$ with a
double zero at $v$, write $b_j,p_k,q_k$ for its corner digits,
remainders and moments at $t=v$, and call $g_s=b_{n+s}$,
$y_s=p_{n+s}$ and $q_s$ (abusing notation) the digits, remainders and
moments of its \emph{tail} at the \emph{levels} $s\ge1$. Put
$\kappa_v=f_v''(v)/2$ and $u_P=v^{n-1}/(2\kappa_v)$. The \emph{free
radius} $F(P)$ is the supremum of the $F\ge0$ such that
$(v,v+Fu_P)\subset W_h$ and every series of $\Fh$ with a double zero in
$(v,v+Fu_P)$ belongs to $\Fh(P)$. The window property gives
$F(P_0)\ge10^{-9}/u_{P_0}$.

\subsection{Local structure at a relative tip}\label{rt:sec:local}
Let $R=10^{-6}$, $\tau=b_h+R$, $C_3=h/(1-\tau)^4$,
$G_1(i)=2h\sum_{k\ge i}k\tau^{k-1}$, $G_2(i)=2h\sum_{k\ge i}k(k-1)\tau^{k-2}$
and $\bar M=2h\tau/(1-\tau)^2$. Then $|f'''|/6\le C_3$ on
$[v-R/2,v+R/2]$ for every $f\in\Fh$ (this interval lies in $[0,\tau)$, since
$v\le b_h+10^{-80}$); a power series $\sum_{k\ge i}c_kz^k$
with $|c_k|\le2h$ has first and second derivatives bounded by $G_1(i)$
and $G_2(i)$ on $|z|\le\tau$; every moment lies in $[0,\bar M]$. The
series $f_v$ shares $P_0$ with a double-zero series at $t_0$, so
$|\kappa_v-\kappa_0|\le G_2(n_0+1)/2+10^{-70}$, where $\kappa_0$ is
enclosed from the certified digits. With $W(i)=G_1(i)/(2\kappa_v-3C_3R-G_2(i))$,
$\zeta_i=(4C_3W(i)^3+\tfrac32G_2(i)W(i)^2)/v^i$ and
\[
 \bar\eta_i=\frac{(1+10^{-3})\,\bar M^2v^{i-2}}{4\kappa_v},\qquad
 L^*(i)=i-4-c_L,
\]
where $c_L$ is the least integer with $\beta^{c_L}\ge\bar M^2/(4\kappa_v)$
(so $\bar\eta_i\le1.001\,v^{L^*(i)+2}$), Proposition~\ref{cert:rt-constants}
certifies for all $i>n_0$ the inequalities
\begin{gather*}
\text{(C4) } 2\kappa_v-3C_3R-G_2(n_0+1)>0,\quad
\text{(C1) } W(i)\le R,\quad
\text{(C2) } \zeta_i<10^{-3},\\
\text{(C5) } \bar\eta_i\le10^{-7},\quad
\text{(D) } \frac{3C_3W(i)+G_2(i)}{2\kappa_v-3C_3W(i)-G_2(i)}\le10^{-5},
\end{gather*}
and the further conditions (C3), (C6), (C7) and (E1)--(E7) quoted where
they are used. If $f\in\Fh$ has corner digits $b'_j$, and $i$ is the first
index with $b'_i\ne b_i$, then, with $\varepsilon=f(v)$ and
$\varepsilon'=f'(v)$,
\begin{equation}\label{rt:eq:first-difference}
 \varepsilon=v^i\bigl[(b_i-b'_i)+(p_i-p'_i)\bigr],\qquad
 \varepsilon'=v^{i-1}\bigl[i\varepsilon/v^i+q_i-q'_i\bigr],
\end{equation}
since $f_v(v)=f_v'(v)=0$.

\begin{lemma}[Taylor expansion at a tip]\label{rt:lem:taylor}
Let $f$ be as above with $i>n_0$, and suppose $f$ has a double zero at
$v+w$ with $|w|\le R/2$. Then $|w|\le W(i)$,
$|\varepsilon-\kappa_vw^2|\le4C_3|w|^3+\tfrac32G_2(i)w^2$, and
\[
 w=\frac{v^{i-1}\bigl(q'_i-q_i-i\varepsilon/v^i\bigr)}{2\kappa_v+\vartheta},
 \qquad |\vartheta|\le3C_3|w|+G_2(i).
\]
Moreover $0\le\varepsilon/v^i\le\bar\eta_i$, with equality on the left only for $w=0$.
\end{lemma}
\begin{proof}
Write $f=f_v+g$ with $g=\sum_{k\ge i}c_kz^k$, $|c_k|\le2h$. Taylor's
formula for $f_v$ at $v$ and for $g$ at $v$ gives
$f(v+w)=\varepsilon+\varepsilon'w+\kappa_vw^2+E$ and
$f'(v+w)=\varepsilon'+2\kappa_vw+E'$ with $|E|\le C_3|w|^3+\frac12G_2(i)w^2$
and $|E'|\le3C_3w^2+G_2(i)|w|$. At a double zero,
$\varepsilon'=-(2\kappa_v+\vartheta)w$ with $\vartheta=E'/w$ and
$\varepsilon=\kappa_vw^2+E'w-E$, which gives the second bound. Since
$\varepsilon'=g'(v)$ is at most $G_1(i)$ in modulus and
$2\kappa_v+\vartheta\ge2\kappa_v-3C_3R-G_2(i)>0$, we get $|w|\le W(i)$.
Inserting \eqref{rt:eq:first-difference} in $\varepsilon'=-(2\kappa_v+\vartheta)w$ gives the
formula for $w$. By (C1) and (C2), $\varepsilon\ge
w^2(\kappa_v-4C_3W(i)-\frac32G_2(i))\ge0$, with equality only for $w=0$.
For the upper bound, first $\varepsilon/v^i\le(\kappa_v+4C_3W(i)+\frac32G_2(i))W(i)^2/v^i$;
inserting this and $|q'_i-q_i|\le\bar M$ in the formula for $w$ and then in
$\varepsilon\le(\kappa_v+4C_3|w|+\frac32G_2(i))w^2$ gives
$\varepsilon/v^i\le\bar\eta_i$; the product of the correction factors
is at most $1+10^{-3}$ by (C6).
\end{proof}

The next lemma is the Shmerkin--Solomyak criterion
\cite[Lemma~3.1 and Corollary~3.2]{SS06}, written for height $h$.

\begin{lemma}[Good prefixes]\label{rt:lem:good-prefix}
Let $Q=1+\sum_{j\le N}c_jz^j$ with $c_j\in[-h,h]\cap\Z$, put
$H(y)=hy^{N+1}/(1-y)$, and let $1/(2h+1)<a<b<1$ satisfy $Q(a)>H(a)$,
$Q(b)>H(b)$, $Q+H>0$ on $[a,b]$ and $Q(x_*)<H(x_*)$ for some
$x_*\in(a,b)$. Then a series of $\Fh$ with initial part $Q$ has a double
zero in $(a,b)$.
\end{lemma}
\begin{proof}
Choose $d_{N+1},d_{N+2},\ldots\in[-h,h]$ successively, each as small as
possible such that $g_k=Q+\sum_{N<j\le k}d_jz^j+h\sum_{j>k}z^j$ stays
positive on $[a,b]$; the choice $d_k=h$ is allowed because it gives
$g_k=g_{k-1}$, and $g_N=Q+H$. The $g_k$ decrease to
$f=Q+\sum_{j>N}d_jz^j\in\Fh$, so $f\ge0$ on $[a,b]$, and
$f(a)\ge Q(a)-H(a)>0$, $f(b)>0$. Suppose that $m=\min_{[a,b]}f>0$. As
$g_k\ge f$, lowering $d_k$ by one keeps $g_k\ge m-b^k>0$ once $b^k<m$,
so $d_k=-h$ for all large $k$; not all $d_j$ equal $-h$, since otherwise
$f(x_*)=Q(x_*)-H(x_*)<0$. Let $k_0$ be the last index with
$d_{k_0}>-h$. Lowering $d_{k_0}$ was not allowed, so $g_{k_0}(x)\le
x^{k_0}$ for some $x\in[a,b]$, and then $f(x)=g_{k_0}(x)-2h\sum_{j>k_0}x^j
\le x^{k_0}(1-2hx/(1-x))<0$, because $x>1/(2h+1)$. This contradiction
shows $m=0$; the minimum is attained in $(a,b)$, where $f'$ vanishes too.
\end{proof}

\begin{lemma}[Realization]\label{rt:lem:realization}
Let $i>n$, let $b'_i\ne b_i$ be a corner digit and $c=b_i-b'_i$, and let
$b''_1,\ldots,b''_L$, with $0\le L\le L^*(i)+1$, be corner digits whose
remaining weight $r_L=(c+p_i-\sum_{r\le L}b''_rv^r)/v^L$ lies in
$[v+\bar\eta_i/v^L,1.7]$. Then some $f\in\Fh(P)$ that agrees with $f_v$
before position $i$ and has the corner digits $b'_i$ at $i$ and $b''_r$
at $i+r$ has a double zero at $v+w$ with
\begin{equation}\label{rt:eq:realized}
 w=\frac{v^{i-1}\bigl(q''-q_i+\varrho\bigr)}{2\kappa_v}(1+\theta'),\qquad
 q''=\sum_{r\le L}rb''_rv^r,
\end{equation}
where $|\theta'|\le10^{-5}$ and
$|\varrho|\le2hv^{L+1}\bigl((L+5)/(1-v)+1/(1-v)^2\bigr)+i\bar\eta_i$.
\end{lemma}
\begin{proof}
Let $N=i+L+5$ and let $Q$ run over the polynomials with constant term $1$,
the coefficients of $f_v$ below $i$, the corner digits $b'_i$ and
$b''_r$ at $i$ and $i+r$, free coefficients at the five positions
$i+L+1,\ldots,N$, and no further terms. By \eqref{rt:eq:first-difference}
for $Q$,
\[
 Q(v)=v^{i+L}\Bigl(r_L-\sum_{j=1}^5b^Q_{i+L+j}v^j\Bigr)-H(v),\qquad
 Q'(v)=v^{i-1}\bigl(iQ(v)/v^i+q_i-q^Q_i\bigr),
\]
where $b^Q$ are the corner digits of $Q$ and $H$ is as in
Lemma~\ref{rt:lem:good-prefix}. On $[v-R/2,v+R/2]$ one has $|f_v''-2\kappa_v|\le6C_3R$ and
$|(Q-f_v)''|\le G_2(i)$, so $Q''\ge2\kappa_v(1-10^{-4})$ by (E5); hence $Q$ has a
unique minimum point $x_Q$ there, with $|x_Q-v|\le\bar Mv^{i-1}/\kappa_v$ and
$Q(v)-\min Q\le Q'(v)^2/(4\kappa_v(1-10^{-4}))\le(1+10^{-4})v^i\bar\eta_i$, by (C6).
Write $M(c)=\min Q$ for the choice $c$ of the free coefficients. By (E3),
$x^k/v^k\in[1-2\cdot10^{-3},1+2\cdot10^{-3}]$ for $k\le N+1$ and $x$ between $v$ and
any of the minimum points, so raising the coefficient at a free position $k$
by one raises $M$ by an amount in $v^k[1-2\cdot10^{-3},1+2\cdot10^{-3}]$; and,
because $2hv>1.01$ (E7), raising it by one while lowering all later free
coefficients from $h$ to $-h$ does not raise $M$. As in the proof of
Lemma~\ref{lem:ivt-completion}, choosing the free coefficients from left to
right therefore reaches every interval of length $v^N(1+2\cdot10^{-3})$ inside
$[M(-h,\ldots,-h),M(h,\ldots,h)]$. Here $M(-h,\ldots,-h)\le Q(v)<0$, since
$2h\sum_{j\le5}v^j>1.7$ (E6), and $M(h,\ldots,h)\ge v^{i+L+1}(1-3\cdot10^{-4})-H(v)\ge
0.9H(v)+v^N(1+2\cdot10^{-3})$, since $r_L\ge v+\bar\eta_i/v^L$,
$v^i\bar\eta_i\le1.001v^{i+L+1}$ and (E4): $v^4(1.9hv/(1-v)+1)\le0.965$. By (E2$'$),
$0.8H(v)>1.01v^N$, so some choice gives $\min Q\in[0.1H(v),0.9H(v)]$.
For it, put $\rho=c_\rho v^{(N+1)/2}$, with $c_\rho=\frac12,\frac34,\frac45,\frac45$
for $h=1,3,4,7$. By (E3), $(N+1)(|x_Q-v|+\rho)/v\le10^{-3}$, so $H$ varies
by a factor at most $1.002$ on $[x_Q-\rho,x_Q+\rho]$. Then $Q>0$ there,
$Q(x_Q)\le0.9H(v)<H(x_Q)$, and $Q(x_Q\pm\rho)\ge\min Q+\kappa_v\rho^2(1-10^{-4})
>H(x_Q\pm\rho)$ by (E1): $\kappa_vc_\rho^2>1.01h/(1-\tau)$. Lemma~\ref{rt:lem:good-prefix}
gives $f$ with initial part $Q$ and a double zero in $(x_Q-\rho,x_Q+\rho)$.
Its first difference from $f_v$ is at $i>n$, so $f\in\Fh(P)$, and
Lemma~\ref{rt:lem:taylor} gives \eqref{rt:eq:realized}: $q'_i$ exceeds $q''$ by
the moment of the corner digits after $i+L$, at most
$2h\sum_{r>L}rv^r$, the term $i\varepsilon/v^i$ is at most $i\bar\eta_i$,
and $\theta'$ is controlled by (D).
\end{proof}

The realization lemma is used with long prescriptions only at depths
$i\ge159$, where $L^*(i)\ge152$. Put $\delta_k=2v^{L^*(k)}$. The next lemma
uses the slice bounds of Proposition~\ref{cert:rt-constants}:
\begin{gather*}
 \Delta_{\mathrm{low}}=\inf\bigl\{\min S_t(c+p)-q^{\mathrm g}_t(p):\ 1\le c\le\bar d,\ p\in[0,\bar p)\cup[1,1+10^{-6}]\bigr\},\\
 \Delta_1=\inf\bigl\{\min S_t(c+p)-q^{\mathrm g}_t(p):\ 1\le c\le\lfloor2\beta\rfloor,\ p\in[0,1+10^{-6}]\bigr\},\\
 q_1^*=\sup\{q^{\mathrm g}_t(x):\ |x-1|\le10^{-6}\},
\end{gather*}
over all $t$ within $10^{-15}$ of the centre of $[a_h,b_h]$, where $\bar p=1$ for
$h\ge3$ and $\bar p=0.496>\beta-1$ for $h=1$, and the integers $L_{\min}$,
$L_G$ and $L_2$: the least $L$ with
$2hv^{L+1}((L+5)/(1-v)+1/(1-v)^2)\le\Delta_{\mathrm{low}}/2-10^{-6}$, the least $L$ with
$v^L(q_1^*+10^{-6})\le\Delta_{\mathrm{low}}-0.0015$, and the least $L$ with
$2hv^{L+1}((L+5)/(1-v)+1/(1-v)^2)\le10^{-3}$. It also uses two finite
checks: (GR) for every weakly admissible word $F=F_1\cdots F_{L_1}$ over
$\{0,\ldots,\bar d\}$ with $1\le L_1<L_G$, $F_{L_1}\ge1$ and $F(v)=\sum_rF_rv^r<1+10^{-6}$,
and every $1\le c\le\lfloor2\beta\rfloor$,
\[
 \min S_t\bigl(c+F(v)\bigr)-q(F)-v^{L_1}(q_1^*+10^{-5})>0.002,\qquad q(F)=\sum_rrF_rv^r;
\]
and (EX) the greedy expansion of every $x\in[1,1+\xi_0]$, $\xi_0=2v^{L^*(159)}$, is
the same for all $t\in[t_0,t_0+10^{-80}]$ ($10^{-140}$ for $h=1$) through $M_E\ge144$
digits, with positive remainders, and has no run of ten zeros among its first
$\max(L_{\min},L_2)+21$ digits.

\begin{lemma}[Greedy tails]\label{rt:lem:greedy}
Let $v=m(P)$ with $p_n<2$, and suppose that $n\ge158$, or that $P=P_0$
and $p_k<1$ for $n_0<k\le158$. Then $p_k<1+\delta_k$ for every $k>n$.
Consequently the tail of $f_v$ is the greedy expansion of $p_n$, except
at \emph{restarts}: indices $k$ with $p_k\in[1,1+\delta_k)$, where the
corner digit is one less than the greedy one and the tail continues as
the greedy expansion of $p_k$. Restarts occur only at indices $k\ge159$.
\end{lemma}
\begin{proof}
Let $k>n$ be the first index with $p_k\ge1+\delta_k$; by hypothesis $k\ge159$.
Put $c=\lfloor p_k-\delta_k\rfloor\ge1$ and $p=p_k-c\in[\delta_k,1+\delta_k)$.
Then $b_k+c\le2h$: for $k\ge n+2$ because $p_{k-1}<1+\delta_{k-1}$ gives
$b_k,c\le\bar d\le h$, and for $k=n+1$ because $b_k+c\le\beta p_n<2\beta\le2h+1$.
We raise $b_k$ by $c$ and use Lemma~\ref{rt:lem:realization} with digits $b''$
after $k$ to be chosen; its double zero lies left of $v$, contradicting
$v=\min\mathcal D_h(P)$, as soon as $q''+\varrho<q_k$, and
$q_k\ge\min S_v(c+p)$. Put $\Delta=\Delta_{\mathrm{low}}$ if $c\le\bar d$
and $\Delta=\Delta_1$ otherwise (then $k=n+1$); both infima include the
strip $p\in[1,1+10^{-6}]$, and $\Delta>\Delta_{\mathrm{low}}/2$. Let $g$ be the
greedy expansion of $p$, with remainders $r_L$.

(1) If $r_L\ge v+\bar\eta_k/v^L$ for some $L\in[L_{\min},L^*(k)]$, take
$b''=g_1\cdots g_L$. Then $q''\le q^{\mathrm g}_v(p)$, $|\varrho|<\Delta_{\mathrm{low}}/2$,
and $\min S_v(c+p)-q^{\mathrm g}_v(p)\ge\Delta$.

(2) Otherwise, writing $r_L=v(g_{L+1}+r_{L+1})$ shows that among the
positions $L_{\min}+1,\ldots,L^*(k)+1$ the expansion $g$ has at most one
nonzero digit, a $1$ at some position $L_1$, after which the remainder is
below $\bar\eta_k/v^{L_1}$. (2a) If $L_1\le L^*(k)-1$ (for $h\ge3$:
$L_1\le L^*(k)+1$), prescribe $g_1\cdots g_{L_1-1}0$. The remaining weight
lies in $[1,1+\bar\eta_k/v^{L_1}]\subset[v+\bar\eta_k/v^{L_1},1.7]$, because
$1.001v^3\le1-v$ for $h=1$ and $1.001v\le1-v$ for $h\ge3$; we conclude as in (1).
(2b) Otherwise $p=F(v)+\xi'$ with $\xi'<2v^{L^*(k)}$, where $F=g_1\cdots g_{L_1}$ is
empty or has its last nonzero digit at $L_1\le L_{\min}$. If $F$ is empty
and $c=1$, then $p_k<1+\delta_k$, contrary to the choice of $k$. If $F$ is
empty and $c\ge2$, use $c-1$ and the weight $1+p$, prescribing the greedy
expansion of $1+p$ up to an $L\in[L_{\min},L_{\min}+9]$ followed by a nonzero
digit, which exists by (EX), and which leaves a remaining weight in
$[v+\bar\eta_k/v^L,1.7]$; then $q''\le q^{\mathrm g}_v(1+p)$, and
$\min S_v(c+p)-q^{\mathrm g}_v(1+p)\ge\Delta$ by the strip part of the bound. If $F$ is
nonempty, prescribe $F$ with its last digit lowered by one, followed by
the greedy expansion of $1+\xi'/v^{L_1}$ up to an $L$ with $L-L_1\in[L_2,L_2+9]$
followed by a nonzero digit (again by (EX)); its weight is $p$, the room
condition holds, $|\varrho|\le0.0011$, and
$q''\le q(F)+v^{L_1}(q_1^*+10^{-6})$. If $L_1\ge L_G$, then
$\min S_v(c+p)-q(F)\ge\Delta\ge v^{L_1}(q_1^*+10^{-6})+0.0015$, by the definition
of $L_G$ and by (K4) of Proposition~\ref{cert:rt-assembly}; if $L_1<L_G$, (GR)
applies. In every case $q''+\varrho<q_k$.

If $p_k\in[1,1+\delta_k)$ the greedy digit at $k$ is $b_k+1$, and the lemma
applied from $k$ on describes the rest of the tail.
\end{proof}

The tips of the construction satisfy the hypothesis of
Lemma~\ref{rt:lem:greedy}: $t_0$ by Proposition~\ref{cert:rt-assembly},
and the later ones because $|P|\ge158$.

\begin{lemma}[Restarts]\label{rt:lem:restarts}
Let $v$ satisfy the hypothesis of Lemma~\ref{rt:lem:greedy}.
\begin{enumerate}
\item After a restart at $k$ the next $M_E$ digits of the tail are those of
$d_\beta(1)$, and two restarts are more than $M_E$ positions apart.
\item For every level $s'$ let $r$ be the distance from $s'$ to the first
restart after it ($r=\infty$ if there is none). Then
$q^{\mathrm g}_v(y_{s'})\le q_{s'}\le q^{\mathrm g}_v(y_{s'})+v^r(q_1^*+10^{-5})$.
\item For every level $s'$ with $g_{s'}\ge1$, every $1\le c\le g_{s'}$ and
$|e|\le10^{-4}$, $\min S_v(c+y_{s'}+e)-q_{s'}$ is at least
$\Delta-v^r(q_1^*+10^{-5})$, where $\Delta=\Delta_{\mathrm{low}}$ for $s'\ge2$ and $\Delta=\Delta_1$
for $s'=1$, and it always exceeds $10^{-3}$.
\end{enumerate}
\end{lemma}
\begin{proof}
At a restart $k$ the tail continues as the greedy expansion of
$p_k\in[1,1+\xi_0]$, which by (EX) agrees for $M_E$ digits with that of $1$
and keeps its remainders in $[0,1)$ there, so no restart occurs in
between. For (ii), let $k=n+s'+r$ be the first restart after $s'$. The
greedy expansion of $y_{s'}$ is $F=g_{s'+1}\cdots g_{s'+r-1}(b_k+1)$ followed
by the greedy expansion of $\xi=p_k-1<\delta_k$, so
\[
 q_{s'}-q^{\mathrm g}_v(y_{s'})=v^r\bigl(q_k-q^{\mathrm g}_v(\xi)\bigr).
\]
Here $q_k\ge\bar dv\ge v$, since the expansion of $p_k\ge1$ starts with the
digit $\bar d$, while $q^{\mathrm g}_v(\xi)$ is of the order of $\xi$; and, by induction
over the later restarts, $q_k\le q_1^*+10^{-6}+v^{M_E}(q_1^*+10^{-5})$. For
(iii), $y_{s'}\in[0,\bar p)\cup[1,1+10^{-6}]$ (for $h=1$ because
$g_{s'}\ge1$ forces $y_{s'}<\beta-1$), so (ii) gives the first bound,
the table lookups of Proposition~\ref{cert:rt-constants} absorbing
$|e|\le10^{-4}$. If $r\ge L_G$ it exceeds $10^{-3}$ by
Proposition~\ref{cert:rt-assembly}~(K4). If $r<L_G$, the word $F$ of the proof of (ii)
is a greedy prefix of $y_{s'}$, $q_{s'}\le q(F)+v^r(q_1^*+10^{-5})$, and (GR)
gives a margin above $0.002$.
\end{proof}

\begin{lemma}[Only lowered digits lie to the right]\label{rt:lem:right}
Let $f\in\Fh$ have a double zero at $v+w$ with $0<w<F(P)u_P$, and let
$i$ be the first index where its corner digits differ from those of $f_v$.
Then $f\in\Fh(P)$, $i>n$, $b'_i=b_i-c$ with $1\le c\le b_i$, and
$p'_i=c+p_i-\eta$ with $0<\eta=\varepsilon/v^i\le\bar\eta_i$.
\end{lemma}
\begin{proof}
The first assertion is the definition of $F(P)$, so $i>n$. Suppose
$b'_i=b_i+c$ with $c\ge1$. By \eqref{rt:eq:first-difference} and
Lemma~\ref{rt:lem:taylor}, $-\zeta_iv^i\le\varepsilon\le v^i(p_i-1)$; hence
$p_i\ge1-\zeta_i$, while $p'_i\le p_i-1+\zeta_i<\delta_i+\zeta_i=:x_i$ by
Lemma~\ref{rt:lem:greedy}. So the corner digits $b'_{i+r}$ vanish for $r<r_0$,
the least $r$ with $v^r<x_{n_0+1}$, and $q'_i\le2hv^{r_0}(r_0/(1-v)+1/(1-v)^2)$.
Then $\varepsilon'/v^{i-1}\ge q_i-q'_i-i\zeta_i>0$ by (C3), since $q_i\ge p_i$;
but $\varepsilon'$ has the sign of $-w$. Hence $b'_i<b_i$, and the bounds on
$\eta$ are those of Lemma~\ref{rt:lem:taylor}.
\end{proof}

\begin{theorem}[Structure at a relative tip]\label{rt:thm:structure}
Let $x\in(v,v+F(P)u_P)$ be a double zero of $f\in\Fh$. With $i$ and $c$ as
in Lemma~\ref{rt:lem:right}, put $s=i-n$. Then $g_s\ge c\ge1$ and
\[
 x-v=u_Pv^s\,(q'-q_s-\iota)(1+\theta),\qquad q'\in S_v(c+y_s-\eta),
\]
with $0<\eta\le\bar\eta_i\le10^{-7}$, $0\le\iota\le i\bar\eta_i\le10^{-5}$ and
$|\theta|\le10^{-5}$. We call the set of these values, for fixed $s$ and $c$,
the \emph{group} $G_{s,c}$; its points are the double zeros in
$(v,v+F(P)u_P)$ of the series $f\in\Fh(P)$ with $b'_{n+s}=b_{n+s}-c$ that agree
with $f_v$ before $n+s$.
\end{theorem}
\begin{proof}
Combine Lemmas~\ref{rt:lem:right} and~\ref{rt:lem:taylor}: $q'=q'_i\in S_v(p'_i)$,
$\iota=i\varepsilon/v^i$, and $(1+\theta)^{-1}=1+\vartheta/(2\kappa_v)$ is
controlled by (D).
\end{proof}

\begin{figure}[tb]
\centering
\includegraphics[width=\textwidth]{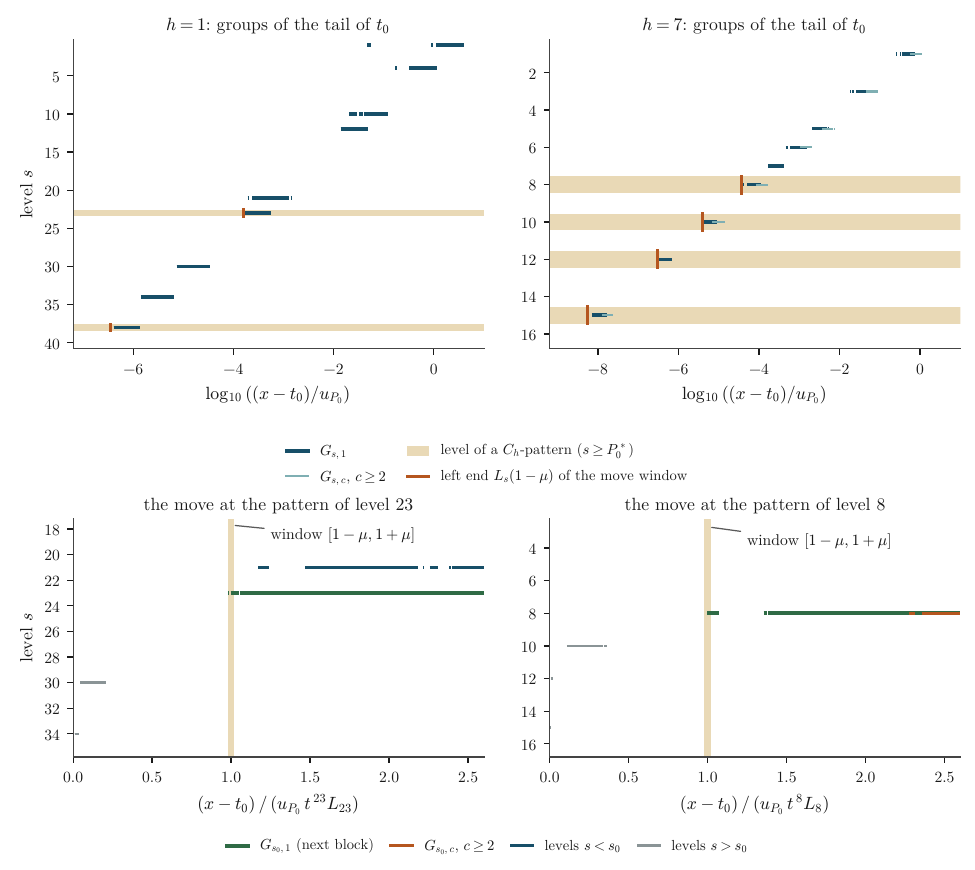}
\caption{Groups at the minimum $t_0$ of $\mathcal D_h$ for $h=1$ (left) and
$h=7$ (right). Top: row $s$ shows, on a logarithmic scale in units of
$u_{P_0}$, outer enclosures of the groups $G_{s,c}$ of the certified tail of
$t_0$, computed from the exact slice tables of
Proposition~\ref{cert:rt-constants}; shaded rows are the $C_h$-patterns at
levels $s\ge P_0^*$, and the tick marks the left end $v^sL_s(1-\mu)$ of the
move window. Bottom: the first such pattern, in units of $u_{P_0}v^sL_s$.
The window $[1-\mu,1+\mu]$ meets no group except $G_{s,1}$, whose minimum
$m(P'_s)$ lies at $1$; the interval from $1-\mu$ to $m(P'_s)$ is a gap of
$\mathcal D_h$. The enclosures contain the true groups; the figure
illustrates Propositions~\ref{cert:rt-move} and~\ref{cert:rt-landing} and
certifies nothing beyond them.}
\label{fig:real-trace-groups}
\end{figure}

\subsection{Moves}\label{rt:sec:moves}
Figure~\ref{fig:real-trace-groups} shows the groups of $t_0$ at $h=1$ and $h=7$. Put $\mu=1/50$. At a level $s$ with $g_s\ge1$ let $P'=Pg_1\cdots g_{s-1}(g_s-1)$,
of length $n+s$, and
\[
 \mathcal L_s=\bigl\{q^{\mathrm g}_t(1+y_s+e)-q_s:\ |e|\le10^{-7}\bigr\},
\]
the range of the first-group position at level $s$, taken over $t$ within
$10^{-15}$ of the centre of $[a_h,b_h]$ and over all greedy branches.

\begin{corollary}[Moves]\label{rt:cor:move}
Let $s\ge2$ with $g_s\ge1$, $y_s<1$ and $n+s\ge158$, let $L_-\le L_+$ bound
$\mathcal L_s$, and assume:
\begin{enumerate}[label=\textup{(\alph*)}]
\item every group of a level $s'>s$ lies below $v^s\bigl(L_-(1-\mu)-5\cdot10^{-4}L_-\bigr)$;
\item no group $G_{s',c}$ with $s'\le s$ and $(s',c)\ne(s,1)$ meets
 $v^s\bigl[L_-(1-\mu)-5\cdot10^{-4}L_-,\,L_+(1+\mu)+5\cdot10^{-4}L_+\bigr]$;
\item $v^sL_+(1+\mu)<F(P)$.
\end{enumerate}
Then $\mathcal D_h(P')\ne\varnothing$, $m(P')-v\ge u_Pv^sL_-(1-10^{-4})$, and the
interval $\bigl(v+u_Pv^sL_-(1-\mu),m(P')\bigr)$ contains no point of
$\mathcal D_h$. If moreover
\begin{enumerate}[label=\textup{(\alph*)},start=4]
\item the tail of a double-zero series of $P'$ at $m(P')$ has no restart among
 its first $R$ digits, where $v^{R+1}(q_1^*+10^{-5})\le7\cdot10^{-5}L_-$,
\end{enumerate}
then $m(P')-v\le u_Pv^sL_+(1+10^{-4})$ and $F(P')\ge(\mu-10^{-4})L_-$.
\end{corollary}
\begin{proof}
Let $x=1+y_s$. Apply Lemma~\ref{rt:lem:realization} with $c=1$ and the
prescription of the proof of Lemma~\ref{rt:lem:greedy} for the weight $x$, with
$L_2$ in place of $L_{\min}$; in case (2b) $F$ is nonempty because $x\ge1$, and
we prescribe $F$ followed by zeros if $\xi'\ge2\bar\eta_i$. Then $|\varrho|\le0.0011$,
and $q''$ exceeds by less than $10^{-40}$ the moment of the greedy expansion of
a weight in $[x-10^{-7},x+10^{-7}]$: in cases (1) and (2a), and for $F$ followed
by zeros, that of $x$ itself; otherwise that of a weight just below $F(v)$,
whose expansion is $F$ lowered followed by the digits of $d_\beta(1)$ for $M_E$
positions, by (EX). So a member of $\Fh(P')$ has a double zero at distance at
most $u_Pv^s(L_++0.0012)$ from $v$, inside $(v,v+F(P)u_P)$ by (c). Hence
$m(P')$ lies in this interval, and by Theorem~\ref{rt:thm:structure} the points
of $\mathcal D_h$ there belong to groups, those of $\mathcal D_h(P')$ forming
$G_{s,1}$. A double-zero series of $P'$ at $m(P')$ satisfies the hypothesis of
Lemma~\ref{rt:lem:greedy} ($|P'|\ge158$, and its weight after $P'$ is
$x-\eta<2$), so its tail is the greedy expansion of that weight up to
restarts, and its moment is at least the greedy one by
Lemma~\ref{rt:lem:restarts}~(ii). The weight lies within $10^{-7}$ of $x$, and
the parameters $v$ and $m(P')$ differ by less than $10^{-80}$, which changes
moments by less than $10^{-70}$. With the factor $1\pm10^{-5}$ and $\iota\le10^{-5}$
of Theorem~\ref{rt:thm:structure} this gives the lower bound for $m(P')-v$.
Now (a) and (b) leave no point of $\mathcal D_h$ between
$v+u_Pv^sL_-(1-\mu)$ and $m(P')$: points of deeper groups lie below that
bound, and the other groups of levels at most $s$ avoid the whole window.
Under (d), Lemma~\ref{rt:lem:restarts}~(ii) bounds the excess of the moment
over the greedy one by $v^{R+1}(q_1^*+10^{-5})\le7\cdot10^{-5}L_-$, whence the
upper bound for $m(P')-v$. For the free radius, a series with a double zero
at $y>m(P')$, $y-m(P')<(\mu-10^{-4})L_-u_{P'}$, has $y-v<F(P)u_P$ by (c); by
Theorem~\ref{rt:thm:structure} and (b) its first difference from $f_v$ is at
$n+s$ with $c=1$, so it belongs to $\Fh(P')$. Here $u_{P'}=u_Pv^s(1+O(10^{-9}))$.
\end{proof}

At a $C_h$-pattern (defined next) one has $y_s<1$: the digit after the
pattern digit is $0$, or, at $h=7$ where $t<1/3$, at most $2$, and all
later remainders are below $1+10^{-7}$.

\subsection{Patterns, cores and the finite lemmas}\label{rt:sec:patterns}
For each $h\in\{1,3,4,7\}$, a \emph{$C_h$-pattern} at level $s$ of a tail
$g$ is a digit $g_s=d\ge1$ such that $g_{s+1}g_{s+2}\cdots$ begins with one
of the following \emph{start words} $u$ for $d$:
\begin{itemize}[leftmargin=*]
\item $h=1$: $u=0000$ (so $d=1$ is followed by at least four zeros, or ends the tail);
\item $h=3$: $u=00$, or $u=01w$ with $w\in\{0,\,10,\,1100,\,11010,\,110110,\,11011100,\,11011101\}$;
\item $h=4$: $u=0$;
\item $h=7$: $u=0$ or $u=10$; or $u=w_1w_2w_3w_4$ with $dw_1w_2w_3w_4$ in the list
 $11111$, $11112$, $11113$, $11123$, $11211$, $11212$, $11213$, $11223$,
 $11311$, $11312$, $11313$, $12113$, $21111$, $21112$, $21113$, $21123$,
 $21211$, $21212$, $21213$, $21223$, $21311$, $21312$, $21313$, $22113$,
 $31111$, $31112$, $31113$, $31123$, $31211$, $31212$, $31213$, $31223$.
\end{itemize}
A tail with no pattern from some level on eventually lies in the
\emph{gapless core} $K_h$. By Parry's criterion, and the expansions
$d_\beta(1)$ of Table~\ref{rt:tab:data}, these cores are: $K_1$, the
concatenations of the blocks $100$ and $1000$; $K_3$, the tails over $\{0,1\}$
in which every $0$ is followed by at least three $1$s and a block of exactly
three $1$s is followed by a block of at least four; $K_4=\{1^\infty\}$;
and $K_7$, the weakly admissible tails over $\{1,2,3\}$ without any of the
$32$ windows.
(At $h=3$ a digit $2$ forces a pattern within five digits; at $h=4$ a $2$ is
followed by a $0$ within three digits; at $h=1$ the block $101$ forces
$10^5$.)

\begin{proposition}[Move Lemma]\label{rt:prop:move}
Let $h\in\{1,3,4,7\}$, let $v$ be a relative tip within $10^{-15}$ of the
centre of $[a_h,b_h]$ satisfying the hypothesis of
Lemma~\ref{rt:lem:greedy}, and let its tail have a $C_h$-pattern at a level
$s\ge J_h+2$, with $J_h=10,3,1,1$. Then $\mathcal L_s\subset[L_{\min}^\circ,L_{\max}^\circ]$
(Table~\ref{rt:tab:data}), condition (a) of Corollary~\ref{rt:cor:move} holds,
and so does (b) for the levels $s-J_h,\ldots,s$.
\end{proposition}

This is Proposition~\ref{cert:rt-move}; the certified clearances exceed
the margins $5\cdot10^{-4}L_\pm$ of the corollary. The levels below $s-J_h$
are handled by Lemma~\ref{rt:lem:restarts}. Let $M_h$ be the least $M$ with
$v^ML_{\max}^\circ(1+\mu)(1+3\cdot10^{-5})<10^{-3}$ and call a pattern at level
$s$ \emph{clean} if the tail has no restart at the levels $s-M_h+1,\ldots,s$.

\begin{lemma}[Beyond the window]\label{rt:lem:beyond}
At a clean $C_h$-pattern of level $s\ge J_h+2$ of a tip as in
Proposition~\ref{rt:prop:move}, condition (b) of Corollary~\ref{rt:cor:move}
also holds for the levels $s'\le s-J_h-1$.
\end{lemma}
\begin{proof}
By Lemma~\ref{rt:lem:restarts}~(iii) every point of $G_{s',c}$ exceeds
$u_Pv^{s'}(10^{-3}-2\cdot10^{-5})$, which is above the window when
$s-s'\ge M_h$. If $J_h+1\le s-s'<M_h$, cleanness puts the first restart
after $s'$ beyond level $s$, at distance $r\ge J_h+2$, and the group lies
above $u_Pv^{s'}(\Delta-v^{J_h+2}(q_1^*+10^{-5}))(1-2\cdot10^{-5})$, which is
above the window by (K1) and (K2) of Proposition~\ref{cert:rt-assembly}.
\end{proof}

The \emph{landing set} $\Lambda_h^{\mathrm{land}}$ consists of the
remainders $y_s$ at the last pattern $s$ of tails eventually in the core,
evaluated on the whole parameter window:
\begin{itemize}[leftmargin=*]
\item $h=1$: $y_s=t^{k+1}(1+t^{e_1}+t^{e_1+e_2}+\cdots)$ with $k\ge4$ and
 $e_j\in\{3,4\}$, or $y_s=0$;
\item $h=3$: $y_s=t^{k+1}w(t)$, where $w(t)=\sum_{i\ge0}w_it^i$ for an infinite
 word $w=1w_1w_2\cdots$ over $\{0,1\}$ satisfying the core condition, with
 $k\ge2$, or with $k=1$ and $w$ beginning with $10$, $110$ or $11101110$; or $y_s=0$;
\item $h=4$: $y_s=t^{k+1}/(1-t)$ with $k\ge1$, or $y_s=0$;
\item $h=7$: $y_s=t^{k+1}w(t)$ for a core word $w$ ($k\ge1$), or $y_s=0$, or
 $y_s=tw(t)$ for a zero-free, weakly admissible word $w$ beginning with one of
 the $32$ windows
 after $d$ and containing no window at the following levels.
\end{itemize}

\begin{proposition}[Landing Lemma]\label{rt:prop:landing}
Let $h\in\{1,3,4,7\}$, $y\in\Lambda_h^{\mathrm{land}}$, $|x-(1+y)|\le\epsilon_h$ with
$\epsilon_h=10^{-145},10^{-80},10^{-75},10^{-80}$, and $t\in[a_h-10^{-90},b_h+10^{-86}]$
($[a_1-10^{-150},b_1+10^{-140}]$ for $h=1$, $b_4+10^{-85}$ for $h=4$). Then the
greedy expansion of $x$ has a $C_h$-pattern at a level $\ge P_0^*$, with
$P_0^*=13,7,6,4$, decided by digits common to all such $x$ and $t$; its
first $R_h=26,12,11,8$ remainders exceed $10^{-20}$.
\end{proposition}

This is Proposition~\ref{cert:rt-landing}. Since $v^{R_h+1}(q_1^*+10^{-5})\le7\cdot10^{-5}L_{\min}^\circ$,
condition (d) of Corollary~\ref{rt:cor:move} then holds at the next tip.

\begin{table}[tbp]\centering\small
\caption{Data of the real-trace argument (Propositions~\ref{cert:rt-minimum}--\ref{cert:rt-landing}).
Lower bounds are rounded down and upper bounds up.}\label{rt:tab:data}
\begin{tabular}{@{}lllll@{}}\toprule
 & $h=1$ & $h=3$ & $h=4$ & $h=7$\\\midrule
$t_0=\min\mathcal D_h$ & $0.6684755322100606$ & $0.4490148201308485$ & $0.3870518533352138$ & $0.2835270145347340$\\
$\beta=1/t_0$, $\bar d$ & $1.49594$, $1$ & $2.22710$, $2$ & $2.58363$, $2$ & $3.52700$, $3$\\
$d_\beta(1)$ & $1010000010\cdots$ & $2010010120\cdots$ & $2110200110\cdots$ & $3130010303\cdots$\\
$n_0$, $N_h$ & $49$, $850$ & $24$, $257$ & $19$, $216$ & $15$, $162$\\
$\kappa_0$ & $13.6617$ & $10.9198$ & $11.8215$ & $17.5228$\\
$\Delta_{\mathrm{low}}$, $\Delta_1$ & $0.0452$, $0.0409$ & $0.4066$, $0.4066$ & $0.7928$, $0.7928$ & $0.7959$, $0.7959$\\
$q_1^*$ & $1.8896$ & $1.2468$ & $1.3407$ & $1.2208$\\
$L_{\min}$, $L_G$, $L_2$ & $22$, $10$, $30$ & $8$, $2$, $15$ & $6$, $1$, $13$ & $4$, $1$, $10$\\
$[L^\circ_{\min},L^\circ_{\max}]$ & $[0.609,2.046]$ & $[0.638,1.547]$ & $[0.519,1.542]$ & $[0.706,1.291]$\\
$J_h$, $P_0^*$, $M_h$, $R_h$ & $10$, $13$, $19$, $26$ & $3$, $7$, $10$, $12$ & $1$, $6$, $8$, $11$ & $1$, $4$, $6$, $8$\\
$F(P_0)\ge$ & $6.79$ & $2.17$ & $0.622$ & $1.61$\\
$i_1$, $X_h\le$ & $841$, $9\cdot10^{-149}$ & $255$, $3.4\cdot10^{-90}$ & $210$, $4.6\cdot10^{-88}$ & $158$, $4.3\cdot10^{-88}$\\
\bottomrule
\end{tabular}
\end{table}

\subsection{Proof of Theorem~\ref{thm:real-traces}}\label{rt:sec:proof}
Fix $h\in\{1,3,4,7\}$ and the constants of Table~\ref{rt:tab:data}. Put
$F_{\min}=(\mu-10^{-4})L^\circ_{\min}$; Proposition~\ref{cert:rt-assembly}
checks $v^{P_0^*}L^\circ_{\max}(1+\mu)<F_{\min}$ and $P_0^*\ge J_h+2$. Call a
relative tip $v=m(P)$ \emph{good} if $v\in[t_0,t_0+X_h]$, $F(P)\ge F_{\min}$,
$P=P_0$ or $|P|\ge158$, and the tail of $f_v$ has a $C_h$-pattern at a level
$\ge P_0^*$.

\smallskip\noindent\emph{Claim.} A good tip $v$ either carries infinitely many
gaps of $\mathcal D_h$ accumulating at $v$, or has a gap of $\mathcal D_h$
on its right followed by a good tip $v'$ with
$v<v'\le v+u_Pv^{P_0^*}L^\circ_{\max}(1+10^{-4})$.

\smallskip\noindent\emph{Proof of the claim.} At a clean pattern of level
$s\ge P_0^*$ with $n+s\ge158$, the hypotheses (a)--(c) of
Corollary~\ref{rt:cor:move} hold: (a) and (b) by Proposition~\ref{rt:prop:move}
and Lemma~\ref{rt:lem:beyond}, and (c) because
$v^sL^\circ_{\max}(1+\mu)<F_{\min}\le F(P)$. (i) Suppose that the tail has
infinitely many such patterns. Each gives the gap
$\bigl(v+u_Pv^sL_-(1-\mu),m(P'_s)\bigr)$, and for two of them, $s<s''$, the gap
and the tip of level $s''$ lie in groups of levels $>s$, hence below the gap
of level $s$ by (a). The gaps are disjoint, separated by the points
$m(P'_s)\in\mathcal D_h$, and accumulate at $v$. (ii) Otherwise the tail has
only finitely many patterns at levels $\ge P_0^*$: if it had infinitely many,
those beyond the last restart, or, if there are infinitely many restarts,
those at relative levels in $[M_h,M_E-10]$ after each restart (which exist by
(EX$'$) of Proposition~\ref{cert:rt-assembly}), would be clean, and all but
finitely many have $n+s\ge158$. Let $s$ be the last pattern at a level
$\ge P_0^*$. Then $n+s\ge158$, by the choice of $P_0$ and $i_1$ if $P=P_0$. The
pattern is clean: a restart at a level in $(s-M_h,s]$ would be followed by
$M_E$ digits of $d_\beta(1)$, which contain a pattern at a relative level at
least $M_h$, hence after $s$. After $s$ the tail lies in the core, and by
the case analysis that defines $\Lambda_h^{\mathrm{land}}$ its remainder $y_s$
belongs to the landing set. By Corollary~\ref{rt:cor:move} and
Lemma~\ref{rt:lem:greedy}, $v'=m(P')$ has the greedy tail, up to restarts,
of a weight $x$ with $|x-(1+y_s)|\le\bar\eta_{n+s}+X_h(2h+\bar d)/(1-\tau)^2\le\epsilon_h$
(Proposition~\ref{cert:rt-assembly}), where $y_s$ may be evaluated at the
parameter $v'$. By Proposition~\ref{rt:prop:landing} the first $R_h$ remainders of
the greedy expansion of $x$ exceed $10^{-20}$, which excludes restarts there
because $\delta_k<10^{-20}$; this is (d), so $F(P')\ge F_{\min}$. The greedy
expansion of $x$ also has a pattern at a level $\ge P_0^*$; it survives in the
tail of $v'$ unless a restart precedes it, in which case the following $M_E$
digits of $d_\beta(1)$ contain a pattern at a level $\ge P_0^*$. Finally
$|P'|=n+s\ge158$ and $v'\le v+u_Pv^sL^\circ_{\max}(1+10^{-4})$. \hfill$\square$

\smallskip\noindent\emph{Start.} By Propositions~\ref{cert:rt-minimum}
and~\ref{cert:rt-assembly}, $t_0=m(P_0)$ has $F(P_0)\ge F_{\min}$, its
remainders satisfy $p_{n_0}<2$ and $p_k<1$ for $n_0<k\le\min(200,N_h-2)$, and
its certified tail has $C_h$-patterns at levels $\ge P_0^*$, the deepest at
absolute position $i_1$. So $t_0$ is good. Every jump of the iteration is
made at an absolute depth at least $i_1$, the first by at most
$v^{i_1-1}L^\circ_{\max}(1+10^{-4})/(2\kappa_v)$ and each later one by a factor
$v^{P_0^*}$ less, so every tip lies in $[t_0,t_0+X_h]$. This interval is
inside the windows of all the finite lemmas, and the weight shifts and
parameter changes between consecutive tips are below $\epsilon_h$.
Hence every tip produced is good.

\smallskip\noindent\emph{Conclusion.} Iterating the claim from $t_0$, either
some tip carries infinitely many gaps accumulating at it, or there are
infinitely many tips $v_0<v_1<\cdots$ with a gap of $\mathcal D_h$ between
consecutive ones. In both cases the gaps are bounded components of
$(0,1)\setminus\mathcal D_h$ separated by points of $\mathcal D_h$, so
$\mathcal D_h$ and $(0,1)\setminus\mathcal D_h$ have infinitely many
components. The accumulation point $t_\infty$ of the gaps lies in $\mathcal A_h$:
in the first case $t_\infty=v$ has a tail with infinitely many nonzero corner
digits and $f_v''(v)>0$, so $f_v+z^N\in\Fh$ for infinitely many $N$ and the
argument of Lemma~\ref{low:lem:countable} gives nonreal zeros converging to
$t_\infty$; in the second case $t_\infty=\lim v_j$ is not isolated in
$\mathcal D_h$, and Lemma~\ref{low:lem:countable} applies.

It remains to treat the diagonal trace. The gaps are separated by tips to
which Lemma~\ref{rt:lem:greedy} applies: in the second case the $v_j$, and in
the first case the tips $m(P'_s)$ of the clean patterns with $n+s\ge158$
(at a pattern $y_s<1$, as noted after Corollary~\ref{rt:cor:move}, so the
weight $1+y_s$ is below $2$). Such a
tip $v$ has $f_v''(v)=2\kappa_v>0$ and corner digits at most $\bar d<2h$ after its
first tail digit, that is, coefficients $a_N>-h$. For such $N$ the series
$f_v-z^N$ is admissible, negative at $v$ and positive at $v\pm r_N$ for
suitable $r_N\to0$, so it has two real zeros in $(v-r_N,v+r_N)$, one on each
side of $v$. Hence $(v,v)$ is a limit of pairs $(\gamma,\lambda)$,
$\gamma\ne\lambda$, of zeros of one member of $\Fh$. Conversely a diagonal
limit of such pairs is a double zero, by coefficient compactness and
Rolle's theorem. So $\Delta_h\cap(0,1)\subseteq\mathcal D_h$ contains all the
separating tips, and the gaps of $\mathcal D_h$ separate them in $\Delta_h$.
This proves Theorem~\ref{thm:real-traces}. \qed

\begin{remark}\label{rt:rem:height-eight}
Height $8$ has $\beta=3.52706$, almost the base of height $7$, but its
slices are fatter, and a floating-point exploration (not certified) finds
moves only at three zero-free windows. The resulting core then has dimension
about $0.63$, and so does its landing set, so a landing branch and bound
cannot be expected to close; the planar core $\mathcal C_8$ is settled
by Proposition~\ref{prop:heights-seven-eight}. The argument of this section
also does not decide whether $\mathcal A_h$ has infinitely many components: a
tip whose tail terminates may be an isolated, inactive double zero.
\end{remark}

\section{The remaining question}\label{sec:remaining-scope}
Theorem~\ref{root:thm:main} determines every connected height.
Theorems~\ref{root:thm:infinite} and~\ref{root:thm:gap} show that
infinitely many components are the typical alternative. They leave the
following precise question.

\begin{conjecture}\label{conj:universal-infinitude}
For every $h\ge1$ that is not a product of consecutive positive
integers, $\mathcal C_h$ has infinitely many connected components.
\end{conjecture}

Let $E$ be the set of nonpronic heights whose cores have only finitely
many components, and put $E_q=E\cap\{q^2-2,q^2-1,q^2\}$. Then
\begin{equation}\label{eq:submission-residual}
 E\subseteq \{1,3,4,9\}\ \cup\ \{25,\,63,\,81\}\ \cup
 \bigcup_{q\ge13}E_q,\qquad |E_q|\le1,
\end{equation}
where $E_{14}\subset\{196\}$, $E_{15}=\varnothing$ and
$E_{16}\subset\{256\}$, by Propositions~\ref{prop:explicit-witnesses}
and~\ref{prop:square-witnesses}.

\subsection{A Diophantine reduction}
At a square-adjacent height with $q\ge4$, a finite anchor set forces a
witness that is eventually constant with the middle exceptional digit
(Corollary~\ref{cor:sa-criterion}).

\begin{proposition}\label{prop:diophantine-reduction}
Let $q\ge4$, $k\in\{0,1,2\}$ and $h=q^2-2+k$. If $\mathcal C_h$ has only
finitely many connected components, then there is a polynomial
$E(z)=\sum_{j=2}^Ne_jz^j$ of degree $N\ge4$ with integer coefficients
$|e_j|\le q-1$ such that
\begin{equation}\label{eq:channel-polynomial}
 R_E(z)=(1-qz)\bigl(1-(q+1)z\bigr)+(k-1)z^2+(1-z)E(z)
\end{equation}
has a double zero $1/\beta$ in $(1/(q+1),1/q)$. Moreover $g(1)^2$
divides $q^2-q-1+k$, where $g$ is the minimal polynomial of $\beta$.
\end{proposition}
\begin{proof}
By Corollary~\ref{cor:sa-criterion}, some square-adjacent witness $F$ is
eventually constant with digit $c=m+q-1=q^2-q-1+k$. Let $N$ be its last
level with $a_N(F)\ne c$; it exists, because the numerator
in~\eqref{eq:channel-series} has discriminant $5-4k\ne0$. Put
$e_j=a_j(F)-c$, so that $|e_j|\le q-1$. Then $F$ is the channel
series~\eqref{eq:channel-series} plus $E$, so $(1-z)F=R_E$, and the
anchor of $F$ is a double zero of $R_E$ in $(1/(q+1),1/q)$ by
Lemma~\ref{lem:sa-anchors}. The reversed polynomial
$x^{N+1}R_E(1/x)$ is monic of degree $N+1$ with integer coefficients,
and $\beta\in(q,q+1)$ is a double root of it. As $\beta$ is not an
integer, $g$ has degree at least two and $g^2$ divides the reversed
polynomial in $\mathbb Z[x]$. Hence $N+1\ge4$, and equality would make
the reversed polynomial equal to $g^2$, whose coefficient of $x^3$ is
even, whereas it is $-(2q+1)$. Evaluating at $x=1$ gives
$g(1)^2\mid R_E(1)=q^2-q-1+k$.
\end{proof}

Conversely, if $R_E$ has a double zero in $(0,1)$ for some $E$ with
$|e_j|\le q-1$, then $R_E/(1-z)$ is a square-adjacent witness that is
eventually constant with digit $m+q-1$. Thus
Conjecture~\ref{conj:universal-infinitude} holds at every height other
than $1$, $3$, $4$ and $9$ if none of the
polynomials~\eqref{eq:channel-polynomial}, for $q\ge4$ and
$k\in\{0,1,2\}$, has a double zero in $(0,1)$.

\begin{question}\label{q:constant-witness}
Can a polynomial~\eqref{eq:channel-polynomial} have a double zero in
$(0,1)$? More generally, does an eventually constant square-adjacent
witness exist?
\end{question}

We know of none. The searches of the companion~\cite{Companion26}
tested every constant digit after each retained prefix, with prefix
lengths $16,14,12,11,10,10$ for $q=4,\ldots,9$, and found no eventually
constant witness; this is a search result, not an exclusion. The
criterion of
Corollary~\ref{cor:sa-criterion} is only sufficient: a finite anchor set
would not rule out infinitely many components produced by a different
mechanism.

The anchor sets at the open heights are large: $\mathcal T_{25}$,
$\mathcal T_{63}$ and $\mathcal T_{81}$ have at least $384$, $332$ and
$369$ points, the numbers of prefixes of length $15$, $12$ and $12$
certified to extend to a square-adjacent witness by
Lemma~\ref{lem:ivt-completion}
(Proposition~\ref{cert:surviving-prefixes} of the
supplement). The searches of the companion found no witness
with a Pisot reciprocal anchor of degree at most six at $25$, or of
degree at most five at $63$ and $81$; they use numerical filters, so
these negative results are not nonexistence proofs. The search behind
Proposition~\ref{prop:square-witnesses} is exhaustive within its bounds
and uses integer arithmetic only. At every other square-adjacent height
with $4\le q\le12$ it finds a series $g^2/(1-z^r)$ with $\deg g\le9$ and
$r\le4$. At the heights $25$, $63$ and $81$ it finds none with $r\le6$
and $\deg g\le16$, $12$ and $11$ respectively, and none with $r\le6$ and
$\deg g\le9$ at the three heights with $q=13$.

\subsection{The heights \texorpdfstring{$1$, $3$, $4$ and $9$}{1, 3, 4 and 9}}
At the heights $1$ and $3$,
Proposition~\ref{prop:local-criterion} never applies. Indeed, let
$F=1+\sum_{j\ge1}a_jz^j$ and $t=1/\beta$ be as in that proposition,
where we may take $h_*=H=h$. The
equations $F(t)=tF'(t)=0$ give
\[
 \sum_{j\ge2}(j-1)a_jt^j=1,\qquad\sum_{j\ge2}(j-2)a_jt^j=2+a_1t,
\]
and the bounds $m\le a_j\le h$ and $a_1\ge-h$ turn these into
$(\beta-1)^2\ge m$ and $(2\beta-h)(\beta-1)^2\le h$. With the second
inequality of~\eqref{eq:compact-comparison-conditions},
$(\beta-1)^2\ge\max\{m,\ell\}$, where $m+\ell=h$. For $h=1$ this gives
$\beta\ge2$ and $(2\beta-1)(\beta-1)^2\ge3$; for $h=3$ it gives
$(\beta-1)^2\ge2$ and $(2\beta-3)(\beta-1)^2>3$. Both contradict
$(2\beta-h)(\beta-1)^2\le h$. At height $4$ the same inequalities,
together with $\Delta>0$, leave only the series
$1-4z+\sum_{j\ge2}a_jz^j$ with $2\le a_j\le4$ and $2.58<\beta<2.70$,
and at height $9$ the series $1-6z+\sum_{j\ge2}a_jz^j$ with
$5\le a_j\le9$ are candidates when $\beta>3.48$. We know no series
satisfying the hypotheses of the proposition at these two heights. The
search program behind Table~\ref{tab:explicit-witnesses}, run for $q=2$
and $q=3$, finds the two series of
Proposition~\ref{prop:heights-seven-eight}, but none in these two
families with a Pisot reciprocal anchor of degree at most six; it uses
numerical filters, so this is not a nonexistence proof. The search of
Proposition~\ref{prop:square-witnesses} finds no series $g^2/(1-z^r)$
in them with $r\le12$ and $\deg g\le36$ at height $4$, or
$\deg g\le24$ at height $9$. At the heights $1$, $3$ and $4$ the real
gaps of Theorem~\ref{thm:real-traces} do not settle the planar question,
as noted in the introduction.

\section*{Proofs and computation}
The proofs of Theorems~\ref{root:thm:main} and~\ref{root:thm:infinite} are
in Sections~\ref{sec:foundations}--\ref{sec:classification}, the proof of
Theorem~\ref{root:thm:gap} in Section~\ref{sec:gap-heights}, and the proof of
Theorem~\ref{thm:real-traces} in Section~\ref{sec:real-traces}, whose finite
inputs are the last four records of the supplement. The continuation
estimates behind Proposition~\ref{prop:pronic-core}, the planar separators
behind Proposition~\ref{prop:residual-separators} and the acceptance rules
of all finite checks are in the supplement. Theorem~\ref{root:thm:infinite}
and the separation argument of Section~\ref{sec:component-cascades} use no
computation. Every finite check is accepted by exact rational, symbolic or
bounded-integer arithmetic in the computational
companion~\cite{Companion26}. The search programs that proposed the
explicit witnesses of Proposition~\ref{prop:explicit-witnesses} use
numerical filters; the certificates that accept these witnesses do not.
The witnesses of Proposition~\ref{prop:square-witnesses} can be checked
by hand: each needs the square of one integer polynomial and its signs
at two integers, and Proposition~\ref{prop:heights-seven-eight} a few
decimal bounds. No negative search result is used as a proof.
Likewise, the programs that
proposed the certificates of the minimum of $\mathcal D_h$ use
floating-point arithmetic with safety margins, and only the exact checks
accept them. Figure~\ref{fig:height-status} is drawn from the statements of
Theorems~\ref{root:thm:main}--\ref{thm:real-traces} and of
Propositions~\ref{prop:explicit-witnesses},
\ref{prop:square-witnesses} and~\ref{prop:heights-seven-eight}; the
other figures are drawn by
Python scripts from exact data and from numerical illustrations identified
as such in their captions. No figure is an image generated by an AI model.

The author conceived the project, directed the research and checked every
argument. AI assistants (Anthropic Claude and OpenAI ChatGPT) were used under
the author's direction to help develop proofs, write and test code, and edit
the text. Every finite verification can be replayed from the companion
without any AI service.

\section*{Data and code availability}
The computational companion~\cite{Companion26}%
\ifx\CompanionDOI\empty
, supplied as supplementary material, contains
\else
\ is deposited at Zenodo
(\href{https://doi.org/\CompanionDOI}{doi:\CompanionDOI}) and supplied as
supplementary material; it contains
\fi
the certificate inputs, the checking programs and the replay logs, the
search programs with their logs, and the data and scripts for the figures
drawn from computed data. An addendum of 10~November 2026 accompanies
it: exact checks of Propositions~\ref{prop:square-witnesses}
and~\ref{prop:heights-seven-eight}, and the program and logs of the
searches of Sections~\ref{sec:heights-seven-eight}
and~\ref{sec:remaining-scope}. The README of the companion maps each
record of
Section~\ref{app:certified} of the supplement to its verifier, under the
numbering and the working names explained in Section~\ref{cert:sec:supplement} of the supplement.

\section*{Acknowledgements}
The author thanks David Juher and Joan Salda\~na for discussions on
collinear connectedness loci and restricted series, Christoph Bandt for
suggestions on exposition and terminology, and Pablo Shmerkin for
comments on related exposition. This work was supported by grant
PID2023-146424NB-I00 (Spain) and the IFUdG 2022--2024 fellowship programme
(Universitat de Girona and Banco Santander).

\section*{Declaration of generative AI and AI-assisted technologies in the manuscript preparation process}
During the preparation of this work the author used Anthropic Claude and
OpenAI ChatGPT in order to assist with developing proofs, writing code and
editing the text. After using these tools, the author reviewed and edited
the content as needed and takes full responsibility for the content of the
published article.

\clearpage
\setcounter{section}{0}
\renewcommand{\thesection}{S\arabic{section}}
\renewcommand{\theHsection}{S\arabic{section}}
\begin{center}
{\large\bfseries SUPPLEMENT}
\end{center}
\medskip
\noindent This supplement contains the technical parts of the proofs of
Theorems~\ref{root:thm:main}, \ref{root:thm:gap} and~\ref{thm:real-traces}
of the paper: the trace coordinates and continuation tests
(Section~\ref{sec:interfaces}), the complete real traces at the pronic heights
behind Proposition~\ref{prop:pronic-core} (Section~\ref{sec:complete-traces}),
the residual planar separators behind
Proposition~\ref{prop:residual-separators} (Section~\ref{sec:sep:framework}),
and the exact finite certificates with their acceptance rules, including
those used in Sections~\ref{sec:gap-heights} and~\ref{sec:real-traces}
(Section~\ref{app:certified}). Numbers without the prefix S refer to the
paper. Theorem~\ref{root:thm:infinite}, the density-one statement and the
separation argument of Section~\ref{sec:component-cascades} do not depend on
this supplement.

\section{Trace coordinates and continuation tests}\label{sec:interfaces}
The continuation proofs in Section~\ref{sec:complete-traces} use two
finite-word tests developed here. Section~\ref{sec:sep:framework}
gives the remaining planar separators, and Section~\ref{app:certified}
states and verifies their finite arithmetic obligations.

Retain $B_h,\mathcal F_h,\mathcal D_h,\mathcal A_h$ and $X$ from
Section~\ref{sec:foundations}, and write
$\mathcal Z_h^{\mathrm{ext}}=\mathcal Z_h\cap X=\Sigma_h$.
Proposition~\ref{prop:global-monic-reduction} supplies the root correspondence,
real slice, symmetries and annulus. In particular,
\begin{equation}\label{eq:root-ifs}
 \mathcal C_h^{\mathrm{ext}}=\operatorname{cl}_X
       (\mathcal Z_h^{\mathrm{ext}}\setminus\mathbb R).
\end{equation}
The other inputs are the inclusion $\mathcal A_h\subseteq\mathcal D_h$
from Lemma~\ref{prop:double-zero}, persistence of double-zero intervals
from Lemma~\ref{low:lem:countable}, and the pronic bound of
Lemma~\ref{lem:pronic-double-zero-bound}. These inputs precede the
connectedness classification and do not assume it.

\begin{lemma}[The pronic endpoint]\label{cor:min-double}
For $q\ge3$ and $h=q(q-1)$, $1/q\in\mathcal D_h$.
\end{lemma}
\begin{proof}
The series $(1-qz)^2/(1-z)=1-(2q-1)z+(q-1)^2\sum_{j\ge2}z^j$
belongs to $\mathcal F_h$ and has a double zero at $1/q$.
\end{proof}

\begin{lemma}[Real collision]\label{lem:real-collision}
A real limit of qualifying zeros of coefficientwise convergent bounded
real series is a multiple zero of the limiting series. Here qualifying
means nonreal or real multiple. If the limiting zero is exactly double
with positive second derivative, nearby qualifying witnesses $g$ at
zeros $z\ne t$ satisfy $g(t)>0$.
\end{lemma}
\begin{proof}
At a nonreal zero the conjugate zero is also present; the vanishing first
divided difference converges to the derivative at the real limit. For real
multiple zeros, pass to the derivative directly. For the last assertion,
$g(t)=|z-t|^2g[t,z,\bar z]$, with confluent arguments for a real double
zero. The second divided difference converges to $f''(t)/2>0$.
\end{proof}

\subsection{Trace coordinates and two entry tests}
Put $N=2h+1$. For a finite word $u=(d_1,\ldots,d_m)\in B_h^m$, $m\ge0$, define
\[
 P_u(z)=z^m-d_1z^{m-1}-\cdots-d_m,\qquad z_u(c)=2cP_u(c),
\]
with $P_\varnothing=1$ for the empty word. Appending $d$ gives $z_{u^\frown d}=c(z_u-2d)$, so
$z_u(c)$ is obtained from the initial point $2c$ by successively applying
the inverse maps $z\mapsto c(z-2d_1),\ldots,z\mapsto c(z-2d_m)$.

At $x>1$, set
\begin{equation}\label{eq:VS-def}
 V_u(x)=P_u(x)+xP'_u(x),\qquad
 S_u(x)=2P_u(x)+xP'_u(x).
\end{equation}
The empty word gives $(V,S)=(1,2)$, and appending $d$ acts by
\begin{equation}\label{eq:VS-rec-real}
 V_{u^\frown d}(x)=xS_u(x)-d,\qquad
 S_{u^\frown d}(x)=2V_{u^\frown d}(x)-xV_u(x).
\end{equation}
Put
\begin{equation}\label{eq:tau-sigma}
 \tau_h(x)=\frac{N-2x}{2x^2},\qquad
 \sigma_h(x)=\frac N{2x}.
\end{equation}
Then
\begin{equation}\label{eq:tau-identities}
 1+x\tau_h(x)=\sigma_h(x),\qquad
 \tau_h(x)>\tfrac12\quad\Longleftrightarrow\quad x^2+2x<N.
\end{equation}
For $c=a+ib$ with $b\ne0$, the complex coordinates and their inverse are
\begin{equation}\label{eq:reim-vs}
 V_c(z)=\frac{\operatorname{Im}z}{2b},\quad
 S_c(z)=\frac{\operatorname{Im}(cz)}{2|c|b};\qquad
 \operatorname{Im}z=2bV_c(z),\quad
 \operatorname{Re}z=2\bigl(|c|S_c(z)-aV_c(z)\bigr).
\end{equation}
Under an inverse digit map $z\mapsto z'=cz-2cd$, the coordinates
$V'_c=V_c(z')$ and $S'_c=S_c(z')$ satisfy
\begin{equation}\label{eq:complex-rec}
 V'_c=|c|S_c-d,\qquad
 S'_c=\frac{2a}{|c|}V'_c-|c|V_c.
\end{equation}

\begin{lemma}[Real-axis continuation]\label{lem:continuity}

For a fixed word $u$, the functions $V_c(z_u(c))$ and $S_c(z_u(c))$
extend continuously to $c=x>1$, with values $V_u(x)$ and $S_u(x)$.
\end{lemma}
\begin{proof}
For a real-coefficient polynomial $g$, the quotient
$\operatorname{Im}g(a+ib)/b$ tends to $g'(x)$ as $a+ib\to x$.
Apply this to $g(c)=2cP_u(c)$ and $g(c)=2c^2P_u(c)$, and use
$|c|\to x$ in the second coordinate.
\end{proof}

\begin{proposition}[Interior criterion]\label{prop:strict-trace}

Let $c=a+ib\in X$ with $a>0$, $b\ne0$, and $|c|^2+2a<N$.
If some finite word $u$ satisfies
\[
 |V_c(z_u(c))|<\frac{N-2a}{2|c|^2},\qquad
 |S_c(z_u(c))|<\frac N{2|c|},
\]
then $c\in\intr\mathcal Z_h^{\mathrm{ext}}$.
For a real parameter $1<x\le h+1$, the same conclusion follows from
\[
 x^2+2x<N,\qquad |V_u(x)|<\tau_h(x),\qquad
 |S_u(x)|<\sigma_h(x).
\]
\end{proposition}
\begin{proof}
\textit{Idea.} Nearest digits keep the inverse orbit of the marked point in
the open \emph{trace rectangle} given by the two inequalities, and a bounded
orbit is the series criterion.
Choose $d\in B_h$ nearest to $|c|S_c$; by the second inequality,
$|c|S_c\in(-h-1/2,h+1/2)$, so $\bigl||c|S_c-d\bigr|\le1/2$, and
\eqref{eq:complex-rec} gives
\[
 |V'_c|\le\tfrac12<\frac{N-2a}{2|c|^2},\qquad
 |S'_c|\le\frac a{|c|}+|c||V_c|<\frac N{2|c|}.
\]
So, along suitable inverse digits, the coordinates stay in the rectangle and,
by \eqref{eq:reim-vs}, the remainders stay bounded; hence
$1-\sum_{k=1}^md_kc^{-k}=z_{(d_1,\ldots,d_m)}(c)/(2c^{m+1})\to0$. Negating
digits gives the criterion of Proposition~\ref{prop:global-monic-reduction}, so
$c\in\mathcal Z_h^{\mathrm{ext}}$. The hypotheses persist
near $c$, and the argument applies at nearby parameters, giving interiority.
At a real $x$, \Cref{lem:continuity} gives them at all nearby nonreal
parameters, and nearby real ones lie in the real slice, since
$x^2+2x<2h+1$ forces $x<h+1$.
\end{proof}

\begin{corollary}[Nearest-digit lift]\label{cor:nearest}

If $1<x\le h+1$, $x^2+2x<N$, and
$|V_u(x)|\le\tau_h(x)$, $|S_u(x)|\le\sigma_h(x)$, then
$x$ is an interior point of $\mathcal Z_h^{\mathrm{ext}}$.
\end{corollary}
\begin{proof}
Choose $d\in B_h$ nearest to $xS_u$. Then
\[
 |V_{u^\frown d}|\le\tfrac12<\tau_h,\qquad
 |S_{u^\frown d}|\le1+x\tau_h=\sigma_h.
\]
A second nearest digit makes both inequalities strict, the second because
$|V_{u^\frown d}|<\tau_h$. Apply \Cref{prop:strict-trace}.
\end{proof}

\begin{corollary}\label{cor:interior-double-zero}

If $x>1$ is an interior point of $\mathcal Z_h^{\mathrm{ext}}$, then
$1/x\in\mathcal D_h$. In particular,
$\bigl(1/\sqrt{h+1},1\bigr)\subset\mathcal D_h$.
\end{corollary}
\begin{proof}
An interior point $x$ is a limit of nonreal points of
$\mathcal Z_h^{\mathrm{ext}}$, which lie in $\mathcal C_h$; so
$x\in\mathcal C_h$. Apply
Lemma~\ref{prop:double-zero} and Proposition~\ref{prop:annulus}.
\end{proof}

The same argument shows that $\mathcal Z_h^{\mathrm{ext}}$ and
$\mathcal C_h^{\mathrm{ext}}$ have the same interior points: an open set inside
$\mathcal Z_h^{\mathrm{ext}}$ has dense nonreal points, which lie in
$\mathcal C_h^{\mathrm{ext}}$ by \eqref{eq:root-ifs}, so it lies in
$\mathcal C_h^{\mathrm{ext}}$; the converse holds as
$\mathcal C_h^{\mathrm{ext}}\subset\mathcal Z_h^{\mathrm{ext}}$. Statements
about interior points below may therefore refer to either set. In particular,
by Proposition~\ref{prop:annulus}, the \emph{central annulus} $\{1<|c|<\sqrt{h+1}\}$ lies
in $\mathcal C_h^{\mathrm{ext}}$.

The second test uses $P=P_u(x)$ and $V=V_u(x)$, so that
$S_u(x)=P+V$; the empty word gives $(P,V)=(1,1)$. By
\eqref{eq:VS-rec-real}, appending a digit $d$ acts by
\begin{equation}\label{eq:PV-map}
 (P,V)\longmapsto(xP-d,\ xV+xP-d).
\end{equation}

\begin{lemma}[Bounded-orbit test]\label{lem:bounded-orbit}

Let $h\ge1$ and $x>1$, and let $d_1,d_2,\ldots\in B_h$. Put $(P_0,V_0)=(1,1)$ and
\[
 (P_m,V_m)=(xP_{m-1}-d_m,\ xV_{m-1}+xP_{m-1}-d_m)\qquad(m\ge1).
\]
If the orbit $(P_m,V_m)_{m\ge0}$ is bounded, then $1/x\in\mathcal D_h$; more
precisely, the admissible series $1-\sum_{k\ge1}d_kw^k$ has a double zero at
$w=1/x$.
\end{lemma}
\begin{proof}
Unrolling \eqref{eq:PV-map} gives
\[
 x^{-m}P_m=1-\sum_{k=1}^md_kx^{-k},\qquad
 x^{-m}\bigl(V_m-(m+1)P_m\bigr)=\sum_{k=1}^mkd_kx^{-k}.
\]
Boundedness makes both left-hand sides tend to zero. Thus the admissible
series $f(w)=1-\sum_{k\ge1}d_kw^k$ satisfies $f(1/x)=f'(1/x)=0$.
\end{proof}

The all-scale arguments below compare a finite-word state with its limit
at an endpoint. The following estimate applies to each choice of trace
coordinates and keeps this propagation separate from the estimates for
the initial states.

\begin{lemma}[Perturbation along a finite word]\label{lem:word-perturbation}
Fix a real square matrix $J$ and a vector $e$ of matching dimension, and let
$T_{x,d}Z=xJZ-de$. For a word $w=(d_1,\ldots,d_n)$, $n\ge1$, put
$T_{x,w}=T_{x,d_n}\circ\cdots\circ T_{x,d_1}$. Suppose
$0<x\le q$, $|d_i|\le h$, $\delta=q-x$,
$\|Z-Z^0\|_\infty\le E$ and $\|Z^0\|_\infty\le G$.
Write $C_j=\|J^j\|_\infty$ for the induced matrix norm and
$B_j=\|J^je\|_\infty$. Then
\begin{equation}\label{eq:word-unrolling}
 T_{x,w}Z=x^nJ^nZ-\sum_{i=1}^nd_i x^{n-i}J^{n-i}e,
\end{equation}
and
\begin{equation}\label{eq:word-parameter-error}
 \|T_{x,w}Z-T_{q,w}Z^0\|_\infty
 \le x^nC_nE+\delta\left(nq^{n-1}C_nG+
 h\sum_{j=1}^{n-1}jq^{j-1}B_j\right).
\end{equation}
Alternatively, put $Z_0^0=Z^0$ and
$Z_i^0=T_{q,(d_1,\ldots,d_i)}Z^0$ for $1\le i<n$. If
$\|JZ_i^0\|_\infty\le D_i$ for $0\le i<n$, then
\begin{equation}\label{eq:word-partial-state-error}
 \|T_{x,w}Z-T_{q,w}Z^0\|_\infty
 \le x^nC_nE+\delta\sum_{i=0}^{n-1}x^{n-1-i}C_{n-1-i}D_i.
\end{equation}
\end{lemma}
\begin{proof}
Iteration gives \eqref{eq:word-unrolling}. Subtract its version at
$(q,Z^0)$ and use $|q^j-x^j|\le jq^{j-1}\delta$ to obtain
\eqref{eq:word-parameter-error}. For the alternative, let $E_i$ be the
difference of the two states after $i$ digits. Then
$E_{i+1}=xJE_i-\delta JZ_i^0$, so iteration and the induced norm give
\eqref{eq:word-partial-state-error}.
\end{proof}

\section{Complete real traces at pronic heights}\label{sec:complete-traces}

The retraction argument in Section~\ref{sec:connected-heights} reduces
connectedness at a pronic height to a complete active real trace. It suffices to show that
$\mathcal C_h^{\mathrm{ext}}\cap(1,\infty)$ is a single interval $(1,a]$, that
is, that every point of $[1/a,1)$ is active and no smaller point is. For
$q\ge3$, Lemma~\ref{lem:pronic-double-zero-bound} confines the double zeros, hence the
active points, to $[1/q,1)$, and an interval of double zeros is active
(Lemma~\ref{low:lem:countable}); so it suffices that every point of $[1/q,1)$ is a
double zero. This section proves this at $h=q(q-1)$, $q\ge3$, and the
corresponding statement with $a=\rho$ at $q=2$. In the bulk, two digits suffice. Near the endpoint the
orbits we construct follow two parabolic channels, and no word of bounded
length works there; the terminal window is
completed in the all-scale proofs below. At height two, a comparison with two explicit
witnesses identifies the endpoint as $1/\rho$.

\subsection{Bulk capture}\label{sec:interior}
In \Cref{sec:interior,sec:residual,sec:pronic-completion}, $q\ge3$ is an integer,
$h=q(q-1)$ and $N=2h+1$, and $\delta=q-x$ for a real parameter $x$. By
Proposition~\ref{prop:annulus}, every $x\in(1,\sqrt{h+1})$ is an interior point of
$\mathcal Z_h^{\mathrm{ext}}$. One of $q$ two-digit words extends this to
$q-1/q$.

\begin{theorem}[Interior below $q-1/q$]\label{thm:bulk}

Let $q\ge3$ be an integer and $h=q(q-1)$. Every
$x\in[\sqrt{h+1},\,q-1/q]$ is an interior point of
$\mathcal Z_h^{\mathrm{ext}}$.
\end{theorem}

\begin{proof}
\textit{Idea.} The first digit $2q-1$ cancels the growth of the empty word;
the second, chosen among $q$ consecutive values, centres the value
coordinate.

Since $(q-\tfrac12)^2=h+\tfrac14<h+1$, we have $0<\delta<1/2$, and
$q\delta\ge1$. For the word $u_r=(2q-1,\,-h+r)$, two steps of
\eqref{eq:VS-rec-real} give
\[
 V_{u_r}=\theta-r,\qquad S_{u_r}=2(\theta-r)-g,
\]
where
\[
 \theta=q-2(q+1)\delta+3\delta^2
 =\bigl(x^2-(h+1)\bigr)+2\bigl(x-q+\tfrac12\bigr)^2+\tfrac12,\qquad
 g=x(2x-2q+1)=(q-\delta)(1-2\delta).
\]
The second expression for $\theta$ gives $\theta\ge1/2$, and $g>0$. Since
$\theta-q+\frac{4q+1}2\delta=3\delta(\delta-\frac12)<0$,
$g=q-2q\delta-\delta(1-2\delta)$ and $q\delta\ge1$,
\[
 \theta<q-\frac{4q+1}2\,\delta\le q-\frac{4q+1}{2q}=q-2-\frac1{2q},\qquad
 g<q-2q\delta\le q-2.
\]
Choose an integer $0\le r\le q-1$ with $|\theta-r|\le1/2$. Both digits of
$u_r$ lie in $B_h$, because $2q-1\le q(q-1)$ for $q\ge3$, and
\[
 |V_{u_r}|\le\tfrac12,\qquad |S_{u_r}|\le1+g<q-1.
\]
Finally,
\[
 N-x^2-2x=(q-1)(q-3)+(2q\delta-2)+\delta(2-\delta)>0,
\]
so $\tau_h(x)>1/2$ by \eqref{eq:tau-identities}, while
$\sigma_h(x)=N/(2x)>N/(2q)=q-1+1/(2q)$. The real case of
\Cref{prop:strict-trace} applies.
\end{proof}

\subsection{Two parabolic channels}\label{sec:residual}
For $q\ge4$ the strip $q-1/q<x<q$ remains. There $x^2+2x<N$, so the interior
criterion still applies, but the words it needs become long as $x\to q$. To see
why, use the real-jet coordinates $(P,V)$ of \eqref{eq:PV-map}, write
$F_{x,d}(P,V)=(xP-d,\,xV+xP-d)$ for the action of a digit $d$, and put
$r=x-1$. The map $F_{x,d}$ fixes $z_d=(d/r,\,-d/r^2)$, and repeating the digit
acts on deviations from $z_d$ by a Jordan block:
\begin{equation}\label{eq:jordan}
 F_{x,d}^{\,\ell}\bigl(z_d+(\xi,\zeta)\bigr)
 =z_d+x^\ell\bigl(\xi,\ \zeta+\ell\xi\bigr)\qquad(\ell\ge0).
\end{equation}
Indeed, one step gives $z_d+x(\xi,\zeta+\xi)$, because $-x+xr-r^2=-1$. We
call the resulting channels \emph{parabolic} because the linear part is $x$
times a unipotent, that is parabolic, matrix; the fixed points themselves are
repelling, with multiplier $x>1$.
Put
\[
 a=q-1,\qquad b=(q-1)(q-2).
\]
At $x=q$ the digits $-a^2$ and $-b$ have the fixed points
\[
 z_A=(-a,\,1),\qquad z_B=\Bigl(-(q-2),\,\frac{q-2}{q-1}\Bigr).
\]
Both lie just above the trace rectangle. Their value coordinates exceed
$\tau_h(q)=1-2/q+1/(2q^2)$ by $(4q-1)/(2q^2)$ and
$(2q^2-5q+1)/(2q^2(q-1))$, while their coordinates $S=P+V$, namely $2-q$ and
$-(q-2)^2/(q-1)$, lie inside $|S|<\sigma_h(q)=q-1+1/(2q)$. The two channels
are joined: the digit $-(a^2+1)$ maps the point $z_A-(0,1/a)$ exactly onto
$z_B$, because
\[
 F_{q,-(a^2+1)}\Bigl(-a,\ 1-\frac1a\Bigr)
 =\Bigl(-qa+a^2+1,\ q-\frac qa-qa+a^2+1\Bigr)
 =\Bigl(2-q,\ 1-\frac1a\Bigr)
\]
when $q=a+1$.

The words that follow the channels are the \emph{parabolic blocks}
\[
 \mathsf A_m=\bigl(2q-1,(-a^2)^{m-1}\bigr),\qquad
 \mathsf B_m=\bigl(2q-1,(-a^2)^{m-2},-(a^2+1)\bigr)\qquad(m\ge2),
\]
where $(d)^k$ denotes $k$ repetitions of the digit $d$ (words are written in
sans-serif), and their extensions
\[
 \mathsf U_{m,\ell}=\mathsf B_m(-b)^\ell\qquad(\ell\ge0).
\]
The first digit sends $(1,1)$ to a state whose $P$-coordinate
$(\delta^2-a^2)/r$ exceeds that of $z_{-a^2}=(-a^2/r,\,a^2/r^2)$ by
$\delta^2/r$. By \eqref{eq:jordan}, each repetition of $-a^2$ multiplies this
deviation by $x$, so the orbit stays near the first channel for about $m$
steps, where $\delta x^m$ is of order $1$. The digit $-(a^2+1)$ then moves it
to the second channel, along which \eqref{eq:jordan} again governs the
deviation from $z_{-b}$. The trace coordinates of the whole family are
explicit. Put
\[
 V_b^\sharp(x):=\frac{b}{(x-1)^2},
 \qquad
 S_b^\sharp(x):=\frac{(2-x)b}{(x-1)^2};
\]
these are the coordinates $(V,S)$ of $z_{-b}$, the fixed state of the digit
$-b$ in \eqref{eq:VS-rec-real}. Set
\begin{equation}\label{eq:u-m}
 u_m(x):=V_{\mathsf B_m}(x)-V_b^\sharp(x),
\end{equation}
\begin{equation}\label{eq:c-m-def}
 c_m(x):=\bigl(S_{\mathsf B_m}(x)-S_b^\sharp(x)\bigr)-\bigl(V_{\mathsf B_m}(x)-V_b^\sharp(x)\bigr).
\end{equation}
Then
\begin{equation}\label{eq:c-m}
 c_m(x)=\frac{\delta(\delta x^{m-1}-1)}{x-1}.
\end{equation}

\begin{proposition}[Parabolic return formulae]\label{prop:branched}

For every $m\ge2$ and $\ell\ge0$,
\begin{align}
 V_{\mathsf U_{m,\ell}}(x)
 &=V_b^\sharp(x)+x^\ell\bigl(u_m(x)+\ell c_m(x)\bigr),
 \label{eq:Uml-V}\\
 S_{\mathsf U_{m,\ell}}(x)
 &=S_b^\sharp(x)+x^\ell\bigl(u_m(x)+(\ell+1)c_m(x)\bigr).
 \label{eq:Uml-S}
\end{align}
\end{proposition}
\begin{proof}
Centred at $(V_b^\sharp,S_b^\sharp)$, appending $-b$ acts by $xJ$, where
\[
 J=\begin{pmatrix}0&1\\-1&2\end{pmatrix},\qquad
 J^\ell=I+\ell(J-I),
\]
because $(J-I)^2=0$. Apply this to the centred state $(u_m,u_m+c_m)$ of
$\mathsf B_m$. For \eqref{eq:c-m}, note that $c_m=P_{\mathsf B_m}(x)+b/r$, since
$S_u-V_u=P_u$ by \eqref{eq:VS-def}. Along $\mathsf B_m$, the deviation $\delta^2/r$
from $-a^2/r$ is multiplied by $x$ at each of the $m-2$ digits $-a^2$, and the
final digit $-(a^2+1)$ adds $1$. Hence
\[
 c_m=\frac{-a^2+x^{m-1}\delta^2}{r}+1+\frac br
 =\frac{x^{m-1}\delta^2-\delta}{r},
\]
because $b-a^2+r=-(q-1)+x-1=-\delta$.
\end{proof}

\subsection{Completion at heights \texorpdfstring{$q(q-1)$}{q(q-1)}}\label{sec:pronic-completion}
For $q\ge4$ the strip $(q-1/q,q)$ splits into intervals captured by six short
words, two compact windows and a terminal window. The three statements below
are proved in \Cref{app:terminal-high}; each uses finitely many exact checks.

\begin{theorem}[Refinement at scale $q^{-3}$\certmark]\label{thm:qminus3}

Let $q\ge4$ and $h=q(q-1)$. Every point of $(1,q)\setminus R_q^{(3)}$ is an
interior point of $\mathcal Z_h^{\mathrm{ext}}$, where
\[
\begin{aligned}
R_q^{(3)}={}&\left(q-\frac1{q^3},q\right)
\cup
\left(q-\frac1{2q^2}-\frac1{2q^3},\,q-\frac1{2q^2}+\frac1{3q^3}\right)\\
&\cup
\left(q-\frac1{2q}-\frac3{4q^3},\,q-\frac1{2q}+\frac3{4q^3}\right).
\end{aligned}
\]
In particular,
\[
        |R_q^{(3)}|=\frac{10}{3q^3}.
\]
\end{theorem}

\begin{proposition}[The compact windows\certmark]\label{prop:compact-residual-closure}

Let $q\ge4$ and $h=q(q-1)$, and put
\[
\begin{aligned}
J_{q,2}
&=\left[
q-\frac1{2q}-\frac3{4q^3},
q-\frac1{2q}+\frac3{4q^3}
\right],\\
J_{q,3}
&=\left[
q-\frac1{2q^2}-\frac1{2q^3},
q-\frac1{2q^2}+\frac1{3q^3}
\right].
\end{aligned}
\]
Every point of $J_{q,2}\cup J_{q,3}$ is an interior point of
$\mathcal Z_h^{\mathrm{ext}}$.
\end{proposition}

\begin{theorem}[Terminal window\certmark]\label{thm:terminal-completion}

For every integer $q\ge4$ and $h=q(q-1)$, every point of
$\left[q-\frac1{q^3},q\right)$ is an interior point of
$\mathcal Z_h^{\mathrm{ext}}$.
\end{theorem}

The short words are $\mathsf A_2,\mathsf B_2,\mathsf A_3,\mathsf B_3$,
$\mathsf A_2(-h)$ and $\mathsf B_2(-b)$; a single coefficient test proves all
their bounds. In the terminal
window, the parabolic blocks supply words of unbounded length. The window is cut into the bands
$\mathcal B_{q,m}=\{x:(q-x)x^{m-1}\le\frac12\le(q-x)x^m\}$, which accumulate
at $q$. On a band, most phases are captured by two families
of overlapping layers, $\mathsf B_m(-b)^\ell(-(b-1))$ and $\mathsf B_{m+1}(-b)^\ell$, whose
trace coordinates are given exactly by \Cref{prop:branched}. In a central
window of phases the capturing words leave the second channel and follow the
fixed state of the top digit $h$ before returning; there are about $q/4$ of
them, and their final digits depend on $q$ modulo~$4$. The limiting overlaps
are rational in $q$ and are checked symbolically on residue classes of $q$
modulo~$4$; explicit error bounds carry them to finite $x$, except for
finitely many pairs $(q,m)$, which are covered by exact rational cells.

One way to see why no gap is expected near the endpoint is that the witness
of the endpoint itself, $(1-qw)^2/(1-w)=1-(2q-1)w+a^2\sum_{j\ge2}w^j$, has the
digits $2q-1,-a^2,-a^2,\ldots$: it is the first channel. Orbits for $x$ near
$q$ shadow it and then correct with the digits $-a^2-k$, $0\le k\le a$, which
lie in $B_h$ exactly because $h-a^2=a$. At $h=2$ the analogous series
$(1-2w)^2/(1-w)$ has the coefficient $-3$, outside the alphabet, and indeed a
gap $(1/2,1/\rho)$ appears (\Cref{sec:height-two}).

At $q=3$ the trace rectangle does not reach the endpoint: $x^2+2x<13$ only
for $x<\sqrt{14}-1<2.75$. The height-six completion below therefore uses the bounded-orbit
test of \Cref{lem:bounded-orbit} instead, in three steps, for
$2999/1000\le x<3$. First, entry: a finite word brings $(1,1)$ into the
rectangle $R_x=[-3/r,3/r]\times[-2/r^2,2/r^2]$, with $r=x-1$, which moves with
$x$; the word is selected by the first crossing $(3-x)x^{M-1}\ge1/5$
(\Cref{seven:thm:terminal-orbit}). Second, renewal: from every point of $R_x$
a finite word returns to $R_x$; when a one-step return fails, two adjacent
digit choices remain, whose paired excursions either return or enter a fixed
compact cone, and all intermediate states stay in one bounded set
(\Cref{seven:thm:continuation}). Third, concatenating the returns gives an
infinite bounded orbit, hence a double zero at $1/x$ by
\Cref{lem:bounded-orbit}. On $(1,2999/1000]$ a finite buffered cover gives
interior points instead (\Cref{seven:thm:compact}).

\begin{theorem}[Pronic trace completion]\label{thm:complete-active-high}\label{thm:pronic-trace-completion}

Let $q\ge3$ be an integer and $h=q(q-1)$. Then
\[
 \mathcal D_h=\mathcal A_h=[1/q,1),\qquad
 \mathcal C_h^{\mathrm{ext}}\cap(1,\infty)=(1,q].
\]
\end{theorem}
\begin{proof}
By Lemma~\ref{lem:pronic-double-zero-bound}, $\mathcal D_h\subseteq[1/q,1)$. Let
$1<x<q$. If $q\ge4$, then $x$ is an interior point of
$\mathcal Z_h^{\mathrm{ext}}$: by Proposition~\ref{prop:annulus} if $x<\sqrt{h+1}$, by
\Cref{thm:bulk} if $\sqrt{h+1}\le x\le q-1/q$, and otherwise by
\Cref{thm:qminus3,prop:compact-residual-closure,thm:terminal-completion},
whose sets cover $(q-1/q,q)$. Then $1/x\in\mathcal D_h$ by
\Cref{cor:interior-double-zero}. If $q=3$, then $x$ is an interior point for
$x\le2999/1000$ (\Cref{seven:thm:compact}), so again $1/x\in\mathcal D_6$, and
$1/x\in\mathcal D_6$ for $2999/1000\le x<3$ by
\Cref{seven:thm:terminal-orbit}. Thus $(1/q,1)\subseteq\mathcal D_h$, and
$1/q\in\mathcal D_h$ by \Cref{cor:min-double}, so $\mathcal D_h=[1/q,1)$. By
Lemma~\ref{low:lem:countable} the whole interval, endpoint included, is active,
so $\mathcal A_h=\mathcal D_h$. Finally,
$\mathcal A_h=\{t\in(0,1):1/t\in\mathcal C_h\}$ gives the positive trace.
\end{proof}

\subsection{Height-two exclusion}\label{sec:height-two}
At $h=2$ the pronic endpoint $q=2$ is not attained. The attachment point is the
tribonacci constant $\rho=1.8392\ldots$, the real root of
$\rho^3=\rho^2+\rho+1$. Two admissible series explain the answer. Put
$G(w)=1-w-w^2-w^3$, so that $G(1/\rho)=0$. The series
\begin{equation}\label{tri:eq:endpoint-series}
 f_\rho(w)=\frac{G(w)^2}{1-w^3}
 =1-2w-w^2+w^3+\frac{w^4+w^5+2w^6}{1-w^3}
\end{equation}
has coefficients $1,-2,-1,1$ followed by the periodic block $1,1,2$, so it
belongs to $\mathcal F_2$, and it has a double zero at $1/\rho$. The series
\[
 f_{1/2}(w)=\frac{(1-2w)^2(1+w)}{1-w}=1-2w-2w^2+\frac{2w^3}{1-w}
\]
belongs to $\mathcal F_2$ and has a double zero at $1/2$, with
$f_{1/2}''(1/2)=24$. Shmerkin and Solomyak proved that $1/2$ is an isolated
double zero at this height and that the next one lies in
$(0.5436,0.5438)$~\cite[Proposition~6.3]{SS06}. The next theorem identifies it
as $1/\rho=0.54368\ldots$ and shows that all larger parameters are double
zeros.

\begin{theorem}[Height-two trace completion]\label{tri:thm:double-zeros}

One has
\[
 \mathcal D_2=\{1/2\}\cup[1/\rho,1),\qquad
 \mathcal A_2=[1/\rho,1).
\]
\end{theorem}

The inclusion $[1/\rho,1)\subseteq\mathcal A_2$ is proved in
\Cref{app:height-two}: a return set around the periodic configuration of
$f_\rho$ is entered along its digit block. The exclusion of all other double
zeros is elementary.

\begin{lemma}[Below one half]\label{tri:lem:small-five-symbol-exclusion}

If $f\in\mathcal F_2$ has a double zero $\lambda\in(0,1/2]$, then
$\lambda=1/2$ and $f=f_{1/2}$.
\end{lemma}
\begin{proof}
Write $f(w)=1+\sum_{k\ge1}a_kw^k$. Since $a_1\ge-2$, and $a_k\le2$ with
$1-k/2<0$ for $k\ge3$,
\[
 0=f(\lambda)-\frac\lambda2f'(\lambda)
 =1+\frac{a_1\lambda}2+\sum_{k\ge3}\Bigl(1-\frac k2\Bigr)a_k\lambda^k
 \ge1-\lambda-\frac{\lambda^3}{(1-\lambda)^2}
 =\frac{(1-2\lambda)(1-\lambda+\lambda^2)}{(1-\lambda)^2}.
\]
The right-hand side is positive for $\lambda<1/2$. At $\lambda=1/2$, equality
forces $a_1=-2$ and $a_k=2$ for $k\ge3$, and then $f(1/2)=0$ forces
$a_2=-2$.
\end{proof}

\begin{lemma}[Between one half and $1/\rho$]\label{lem:two-exclusion}

No $f\in\mathcal F_2$ has a double zero in $(1/2,1/\rho)$.
\end{lemma}
\begin{proof}
\textit{Idea.} Compare $f$ with the two witnesses. If $f$ does not begin like
$f_{1/2}$, its first deviation from $f_\rho$ must be negative at $\lambda$ and
increasing there, which the coefficient bounds forbid. If $f$ begins like
$f_{1/2}$, its first deviation from $f_{1/2}$ must come so late that the value
equation fails.

Fix $\lambda\in(1/2,1/\rho)$ and $f(w)=1+\sum_{k\ge1}a_kw^k\in\mathcal F_2$
with $f(\lambda)=f'(\lambda)=0$, and put $x=1/\lambda$ and
$\lambda_+=5437/10000$. Since $G$ decreases on $(0,1)$ and
$G(\lambda_+)=-32678453\cdot10^{-12}<0=G(1/\rho)$, we have
$\lambda<\lambda_+$ and $G(\lambda)>0$.
We use the series
\[
 \Psi_0=\frac{w+w^2}{1-w^3},\qquad
 \Psi_1=\frac{w+w^3}{1-w^3},\qquad
 \Psi_2=\frac{w^2+w^3}{1-w^3},
\]
whose coefficients lie in $\{0,1\}$, and four facts.
\begin{enumerate}
\item $1-\Psi_0=G/(1-w^3)$.
\item The coefficients $e_n$ of $f_\rho$ satisfy $e_n\ge1$ for $n\ge3$, and for
 $n\ge3$ the series $\sum_{k\ge1}(2-e_{n+k})w^k$ equals $\Psi_s$ with
 $s\equiv n\pmod3$.
\item On $(0,1)$, $\Psi_s\le\Psi_0$ and $\Psi_s'\le\Psi_0'$ for $s=1,2$, and
 $(\Psi_0-w)'\le\Psi_2'$ and $(\Psi_1-w)'\le\Psi_2'$. Indeed, the differences
 \[
  \frac{w^2}{1+w+w^2},\quad\frac{w(1+w)}{1+w+w^2},\quad
  \frac{w^3}{1+w+w^2},\quad\frac{w^2(1+w)}{1+w+w^2}
 \]
 are nonnegative and increasing.
\item For $0<\lambda\le\lambda_+$,
 \[
  \lambda^3<2(1-\lambda)^3,\qquad
  \Psi_2'(\lambda)=\frac{2\lambda+3\lambda^2+\lambda^4}{(1-\lambda^3)^2}<3.
 \]
 The two sides of the first are monotone in opposite directions, and
 $\Psi_2'$ increases, so it is enough to check both at $\lambda_+$, where
 $\lambda_+^3-2(1-\lambda_+)^3=-29289174641\cdot10^{-12}$ and
 $\Psi_2'(\lambda_+)<2.927$.
\end{enumerate}

\emph{Case $(a_1,a_2)\ne(-2,-2)$.} Since $f_\rho(\lambda)>0$, $f\ne f_\rho$.
Let $\nu$ be the first index with $a_\nu\ne e_\nu$, and write
$f=f_\rho+w^\nu C$ with $C=\sum_{k\ge0}c_kw^k$ and $c_0\ne0$. Then
$c_k\le2-e_{\nu+k}\le1$ for $k\ge2$. Put $\varpi=x^\nu G(\lambda)>0$ and
$\zeta=-f_\rho'(\lambda)/G(\lambda)=2\Psi_0'(\lambda)+3\lambda^2G(\lambda)/(1-\lambda^3)^2$,
so that $\zeta>2\Psi_0'(\lambda)$. The equations $f(\lambda)=f'(\lambda)=0$
give
\begin{equation}\label{eq:two-C}
 C(\lambda)=-\frac{\varpi G(\lambda)}{1-\lambda^3}<0,\qquad
 C'(\lambda)=x^\nu\bigl(\nu xf_\rho(\lambda)-f_\rho'(\lambda)\bigr)\ge\varpi\zeta>0.
\end{equation}
If $\nu=1$, then $c_0\ge1$ and $c_1\ge-1$. If $\nu=2$, then $c_0=a_2+1$, which
is $-1$ only when $a_1=a_2=-2$. If $\nu\ge3$, then $c_0\le1$. So one of three
cases occurs.

\emph{(A) $c_0\le-1$ and $\nu\ge3$.} By facts~(ii) and~(iii),
$C(\lambda)\le-1+\Psi_s(\lambda)\le-1+\Psi_0(\lambda)$ and
$C'(\lambda)\le\Psi_0'(\lambda)$. With fact~(i) and \eqref{eq:two-C}, the first
gives $\varpi\ge1$, and the second gives
$2\varpi\Psi_0'(\lambda)<\varpi\zeta\le\Psi_0'(\lambda)$, so $\varpi<1/2$.

\emph{(B) $c_0\ge2$, or $c_0=1$ and $c_1\ge-2$.} Then
$2c_0/\lambda+c_1\ge2/\lambda-2$; for $c_0\ge2$ use $c_1\ge-4$ and
$\lambda<1$. Since $c_k\le1$ for $k\ge2$ and $C(\lambda)<0$,
\[
 C'(\lambda)-\frac{2C(\lambda)}\lambda=\sum_{k\ge0}(k-2)c_k\lambda^{k-1}
 \le-\frac{2c_0}\lambda-c_1+\sum_{k\ge2}(k-2)\lambda^{k-1},
\]
and the last sum is $\lambda^2/(1-\lambda)^2$. Hence
\[
 C'(\lambda)<\frac{\lambda^2}{(1-\lambda)^2}-\frac2\lambda+2
 =\frac{\lambda^3-2(1-\lambda)^3}{\lambda(1-\lambda)^2}<0
\]
by fact~(iv), contrary to \eqref{eq:two-C}.

\emph{(C) $c_0=1$ and $c_1\le-3$.} Here $\nu\ge2$. The series
$\Xi=\sum_{k\ge2}(2-e_{\nu+k})w^k$ dominates $\sum_{k\ge2}c_kw^k$
coefficientwise. It equals $\Psi_2$ if $\nu=2$, and if $\nu\ge3$ it is
$\Psi_s$ minus its linear term, namely $\Psi_0-w$, $\Psi_1-w$ or $\Psi_2$. By
fact~(iii), $\Xi'\le\Psi_2'$, so $C'(\lambda)\le-3+\Psi_2'(\lambda)<0$ by fact~(iv),
contrary to \eqref{eq:two-C}.

\emph{Case $a_1=a_2=-2$.} Since $f_{1/2}(\lambda)>0$, $f\ne f_{1/2}$. Let
$n\ge3$ be the first index with $a_n\ne2$, and write
$f_{1/2}-f=w^n D$ with $D=\sum_{k\ge0}(2-a_{n+k})w^k$. Its coefficients
lie in $\{0,\ldots,4\}$ and the first is positive, so
\[
 1\le D(\lambda)\le\frac4{1-\lambda},\qquad
 0\le\lambda D'(\lambda)\le\frac{4\lambda}{(1-\lambda)^2}.
\]
The equations become $\lambda^n D(\lambda)=f_{1/2}(\lambda)$ and, after
division,
\[
 n+\frac{\lambda D'(\lambda)}{D(\lambda)}
 =\frac{\lambda f_{1/2}'(\lambda)}{f_{1/2}(\lambda)}
 =\frac{4\lambda}{2\lambda-1}+\frac{2\lambda}{1-\lambda^2}.
\]
Put $\varepsilon=2\lambda-1\in(0,2\lambda_+-1)$. Since
$4\lambda/(2\lambda-1)=2+2/\varepsilon$, $2\lambda/(1-\lambda^2)>0$ and
$4\lambda/(1-\lambda)^2<4\lambda_+/(1-\lambda_+)^2=\frac{217480000}{20820969}<\frac{21}2$,
\[
 n>\frac2\varepsilon-\frac{17}2>\frac2{2\lambda_+-1}-\frac{17}2
 =\frac{12571}{874}>14,
\]
so $n\ge15$. On the other hand, since $2/\varepsilon<n+17/2$,
\[
 x^n=\frac{D(\lambda)}{f_{1/2}(\lambda)}\le\frac4{\varepsilon^2(1+\lambda)}
 <\Bigl(\frac2\varepsilon\Bigr)^2<\Bigl(n+\frac{17}2\Bigr)^2.
\]
But $x>1/\lambda_+>9/5$, and $(9/5)^n>(n+\frac{17}2)^2$ for every
$n\ge15$: it holds at $n=15$, as $(9/5)^{15}>6746>\frac{2209}4$, and each
step multiplies the left side by
$9/5$ and the right side by at most $(49/47)^2<9/5$. This contradiction
completes the proof.
\end{proof}

\begin{lemma}[One half is inactive]\label{lem:half-inactive}

One has $1/2\notin\mathcal A_2$; equivalently, $2\notin\mathcal C_2$.
\end{lemma}
\begin{proof}
Suppose that $f_j\in\mathcal F_2$ have nonreal zeros
$z_j\to1/2$. By \Cref{lem:real-collision}, every coefficientwise limit has a
double zero at $1/2$. It is $f_{1/2}$ by
\Cref{tri:lem:small-five-symbol-exclusion}, so $f_j\to f_{1/2}$, and since
$f_{1/2}''(1/2)>0$, \Cref{lem:real-collision} gives $f_j(1/2)>0$ for large $j$. But for
large $j$ the first two coefficients of $f_j$ equal $-2$ and all others are at
most $2$, so
\[
 f_j(1/2)\le1-1-\frac12+\sum_{k\ge3}\frac2{2^k}=0,
\]
a contradiction.
\end{proof}

\begin{proof}[Proof of \Cref{tri:thm:double-zeros}]
By \Cref{tri:lem:small-five-symbol-exclusion,lem:two-exclusion},
$\mathcal D_2\cap(0,1/\rho)=\{1/2\}$, attained by $f_{1/2}$.
\Cref{low:prop:ternary} gives $[1/\rho,1)\subseteq\mathcal A_2\subseteq\mathcal D_2$.
This proves the formula for $\mathcal D_2$, and \Cref{lem:half-inactive}
gives $1/2\notin\mathcal A_2$.
\end{proof}

\medskip\noindent\textbf{All-scale completion.}\label{app:terminal}
Every estimate used at infinitely many scales is proved below. Only the
exact finite obligations in Section~\ref{app:certified} are delegated to code.

\subsection{Heights \texorpdfstring{$q(q-1)$ with $q\ge4$}{q(q-1) with q >= 4}}\label{app:terminal-high}

Throughout, $q\ge4$, $h=q(q-1)$, $N=2h+1$, and $a$, $b$, $\delta=q-x$, $r=x-1$
are as in \Cref{sec:residual}; \emph{interior} means an interior point of
$\mathcal Z_h^{\mathrm{ext}}$, and $\tau_*=\tau_h(q)=1-\frac2q+\frac1{2q^2}$,
$\sigma_*=\sigma_h(q)=a+\frac1{2q}$. We also write $\bar q=q-\frac12$, and
$j$ for the offset of an outer layer; words are named by their digits or by
the blocks of \Cref{sec:residual}.

This subsection proves
\Cref{thm:qminus3,prop:compact-residual-closure,thm:terminal-completion},
which with Proposition~\ref{prop:annulus} and \Cref{thm:bulk} make every $x\in(1,q)$ interior, as
\Cref{thm:complete-active-high} requires. Six short words leave three
windows; finitely many words cover the two compact ones. In the terminal
window, phase bands accumulate at $q$; on each band two families of outer
layers and a central hierarchy of returning words overlap with explicit error
bounds, and finitely many pairs $(q,m)$ are checked exactly.

\subsubsection{Short words}
For $0<\delta<1/q$ we use $\mathsf A_2,\mathsf B_2,\mathsf A_3,\mathsf B_3$ and
$\mathsf A_2(-h)=(2q-1,-a^2,-h)$ and $\mathsf B_2(-b)=(2q-1,-(a^2+1),-b)
=\mathsf U_{2,1}$. They capture the following
$\delta$-intervals, all nondegenerate for $q\ge4$; the degree columns are
explained in the proof below.
\begin{center}\small
\renewcommand{\arraystretch}{1.05}
\begin{tabular}{@{}cllcc@{}}
\toprule
Word & $\delta_-(q)$ & $\delta_+(q)$ & $V$ degrees & $S$ degrees\\
\midrule
$\mathsf A_2$ & $q^{-2}-q^{-3}$ & $(2q)^{-1}-(4q^2)^{-1}$ & $(14,4)$ & $(11,3)$\\
$\mathsf B_2$ & $(2q)^{-1}+(2q^2)^{-1}$ & $q^{-1}$ & $(10,4)$ & $(8,3)$\\
$\mathsf A_3$ & $q^{-3}$ & $(2q^2)^{-1}-(3q^3)^{-1}$ & $(17,5)$ & $(14,4)$\\
$\mathsf B_3$ & $(2q^2)^{-1}+(2q^3)^{-1}$ & $q^{-2}$ & $(17,5)$ & $(14,4)$\\
$\mathsf A_2(-h)$ & $(2q)^{-1}-(4q^2)^{-1}$ & $(2q)^{-1}-3/(4q^3)$ & $(17,5)$ & $(14,4)$\\
$\mathsf B_2(-b)$ & $(2q)^{-1}+3/(4q^3)$ & $(2q)^{-1}+(2q^2)^{-1}$ & $(17,5)$ & $(14,4)$\\
\bottomrule
\end{tabular}
\end{center}

\begin{lemma}[Direct interval positivity\certmark]\label{lem:short-word-direct}

For every word and interval in the table, every real $q\ge4$, and
$\delta\in[\delta_-(q),\delta_+(q)]$, the four margins
\begin{equation}\label{eq:short-four-margins}
 N-2x\pm2x^2V_u(x),\qquad N\pm2xS_u(x),
 \qquad N=2q(q-1)+1,
\end{equation}
are strictly positive, giving the interval captures used in
\Cref{thm:qminus3}.
\end{lemma}
\begin{proof}
\textit{Idea.} After rescaling, each margin is a polynomial with nonnegative
Bernstein coefficients on a half-line times an interval.
Set $\delta=\delta_-+(\delta_+-\delta_-)t$ with $0\le t\le1$, form the margins
from the literal words by \eqref{eq:VS-rec-real}, multiply by a positive
constant times a power of $q$, and put $q=4+s$. The result $p(s,t)$ is a
polynomial whose degrees in $(s,t)$ are at most the degree pair of its
coordinate. If $p=\sum a_{ke}s^kt^e$, then
\begin{equation}\label{eq:short-normalized-certificate}
 p(s,t)=\sum_{k=0}^{d_s}s^k\sum_{i=0}^{d_t}
 b_{ki}\binom{d_t}{i}t^i(1-t)^{d_t-i},\qquad
 b_{ki}=\sum_{e\le i}a_{ke}\binom ie\Big/\binom{d_t}e .
\end{equation}
If all $b_{ki}\ge0$ and all $b_{0i}>0$, then $p>0$ for $s\ge0$ and
$0\le t\le1$, since the Bernstein basis is nonnegative and sums to one. The
positivity of the twenty-four margins is \Cref{cert:short-word-strip}\certmark;
its certificate is this coefficient test.
\end{proof}

\begin{proof}[Proof of \Cref{thm:qminus3}]
\textit{Idea.} The six capture intervals tile $(0,1/q]$ in $\delta$ except for
three windows. Proposition~\ref{prop:annulus} and \Cref{thm:bulk} cover $(1,q-1/q]$. For
$0<\delta\le1/q$, \Cref{lem:short-word-direct} gives strict capture on the six
closed $\delta$-intervals of the table. In increasing $\delta$ they belong to
$\mathsf A_3,\mathsf B_3,\mathsf A_2,\mathsf A_2(-h),\mathsf B_2(-b),\mathsf B_2$;
consecutive ones overlap or share an endpoint except for two gaps, and the
first starts at $q^{-3}$. The omitted parts are exactly the three windows of
$R_q^{(3)}$, of lengths $q^{-3}$, $5/(6q^3)$ and $3/(2q^3)$, whose sum is
$10/(3q^3)$. The six words have digits in $B_h$ for every integer $q\ge4$,
and $x<q$ gives $x^2+2x<q^2+2q<N$, so \Cref{prop:strict-trace} gives
interiority on the union.
\end{proof}

\subsubsection{The compact windows}
\begin{proof}[Proof of \Cref{prop:compact-residual-closure}]
\textit{Idea.} Finitely many words per window, with margins positive
uniformly in $q$, as for the short words. A word with digits in $B_h$ whose
four margins \eqref{eq:short-four-margins} are positive at $x$ captures $x$,
and then $x$ is interior by \Cref{prop:strict-trace}, as $x^2+2x<N$. The
statement is \Cref{cert:compact-windows}\certmark{} (interior points of
$\mathcal C_h^{\mathrm{ext}}$ and of $\mathcal Z_h^{\mathrm{ext}}$ coincide, by
the remark after \Cref{cor:interior-double-zero}). Its certificate covers $J_{q,2}$ and $J_{q,3}$ by
consecutive closed cells, each with such a word: rational cells for finitely
many $q$ and, on residue classes of $q$ modulo~$4$, cells whose endpoints and
words depend polynomially on $q$, with margins positive on each cell by
\Cref{sep:lem:bernstein}.
\end{proof}

\subsubsection{Terminal capture at every scale}\label{sec:terminal-completion}
Write $Z_u=(V_u,S_u)$; a word $u$ \emph{captures} $x$ if
$|V_u(x)|<\tau_h(x)$ and $|S_u(x)|<\sigma_h(x)$. For $x<q$, $x^2+2x<N$,
$\tau_h(x)>\tau_*$ and $\sigma_h(x)>\sigma_*$, so capture gives interiority
(\Cref{prop:strict-trace}). The estimates below cover all large scales and, at
each remaining scale, all large $q$; \eqref{eq:terminal-finite-range} lists the
finite complement.

\begin{definition}[Trace cells and linked covers]\label{def:terminal-trace-cover}

In a common real coordinate $T$, a trace cell has four strict affine
endpoint inequalities
\[
 L_V<T<U_V,\qquad L_S<T<U_S.
\]
They express the two strict trace bounds of its assigned word.
An ordered list of cells is \emph{linked} if each lower endpoint of
one cell is strictly below both upper endpoints of its predecessor.
\end{definition}

\begin{lemma}[Strict chain principle]\label{lem:terminal-cover-principle}

Let $C_0,\ldots,C_k$ be a linked list of trace cells in a common
coordinate $T$. If $T$ exceeds both lower endpoints of $C_0$ and is
less than both upper endpoints of $C_k$, then one of the cells
contains $T$. Whenever its four inequalities imply strict capture
at a parameter $x$, the point $x$ is an interior point of
$\mathcal Z_h^{\mathrm{ext}}$.
\end{lemma}
\begin{proof}
Suppose both lower inequalities hold at one cell. Either its two upper
inequalities also hold, or a failed upper inequality, together with the link
inequalities, implies both lower inequalities of the next cell. Induction
reaches a capturing cell, since both upper inequalities hold at the final
cell. Apply \Cref{prop:strict-trace}.
\end{proof}

A cell whose estimates hold only for $T$ above a constant carries that
constant as a fifth, constant lower endpoint; linked lists and
\Cref{lem:terminal-cover-principle}, with the same proof, apply to such
cells.

All chains below are evaluated at a fixed parameter $x$, so their endpoints
are numbers and the comparisons give pointwise capture.

\paragraph{Phase bands.} For $m\ge4$ put
\begin{equation}\label{eq:terminal-centre}
 \mathcal B_{q,m}=\Bigl\{x\in[q-q^{-3},q):
 (q-x)x^{m-1}\le\tfrac12\le(q-x)x^m\Bigr\}.
\end{equation}
For $x\in\mathcal B_{q,m}$ define the phases $s=1-2\delta x^{m-1}$ and
$t=2\delta x^m-1$ and the centred coordinates $T_m=-x^ms\le0$ and
$T_{m+1}=x^{m+1}t\ge0$, so that
\begin{equation}\label{eq:terminal-transition}
 t=x(1-s)-1=r-xs.
\end{equation}
With $\bar q=q-\frac12$ and $D_m=1/(2\bar q^{\,m-1})$ we shall repeatedly use
$0<\delta\le D_m$, $x\ge\bar q$ and $r\ge q-\frac32$.
The two outer families are
$\mathsf W_{m,\ell}=\mathsf B_m(-b)^\ell(-(b-1))$ and
$\mathsf U_{m+1,\ell}=\mathsf B_{m+1}(-b)^\ell$, with
$\mathsf W_{m,-1}=\mathsf A_m$. A $\mathsf W$-layer of offset $j$ has
$\ell=m-j-2$ and coordinate $z_W=x^{\ell+2}s=-T_m/x^j$; a $\mathsf U$-layer
has $\ell=m-j$ and $z_U=x^{\ell+1}t=T_{m+1}/x^j$. The variable layers are
$1\le j\le m-1$ ($\mathsf W$) and $1\le j\le m$ ($\mathsf U$); the fixed cells
have $j\in\{0,-1\}$. \Cref{tab:terminal-roadmap} lists the cells of one band
in the order in which they are linked below, with the coordinate in which each
is used and the results that control it.

\begin{table}[tbp]
\centering\small
\caption{The linked chain on one band $\mathcal B_{q,m}$, listed from the
centre $T_m=0$ to the centre $T_{m+1}=0$; read in the opposite order, it is
linked in the sense of \Cref{def:terminal-trace-cover} (proof of
\Cref{thm:terminal-completion}). For $4\le q\le7$ the chains of
\Cref{tab:low-q-chains} replace the central hierarchy from $\mathsf P_m$ to
$\mathsf R_{m,1,0}$, and the proof of \Cref{prop:terminal-central-hierarchy}
bounds the fixed cells. The pairs $(q,m)$ of \eqref{eq:terminal-finite-range}
and the shoulder of \Cref{prop:terminal-shoulder} are covered instead by the
rational and symbolic cells of \Cref{cert:terminal-finite}.}
\label{tab:terminal-roadmap}
\begin{tabular}{@{}p{.17\linewidth}p{.33\linewidth}p{.19\linewidth}p{.23\linewidth}@{}}
\toprule
Stage & Words & Coordinate & Controlled by\\
\midrule
central hierarchy & $\mathsf P_m$, residue chains,
 $\Phi_{K-1},\ldots,\Phi_0=\mathsf R_{m,0,0}$ & $T_m\in[-\frac12,0]$ &
 \Cref{lem:terminal-base-errors,prop:terminal-central-hierarchy}\\
return & $\mathsf R_{m,1,0}$ & $T_m\in[-\frac23,0]$ & the same\\
fixed cells & $\mathsf W_{m,m-1}$, $\mathsf W_{m,m-2}$ ($j=-1,0$) &
 $z_W=-T_m/x^j$ &
 \Cref{lem:terminal-fixed-offset,prop:terminal-central-hierarchy}\\
variable layers & $\mathsf W_{m,\ell}$ ($j=1,\ldots,m-1$), the seam
 \eqref{eq:terminal-transition}, $\mathsf U_{m+1,\ell}$ ($j=m,\ldots,1$) &
 $z_W$, then $z_U$ &
 \Cref{lem:terminal-outer-formulas,prop:terminal-outer-annuli}\\
fixed cells & $\mathsf U_{m+1,m}$, $\mathsf P_{m+1}$ ($j=0,-1$) &
 $z_U=T_{m+1}/x^j$ &
 \Cref{lem:terminal-fixed-offset,prop:terminal-central-hierarchy}; end of the
 proof of \Cref{thm:terminal-completion}\\
\bottomrule
\end{tabular}
\end{table}

\begin{lemma}[Phase selection]\label{lem:terminal-coordinate-monotonicity}

Every $x\in[q-q^{-3},q)$ satisfies
\[
 x\le q-\left(\frac12+\frac{3}{4q^4}\right)q^{-3}
 \quad\text{or}\quad x\in\mathcal B_{q,m}\text{ for some }m\ge4.
\]
Throughout this interval, $x\ge q-1/2>1$ and
$x\in[q-1/q,q]$. On each selected band,
$0\le s<1$, $0\le t\le x-1$, and $0<\delta\le D_m$.
\end{lemma}
\begin{proof}
Choose the least $m$ with $\delta x^m\ge1/2$. Since $\delta x^2<q^{-1}<1/2$,
$m\ge3$, and minimality gives $\delta x^{m-1}<1/2$. If $m=3$, then
$\delta\ge1/(2x^3)\ge1/(2q^3)$, and $x=q(1-\delta/q)$ gives
$\delta\ge\frac{(1-\delta/q)^{-3}}{2q^3}\ge\frac{1+3\delta/q}{2q^3}
\ge\frac1{2q^3}\bigl(1+\frac3{2q^4}\bigr)$, which is the first alternative. Otherwise $x\in\mathcal B_{q,m}$ with $m\ge4$.
Then $s\in[0,1)$ and $t=x(1-s)-1\in[0,x-1]$. Finally $x\ge q-q^{-3}\ge\bar q$,
and $\delta x^{m-1}\le1/2$ gives $\delta\le D_m$.
\end{proof}

\begin{lemma}[Exact outer trace formulae]\label{lem:terminal-outer-formulas}

For the words above,
\begin{align}
 V_{\mathsf W}={}&\frac b{r^2}-1+\frac{z_W}{r}
 +\frac{\delta x^{\ell+1}}2
 \Bigl(\frac{1+s}{r^2}+\frac{j-(2m-j-2)s}{r}\Bigr),\label{eq:terminal-WV}\\
 S_{\mathsf W}={}&\frac{(2-x)b}{r^2}-2+\frac{z_W}{r}
 +\frac{\delta x^{\ell+1}}2
 \Bigl(\frac{1+s}{r^2}+\frac{j-1-(2m-j-1)s}{r}\Bigr),\label{eq:terminal-WS}\\
 V_{\mathsf U}={}&\frac b{r^2}-\frac{z_U}{r}
 +\frac{\delta x^\ell}2
 \Bigl(\frac{1-t}{r^2}+\frac{j+(2m-j)t}{r}\Bigr),\label{eq:terminal-UV}\\
 S_{\mathsf U}={}&\frac{(2-x)b}{r^2}-\frac{z_U}{r}
 +\frac{\delta x^\ell}2
 \Bigl(\frac{1-t}{r^2}+\frac{j-1+(2m+1-j)t}{r}\Bigr).\label{eq:terminal-US}
\end{align}
The same identities hold at the fixed offsets $j=0,-1$ whenever
the corresponding word is defined.
\end{lemma}
\begin{proof}
Use \eqref{eq:Uml-V}--\eqref{eq:Uml-S} for the block $(-b)^\ell$,
\eqref{eq:VS-rec-real} for the final digit of $\mathsf W$, and the formula
$u_m=1+a/r^2-2xy/r-\delta y/r^2+(m-1)\delta y/r$ with $y=\delta x^{m-1}$,
which follows from \eqref{eq:jordan} along $\mathsf B_m$ as \eqref{eq:c-m}
does: explicitly,
$V_{\mathsf B_m}=1+\frac{a^2}{r^2}+x^{m-1}\bigl(1-2\delta-\frac{a^2}{r^2}+(m-1)\frac{\delta^2}r\bigr)$,
where $1-2\delta-\frac{a^2}{r^2}=-\delta\bigl(\frac{2x}r+\frac\delta{r^2}\bigr)$
and $a^2-b=a$. Substituting $s$ and $t$ gives the four identities.
\end{proof}
Since $\delta x^{m-1}=\frac{1-s}2$ and $\delta x^m=\frac{1+t}2$, the prefactors
are $\delta x^{\ell+1}=\frac{1-s}{2x^j}$ and $\delta x^\ell=\frac{1+t}{2x^j}$.
Hence the last line of \eqref{eq:terminal-WV} is
$L_j(x)-\frac{s}{4x^j}\bigl(\frac s{r^2}+\frac{2m-2-(2m-j-2)s}{r}\bigr)$
and that of \eqref{eq:terminal-UV} is
$L_j(x)+\frac{t}{4x^j}\bigl(\frac{2m+(2m-j)t}{r}-\frac t{r^2}\bigr)$, where
$L_j(x)=\frac1{4x^j}\bigl(\frac1{r^2}+\frac jr\bigr)$; for $S$ replace the
numerator $j$ of $L_j$ by $j-1$ (keeping the factor $x^{-j}$), and $2m-j-2$,
$2m-j$ by $2m-j-1$, $2m+1-j$.

\paragraph{Limiting states.} Dropping the last lines and putting $x=q$ gives
the affine states
\[
 Z_W^0(z)=\Bigl(\frac{z-1}a,\frac{z-(a^2+1)}a\Bigr),\qquad
 Z_U^0(z)=\Bigl(\frac{q-2-z}a,-\frac{(q-2)^2+z}a\Bigr).
\]
Since $a\sigma_*=a^2+\frac a{2q}$, their capture intervals are
\begin{equation}\label{eq:terminal-limit-intervals}
 \lambda_W^-<z<\lambda_W^+\ \ (\mathsf W),\qquad
 \lambda_U^-<z<\lambda_U^+\ \ (\mathsf U),\qquad
 \lambda_W^-=\frac{q+1}{2q},\ \ \lambda_W^+=1+a\tau_*,\ \
 \lambda_U^\mp=q-2\mp a\tau_*.
\end{equation}
At $x=q$ the zero-phase layer terms are $\alpha_j=L_j(q)=\frac{ja+1}{4a^2q^j}$
($V$) and $\beta_j=\frac{(j-1)a+1}{4a^2q^j}$ ($S$). The corrected states
$\overline Z_{W,j}=Z_W^0+(\alpha_j,\beta_j)$ and
$\overline Z_{U,j}=Z_U^0+(\alpha_j,\beta_j)$ have the capture intervals
$\overline I_{W,j}=(\lambda_W^--a\beta_j,\lambda_W^+-a\alpha_j)$ and, for
$j\ge0$, $\overline I_{U,j}=(\lambda_U^-+a\alpha_j,\lambda_U^++a\alpha_j)$; for
$j=-1$ the capture interval of $\overline Z_{U,-1}$ lies in
$\overline I_{U,-1}$, the $S$-bound being binding above when $q\ge9$. All these
states have $z$-slope $\pm1/a$ in both coordinates, so a state error
$\varepsilon$ is absorbed by shrinking a capture interval by $a\varepsilon$.

\begin{proposition}[Outer-annulus cover]\label{prop:terminal-outer-annuli}

Suppose either $q\ge8$, $m\ge4$, or $q\in\{4,5,6,7\}$, $m\ge11$.
Then all variable $\mathsf W$- and $\mathsf U$-layers in
$\mathcal B_{q,m}$ form a linked strict trace chain between
their fixed handoff cells.
The final $\mathsf W$-layer and first $\mathsf U$-layer overlap under the
transition \eqref{eq:terminal-transition}.
\end{proposition}
\begin{proof}
\textit{Idea.} At $x=q$ the layer states are affine in $z$, and consecutive
layers differ by the factor $x$ in $z$; the finite-$x$ errors are smaller
than the overlaps.

Let $q\ge8$. On \eqref{eq:terminal-limit-intervals}, $z_W<\lambda_W^+<q$ and
$z_U<\lambda_U^+<2q$, so $s/x^j=z_Wx^{-m}<2qD_m/x$,
$t/x^j=z_Ux^{-m-1}<4qD_m/x^2$ and $t\le z_U/x<3$. With
$|b(r^{-2}-a^{-2})|\le2a\delta/r^2$, $b(r^{-1}-a^{-1})\le a\delta/r$,
$z(r^{-1}-a^{-1})=z\delta/(ar)$, $0\le2m-2-(2m-j-2)s\le2m$, and $j/x^j\le1/x$,
$(j-1)/x^j\le1/x^2$ for $j\ge1$, the four coordinate errors of
$Z_{\mathsf W}-Z_W^0(z_W)$ and $Z_{\mathsf U}-Z_U^0(z_U)$ are at most
\begin{align*}
 \mathcal E_{W,V}&=D_m\Bigl(\frac{2a}{r^2}+\frac q{ar}\Bigr)
 +\frac1{4x}\Bigl(\frac1{r^2}+\frac1r\Bigr)+\frac{qD_m}{2x}\Bigl(\frac1{r^2}+\frac{2m}r\Bigr),\\
 \mathcal E_{W,S}&=D_m\Bigl(\frac{2a}{r^2}+\frac ar+\frac q{ar}\Bigr)
 +\frac14\Bigl(\frac1{xr^2}+\frac1{x^2r}\Bigr)+\frac{qD_m}{2x}\Bigl(\frac1{r^2}+\frac{2m}r\Bigr),\\
 \mathcal E_{U,V}&=D_m\Bigl(\frac{2a}{r^2}+\frac{2q}{ar}\Bigr)
 +\frac1{4x}\Bigl(\frac1{r^2}+\frac1r\Bigr)+\frac{qD_m}{x^2}\Bigl(\frac3{r^2}+\frac{8m}r\Bigr),\\
 \mathcal E_{U,S}&=D_m\Bigl(\frac{2a}{r^2}+\frac ar+\frac{2q}{ar}\Bigr)
 +\frac14\Bigl(\frac1{xr^2}+\frac1{x^2r}\Bigr)+\frac{qD_m}{x^2}\Bigl(\frac3{r^2}+\frac{8m}r\Bigr).
\end{align*}
Each term decreases in $x$ and, through $D_m$ and $mD_m$, in $m$; so put
$x=\bar q$, $r=q-\frac32$ and $m=4$. Then each term of $q\mathcal E$ is a
product of the factors $q/\bar q$, $a/r$, $q/r$, $1/a$, $1/r$ and
$1/\bar q$, all decreasing in $q$. At $q=8$ the four values of $\mathcal E$
are $\frac{437863}{59889375}$, $\frac{248128}{59889375}$,
$\frac{6799049}{898340625}$ and $\frac{3953024}{898340625}$, all below
$\frac1{40}$, so $q\mathcal E\le8\,\mathcal E|_{q=8}<\frac15$. Hence
\begin{equation}\label{eq:terminal-uniform-error}
 \|Z_{\mathsf W}-Z_W^0(z_W)\|_\infty<\frac1{5q},\qquad
 \|Z_{\mathsf U}-Z_U^0(z_U)\|_\infty<\frac1{5q},
\end{equation}
and the intervals $I_W=(\lambda_W^-+\frac a{5q},\lambda_W^+-\frac a{5q})$ and
$I_U=(\lambda_U^-+\frac a{5q},\lambda_U^+-\frac a{5q})$ suffice. Layer $j$
covers $-T_m\in x^jI_W$, resp.\ $T_{m+1}\in x^jI_U$; since $x<q$, successive
layers overlap by at least
\begin{align}
 \Bigl(\lambda_W^+-\frac a{5q}\Bigr)-q\Bigl(\lambda_W^-+\frac a{5q}\Bigr)
 &=\frac{a(3q^2-22q+5)}{10q^2}>0,\label{eq:terminal-W-overlap}\\
 \Bigl(\lambda_U^+-\frac a{5q}\Bigr)-q\Bigl(\lambda_U^-+\frac a{5q}\Bigr)
 &=\frac{a(8q^2-17q+5)}{10q^2}>0.\label{eq:terminal-U-overlap}
\end{align}
At the middle seam ($\mathsf A_m$, $j=m-1$, to $\mathsf B_{m+1}$, $j=m$),
\eqref{eq:terminal-transition} reads $z_U=xr-xz_W$, and the overlap is at
least $\lambda_U^+-\frac a{5q}+x(\lambda_W^+-\frac a{5q})-xr$. This decreases
in $x$ and equals \eqref{eq:terminal-U-overlap} at $x=q$, because
$\lambda_W^++\lambda_U^-=a$.

Now let $4\le q\le7$ and $m\ge11$, and compare with the corrected states. On
$\overline I_{W,j}$ and $\overline I_{U,j}$ ($j\ge1$),
\begin{equation}\label{eq:terminal-lowq-error}
 \|Z_{\mathsf W}-\overline Z_{W,j}(z_W)\|_\infty<2mD_m,\qquad
 \|Z_{\mathsf U}-\overline Z_{U,j}(z_U)\|_\infty<2mD_m.
\end{equation}
Indeed, the fixed rational part costs at most
$(\frac{2a}{r^2}+\frac ar+\frac{2q}{ar})\delta\le\frac{242}{75}\delta<4D_m$
(each summand decreases in $q$, and $\frac{242}{75}$ is the value at $q=4$);
moving $L_j$ from $x$ to $q$ costs less than $D_m$, as
$|L_j'|\le x^{-j}(j^2+2j+2)/(4r)\le\frac17$ for $x\ge\frac72$; and the remainders are at most
$\frac{qD_m}{2x}(\frac1{r^2}+\frac{2m}r)<mD_m$ ($\mathsf W$) and, as
$t\le z_U/x<(\lambda_U^++a\alpha_1)/\bar q\le\frac{11041}{7644}<\frac32$,
at most $\frac{qD_m}{x^2}(\frac3{2r^2}+\frac{5m}r)<mD_m$ ($\mathsf U$), both
largest at $q=4$. The sum is $(m+5)D_m<2mD_m$. With $\varepsilon_m=2mD_m$
the shrunken intervals are
$(\lambda_W^--a\beta_j+a\varepsilon_m,\lambda_W^+-a\alpha_j-a\varepsilon_m)$
and $(\lambda_U^-+a\alpha_j+a\varepsilon_m,\lambda_U^++a\alpha_j-a\varepsilon_m)$.
By $\alpha_j=q\beta_{j+1}$ and $\alpha_j-q\alpha_{j+1}=-1/(4aq^j)$, the
successive overlaps and the middle seam are at least
\[
 \lambda_W^+-q\lambda_W^--a(q+1)\varepsilon_m,\qquad
 \lambda_U^+-q\lambda_U^--\frac1{4q^j}-a(q+1)\varepsilon_m,\qquad
 \lambda_U^+-q\lambda_U^--aq\alpha_{m-1}-a(q+1)\varepsilon_m.
\]
Here $\lambda_W^+-q\lambda_W^-=\frac{q-5}2+\frac5{2q}-\frac1{2q^2}\ge\frac3{32}$
and $\lambda_U^+-q\lambda_U^-=a(1-\frac3{2q}+\frac1{2q^2})\ge\frac{63}{32}$
increase in $q$, while $\frac1{4q^j}\le\frac1{16}$, $aq\alpha_{m-1}<10^{-3}$
and $a(q+1)\varepsilon_m\le\frac{168960}{282475249}<10^{-3}$, the value at
$q=4$, $m=11$. (For the seam, as for $q\ge8$, the bound decreases in $x$; we
put $x=q$ and drop the term $a\alpha_m>0$.) All three bounds are positive.
\end{proof}

\paragraph{Fixed cells.} With $\mathsf P_m=\mathsf U_{m,m}$, the four fixed
cells are $\mathsf W_{m,m-2}$ ($j=0$, $z=-T_m$), $\mathsf W_{m,m-1}$ ($j=-1$,
$z=-xT_m$), $\mathsf U_{m+1,m}$ ($j=0$, $z=T_{m+1}$) and $\mathsf P_{m+1}$
($j=-1$, $z=xT_{m+1}$). The limiting $T$-slope of offset $j$ has
absolute value $1/(aq^j)$.

\begin{lemma}[Fixed-offset estimate]\label{lem:terminal-fixed-offset}

If $q\ge8$ and $m\ge4$, restrict each of the four fixed cells listed
above to the following domain, which contains its corrected limiting capture
interval: for $j\in\{-1,0\}$, require
\[
 \lambda_W^--a\beta_j<z_W<\lambda_W^+-a\alpha_j
 \quad\text{or}\quad
 \lambda_U^-+a\alpha_j<z_U<\lambda_U^++a\alpha_j,
\]
according to the side, $\mathsf W$ or $\mathsf U$. On these domains, in particular,
$z_W<5q/4$ and $z_U<2q$. At the same $T$, the finite state then
differs from its corrected limiting state by less than $mD_m$ in
sup norm.
\end{lemma}
\begin{proof}
\textit{Idea.} As in \Cref{prop:terminal-outer-annuli}, but the corrected
limits absorb the layer terms, so only $O(m\delta)$ remains.
Put $i=-j\in\{0,1\}$. The bounds on $z$ follow from the upper endpoints, for
instance $\frac{5q}4-(\lambda_W^+-a\alpha_{-1})=\frac{q^3+8a^2q-2a^2}{4aq^2}$.
They give $0\le s,t\le\frac12$, since the corresponding powers of $x\ge\bar q$
have exponent at least four. In \eqref{eq:terminal-WV}--\eqref{eq:terminal-US}
the rational part differs from its value at $q$ by less than $2\delta$, its
larger coefficient being $\frac{2a}{r^2}+\frac ar+\frac{2q}{ar}<2$. At zero
phase the layer terms are $L_i^V=\frac1{4r^2}-\frac i4$ and
$L_i^S=\frac1{4r^2}-\frac1{4r}-\frac i2$, with derivatives $-\frac1{2r^3}$ and
$\frac{r-2}{4r^3}$, so their drift from $x$ to $q$ is at most $\delta/16$.
For the remainders put $u=(1-s)x^i/4$ ($\mathsf W$) or $u=(1+t)x^i/4$
($\mathsf U$), so that $us=\delta z_W/(2x)$ and $ut=\delta z_U/(2x)$. The
remainders then equal
$-us\bigl(\frac{s}{(1-s)r^2}+\frac{2m-2}r-\frac{is}{(1-s)r}\bigr)$ and
$ut\bigl(\frac{2m}r+\frac{it}{(1+t)r}-\frac{t}{(1+t)r^2}\bigr)$ for the $V$
coordinates; for the $S$ coordinates the middle terms carry $i+1$ in place of
$i$. By $s,t\le\frac12$ all of them are at most
$\frac{5q\delta}{8x}(\frac1{r^2}+\frac{2m-2}r)$ and
$\frac{q\delta}x(\frac1{3r^2}+\frac{2m+2/3}r)$, and with $q/x\le\frac{16}{15}$
and $1/r\le\frac2{13}$ each is at most $(\frac m3+\frac18)\delta$. Finally,
comparison at the same $T$ rather than the same $z$ changes $z$ by $\delta|T|$
when $i=1$ and costs at most $\delta|T|/a<\delta/3$; it vanishes when $i=0$.
Thus each coordinate error is less than
\[
 \Bigl(2+\frac1{16}+\frac m3+\frac18+\frac13\Bigr)\delta
 =\Bigl(\frac m3+\frac{121}{48}\Bigr)\delta<mD_m.\qedhere
\]
\end{proof}

\paragraph{Central words.} The central cells use $\mathsf P_m$ and the returns
\[
 \mathsf R_{m,\ell,k}=\mathsf B_m(-b)^\ell\bigl(-(b-1),-h,-q\bigr)(h)^{m-\ell-3}(h-k),
 \qquad \Phi_k=\mathsf R_{m,0,k},
\]
on a central window $-\vartheta\le T\le0$ with $T=T_m$ and
$\vartheta\in\{\frac12,\frac23\}$. Put $y=\delta x^{m-1}=\frac12(1+T/x^m)$ and
$\kappa=y-\frac12=T/(2x^m)$. For $\ell\in\{0,1\}$,
\eqref{eq:VS-rec-real} for the digits $-(b-1),-h,-q$ after
$\mathsf U_{m,\ell}$, and \eqref{eq:jordan} for the run $(h)^{m-\ell-3}$
about $z_h$ (the last digit adds $(k,k)$ to $(P,V)$), give the closed formulae
\begin{equation}\label{eq:terminal-R-closed}
\begin{aligned}
 V_{\mathsf R_{m,\ell,k}}&=k-\frac{2\ell y(x+1)}x
 +\frac{-h-rx^2T+yg_0+y^2g_1+\delta yf_\ell}{r^2},\\
 S_{\mathsf R_{m,\ell,k}}&=V_{\mathsf R_{m,\ell,k}}+\frac hr+k
 -\frac{\delta y(x+1)}{x^{\ell+1}r}+\frac{x^2y(y-1)}r,
\end{aligned}
\end{equation}
where $g_0=-2+(m+4)x^2-(m+1)x^3$, $g_1=x^2(2mr-1)$ and
$f_\ell=\bigl((m-1-\ell)(1-x^2)+2x\bigr)/x^{\ell+1}$. At $x=q$, $y=\frac12$,
the limiting states are
\begin{equation}\label{eq:terminal-central-limits}
\begin{aligned}
 Z^0_{\mathsf P}(T)&=\Bigl(\frac{(q-2)(3q-4)}{4a^2},\
 -\frac{(q-2)^2}a-\frac{q(2q-3)}{4a^2}\Bigr)-\frac qaT\,(1,1),\\
 Z^0_{\ell,k}(T)&=\bigl(\varphi_{\ell,k},\ \varphi_{\ell,k}+q+k-\tfrac{q^2}{4a}\bigr)
 -\frac{q^2}aT\,(1,1),\qquad
 \varphi_{\ell,k}=k-\ell\Bigl(1+\frac1q\Bigr)+\frac1a-\frac{q^2(2q-3)}{4a^2}.
\end{aligned}
\end{equation}
Appending a digit $d$ acts by $Z\mapsto xJZ-d\,(1,2)$, with $J$ as in the proof
of \Cref{prop:branched}; since $J(1,1)=(1,1)$, each appended digit multiplies
the limiting $T$-slopes $q/a$ ($\mathsf P$) and $q^2/a$ ($\mathsf R$) by $q$.

\Needspace{6\baselineskip}
\begin{lemma}[Base errors]\label{lem:terminal-base-errors}

Let $m\ge4$, $\ell\in\{0,1\}$, $0\le k\le q/4$ and $x\in\mathcal B_{q,m}$
with $-\vartheta\le T_m\le0$. The limiting states
\eqref{eq:terminal-central-limits} have sup norm at most $2q$, and:
\begin{enumerate}
\item $\|Z_{\mathsf R_{m,\ell,k}}-Z^0_{\ell,k}(T_m)\|_\infty<(m+3)\delta$ if
 $\vartheta=\frac12$ or $q\ge8$, and $<(2m+8)\delta$ if $\vartheta=\frac23$;
\item $\|Z_{\mathsf P_m}-Z^0_{\mathsf P}(T_m)\|_\infty<(\frac m5+\frac65)\delta$
 if $\vartheta=\frac12$, and $<2D_m$ if moreover $q\ge8$.
\end{enumerate}
\end{lemma}
\begin{proof}
\textit{Idea.} After freezing $y=\frac12$, the closed formulae are smooth in
$x$; the rest is of order $\delta$ with a coefficient linear in $m$, because
$|\kappa|=|T|\delta/(2xy)$.

Since $2y\ge1-\frac23x^{-4}>\frac{99}{100}$, we have
$|\kappa|\le\vartheta'\delta/x$ with $\vartheta'=\frac{100}{99}\vartheta$.
The coefficient bounds below decrease in $x$, and then in $q$; we evaluate
them at $x=\bar q$, $r=q-\frac32$ and $q=4$ or $q=8$.

(i) Write $yg_0+y^2g_1=(\frac{g_0}2+\frac{g_1}4)+\kappa H$, where
$H=g_0+(y+\frac12)g_1$ and $\frac{g_0}2+\frac{g_1}4=-1+\frac74x^2-\frac12x^3$.
Then $V_{\mathsf R}=F_\ell(x)+\Delta$, where
$F_\ell(x)=k-\ell\frac{x+1}x+\frac{-h-1+\frac74x^2-\frac12x^3-rx^2T}{r^2}$
satisfies $F_\ell(q)=\varphi_{\ell,k}-\frac{q^2}aT$, and
$\Delta=\frac{\delta yf_\ell+\kappa H}{r^2}-2\ell\kappa\frac{x+1}x$. For
$t'\in[x,q]$ and $\rho'=t'-1$,
\[
 \partial_{t'}F_\ell=\frac{\mathcal N}{2\rho'^3}+\frac\ell{t'^2},\qquad
 \mathcal N=-\rho'^3+4\rho'^2-1+(q-t')\bigl(8\rho'+4+4(q-t')\bigr)+2|T|(\rho'^3-\rho'),
\]
and $q-t'\le D_4\le\frac4{343}$ gives
$-\rho'^3<\mathcal N<(2\vartheta-1)\rho'^3+4\rho'^2-(2\vartheta-\frac15)\rho'$;
at $\rho'=\frac52$, $t'=\frac72$ this yields $|\partial_{t'}F_\ell|<1$. Next,
$0<-x^{\ell+1}f_\ell<(m-1-\ell)(x^2-1)$ for $m\ge4$, and
$0<H<x^2((m-1)r+2)$, because $g_1>0$, $y+\frac12\in[\frac78,1]$ and
$g_0+\frac78g_1=x^2\bigl((\frac{3m}4-1)x+\frac{25}8-\frac{3m}4\bigr)-2>0$.
Hence
\[
 \frac{|V_{\mathsf R}-V^0_{\ell,k}|}\delta<1+(m-1)\mu
 +2\vartheta'\Bigl(\frac x{r^2}+\frac{x+1}{x^2}\Bigr),
 \qquad \mu=\frac{x+1}{2xr}+\vartheta'\frac xr.
\]
By \eqref{eq:terminal-R-closed}, $S_{\mathsf R}-S^0_{\ell,k}$ adds to this
$h(\frac1r-\frac1a)=\frac{q\delta}r$, at most $\frac{(x+1)\delta}{2xr}$, and
$\frac{x^2y(y-1)}r+\frac{q^2}{4a}=\frac\delta4(1-\frac1{ar})+\frac{x^2\kappa^2}r
<(\frac14+\frac1{100})\delta$. Writing the total as $(m-1)\mu+\nu$, in the
three cases $(\vartheta=\frac12,q\ge4)$, $(\vartheta=\frac23,q\ge4)$ and
$(\vartheta=\frac23,q\ge8)$ the values of $(\mu,\nu)$ are
$(\frac{3341}{3465},\frac{983263}{242550})$,
$(\frac{12473}{10395},\frac{3176989}{727650})$ and
$(\frac{5561}{6435},\frac{6202231}{2053350})$, below
$(0.9643,4.054)$, $(1.2,4.367)$ and $(0.8642,3.021)$. So the error is below
$(0.9643m+3.09)\delta<(m+3)\delta$, $(1.2m+3.17)\delta<(2m+8)\delta$ and
$(0.8642m+2.16)\delta<(m+3)\delta$.

(ii) By \eqref{eq:Uml-V}--\eqref{eq:Uml-S} with $\ell=m$, \eqref{eq:c-m},
the formula for $u_m$ in the proof of \Cref{lem:terminal-outer-formulas}, and
$x^{1-m}=\delta/y$,
\[
\begin{aligned}
 V_{\mathsf P}-V^0_{\mathsf P}&=(b+\tfrac14)\bigl(\tfrac1{r^2}-\tfrac1{a^2}\bigr)
 +\tfrac{T(2-x)\delta}{4r^2y}+T\delta\bigl(\tfrac{m-1/2}r-\tfrac1{ar}-\tfrac1{2r^2}\bigr),\\
 S_{\mathsf P}-S^0_{\mathsf P}&=(b+\tfrac14)\bigl[\bigl(\tfrac1{r^2}-\tfrac1{a^2}\bigr)-\bigl(\tfrac1r-\tfrac1a\bigr)\bigr]
 +\tfrac{T(3-2x)\delta}{4r^2y}+T\delta\bigl(\tfrac mr-\tfrac1{ar}-\tfrac1{2r^2}\bigr).
\end{aligned}
\]
As $b+\frac14=(a-\frac12)^2\le r^2$, the first terms are at most
$\frac{2\delta}a$ and $\max(\frac2a,\frac ra)\delta\le\delta$; for
$|T|\le\frac12$ the second terms are at most $0.07\delta$ and $0.17\delta$, and
the third at most $\frac{m\delta}{2r}\le\frac{m\delta}5$. This gives the first
bound. If $q\ge8$, then $\delta=yx^{1-m}\le D_m\theta^{m-1}$ with
$\theta=\bar q/x\in[\frac{15}{16},1-\frac1{3q}]$, since
$x-\bar q\ge\frac12-D_4>\frac13$. As $\max_n n\theta^n\le\frac1{e(1-\theta)}$,
we get $m\theta^{m-1}\le\frac{16}{15}\cdot\frac{3q}e$, and the third terms are
at most $\frac{16}{15}\cdot\frac{3q}{2er}D_m<0.725D_m$, while the others are
at most $1.08D_m$. So both errors are below $2D_m$.

Finally, the limiting coordinates \eqref{eq:terminal-central-limits} are
affine in $k$ and $T$, so it suffices to bound them at the corners of
$0\le k\le\frac q4$, $\ell\in\{0,1\}$, $-\frac23\le T\le0$. The largest is the
$S$ coordinate of $Z^0_{0,q/4}(-\frac23)$, which equals
$\frac{17}{12}q+\frac16+\frac{17}{12a}+\frac1{4a^2}$, as $\frac{q^2}a=q+1+\frac1a$
and $\frac{2q-3}{4a}=\frac12-\frac1{4a}$; it is below $2q$ by
$\frac7{12}q-\frac16-\frac{17}{12a}-\frac1{4a^2}\ge\frac53$. At each of the
other corners, $2q\mp$ (coordinate) is likewise a rational function of $q$ with
nonnegative coefficients after $q=4+s$.
\end{proof}

\paragraph{Suffixes and endpoint movements.} Let $Z$ and $Z^0$ be a finite
state and its limit at the same $T$, with $\|Z-Z^0\|_\infty\le E_0$ and
$\|Z^0\|_\infty\le G_0$. Apply \Cref{lem:word-perturbation} with the
matrix $J$ of \Cref{prop:branched} and $e=(1,2)$. Since
$\|J\|_\infty=3$, $\|J^2\|_\infty=5$ and $\|Je\|_\infty=3$, one and two
appended digits give, respectively,
\begin{equation}\label{eq:terminal-suffix}
 E_{i+1}\le3qE_i+3\delta G_i,\qquad G_{i+1}\le3qG_i+2h,
\end{equation}
$E_2\le5q^2E_0+\delta(10qG_0+3h)$. Dividing by the $T$-slope gives the
normalized endpoint movement. With the base errors of
\Cref{lem:terminal-fixed-offset,lem:terminal-base-errors}, written as
$E_0\le cD_m$, and with $a<q$, $G_0=4q$, $G_1\le14q^2$, $G_2\le44q^3$ (or
$G_0=2q$, $G_1\le8q^2$ where marked $^\circ$), the movements of cells based on
the listed words are as follows. The fixed cells and the returns carry no
suffix; the central cells carry suffixes of at most two digits, or three
where stated. For $4\le q\le7$ the two $\mathsf R$ rows give the allowances
inserted in the record: on $-\frac23\le T\le0$ the base error proved for these
returns is $(2m+8)\delta$ (row $\mathsf R^\circ$); its movement
$(10m+63)D_m/q$ decreases in $m$ and is at most $173D_{11}/q$, below the values
$\eta_{R_1}(11)=204D_{11}/q$ and $\eta_{R_0}(11)=285D_{11}/q$ that the record
inserts.
\begin{equation}\label{eq:terminal-endpoint-movements}
\renewcommand{\arraystretch}{1.1}
\begin{array}{@{}l@{\ \ }l@{\ \ }l@{\quad}l@{\ \ }l@{\ \ }l@{}}
\multicolumn{3}{@{}l}{q\ge8} & \multicolumn{3}{l@{}}{4\le q\le7,\ m\ge11}\\\hline
\text{cells} & c & \text{movement} & \text{cells} & c & \text{movement}\\
\mathsf W_{m,m-2},\mathsf U_{m+1,m} & m & maD_m &
 \mathsf R_{m,1,0} & m+3 & \eta_{R_1}=(9m+105)D_m/q\\
\mathsf W_{m,m-1},\mathsf P_{m+1} & m & maD_m/q &
 \mathsf R_{m,0,k} & m+12 & \eta_{R_0}=(9m+186)D_m/q\\
\mathsf R_{m,1,0},\ T\ge-\frac23 & m+3 & (m+11)D_m/q &
 \mathsf P & m+10 & \eta_{P,\le2}=(9m+168)D_m\\
\mathsf R_{m,0,0},\ T\ge-\frac12 & m+3 & (m+12)D_m/q &
 \mathsf P,\ \text{three digits} & m+10 & \eta_{P,3}=(27m+636)D_m\\
\mathsf P,\ m\ge11 & m+12 & \eta_P=(9m+186)D_m &
 \mathsf R^\circ,\ T\ge-\frac23 & 2m+8 & (10m+63)D_m/q\\
\Phi,\ m\ge11 & m+12 & \eta_\Phi=\eta_P/q &&&\\
\mathsf P^\circ,\ m\le10 & 2 & 60D_m &&&\\
\Phi^\circ,\ m\le10 & m+10 & (9m+132)D_m/q &&&
\end{array}
\end{equation}
For instance, $E_0=(m+12)D_m$ and $G_0=4q$ give $E_2\le q^2(9m+186)D_m$, and
the two-digit bound with $E_0=(2m+8)D_m$, $G_0=2q$ gives the last row. Each allowance $(\alpha m+\beta)D_m$
in the table decreases for $m\ge11$, since consecutive values have ratio
$(1+\frac\alpha{\alpha m+\beta})/\bar q<1$; thus $\eta_P\le285D_{11}$ and
$\eta_\Phi\le285D_{11}/q$. The return estimates come from the direct closed
formulae \eqref{eq:terminal-R-closed}, not from removing digits from another
state.

\paragraph{The central hierarchy.} If $q=4K+\nu$ with $0\le\nu\le3$, put
$d=(-4q^2+11q-8+\nu)/4$, and write $\Phi_k(\mathbf e)$ and $\mathsf P(\mathbf e)$
for appended suffixes, with $\mathsf P=\mathsf P_m$. After the consecutive
cells $\Phi_0,\ldots,\Phi_{K-1}$ use the residue chains
\[
\begin{aligned}
 \nu=0:&\ \ \Phi_{K-1}(h-q-K+2)\leadsto\Phi_{K-1}(h-q-K+1)\leadsto\mathsf P(d,d)\leadsto\mathsf P(d)\leadsto\mathsf P;\\
 \nu=1:&\ \ \Phi_{K-1}(h-6K)\leadsto\mathsf P(d,-b-2K)\leadsto\mathsf P(d,-b-2K-1)\leadsto\mathsf P(d)\leadsto\mathsf P;\\
 \nu=2:&\ \ \Phi_{K-1}(h-q-3K)\leadsto\mathsf P(d,-b-3K-1)\leadsto\mathsf P(d,-b-3K-2)\leadsto\mathsf P(d)\leadsto\mathsf P(d-1)\leadsto\mathsf P;\\
 \nu=3:&\ \ \Phi_{K-1}(b)\leadsto\Phi_{K-1}(b-1)\leadsto\Phi_K\leadsto\Phi_K(q(q-2)+2)\leadsto\Phi_K(q(q-2)+1)\\
 &\ \ \leadsto\mathsf P(d-1,-b)\leadsto\mathsf P(d-1)\leadsto\mathsf P.
\end{aligned}
\]

\begin{proposition}[Central hierarchy and fixed seams\certmark]\label{prop:terminal-central-hierarchy}

The ordered central hierarchy, together with
\[
 \mathsf W_{m,m-2},\mathsf W_{m,m-1},\mathsf R_{m,1,0},\mathsf R_{m,0,0}
 \quad\text{and}\quad
 \mathsf P_{m+1},\mathsf U_{m+1,m},
\]
has consecutive strict overlaps in each of the following ranges:
\[
\begin{array}{c|ccccc}
m&m\ge11&4&5&6&7,8,9,10\\ \hline
q&q\ge8&q\ge275&q\ge23&q\ge11&q\ge8.
\end{array}
\]
For $q=4,5,6,7$ and $m\ge11$, explicit replacement chains have
strict overlaps. Their ordered words and rational limiting endpoints
are listed in \Cref{tab:low-q-chains}; their
finite-$m$ handoffs are checked in exact arithmetic. In every listed range,
the first central cell strictly contains $T=-1/2$ and the final
$\mathsf P$-based cell strictly contains $T=0$.
All overlaps remain strict after the finite-$x$ endpoint movements.
\end{proposition}

\begin{table}[tbp]
\centering\scriptsize
\setlength{\tabcolsep}{4pt}
\caption{Replacement chains for $q=4,5,6,7$ and $m\ge11$ in
\Cref{prop:terminal-central-hierarchy}, with $m$ omitted:
$\mathsf R_{\ell,k}=\mathsf R_{m,\ell,k}$ and $\mathsf P=\mathsf P_m$, followed by
the displayed suffix digits. The endpoints bound the open interval of $T$ on
which the limiting state satisfies $|V|<\tau_*$ and $|S|<\sigma_*$. In each
chain the first interval contains $T=-1/2$, the last contains $T=0$, and
consecutive intervals overlap, by at least $3/8192$, $24/15625$, $65/46656$
and $103/67228$ for $q=4,5,6,7$; all suffix digits lie in $B_h$.}
\label{tab:low-q-chains}
\begin{tabular}{@{}lrr@{}}
\toprule
cell & lower & upper\\
\midrule
\multicolumn{3}{@{}l}{$q=4$}\\
$\mathsf R_{1,0}$ & $-259/384$ & $-751/1536$\\
$\mathsf R_{1,0}(9)$ & $-3289/6144$ & $-2983/6144$\\
$\mathsf R_{1,0}(9,5)$ & $-12217/24576$ & $-11911/24576$\\
$\mathsf R_{1,0}(8,12)$ & $-11929/24576$ & $-11623/24576$\\
$\mathsf R_{1,0}(8)$ & $-745/1536$ & $-2695/6144$\\
$\mathsf R_{0,0}$ & $-169/384$ & $-391/1536$\\
$\mathsf R_{0,0}(9)$ & $-1849/6144$ & $-1543/6144$\\
$\mathsf R_{0,0}(9,5)$ & $-6457/24576$ & $-6151/24576$\\
$\mathsf R_{0,0}(8)$ & $-385/1536$ & $-1255/6144$\\
$\mathsf P(-7,-7)$ & $-1465/6144$ & $-1159/6144$\\
$\mathsf P(-7)$ & $-313/1536$ & $-7/1536$\\
$\mathsf P$ & $-25/384$ & $281/384$\\
\midrule
\multicolumn{3}{@{}l}{$q=5$}\\
$\mathsf R_{1,0}$ & $-6887/10000$ & $-4903/10000$\\
$\mathsf R_{0,0}$ & $-4967/10000$ & $-2983/10000$\\
$\mathsf R_{0,0}(14)$ & $-15767/50000$ & $-13783/50000$\\
$\mathsf P(-13,-14)$ & $-14967/50000$ & $-2703/10000$\\
$\mathsf P(-13,-14,-20)$ & $-69367/250000$ & $-67383/250000$\\
$\mathsf P(-13,-15,-11)$ & $-67767/250000$ & $-65783/250000$\\
$\mathsf P(-13,-15)$ & $-13367/50000$ & $-11383/50000$\\
$\mathsf P(-13)$ & $-2567/10000$ & $-143/2000$\\
$\mathsf P(-14)$ & $-967/10000$ & $1017/10000$\\
\bottomrule
\end{tabular}\hspace{2.5em}
\begin{tabular}{@{}lrr@{}}
\toprule
cell & lower & upper\\
\midrule
\multicolumn{3}{@{}l}{$q=6$}\\
$\mathsf R_{1,0}$ & $-8797/12960$ & $-6347/12960$\\
$\mathsf R_{0,0}$ & $-6697/12960$ & $-4247/12960$\\
$\mathsf R_{0,0}(21)$ & $-26497/77760$ & $-24047/77760$\\
$\mathsf P(-21,-24)$ & $-25417/77760$ & $-3947/12960$\\
$\mathsf P(-21,-24,-30)$ & $-144217/466560$ & $-141767/466560$\\
$\mathsf P(-21,-25,-19)$ & $-142417/466560$ & $-139967/466560$\\
$\mathsf P(-21,-25)$ & $-23617/77760$ & $-21167/77760$\\
$\mathsf P(-21)$ & $-3817/12960$ & $-347/2160$\\
$\mathsf P(-21,-30)$ & $-14617/77760$ & $-12167/77760$\\
$\mathsf R_{0,1}(26)$ & $-2077/12960$ & $-11447/77760$\\
$\mathsf P(-22)$ & $-2017/12960$ & $433/12960$\\
\midrule
\multicolumn{3}{@{}l}{$q=7$}\\
$\mathsf R_{1,0}$ & $-38411/57624$ & $-28187/57624$\\
$\mathsf R_{0,0}$ & $-30347/57624$ & $-20123/57624$\\
$\mathsf R_{0,0}(30)$ & $-148535/403368$ & $-138311/403368$\\
$\mathsf R_{0,0}(29)$ & $-19847/57624$ & $-131255/403368$\\
$\mathsf R_{0,1}$ & $-2711/8232$ & $-13067/57624$\\
$\mathsf R_{0,1}(37)$ & $-99143/403368$ & $-88919/403368$\\
$\mathsf R_{0,1}(36)$ & $-12791/57624$ & $-81863/403368$\\
$\mathsf P(-32,-30)$ & $-87383/403368$ & $-77159/403368$\\
$\mathsf P(-32)$ & $-11531/57624$ & $-1307/57624$\\
$\mathsf P$ & $-947/8232$ & $9277/8232$\\
\bottomrule
\end{tabular}
\end{table}

\begin{proof}
\textit{Idea.} At $x=q$ every cell is an explicit rational interval in $T$,
and the finite-$x$ movements \eqref{eq:terminal-endpoint-movements} are smaller
than all overlaps.

By \eqref{eq:terminal-central-limits} every limiting trace coordinate in the
chains is affine in $T$, so the four trace inequalities of a cell give a
rational open interval. Consecutive $\Phi_k$-intervals overlap, the residue
chains close the interval from $\Phi_{K-1}$ to $\mathsf P$, and every displayed
digit lies in $B_h$. On each progression $q=q_0+4i$ the overlaps and the two
end containments, after the movements \eqref{eq:terminal-endpoint-movements}
are inserted, are bounded below by positive coefficients
(\Cref{cert:terminal-finite}\certmark). Each cell is used only where its
estimate holds: $\mathsf P$- and $\Phi$-based cells on $-\frac12\le T\le0$,
$\mathsf R_{m,1,0}$ on $-\frac23\le T\le0$; for $q\ge8$ the record shows that
$\mathsf R_{m,1,0}$ and $\mathsf R_{m,0,0}$ both strictly capture $T=-\frac12$
and that the $\mathsf W_{m,m-1}/\mathsf R_{m,1,0}$ overlap meets $T>-\frac23$.

For $q\ge8$ the variable annuli join the fixed cells directly. Put
$L_W=\lambda_W^-+\frac a{5q}$ and $L_U=\lambda_U^-+\frac a{5q}$. The $j=1$
$\mathsf W$-annular interval has boundary $T=-xL_W$, and $\mathsf W_{m,m-2}$
strictly contains the limiting target $-qL_W$ after the additional shrink
$D_m>\delta L_W$. Likewise the $j=1$ $\mathsf U$-annular interval has
boundary $T=xL_U$, and $\mathsf U_{m+1,m}$ contains $qL_U$ after the
additional shrink $qD_m>\delta L_U$. As these containments are strict, each
fixed cell overlaps the adjacent open annular interval; the record checks
them with the movements \eqref{eq:terminal-endpoint-movements}.

For $q=4,5,6,7$ and $m\ge11$ the chains of \Cref{tab:low-q-chains} replace the
central hierarchy from $\mathsf P_m$ to $\mathsf R_{m,1,0}$. Their $\mathsf R$-based cells are used on
$-\frac23\le T\le0$, with movement $(10m+63)D_m/q\le173D_{11}/q$, below
$\eta_{R_1}(11)=204D_{11}/q$ and $\eta_{R_0}(11)=285D_{11}/q$. Next to the
chain lie $\mathsf W_{m,m-3}$ and $\mathsf U_{m+1,m-1}$ ($j=1$) and the four
fixed cells. At the same $T$ their errors are below $2mD_m$ for $j=1$ (the
proof of \eqref{eq:terminal-lowq-error} gives $(m+5)D_m$, passing from $z$ to
$T$ adds less than $D_m$, and $m+6<2m$)
and below $(m+6)D_m<3mD_m$ for $j\in\{0,-1\}$, since with $q/x\le\frac87$,
$1/r\le\frac25$ and $a\ge3$ the four contributions in the proof of
\Cref{lem:terminal-fixed-offset} become $\frac{242}{75}\delta$,
$\frac\delta{16}$, $(\frac{32}{35}m+\frac{64}{175})\delta$ and
$\frac{16}{21}\delta$, whose sum $(\frac{32}{35}m+\frac{12367}{2800})\delta$ is
less than $(m+6)\delta$. Divided by the slopes $1/(aq^j)$, the movements are at
most $2maqD_m$, $3maD_m$ and $3maD_m/q$; they decrease in $m$, and at $m=11$
they are at most $2\cdot11\cdot6\cdot7\,D_{11}=924D_{11}<1000D_{11}$. With
this single allowance for the $\mathsf W$ and $\mathsf U$ cells of
$\mathsf W_{m,m-3}\leadsto\mathsf W_{m,m-2}\leadsto\mathsf W_{m,m-1}$ and
$\mathsf P_{m+1}\leadsto\mathsf U_{m+1,m}\leadsto\mathsf U_{m+1,m-1}$, the
allowance $\eta_{P,\le2}$ for $\mathsf P_{m+1}$, and the decreasing
$\mathsf R$- and $\mathsf P$-allowances of the chains at $m=11$, every
consecutive overlap of the four chains stays strict, and the
$\mathsf W_{m,m-1}/\mathsf R_{m,1,0}$ overlap meets $T>-\frac23$
(\Cref{cert:terminal-finite}\certmark).
\end{proof}

\paragraph{The shoulder and the finite range.}
\begin{proposition}[Initial shoulder\certmark]\label{prop:terminal-shoulder}

For every $q\ge4$, every point of
$\bigl[q-\frac1{q^3},q-\bigl(\frac12+\frac{3}{4q^4}\bigr)q^{-3}\bigr]$
is an interior point of $\mathcal Z_{q(q-1)}^{\mathrm{ext}}$.
\end{proposition}
\begin{proof}
\textit{Idea.} The words $\mathsf B_4(-b)^i$ follow the second channel
across the shoulder. For $q\ge12$, this rational interval is covered by five
closed cells of the phase $q^3(q-x)$ carrying the words $\mathsf B_4$,
$\mathsf B_4(-b)$, $\mathsf B_4(-b)^2$, $\mathsf B_4(-b)^3$,
$\mathsf B_4(-b)^4$, whose four strict trace margins, cleared of positive
factors, are positive for all $q\ge12$ by \Cref{sep:lem:bernstein}. For
$4\le q\le11$, exact univariate Bernstein subdivision gives rational pieces;
only $q=4$ needs the additional word $\mathsf B_4(-(a^2+1))$. Every digit lies
in $B_h$. The statement is part~(ii) of
\Cref{cert:terminal-finite}\certmark, whose certificate consists of these
checks.
\end{proof}

The low-height chains cover $4\le q\le7$ at all depths $m\ge11$. The finite
complement of these ranges and of those in
\Cref{prop:terminal-central-hierarchy} is
\begin{equation}\label{eq:terminal-finite-range}
 m=4,\ 4\le q\le274;\qquad m=5,\ 4\le q\le22;\qquad
 m=6,\ 4\le q\le10;\qquad 7\le m\le10,\ 4\le q\le7.
\end{equation}
For each of these pairs, exact endpoint phase tests enclose
$\mathcal B_{q,m}$ in a rational interval, and consecutive rational cells cover
it by words with digits in $B_h$ whose four strict trace margins are positive
by rational Bernstein certificates (\Cref{cert:terminal-finite}\certmark).

\begin{proof}[Proof of \Cref{thm:terminal-completion}]
\textit{Idea.} On each band the cells above form one linked chain from one
centre to the next. Fix $x\in[q-q^{-3},q)$. By \Cref{lem:terminal-coordinate-monotonicity},
either \Cref{prop:terminal-shoulder} applies or $x\in\mathcal B_{q,m}$ for some
$m\ge4$. If $(q,m)$ lies in \eqref{eq:terminal-finite-range}, the finite
certificate applies. Otherwise, at the fixed $x$, $T_m$ and
$T_{m+1}=x^{m+1}r+x^2T_m$ are increasing affine functions of each other, and
by \Cref{prop:terminal-outer-annuli,lem:terminal-fixed-offset,prop:terminal-central-hierarchy}
the following list, read in the opposite order, is linked after all endpoint
movements (for $4\le q\le7$ the chains of \Cref{tab:low-q-chains} replace the
central hierarchy, joined as in the proof of
\Cref{prop:terminal-central-hierarchy}): the central
hierarchy from $\mathsf P_m$ to $\mathsf R_{m,0,0}$, then $\mathsf R_{m,1,0}$,
$\mathsf W_{m,m-1}$, $\mathsf W_{m,m-2}$, the variable layers across the middle
seam \eqref{eq:terminal-transition}, $\mathsf U_{m+1,m}$ and $\mathsf P_{m+1}$.
Its end cells contain the ends $T_m=0$ and $T_{m+1}=0$ of $\mathcal B_{q,m}$:
for $\mathsf P_m$ by \Cref{prop:terminal-central-hierarchy}; and
$\mathsf P_{m+1}$, in the coordinate $T_{m+1}$, has the same limiting state
$Z^0_{\mathsf P}$ of \eqref{eq:terminal-central-limits}. This affine state
satisfies $|V|<\tau_*$ and $|S|<\sigma_*$ exactly for
$\max(\ell_V,\ell_S)<T<\min(u_V,u_S)$, where
\[
\begin{aligned}
 \ell_V&=-\frac{q^4-6q^3+14q^2-12q+2}{4q^3a}, &
 u_V&=\frac{7q^4-26q^3+30q^2-12q+2}{4q^3a},\\
 \ell_S&=-\frac{8q^4-30q^3+43q^2-24q+2}{4q^2a}, &
 u_S&=\frac{6q^3-15q^2+8q+2}{4q^2a}.
\end{aligned}
\]
The numerators of $-\frac1{16}-\ell_V$, $u_V-\frac1{16}$, $-\frac1{16}-\ell_S$
and $u_S-\frac1{16}$ over the positive denominators $16q^3a$, $16q^3a$,
$16q^2a$, $16q^2a$ are $(q-2)(3q^3-17q^2+22q-4)$, $(q-2)(27q^3-49q^2+22q-4)$,
$32q^4-121q^3+173q^2-96q+8$ and $23q^3-59q^2+32q+8$, whose coefficients in
$q-4$ are positive; the smallest margin is $\frac1{384}$, at $q=4$. So the
limiting interval contains $[-\frac1{16},\frac1{16}]$, while the movement
$maD_m/q$ (or $3maD_m/q$) is below $\frac1{16}$.
\Cref{lem:terminal-cover-principle} gives capture.
\end{proof}

\subsubsection{Finite checks delegated to the records}
The records enter this section only through their statements;
\Cref{cert:tab:records} lists their uses. For
\Cref{cert:short-word-strip} the certificate consists of the coefficient
conditions \eqref{eq:short-normalized-certificate} for the margins of
\Cref{lem:short-word-direct}; for \Cref{cert:compact-windows}, of covers of
the two windows by rational cells for finitely many $q$ and by symbolic cells
on residue classes of $q$ modulo~$4$; for \Cref{cert:terminal-finite}, of
rational cells for the pairs of \eqref{eq:terminal-finite-range}, rational and
symbolic cells for the shoulder, the coefficient positivity, on residue
classes of $q$, of the limiting overlaps and end containments of
\Cref{prop:terminal-central-hierarchy} with the allowances stated in the
record (the fixed seams are checked with the larger movements $3maD_m$ and
$3maD_m/q$), and the rational comparisons of the chains of
\Cref{tab:low-q-chains} at $m=11$. All monotonicity statements are proved above.

\subsection{Height six}\label{seven:app:terminal-proof}
We prove \Cref{seven:thm:compact,seven:thm:terminal-orbit}, the case $q=3$ of
\Cref{thm:complete-active-high}, by the scheme of \Cref{sec:pronic-completion}:
buffered returns give interior points up to $2999/1000$, and bounded orbits in
a moving rectangle give double zeros beyond. We use $F_{x,d}$, $r=x-1$,
$\delta=3-x$, $z_d$ and \eqref{eq:jordan} from \Cref{sec:residual}, with
$h=6$, and $F_{x,d+1}=F_{x,d}-(1,1)$. Words and orbits are \emph{admissible}
if their digits lie in $B_6$; $F_{x,u}=F_{x,d_m}\circ\cdots\circ F_{x,d_1}$
for $u=(d_1,\ldots,d_m)$, so $F_{x,u}(1,1)=(P_u(x),V_u(x))$.

\subsubsection{The compact interval}

\begin{lemma}[Buffered control]\label{seven:lem:control}

Let $x>1$ and $\eta>0$, and let $Y\subset\R^2$ be compact. Suppose that $Y$ is
the union of finitely many closed sets $Y_i$, each with an assigned admissible
nonempty word $u_i$ such that
\begin{equation}\label{seven:eq:return}
 F_{x,u_i}(Y_i)+[-\eta,\eta]^2\subset Y,
\end{equation}
and that some admissible word $v$ satisfies
\begin{equation}\label{seven:eq:entry}
 F_{x,v}(1,1)+[-\eta,\eta]^2\subset Y.
\end{equation}
Then $x$ is an interior point of $\mathcal Z_6^{\mathrm{ext}}$.
\end{lemma}
\begin{proof}
\textit{Idea.} The buffers absorb the change from the real map to the
remainder map at nearby complex parameters.
For $c=x'+iy$, $x'>0$, $y\ne0$, write a remainder $Z$ as $P=\re Z$,
$V=\re Z+x'\im Z/y$; then $Z\mapsto cZ-d$ becomes
\begin{equation}\label{seven:eq:complex-map}
 P^+=x'P-d+\frac{y^2}{x'}(P-V),\qquad V^+=x'V+P^+,
\end{equation}
and $Z=1$ has state $(1,1)$. This map $F_{c,d}$ is $F_{x',d}$ at $y=0$, so
$F_{c,u}(z)$ is continuous in $(c,z)$. By compactness of the $Y_i$ and
finiteness of the words, some open disc $U\ni x$ in $\{|c|>1,\ \re c>0\}$
has, for $c\in U$, $F_{c,u_i}$ within $\eta$ of $F_{x,u_i}$ on $Y_i$ and
$F_{c,v}(1,1)$ within $\eta$ of $F_{x,v}(1,1)$ in the sup norm; so
$F_{c,u_i}(Y_i)\subset Y$ and $F_{c,v}(1,1)\in Y$ by
\eqref{seven:eq:return}--\eqref{seven:eq:entry}. For $c\in U$, apply $v$,
then at each state in $Y$ the word of some $Y_i$ containing it. States ending
a word lie in $Y$, and the others are images of $(1,1)$ or of the $Y_i$ under
finitely many continuous maps; so the remainders $Z_k=P_k+i(y/x')(V_k-P_k)$
($Z_k=P_k$ if $y=0$) are bounded. As $Z_k=c^k\bigl(1-\sum_{j\le k}d_jc^{-j}
\bigr)$ and $|c|>1$, $1-\sum_{j\ge1}d_jc^{-j}=0$, and Proposition~\ref{prop:global-monic-reduction}
(digits $-d_j$) with \eqref{eq:root-ifs} gives $U\subset\mathcal
Z_6^{\mathrm{ext}}$.
\end{proof}

\begin{theorem}[Compact height-six interval\certmark]\label{seven:thm:compact}

One has $(1,2999/1000]\subset\intr\mathcal Z_6^{\mathrm{ext}}$.
\end{theorem}
\begin{proof}
\textit{Idea.} The annulus, one word and buffered control cover three ranges.
By Proposition~\ref{prop:annulus}, $(1,\sqrt7)\subset\intr\mathcal Z_6^{\mathrm{ext}}$.
On $[\sqrt7,137/50]$ the word $(5,-5)$ has $V=3x^2-10x+5$ and
$S=4x^2-15x+10$ by \eqref{eq:VS-rec-real}, increasing for $x\ge15/8$, with
endpoint values $26-10\sqrt7>-\frac12$, $\frac{307}{2500}$ and
$38-15\sqrt7>-2$, $-\frac{2674}{2500}$. As $(137/50+1)^2=34969/2500<14$,
$x^2+2x<13$, so $|V|<\frac12<\tau_6(x)$ by \eqref{eq:tau-identities} and
$|S|<2<13/(2x)=\sigma_6(x)$; \Cref{prop:strict-trace} applies. For $x\in[137/50,2999/1000]$, choose a band $I$ containing $x$.
Part~(i) of \Cref{cert:height-six} supplies the compact set $Y_I$, its finite
closed cell cover, buffered one-digit returns, and a buffered entry from
$(1,1)$.  Hence \Cref{seven:lem:control} applies, proving that $x$ is an
interior point of $\mathcal Z_6^{\mathrm{ext}}$.
\end{proof}

By \Cref{cor:interior-double-zero}, $[1000/2999,1)\subset\mathcal D_6$.

\subsubsection{A moving rectangle and paired renewal}
\label{seven:sec:continuation}
For $2999/1000\le x\le3$ let $R_x=[-\bar p,\bar p]\times[-\bar v,\bar v]$,
with $\bar p=3/r$ and $\bar v=2/r^2$, be the rectangle of
\Cref{sec:pronic-completion}, and $p_0=-2/r$, so $z_{-2}=(p_0,\bar v)$. Put
\begin{equation}\label{seven:eq:renewal-thresholds}
 \lambda_-=(2x-4-2\bar v)/x^2,\quad \kappa=x\lambda_-,\quad
 \lambda_+=(2x-3)/x^2,\quad \xi_*=(x-2-1/r)/x^2,\quad W=1-1/r;
\end{equation}
these are positive, and $1.999\le r\le2$ gives $1\le2\bar v<1.002$. By
\eqref{eq:jordan} and $F_{x,d+1}=F_{x,d}-(1,1)$,
\begin{align}
 F_{x,-6}\bigl(z_{-2}+(\xi-1,\zeta-1)\bigr)
 &=z_{-2}+\bigl(4-x+x\xi,\ 4-2x+x(\xi+\zeta)\bigr),\label{seven:eq:fast}\\
 F_{x,(-1,-j-2)}\bigl(z_{-2}+(\xi,\zeta)\bigr)
 &=z_{-2}+\bigl(x^2\xi-x+j,\ x^2(\zeta+2\xi)-2x+j\bigr),
 \label{seven:eq:renewal-exit}
\end{align}
and $F_{x,-5}$ gives the first image minus $(1,1)$. Two captures recur.
\begin{enumerate}[label=(\alph*)]
\item If $-\delta-1/r\le x\xi\le xW$ and $\kappa\le\xi+\zeta\le(2x-3)/x$, then
$-6$ or $-5$ sends $z_{-2}+(\xi-1,\zeta-1)$ into $R_x$. The vertical
constraints hold for $\xi+\zeta$ in $[\kappa,(2x-4)/x]$ and
$[(2x-3-2\bar v)/x,(2x-3)/x]$, which overlap as $2\bar v\ge1$; the horizontal
ones are $-\delta-1/r\le x\xi\le x-4+5/r$, and $xW\le x-4+5/r$ is $x\le3$.
\item If $-x/(4r)\le x^2\xi\le x-2$ and $\lambda_-\le\zeta+2\xi\le\lambda_+$,
then $(-1,-6)$ or $(-1,-5)$ sends $z_{-2}+(\xi,\zeta)$ into $R_x$. The
vertical constraints hold for $\zeta+2\xi$ in $[\lambda_-,(2x-4)/x^2]$ and
$[(2x-3-2\bar v)/x^2,\lambda_+]$, which overlap likewise; the horizontal ones reduce to
$1-x/4+\delta r\ge0$, $4-x+(1-x/4)/r\ge0$ and $2\le5/r$ (which also gives
the bound $1\le5/r$ for the digit $-5$).
\end{enumerate}

\begin{theorem}[Uniform rectangle continuation\certmark]\label{seven:thm:continuation}

Every point of $R_x$ admits an infinite bounded admissible orbit under
\eqref{eq:PV-map}, for every $2999/1000\le x\le3$.
More precisely, every point has a nonempty finite return to $R_x$.
Return lengths need not be uniformly bounded, but all intermediate
states in the construction lie in one bounded set, independent of $x$
and of the return length.
\end{theorem}

\begin{lemma}[Compact paired return\certmark]\label{seven:lem:paired-cover}

Let
\[
 \mathcal Q=\{(u,w):u\ge1/28,\ 0\le w\le u,\ u+w\le1/3\}.
\]
For every $(u,w)\in\mathcal Q$ and $2999/1000\le x\le3$, at least one of
\[
 (p_0+u,\bar v+w),\qquad(p_0+u-1,\bar v+w-1)
\]
has an admissible word of length at most three whose image satisfies
\[
 |P|\le\bar p-\delta/100,\qquad |V|\le\bar v-\delta/100.
\]
\end{lemma}
\begin{proof}
This is part (ii) of \Cref{cert:height-six}\certmark.
\end{proof}

\begin{proof}[Proof of \Cref{seven:thm:continuation}]
\textit{Idea.} A failed one-step return leaves two adjacent states; a fast
return, the cone $\mathcal Q$ or an excursion along $z_{-2}$ takes over, and
an excursion returns or renews the pair with $u$ multiplied by at least $x^2$.

\emph{Pairs.} By the symmetry $(P,V,d)\mapsto(-P,-V,-d)$ let
$0\le V\le\bar v$. Then $F_{x,d}(P,V)\in R_x$ if and only if $xP-d$ lies in
$\mathcal J=[\max(-\bar p,-xV-\bar v),\,-xV+\bar v]$. As $-x\bar p-\bar v>-6$
and $x\bar p+r\bar v=(3x+2)/r<6$, the digit $d_1=\lceil xP+xV-\bar v\rceil$
is legal, and $xP-d_1$, the largest point of $xP-\Z$ not above $\mathcal J$,
is a return if $\mathcal J$ has length at least one. As $2\bar v\ge1$, a
failure forces $\mathcal J$ to start at $-\bar p$ and $xP-d_1<-\bar p$; then
the digits $d_1-1$ and $d_1$ give
\begin{equation}\label{seven:eq:paired-wedge}
 z=(p_0+u,\bar v+w),\qquad z-(1,1),\qquad 0<w\le u<W,
\end{equation}
with $u=xP-d_1+1-p_0<1-\bar p-p_0=W$ and $w=x(V-\bar v)+u\in(0,u]$, positive
as $xP-d_1+1$ lies above $\mathcal J$. Both digits are legal, since
$d_1-1=xP-p_0-u\in(-x\bar p-p_0-W,\,x\bar p-p_0)=(-4,(3x+2)/r)$. The upper
state alone need not return.

\emph{Fast return and cone.} If $u+w\ge\kappa$, capture (a) applies to the
lower state with $(\xi,\zeta)=(u,w)$, whose horizontal condition holds as
$0<xu<xW$, and as $u+w<2W\le(2x-3)/x$, which is
$x\le3$. If $u+w<\kappa$, then $2\bar v<1.002$, $2x-4\ge1.998$, $x^2\le9$ and
$5x-12\le3\le6\bar v$ give $\lambda_-/3\ge83/2250>1/28$ and $\kappa\le1/3$;
so \Cref{seven:lem:paired-cover} gives a return if $u\ge\lambda_-/3$.

\emph{Excursions.} Let $0<u<\lambda_-/3$. By \eqref{eq:jordan}, $m$ digits
$-2$ take the upper state to $z_{-2}+(\xi,\zeta)$ with
$(\xi,\zeta)=x^m(u,w+mu)$; put $\gamma=\zeta+2\xi=x^m(w+(m+2)u)$, which grows
without bound. Let $m$ be first with $\gamma\ge\lambda_-$. Then $m\ge1$, as
$w+2u\le3u<\lambda_-$, and $\gamma=x\gamma'+\xi$ with
$\gamma'=(\xi/x)(m+1+w/u)<\lambda_-$ the preceding value, so
\begin{equation}\label{seven:eq:crossing-bounds}
 \xi<x\lambda_-/3\le(x-2)/x^2,\qquad
 \lambda_-\le\gamma<x\lambda_-+\xi\le1/x+\xi.
\end{equation}
The first bound is $u<\lambda_-/3$ if $m=1$ and follows from
$\gamma'<\lambda_-$ if $m\ge2$; the second is
$(x-3)(2x^3-5x^2+7x-2)\le0$, and the last is $2x-5\le2\bar v$. If
$\gamma\le\lambda_+$, capture (b) applies, as $0<x^2\xi<x-2$ by
\eqref{seven:eq:crossing-bounds}. If $\gamma>\lambda_+$ and
$\xi\ge\xi_*$, take $j=2$ in \eqref{seven:eq:renewal-exit}: $P\ge-\bar p$ is
exactly $\xi\ge\xi_*$, $P\le p_0$, and $\lambda_+\ge(2x-2-2\bar v)/x^2$,
$\gamma<1/x+\xi\le(2x-2)/x^2$ give the vertical bounds. If
$\gamma>\lambda_+$ and $\xi<\xi_*$, then $j=3$ gives the renewed upper state
$z_{-2}+(u',w')$ with $u'=x^2\xi+\delta$, $w'=x^2\gamma+3-2x$, and $j=2$ the
renewed lower state. Here $w'>0$ is $\gamma>\lambda_+$, $w'\le u'$ is
$\gamma\le\xi+1/x$, and $u'<x^2\xi_*+\delta=W$; so
\eqref{seven:eq:paired-wedge} holds again, with $u'\ge x^2u$. After finitely
many renewals, a fast return, the cone or an excursion returns.

\emph{Boundedness.} Concatenating returns gives an infinite admissible orbit.
Each state is a point of $R_x$, a pair state or an excursion state
$z_{-2}+(\xi',\zeta')$ with $0\le\xi'<1/9$, $0\le\zeta'<1/2$ (by
\eqref{seven:eq:crossing-bounds}, as $(x-2)/x^2\le1/9$ and $1/x+1/9<1/2$),
all of sup norm at most $5/2$, or the image of one of these under at most
three digits; as $\|F_{x,d}(z)\|_\infty\le6\|z\|_\infty+6$, every state has
sup norm at most $6(6(6\cdot\frac52+6)+6)+6=798$.
\end{proof}

\subsubsection{Entry at every terminal scale}
\label{seven:sec:all-scale-entry}

\begin{theorem}[Terminal bounded orbits\certmark]\label{seven:thm:terminal-orbit}

For every $2999/1000\le x<3$, the state $(1,1)$ has a finite admissible
word entering $R_x$, and hence an infinite bounded admissible orbit.
Consequently
\[
 (1/3,1000/2999]\subset\mathcal D_6.
\]
\end{theorem}

First let $0<\delta\le\delta_0:=1/(5\cdot3^{15})$, and let $M$ be the first
positive integer with
\begin{equation}\label{seven:eq:entry-scale}
 \alpha=\delta x^{M-1}\ge\frac15.
\end{equation}
As $\delta x^{15}<\delta_03^{15}=1/5$, $M\ge17$. Minimality gives
$\delta x^{M-2}<1/5$, so, with $x\ge2$, $r\ge1$ and the decrease of $M/2^M$,
\begin{equation}\label{seven:eq:entry-bounds}
 1/5\le\alpha<x/5<3/5,\qquad
 M\delta/r<M/(5x^{M-2})\le\tfrac45M/2^M\le\tfrac45\cdot17/2^{17}<1/1000.
\end{equation}
At $q=3$ the blocks of \Cref{sec:residual} are $\mathsf A_M=5(-4)^{M-1}$,
$\mathsf B_M=5(-4)^{M-2}(-5)$ and $\mathsf U_{M,k}=\mathsf B_M(-2)^k$; the
$\mathsf B_M$ state, $z_{-2}+(\xi,\zeta)$ say, is the $\mathsf A_M$ state plus
$(1,1)$, and $\xi=c_M$ by \eqref{eq:c-m}. As
$F_{x,5}(1,1)=z_{-4}+\bigl(\delta^2/r,\,-\delta(2x^2-3x+3)/r^2\bigr)$ and
$z_{-4}-z_{-2}+(1,1)=(1-2/r,\,1+2/r^2)$, \eqref{eq:jordan} at $z_{-4}$ for
$M-1$ steps, the shift of the last digit and $\delta x^{M-1}=\alpha$ give
\begin{equation}\label{seven:eq:entry-q}
 \xi=\frac{\delta(\alpha-1)}r<0,\qquad
 \zeta=\frac{x(1-2\alpha)}r+
 \delta\Bigl(\frac{(M-1)\alpha}r+\frac{1-\alpha}{r^2}\Bigr).
\end{equation}
By \eqref{eq:jordan}, the $\mathsf U_{M,k}$ state is
$z_{-2}+x^k(\xi,\zeta+k\xi)=:z_{-2}+(\xi_k,\zeta_k)$.

\begin{lemma}[The descending half of the first phase]
\label{seven:lem:descending-entry}

If $\alpha\ge1/2$, an admissible continuation of $\mathsf B_M$ enters
$R_x$ after at most $M$ further digits $-2$.
\end{lemma}
\begin{proof}
\textit{Idea.} Along the $-2$ excursion the vertical displacement turns
negative within $M$ steps, while the horizontal one stays small.
By \eqref{seven:eq:entry-q},
$\zeta+M\xi=x(1-2\alpha)/r+\delta\bigl(((2M-1)\alpha-M)/r+(1-\alpha)/r^2\bigr)$.
As a function of $\alpha$ this has derivative
$[-2x+\delta(2M-1)-\delta/r]/r<0$, as $M\delta<r/1000<1$ by
\eqref{seven:eq:entry-bounds}, and value
$\delta(2-x)/(2r^2)<0$ at $\alpha=\frac12$; so $\zeta_M<0$. If $\zeta\le0$,
$\mathsf B_M$ already enters $R_x$: $\zeta>-x/(5r)\ge-2\bar v$ as
$\alpha<3/5$, and $|\xi|\le\delta/(2r)$ gives $|P|<\bar p$. Otherwise let
$k\le M$ be first with $\zeta_k\le0$. Then $\zeta_k-\xi_k=x\zeta_{k-1}>0$, so
$\xi_k<\zeta_k\le0$, and
$|\xi_k|\le x^M\delta(1-\alpha)/r=x\alpha(1-\alpha)/r\le x/(4r)$. Thus
$-\bar p\le-(2+x/4)/r\le P\le p_0$ and $-\bar v\le\bar v-x/(4r)\le V\le\bar v$,
by $x<4$ and $xr\le6<16$.
\end{proof}

Suppose henceforth that $\alpha\le1/2$, and set
\begin{equation}\label{seven:eq:second-phase}
 \Theta=x^M(2\alpha-1)\le0,\qquad
 \gamma_k=\zeta_k+2\xi_k=x^k\bigl(\zeta+(k+2)\xi\bigr).
\end{equation}
If $k\le M-2$ and $\lambda_-\le\gamma_k\le\lambda_+$, capture (b) applies to
the $\mathsf U_{M,k}$ state, since
$0>x^2\xi_k\ge x^M\xi=-x\alpha(1-\alpha)/r\ge-x/(4r)$. Also
\begin{equation}\label{seven:eq:xl-u}
 (x\lambda_--\lambda_+)\,x^2r^2=(x-3)(2x^3-4x^2+5x-1)\le0,
 \quad\text{so}\quad \kappa=x\lambda_-\le\lambda_+.
\end{equation}

\begin{lemma}[The outer part of the second phase]
\label{seven:lem:outer-entry}

Suppose $\alpha\le1/2$. There is finite entry into $R_x$ if either
$\gamma_0>\lambda_+$, or $\gamma_0\le\lambda_+$ and $\Theta\le-4/5$.
\end{lemma}
\begin{proof}
\textit{Idea.} Either $\mathsf A_M$ is one digit from $R_x$, or the $-2$
excursion meets $[\lambda_-,\lambda_+]$ by step $M-2$.
If $\gamma_0>\lambda_+$, then $\beta_0=\zeta+\xi=\gamma_0-\xi>\lambda_+\ge
\kappa$ by \eqref{seven:eq:xl-u}, and by
\eqref{seven:eq:entry-bounds}--\eqref{seven:eq:entry-q} and $\bar v<0.501$,
\[
 \beta_0=1+\bar v-3\alpha+\delta(M\alpha-1)/r-\alpha\delta x/r^2
 \le\tfrac25+\bar v+M\delta/(2r)<\tfrac{1803}{2000}<\tfrac{99}{100}<
 (2x-3)/x.
\]
As $0>x\xi\ge-\delta-\delta/r\ge-\delta-1/r$, capture (a) applies to the
$\mathsf A_M$ state $z_{-2}+(\xi-1,\zeta-1)$.

Now let $\gamma_0\le\lambda_+$ and $\Theta\le-4/5$. Substituting
$x^{M-2}\delta=\alpha/x$, $x^{M-1}(1-2\alpha)=-\Theta/x$ and
$x\alpha(2\alpha-1)=\delta\Theta$, and using $x/r=\frac32+\delta/(2r)$, into
$\gamma_{M-2}=x^{M-2}(\zeta+M\xi)$ gives
\begin{equation}\label{seven:eq:outer-last-state}
 \gamma_{M-2}=\frac{\alpha(1-\alpha x)}{xr^2}
   -\frac{\Theta}{x^2}\Bigl(\frac32-\frac{\delta(M-\frac12)}r\Bigr)
 \ge-\frac{21}{1000}+\frac45\cdot\frac{1499}{1000\cdot9}
 =\frac{5051}{45000}>\frac19\ge\lambda_-.
\end{equation}
Here the first term decreases for $\alpha\ge1/5>1/(2x)$, so it is at least
its value $(2-x)/(4xr^2)>-21/1000$ at $\alpha=\frac12$; the second is bounded
by $-\Theta\ge4/5$, \eqref{seven:eq:entry-bounds} and $x^2\le9$; and
$\lambda_-\le(2x-5)/x^2\le1/9$ as $x^2-18x+45=(x-3)(x-15)\ge0$. If
$\gamma_0\ge\lambda_-$, capture (b) applies with $k=0$. Otherwise, for the
first $k\le M-2$ with $\gamma_k\ge\lambda_-$, $\xi_{k-1}<0$ and
\eqref{seven:eq:xl-u} give
$\gamma_k=x(\gamma_{k-1}+\xi_{k-1})<x\lambda_-\le\lambda_+$, and capture (b)
applies. This is uniform even when $-\Theta$ is unbounded.
\end{proof}

\subsubsection{The central second phase}\label{seven:sec:central-phase}
Only $-4/5\le\Theta\le0$ remains. For $\ell=0,1$ let
$\mathsf P_M=\mathsf U_{M,M}$, $\mathsf C_M=\mathsf U_{M,M-1}(-1,-6)$ and
$\mathsf R_{M,\ell}=\mathsf U_{M,\ell}(-1,-6,-3)6^{M-\ell-2}$, an excursion
around $z_6=(6/r,-6/r^2)$. With $G_\ell=(x+1)/(rx^{\ell+1})$, their exact
coordinates are
\begin{align}
 P_{\mathsf P}&=\bigl(-2+x\alpha(\alpha-1)\bigr)/r,\quad
 V_{\mathsf P}=\bigl(2+x\alpha(1-\alpha)\bigr)/r^2-x\alpha^2/r
              -(x-M\delta)\Theta/r,
 \label{seven:eq:central-P}\\
 P_{\mathsf R_\ell}&=\bigl(6+x^2\alpha(\alpha-1)\bigr)/r-\delta\alpha G_\ell,
 \label{seven:eq:central-RP}\\
 V_{\mathsf R_\ell}&=-6/r^2+x^2\alpha(\alpha-1)(2x-3)/r^2-x^2\Theta/r
  +2x\alpha G_\ell\notag\\
 &\quad+\delta\bigl[(M-1)x\Theta/r+\alpha G_\ell
 \bigl(\ell+1-M+x/r-x/(x+1)\bigr)\bigr],
 \label{seven:eq:central-RV}\\
 P_{\mathsf C}&=xP_{\mathsf P}-x+6,\qquad
 V_{\mathsf C}=x(V_{\mathsf P}+P_{\mathsf P})-2x+6.
 \label{seven:eq:central-bridge}
\end{align}
For $\mathsf P_M$ take $k=M$ above, with $x^M\xi=x\alpha(\alpha-1)/r$ and
$x\alpha(2\alpha-1)=\delta\Theta$. For $\mathsf R_{M,\ell}$, the map
$F_{x,(-1,-6,-3)}(P,V)=(x^3P+x^2+6x+3,\ x^3(V+3P)+3x^2+12x+3)$ sends $z_{-2}$
to $z_6+\bigl(-\delta^2(x+1)/r,\ \delta(3x^3-5x^2+3x+3)/r^2\bigr)$; apply it
to $z_{-2}+(\xi_\ell,\zeta_\ell)$ and then \eqref{eq:jordan} at $z_6$ for
$M-\ell-2$ steps, using $x^{M-1}\delta=\alpha$ and
$x\alpha(2\alpha-1)=\delta\Theta$. Finally, $\mathsf C_M$ is $\mathsf P_M$ with
its last digit raised to $-1$, which subtracts $(1,1)$, followed by $-6$. At
$x=3$, $\alpha=1/2$, keeping $\Theta$, the limiting states are
\begin{equation}\label{seven:eq:central-limits}
 \mathsf P:\ \bigl(-\tfrac{11}8,\ \tfrac5{16}-\tfrac32\Theta\bigr),\qquad
 \mathsf R_\ell:\ \bigl(\tfrac{15}8,\ \tfrac2{3^\ell}-\tfrac{51}{16}
 -\tfrac92\Theta\bigr),\qquad
 \mathsf C:\ \bigl(-\tfrac98,\ -\tfrac{51}{16}-\tfrac92\Theta\bigr).
\end{equation}

\begin{table}[htbp]
\centering\footnotesize
\caption{Central cells: a closed partition of $[-4/5,0]$ in $\Theta$.}
\label{seven:tab:central-cells}
\begin{tabular}{@{}ll@{\qquad}ll@{}}
\toprule
$[-4/5,-3/5]$ & $\mathsf C_M$ & $[-11/27,-781/1944]$ & $\mathsf P_M(-2,-5,-3,-2)$\\
$[-3/5,-179/324]$ & $\mathsf R_{M,1}(6)$ & $[-781/1944,-107/324]$ & $\mathsf P_M(-2,-5)$\\
$[-179/324,-175/324]$ & $\mathsf R_{M,1}(5,3)$ & $[-107/324,-103/324]$ & $\mathsf P_M(-2,-6,0)$\\
$[-175/324,-155/324]$ & $\mathsf R_{M,1}(5)$ & $[-103/324,-83/324]$ & $\mathsf P_M(-2,-6)$\\
$[-155/324,-37/81]$ & $\mathsf P_M(-2,-4,-6)$ & $[-83/324,-79/324]$ & $\mathsf P_M(-3,-1,-6)$\\
$[-37/81,-110/243]$ & $\mathsf R_{M,1}(4,5,1)$ & $[-79/324,-59/324]$ & $\mathsf R_{M,0}(5)$\\
$[-110/243,-35/81]$ & $\mathsf R_{M,1}(4,5)$ & $[-59/324,-55/324]$ & $\mathsf P_M(-3,-2,-3)$\\
$[-35/81,-104/243]$ & $\mathsf P_M(-2,-5,-2,-5)$ & $[-55/324,-7/54]$ & $\mathsf P_M(-3,-2)$\\
$[-104/243,-11/27]$ & $\mathsf R_{M,1}(4,4)$ & $[-7/54,0]$ & $\mathsf P_M(-3)$\\
\bottomrule
\end{tabular}
\end{table}

\begin{lemma}[Eighteen central entry cells\certmark]\label{seven:lem:central-entry}

For every scale and parameter in \eqref{seven:eq:entry-scale} with
$-4/5\le\Theta\le0$, one of the words in
\Cref{seven:tab:central-cells} sends $(1,1)$ strictly into
$R_0=[-3/2,3/2]\times[-1/2,1/2]\subset R_x$.
\end{lemma}
\begin{proof}
\textit{Idea.} At $x=3$, $\alpha=\frac12$ the images are affine in $\Theta$ and
lie in $R_0$ with margins; the actual images differ by $O(M\delta)$.
By \Cref{cert:height-six}\certmark(iii), with suffixes of length at most
four, every limiting image lies in
$[-\frac{11}8,\frac{11}8]\times[-\frac{39}{80},\frac{39}{80}]$, so with
margins $\frac18$ and $\frac1{80}$ in $R_0$, and every limiting partial state,
the base included, lies in $[-3,3]\times[-\frac32,\frac32]$.

To transfer this uniformly in $M$, note that, by $\alpha\le\frac12$,
$|\Theta|\le\frac45$ and $x^M\ge(2999/1000)^{17}$,
$|\alpha-\frac12|=|\Theta|/(2x^M)$ gives $\alpha\ge\frac{49}{100}$, and then
$|\alpha-\frac12|=|\Theta|\delta/(2\alpha x)<\delta/3$. Write each base
coordinate as $f(x,\alpha,\Theta)+\delta g$, where $\delta g$ is the displayed
term with the factor $\delta$: none in $P_{\mathsf P}$, $M\delta\Theta/r$ in
$V_{\mathsf P}$, $-\delta\alpha G_\ell$ in $P_{\mathsf R_\ell}$, the bracket in
$V_{\mathsf R_\ell}$; treat $x,\alpha,\Theta$ as independent. On the box
$2999/1000\le x\le3$, $\frac{49}{100}\le\alpha\le\frac12$ we have
$1.999\le r\le2$, $x(x-2)\le3$, $\alpha(1-\alpha)\le\frac14$,
$|2\alpha-1|\le\frac1{50}$, and $|\Theta|\le\frac45$, or $|\Theta|\le\frac35$
for the $\mathsf R$ words. Hence
\begin{align*}
 |\partial_xP_{\mathsf P}|&=\bigl(2+\alpha(1-\alpha)\bigr)/r^2<1,&
 |\partial_\alpha P_{\mathsf P}|&=x|2\alpha-1|/r<1,\\
 |\partial_xV_{\mathsf P}|&=\bigl|(\alpha+\Theta)/r^2
   -2\bigl(2+x\alpha(1-\alpha)\bigr)/r^3\bigr|<2,&
 |\partial_\alpha V_{\mathsf P}|&=\bigl|x(1-2\alpha)/r^2-2x\alpha/r\bigr|<2,\\
 |\partial_xP_{\mathsf R_\ell}|&=\bigl(6+\alpha(1-\alpha)x(x-2)\bigr)/r^2<2,&
 |\partial_\alpha P_{\mathsf R_\ell}|&=x^2|2\alpha-1|/r<1,\\
 |\partial_xV_{\mathsf R_\ell}|&\le12/r^3+\tfrac34+\tfrac95/r^2+1<5,&
 |\partial_\alpha V_{\mathsf R_\ell}|&\le\tfrac7{50}+8/r<6.
\end{align*}
The last row uses $\partial_xV_{\mathsf R_\ell}=12/r^3+\alpha(\alpha-1)Q'
-\Theta x(x-2)/r^2+2\alpha J'_\ell$ and $\partial_\alpha V_{\mathsf R_\ell}=
(2\alpha-1)Q+2J_\ell$, where $Q=x^2(2x-3)/r^2\in(0,7]$,
$Q'=2x(x^2-3x+3)/r^3\in(0,3]$, $J_\ell=xG_\ell\in(0,4/r]$ and
$|J'_\ell|\le1$. Also $|g|\le\frac M2,\frac12,\frac32M+3$ for $V_{\mathsf P}$,
$P_{\mathsf R_\ell}$, $V_{\mathsf R_\ell}$, by $|\Theta|/r\le\frac12$,
$\alpha G_\ell<\frac12$, $x|\Theta|/r<1$ for $\mathsf R$, and
$-M\le\ell+1-M+x/r-x/(x+1)\le0$. By the mean value theorem, a non-bridge base
coordinate is within $(5+\frac63+\frac32M+3)\delta\le\epsilon_M:=(2M+10)\delta$
of its limit, the $V_{\mathsf R_\ell}$ row being the worst; for $\mathsf P$
the errors are $e_P\le\frac43\delta$ and $e_V\le(\frac83+\frac M2)\delta$,
also when $|\Theta|\le\frac45$. By \eqref{seven:eq:central-bridge},
$P_{\mathsf C}$ and $V_{\mathsf C}$ differ from their limits by
$x(P_{\mathsf P}+\frac{11}8)+\frac{19}8\delta$ and
$x\bigl(V_{\mathsf P}+P_{\mathsf P}-\frac5{16}+\frac32\Theta+\frac{11}8\bigr)
+\delta\bigl(\frac{49}{16}+\frac32\Theta\bigr)$, so by at most
$3e_P+\frac{19}8\delta$ and $3(e_P+e_V)+\frac{49}{16}\delta$, both below
$(2M+25)\delta$.

For a nonempty suffix of length $k\le4$, apply
\eqref{eq:word-partial-state-error} with
$J=\bigl(\begin{smallmatrix}1&0\\1&1\end{smallmatrix}\bigr)$ and $q=3$.
Here $\|J^j\|_\infty=j+1$, while each limiting partial state has
$\|J(P,V)\|_\infty=\max\{|P|,|P+V|\}\le9/2$.
Thus the propagated base error is at most $405\epsilon_M$, and the
parameter error is at most
$\frac92\sum_{j=0}^3 3^j(j+1)\delta=639\delta<1166\delta$.
The bridge cell has an empty suffix. Keeping the common allowance,
every final coordinate is therefore within
$405(2M+10)\delta+1166\delta=(810M+5216)\delta$ of its limit. As
$\alpha\le\frac12$, $\delta\le1/(2x^{M-1})\le\delta_M:=\frac12(2999/1000)^{1-M}$,
and $(810M+5216)\delta_M$ decreases in $M$, with ratio
$\frac{810M+6026}{810M+5216}\cdot\frac{1000}{2999}<\frac{6026}{5216}\cdot
\frac{1000}{2999}<1$, from $18986\,\delta_{17}<2.3\cdot10^{-4}<\frac1{80}$. So
all eighteen images lie strictly inside $R_0$, and $R_0\subset R_x$ as
$\bar p\ge\frac32$ and $\bar v\ge\frac12$.
\end{proof}

\subsubsection{The finite join}\label{seven:sec:finite-join}

\begin{lemma}[Compact entry into the moving rectangle\certmark]
\label{seven:lem:compact-entry}

For every
\[
 \frac{2999}{1000}\le x\le
 3-\frac1{24\cdot3^{14}}=\frac{344373767}{114791256},
\]
a finite admissible word sends $(1,1)$ strictly into $R_x$.
\end{lemma}
\begin{proof}
This is part (iv) of \Cref{cert:height-six}\certmark.
\end{proof}

\begin{proof}[Proof of \Cref{seven:thm:terminal-orbit}]
For $0<3-x\le\delta_0$ choose $M$ by \eqref{seven:eq:entry-scale}.
\Cref{seven:lem:descending-entry} treats $\alpha\ge1/2$,
\Cref{seven:lem:outer-entry} treats $\alpha\le1/2$ with $\gamma_0>\lambda_+$
or $\Theta\le-4/5$, and \Cref{seven:lem:central-entry} treats
$-4/5\le\Theta\le0$; so $(1,1)$ enters $R_x$. \Cref{seven:lem:compact-entry}
supplies the remaining parameters, with overlap since
$1/(24\cdot3^{14})=\frac58\delta_0<\delta_0$. After entry,
\Cref{seven:thm:continuation} continues the orbit; with its finite prefix it
is infinite and bounded. By \Cref{lem:bounded-orbit}, $1/x\in\mathcal D_6$.
\end{proof}

\subsection{Height two}\label{app:height-two}

This subsection proves the inclusion $[1/\rho,1)\subseteq\mathcal A_2$ of
\Cref{tri:thm:double-zeros}; the exclusion of the other double zeros is proved
in \Cref{sec:height-two}. For $x\in[433/250,\rho)$ we build a bounded real
orbit and apply \Cref{lem:bounded-orbit}. The orbit enters a fixed diamond $Q$
along the periodic digit block of the tribonacci witness $f_\rho$ of
\eqref{tri:eq:endpoint-series}, followed by a short suffix, and a finite return
cover then brings it back to $Q$ again and again.

\subsubsection{Trace states and the return diamond}
Fix $x>1$. A \emph{state} is a pair $Y=(S,V)\in\R^2$. For a digit $d\in B_2$
and a word $w=(d_1,\ldots,d_n)$ over $B_2$ put
\begin{equation}\label{tri:eq:trace-map}
 \Theta_{x,d}(S,V)=(2xS-xV-2d,\ xS-d),\qquad
 \Theta_{x,w}=\Theta_{x,d_n}\circ\cdots\circ\Theta_{x,d_1}.
\end{equation}
By \eqref{eq:VS-rec-real}, $\Theta_{x,w}(2,1)=(S_w(x),V_w(x))$. In the
coordinates $(P,V)=(S-V,V)$ the map $\Theta_{x,d}$ is \eqref{eq:PV-map}, and
$(2,1)$ becomes $(1,1)$. Hence, by \Cref{lem:bounded-orbit}, if a digit
sequence $d_1,d_2,\ldots$ in $B_2$ moves $(2,1)$ along a bounded orbit, then
$1-\sum_kd_kz^k\in\mathcal F_2$ has a double zero at $z=1/x$. The linear part
$L(P,V)=(xP,xV+xP)$ of \eqref{eq:PV-map} satisfies $L^n(P,V)=x^n(P,V+nP)$.
Rewritten in $(S,V)$ and unrolled, this gives
\begin{equation}\label{tri:eq:unrolled}
 \Theta_{x,w}Y=x^nJ^nY-\sum_{k=1}^nd_k\,x^{n-k}\,(n-k+2,\,n-k+1),\qquad
 J^n=\begin{pmatrix}n+1&-n\\ n&1-n\end{pmatrix},
\end{equation}
with states read as column vectors. For $n\ge1$ the absolute row sums of $J^n$
are $2n+1$ and $2n-1$.

In the coordinates
\[
 u_1=\frac{5V}7,\qquad u_2=\frac{2S}3-\frac V5=\frac{2P}3+\frac{7V}{15},
\]
put $\|Y\|_Q=|u_1|+|u_2|$ and
\begin{equation}\label{tri:eq:diamond}
 Q=\{Y\in\R^2:\|Y\|_Q\le1\}.
\end{equation}
Thus $Y\in Q$ if and only if its four \emph{signed margins}
$1-\epsilon_1u_1-\epsilon_2u_2$, $\epsilon_1,\epsilon_2\in\{-1,1\}$, are
nonnegative, and $Y\in\interior Q$ if and only if they are positive. Since
\begin{equation}\label{tri:eq:diamond-inverse}
 (S,V)=\Bigl(\frac{21}{50}u_1+\frac32u_2,\ \frac75u_1\Bigr),
\end{equation}
every $Y\in Q$ has $\|Y\|_\infty\le3/2$, where $\|Y\|_\infty=\max\{|S|,|V|\}$.
Conversely,
\begin{equation}\label{tri:eq:norm-factor}
 \|Y\|_Q\le\tfrac23|S|+\bigl(\tfrac57+\tfrac15\bigr)|V|
 \le\tfrac{166}{105}\|Y\|_\infty .
\end{equation}
Finally, with $G$ as in \eqref{tri:eq:endpoint-series}, the cubic
$\chi(x)=x^3G(1/x)=x^3-x^2-x-1$ increases on $[1,\infty)$, and
\[
 \chi\bigl(\tfrac{183}{100}\bigr)=-\frac{50413}{10^6}<0<
 \frac{3828037384313}{15625\cdot10^{15}}
 =\chi\bigl(\tfrac{4598217}{2500000}\bigr),
\]
so that $183/100<\rho<4598217/2500000<46/25$.

\begin{proposition}[Closed return cover\certmark]\label{tri:prop:core}

For every $x\in[433/250,4598217/2500000]$ and $Y\in Q$, some digit word $w$,
with $1\le|w|\le10$, satisfies $\Theta_{x,w}Y\in Q$. The parameter interval
contains $\rho$.
\end{proposition}
\begin{proof}
For a fixed word $w$, \eqref{tri:eq:unrolled} shows that the four signed
margins of $\Theta_{x,w}Y$ are affine in $(u_1,u_2)$, with coefficients that
are polynomials in $x$. On a closed rational parameter band, exact rational
enclosures of these polynomials over the band define a polygon of states, the
\emph{word region} of $w$, on which all four margins are nonnegative for every
$x$ in the band; thus $\Theta_{x,w}$ maps the word region into $Q$. Part~(i)
of \Cref{cert:height-two-return}\certmark{} splits $[433/250,4598217/2500000]$ into
consecutive closed rational bands and covers $Q$, on each band, by word regions
of words of length between $1$ and $10$. This proves the first assertion. The
second follows from the bracket for $\rho$ displayed above.
\end{proof}

\subsubsection{Tribonacci entry}\label{tri:sec:phase}
Let
\[
 a=(2,1,-1),\qquad b=(-1,-1,-2),\qquad Y_m(x)=\Theta_{x,ab^m}(2,1)\quad(m\ge0),
\]
where $b^m$ is the $m$-fold concatenation. In the convention of
\Cref{lem:bounded-orbit}, the digits of $f_\rho$ in \eqref{tri:eq:endpoint-series}
are $a$ followed by $b,b,\ldots$. Appending a digit $d$ to a word $u$ sends
$P_u(x)$ to $xP_u(x)-d$, so appending $b$ acts by $P\mapsto x^3P+x^2+x+2$.
This affine map fixes $\bar P(x)=-(x^2+x+2)/(x^3-1)=-1/(x-1)-1/(x^3-1)$, and
\[
 \Delta(x):=P_a(x)-\bar P(x)=x^3-2x^2-x+1+\frac{x^2+x+2}{x^3-1}
 =\frac{\chi(x)^2}{x^3-1}.
\]
Hence $P_{ab^m}=\bar P+x^{3m}\Delta$ as rational functions of $x$.
Differentiating, and using $S=2P+xP'$ and $V=P+xP'$ from \eqref{eq:VS-def},
\begin{equation}\label{tri:eq:centered-return}
 Y_m(x)=\bar Y(x)+x^{3m}\bigl(\omega+(3m+1)\Delta,\ \omega+3m\Delta\bigr)(x),
 \qquad \omega=\Delta+x\Delta',
\end{equation}
where $\bar Y=(\bar S,\bar V)=(2\bar P+x\bar P',\ \bar P+x\bar P')$ is the
state fixed by the block:
\begin{equation}\label{tri:eq:fixed-state}
 \bar S(x)=\frac{2-x}{(x-1)^2}+\frac{x^3+2}{(x^3-1)^2},\qquad
 \bar V(x)=\frac1{(x-1)^2}+\frac{2x^3+1}{(x^3-1)^2}.
\end{equation}

Put $\varepsilon=\rho-x$. As $\rho^3=\rho^2+\rho+1$, we have
$\chi(x)=-\varepsilon\psi(x)$ with $\psi(x)=x^2+(\rho-1)x+\rho^2-\rho-1$, so
that $\psi(\rho)=\chi'(\rho)$. With $\eta=\psi^2/(x^3-1)$,
\begin{equation}\label{tri:eq:eta}
 \Delta=\varepsilon^2\eta,\qquad
 \omega=-2x\varepsilon\eta+\varepsilon^2(\eta+x\eta').
\end{equation}
In particular $\Delta(\rho)=\omega(\rho)=0$, so $Y_m(\rho)=\bar Y(\rho)$ for
every $m$. Put $\kappa=2\eta(\rho)/\rho^2$, so that
$\omega(x)=-\kappa\rho^3\varepsilon+O(\varepsilon^2)$. Multiply
\eqref{tri:eq:fixed-state} by $(x^3-1)^2=(x-1)^2(x^2+x+1)^2$; the polynomial
identities
\begin{align*}
 (2-x)(x^2+x+1)^2+x^3+2-(7+2x-3x^2)(x^3-1)^2
   &=(3x^5+x^4-3x^3-6x^2-4x+3)\,\chi(x),\\
 (x^2+x+1)^2+2x^3+1-(5+2x-2x^2)(x^3-1)^2
   &=(2x^5-3x^3-5x^2-3x+3)\,\chi(x),\\
 2\chi'(x)^2-(16+52x-32x^2)\,x^2(x^3-1)
   &=(32x^4-20x^3-4x^2-6x-2)\,\chi(x),
\end{align*}
taken at $x=\rho$, give
\begin{equation}\label{tri:eq:Ystar}
 \bar Y(\rho)=Y_*:=(7+2\rho-3\rho^2,\ 5+2\rho-2\rho^2),\qquad
 \kappa=16+52\rho-32\rho^2.
\end{equation}

For $x<\rho$ and $m\ge0$ let $\tau_m(x)=(\rho-x)\rho^{3m+3}$ be the
\emph{phase} of $x$ at scale $m$ (it is unrelated to $\tau_h$ in
\eqref{eq:tau-sigma}), and put
\[
 \tau_-=\frac9{80},\qquad \tau_+=\tau_-\rho^3,\qquad
 Y_\infty(\tau)=Y_*-\kappa\tau\,(1,1).
\]
For $x\in[\rho-\tau_-,\rho)$, the least $m\ge0$ with $\tau_m(x)\ge\tau_-$ is
the \emph{entry scale} of $x$; since $\tau_0(x)\le\tau_+$ and
$\tau_{m+1}=\rho^3\tau_m$, it satisfies $\tau_m(x)\in[\tau_-,\tau_+]$. By
\eqref{tri:eq:entry-error} below, $Y_m(\rho-\tau\rho^{-3m-3})\to
Y_\infty(\tau)$ as $m\to\infty$, for each $\tau\in[\tau_-,\tau_+]$: after the
blocks, the state lies near the segment $Y_\infty([\tau_-,\tau_+])$, at a
point fixed by the phase. One of the eight suffixes
\[
\begin{gathered}
 s_1=(-1,-2,0,-2,0,-2),\quad s_2=(-1,-2,0,-2),\quad s_3=(-1,-2),\quad
 s_4=(-2,0),\\
 s_5=(-2),\quad s_6=(-2,-2),\quad s_7=(-2,-2,0,2,2,2),\quad
 s_8=(-2,-2,0,2,2)
\end{gathered}
\]
moves each point of this segment well inside $Q$: by part~(iii) of
\Cref{cert:height-two-return}\certmark, for every $\tau\in[\tau_-,\tau_+]$ some $j$
satisfies
\begin{equation}\label{tri:eq:phase-margin}
 \bigl\|\Theta_{\rho,s_j}Y_\infty(\tau)\bigr\|_Q\le1-\mu,\qquad
 \mu=\frac{223}{46500}.
\end{equation}
The record assigns $s_j$ to the $j$th of eight consecutive closed phase cells.
As $\Theta_{\rho,s_j}Y_\infty(\tau)$ is affine in $\tau$, so are its signed
margins, and it suffices to check them at the two ends of each cell. In the
coordinates $(P,V)$, $Y_\infty(\tau)=(2-\rho^2,\ 5+2\rho-2\rho^2-\kappa\tau)$,
the form used in the record.

\subsubsection{Uniform control of all entry scales}\label{tri:app:finite-phase}
\begin{lemma}[Entry at every scale\certmark]\label{tri:lem:finite-phase-estimates}

Let $m\ge0$, $\tau\in[\tau_-,\tau_+]$ and $x=\rho-\tau\rho^{-3m-3}$.
\begin{enumerate}
\item If $m\le8$, some word $s$ over $B_2$ satisfies
 $\Theta_{x,s}Y_m(x)\in\interior Q$.
\item If $m\ge9$, then for $j=1,\ldots,8$,
 \[
  \bigl\|\Theta_{x,s_j}Y_m(x)-\Theta_{\rho,s_j}Y_\infty(\tau)\bigr\|_Q
  <\frac9{2000}<\mu,
 \]
 so $\Theta_{x,s_j}Y_m(x)\in\interior Q$ whenever $s_j$ satisfies
 \eqref{tri:eq:phase-margin}.
\end{enumerate}
\end{lemma}
\begin{proof}
Part~(i) is part~(ii) of \Cref{cert:height-two-return}\certmark.

\textit{Idea for (ii).} By \eqref{tri:eq:centered-return}, the entry state
lies within $O(m\rho^{-3m})$ of the limiting segment, a suffix of length at
most six enlarges errors by a bounded factor, and the gain $\rho^{-3}$ per
scale beats the factor $m$. The constants below are crude.

Write $\varepsilon=\rho-x$, so that $\varepsilon\rho^{3m}=\tau\rho^{-3}\le9/80$.
As $183/100<\rho<46/25$, we have $\rho^3<(46/25)^3<25/4$; as $m\ge9$, also
$\varepsilon\le\frac9{80}\rho^{-27}<10^{-6}$, and thus $9/5\le x\le\rho$.

\emph{The block.} For $9/5\le\xi\le46/25$,
\[
 0\le\psi(\xi)\le\tfrac{11}2,\qquad
 0\le\psi'(\xi)=2\xi+\rho-1\le\tfrac{23}5,\qquad
 \tfrac{24}5\le\xi^3-1\le\tfrac{21}4,\qquad
 3\xi^2\le\tfrac{21}2 .
\]
Hence $0\le\eta\le(11/2)^2/(24/5)<7$, and
$\eta'=2\psi\psi'/(\xi^3-1)-3\xi^2\psi^2/(\xi^3-1)^2$ is a difference of
numbers in $[0,11]$ and $[0,14]$, so $|\eta'|\le14$. By the mean value theorem
$|\eta(\rho)-\eta(x)|\le14\varepsilon$, so
$|\rho\eta(\rho)-x\eta(x)|\le\rho|\eta(\rho)-\eta(x)|+\varepsilon\eta(x)
\le33\varepsilon$; also $|\eta+x\eta'|\le7+\frac{46}{25}\cdot14<33$. As
$\kappa\rho^3=2\rho\eta(\rho)$, \eqref{tri:eq:eta} gives
$\omega+\kappa\rho^3\varepsilon=2\varepsilon\bigl(\rho\eta(\rho)-x\eta(x)\bigr)
+\varepsilon^2(\eta+x\eta')(x)$, and therefore
\begin{equation}\label{tri:eq:delta-omega}
 |\Delta(x)|\le7\varepsilon^2,\qquad
 |\omega(x)+\kappa\rho^3\varepsilon|\le99\varepsilon^2 .
\end{equation}
Differentiating \eqref{tri:eq:fixed-state},
\[
 \bar S'=\frac{x-3}{(x-1)^3}-\frac{3x^2(x^3+5)}{(x^3-1)^3},\qquad
 \bar V'=-\frac2{(x-1)^3}-\frac{6x^2(x^3+2)}{(x^3-1)^3}.
\]
On $[9/5,46/25]$ we have $|x-3|\le6/5$, $(x-1)^3\ge64/125$,
$(x^3-1)^3>110$, $3x^2\le11$, $x^3+5\le12$, $6x^2\le21$ and $x^3+2\le9$, so
$|\bar S'|<\frac{75}{32}+\frac{132}{110}<4$ and
$|\bar V'|<\frac{125}{32}+\frac{189}{110}<6$. By the mean value theorem,
$\|\bar Y(x)-\bar Y(\rho)\|_\infty\le6\varepsilon$.

\emph{The entry error.} Since $\kappa\tau=\kappa\rho^3\varepsilon\rho^{3m}$,
for $i\in\{3m,3m+1\}$
\[
 x^{3m}(\omega+i\Delta)+\kappa\tau
 =x^{3m}(\omega+\kappa\rho^3\varepsilon)
 +\kappa\rho^3\varepsilon(\rho^{3m}-x^{3m})+i\,x^{3m}\Delta .
\]
Here $0<x^{3m}\le\rho^{3m}$, $0\le\rho^{3m}-x^{3m}\le3m\rho^{3m}\varepsilon$,
and $0<\kappa<4$, because $\kappa=16+52\rho-32\rho^2$ decreases in $\rho$ from
$2497/625$ at $183/100$ to $2088/625$ at $46/25$. By
\eqref{tri:eq:centered-return}, \eqref{tri:eq:delta-omega},
$\varepsilon\rho^{3m}\le9/80$ and $\rho^3<7$,
\begin{equation}\label{tri:eq:entry-error}
 \|Y_m(x)-Y_\infty(\tau)\|_\infty
 \le6\varepsilon+\rho^{3m}(106+21m+12m\rho^3)\varepsilon^2
 \le\frac{3+2m}{\rho^{3m}},
\end{equation}
because $6\cdot\frac9{80}+(106+105m)\bigl(\frac9{80}\bigr)^2<3+2m$.

\emph{The suffix.} Let $s=s_j$ and $1\le\ell=|s|\le6$. We apply
\Cref{lem:word-perturbation} with the matrix $J$ of
\eqref{tri:eq:unrolled}, $e=(2,1)$, $h=2$ and $q=\rho$.
The row sums give $C_\ell\le13$ and $B_j=j+2$, hence
$x^\ell C_\ell\le13(46/25)^6<505$.
For the limiting initial state put $Y=Y_\infty(\tau)$.
On $[183/100,46/25]$ both coordinates of $Y_*$ decrease in
$\rho$, from $6133/10^4$ and $19622/10^4$ to $327/625$ and $1193/625$, so
$Y_*\in(0,1)\times(0,2)$; with $0\le\kappa\tau\le4\tau_+<3$ this gives
$\|Y\|_\infty\le3$. Since $j\rho^{j-1}<130$ for $1\le j\le6$,
the parameter term in \eqref{eq:word-parameter-error} is bounded by
\[
 130\left(13\cdot3+2\sum_{j=1}^{5}(j+2)\right)\varepsilon
 \le93\cdot130\,\varepsilon
 =12090\,\varepsilon<\frac{1361}{\rho^{3m}} .
\]

\emph{Conclusion.} Substituting these bounds and
\eqref{tri:eq:entry-error} into \eqref{eq:word-parameter-error}, then
using \eqref{tri:eq:norm-factor}, gives
\[
 \bigl\|\Theta_{x,s}Y_m(x)-\Theta_{\rho,s}Y_\infty(\tau)\bigr\|_Q
 \le\frac{166}{105}\cdot\frac{505(3+2m)+1361}{\rho^{3m}} .
\]
At $m=9$, as $\rho>9/5$, this is less than
$\frac{166}{105}\cdot11966\cdot(5/9)^{27}<0.0025<9/2000$. Passing from $m$ to
$m+1$ multiplies the numerator by less than $2$ and the denominator by
$\rho^3>2$, so the bound holds for every $m\ge9$. Finally $9/2000<\mu$, so
under \eqref{tri:eq:phase-margin} we get
$\|\Theta_{x,s_j}Y_m(x)\|_Q<1-\mu+\mu=1$.
\end{proof}

\begin{proposition}[Tribonacci entry\certmark]\label{tri:prop:phase-entry}

For every
\[
 x\in\left[\rho-\frac9{80},\rho\right),
\]
there is a finite digit word $w$ such that $\Theta_{x,w}(2,1)\in Q$.
\end{proposition}
\begin{proof}
Let $m$ be the entry scale of $x$ and $\tau=\tau_m(x)\in[\tau_-,\tau_+]$, so
that $x=\rho-\tau\rho^{-3m-3}$. If $m\le8$, take $w=ab^ms$ with $s$ from
\Cref{tri:lem:finite-phase-estimates}(i). If $m\ge9$, choose $j$ by
\eqref{tri:eq:phase-margin} and take $w=ab^ms_j$; then
\Cref{tri:lem:finite-phase-estimates}(ii) gives $\Theta_{x,w}(2,1)\in Q$.
\end{proof}

\subsubsection{Completion}
\begin{proposition}[The height-two active trace\certmark]\label{low:prop:ternary}

One has $[1/\rho,1)\subseteq\mathcal A_2\subseteq\mathcal D_2$. Equivalently,
$(1,\rho]\subseteq\mathcal C_2^{\mathrm{ext}}$.
\end{proposition}
\begin{proof}
\textit{Idea.} An entry word followed by return words gives a bounded orbit,
hence a double zero; an interval of double zeros is active.

Let $433/250\le x<\rho$. Then $x\ge\rho-9/80$, as $\rho<46/25$, and
\Cref{tri:prop:phase-entry} gives a word $w_0$ with
$Y^{(0)}=\Theta_{x,w_0}(2,1)\in Q$. As $x<\rho<4598217/2500000$,
\Cref{tri:prop:core} assigns to each $Y\in Q$ a word $w(Y)$ with
$1\le|w(Y)|\le10$ and $\Theta_{x,w(Y)}Y\in Q$. Put
$Y^{(k+1)}=\Theta_{x,w(Y^{(k)})}Y^{(k)}$ and follow the infinite digit sequence
$w_0\,w(Y^{(0)})\,w(Y^{(1)})\cdots$. Its orbit from $(2,1)$ consists of the
states along $w_0$, the return states $Y^{(k)}\in Q$, and states
$\Theta_{x,v}Y^{(k)}$ with $|v|\le9$. On $Q$, $\|Y\|_\infty\le3/2$, and if
$\|Y\|_\infty\le M$ with $M\ge1$, then
$\|\Theta_{x,d}Y\|_\infty\le3xM+4\le10M$, as $x<2$. So every state after
$w_0$ has $\|\cdot\|_\infty\le10^9\cdot3/2$. The orbit is bounded, and so is
its $(P,V)$-orbit from $(1,1)$; by \Cref{lem:bounded-orbit},
$1/x\in\mathcal D_2$.

At $x=\rho$, the digits $a,b,b,\ldots$ of $f_\rho$ give the state
$Y_m(\rho)=\bar Y(\rho)$ after every block, and two further states inside
each block. This orbit is bounded, so $1/\rho\in\mathcal D_2$, with witness $f_\rho$.
For $1<x<433/250$ we have $x<\sqrt3$, since $433^2=187489<187500=3\cdot250^2$,
and $1/x\in\mathcal D_2$ by \Cref{cor:interior-double-zero}.

Thus $[1/\rho,1)\subseteq\mathcal D_2$. This interval is nondegenerate, so
Lemma~\ref{low:lem:countable} gives $[1/\rho,1)\subseteq\mathcal A_2$, including
the endpoint $1/\rho$, and $\mathcal A_2\subseteq\mathcal D_2$. By the
definition of $\mathcal A_2$, $(1,\rho]\subseteq\mathcal C_2$, and these
points lie in $X$.
\end{proof}

\Cref{tri:thm:double-zeros} combines this inclusion with the exclusion proved
in \Cref{sec:height-two}.

\paragraph{Finite checks delegated to \Cref{cert:height-two-return}.}
The closed return cover, finite entry scales $0\le m\le8$, limiting
phase cover and algebraic bracket for $\rho$ are the four exact obligations
in that record. The entry and phase checks are made at $r=\rho$ over
$\mathbb Q(\rho)$. The estimates for all scales $m\ge9$ are proved above.

\section{Residual planar separators}\label{sec:sep:framework}
This section proves Proposition~\ref{prop:residual-separators}, which
supplies the disconnectedness cases outside the ranges of
Theorem~\ref{root:thm:infinite}. The required heights are
$q^2-2\le h\le q^2$ for $q\ge3$, together with $h=1,3,4$.

\subsection{Necessary bounds and exact exclusion}
By Proposition~\ref{prop:global-monic-reduction}, $c\in X$ lies in
$\mathcal Z_h^{\mathrm{ext}}$ if and only if
$c\in T_h(c)=\{\sum_{j\ge0}b_jc^{-j}:b_j\in B_h\}$, a bounded set. An
\emph{inverse orbit} $z_0=c$, $z_j=c(z_{j-1}-d_j)$, $d_j\in B_h$, is
\emph{admissible} if all $z_j$ lie in $T_h(c)$. As
$c(z-b_0)=\sum_{j\ge0}b_{j+1}c^{-j}$ for $z=\sum_{j\ge0}b_jc^{-j}$, and
$c=\sum_{j=1}^md_jc^{1-j}+c^{-m}z_m$, the point $c$ lies in
$\mathcal Z_h^{\mathrm{ext}}$ exactly when it has an admissible inverse orbit.

Let $c=x+iy$ with $x>1$, $y\ne0$ and $R=|c|^2$, and put
$\mathsf V_c(z)=\im z/y$ and $\mathsf S_c(z)=\im(cz)/y$ in the dynamical
plane; the marked point $z=c$ has $(\mathsf S_0,\mathsf V_0)=(2x,1)$. The
orbit point $z_j$ is half the point $z_u(c)$ defined in Section~\ref{sec:interfaces} for
$u=(d_1,\ldots,d_j)$, so $\mathsf V_c(z_j)=V_c(z_u(c))$ and
$\mathsf S_c(z_j)=|c|\,S_c(z_u(c))$ in the notation of \eqref{eq:reim-vs}. If
$z'=c(z-d)$, then $\im z'=\im(cz)-dy$ and $cz'=c^2(z-d)=2xz'-R(z-d)$, as
$c^2=2xc-R$. Hence the inverse digit map acts by
\begin{equation}\label{sep:eq:recurrence}
 \mathsf V_{j+1}=\mathsf S_j-d_{j+1},\qquad
 \mathsf S_{j+1}=2x\mathsf V_{j+1}-R\mathsf V_j.
\end{equation}
For a fixed word these are polynomials in $(x,R)$, which extend to $y=0$ with
$R=x^2$.

\begin{lemma}[Outer support bounds]\label{sep:lem:outer}

For $z\in T_h(c)$,
\[
 |\mathsf V_c(z)|\le\beta_c,\qquad |\mathsf S_c(z)|\le h+\beta_c,\qquad
 \beta_c=\frac h{(|c|-1)^2}.
\]
Consequently an admissible inverse orbit satisfies
\begin{equation}\label{sep:eq:coarse}
 |\mathsf V_j|\le\beta:=\frac h{(x-1)^2},\qquad |\mathsf S_j|\le h+\beta.
\end{equation}
Moreover, $c\in\mathcal Z_h^{\mathrm{ext}}\setminus\R$ implies
$|c|\le1+\sqrt h$.
\end{lemma}
\begin{proof}
With $w=1/c$, the identity $w^k-\bar w^k=(w-\bar w)\sum_{\ell=0}^{k-1}w^\ell
\bar w^{k-1-\ell}$ gives $|\im(c^{-k})|\le k|y||c|^{-k-1}$. Summing over the
digits of $z$, with $\sum_{k\ge1}k|c|^{-k-1}=(|c|-1)^{-2}$, gives the
$\mathsf V$ bound and the $\mathsf S$ bound, in which $b_0$ adds
$b_0\im c/y=b_0$. As $|c|\ge x$, \eqref{sep:eq:coarse} follows. Finally,
$\mathsf V_0=1$ in an admissible orbit, so $1\le h/(|c|-1)^2$.
\end{proof}

Strict inequalities \emph{exclude all admissible continuations through $m$
steps} at $c_0$ if each $(d_1,\ldots,d_m)\in B_h^m$ has some $0\le j\le m$
with $|\mathsf V_j|>\beta_c$ or $|\mathsf S_j|>h+\beta$, the coordinates
\eqref{sep:eq:recurrence} being taken at $(x,R)=(\re c_0,|c_0|^2)$; $\beta$
may replace $\beta_c$ in the first test.

\begin{lemma}[Strict finite exclusion and real contacts]\label{sep:lem:strict}

Suppose strict inequalities in the continuous variables $(x,R)$ exclude all
admissible continuations of \eqref{sep:eq:recurrence} through finitely many
steps at a parameter $c_0$, with $\re c_0>1$. Then a neighbourhood of $c_0$
contains no nonreal points of $\mathcal Z_h^{\mathrm{ext}}$. In particular,
$c_0\notin\mathcal C_h^{\mathrm{ext}}$. This remains true when $c_0$ is real.
\end{lemma}
\begin{proof}
Finitely many words are involved, and the excluding inequalities are strict
and continuous in $(x,R)$ for $x>1$, so they persist at $(\re c,|c|^2)$ for
$c$ near $c_0$, with $R=x^2$ if $c_0$ is real. At nonreal such $c$,
\Cref{sep:lem:outer} is necessary for an admissible inverse orbit, so none
exists. As
$\mathcal C_h^{\mathrm{ext}}=\operatorname{cl}_X(\mathcal Z_h^{\mathrm{ext}}\setminus\R)$ by
\eqref{eq:root-ifs}, $c_0\notin\mathcal C_h^{\mathrm{ext}}$.
\end{proof}

We call $R+2x<N$ the \emph{lens condition}; in these coordinates
\Cref{prop:strict-trace} reads as follows.

\begin{lemma}[Strict inner trap]\label{sep:lem:inner}

Put $N=2h+1$. Assume $R+2x<N$, and define
\[
 \tau_h(c)=\frac{N-2x}{2R}>\frac12.
\]
If a finite word brings the initial state $(2x,1)$ into
\begin{equation}\label{sep:eq:inner}
 |\mathsf V|<\tau_h(c),\qquad |\mathsf S|<N/2,
\end{equation}
then $c$ is an interior point of $\mathcal Z_h^{\mathrm{ext}}$. The
conclusion also applies to a real $c=x>1$ when the continued finite-word
coordinates satisfy the strict inequalities.
\end{lemma}
\begin{proof}
After a word $u$ the orbit point is $cP_u(c)=\frac12z_u(c)$, so, in the
coordinates of \Cref{lem:continuity}, $\mathsf V=V_c(z_u(c))$ and
$\mathsf S=|c|S_c(z_u(c))$, and these are the inequalities of
\Cref{prop:strict-trace}. At a real $x$ the continued values are $V_u(x)$ and
$xS_u(x)$ (\Cref{lem:continuity}), $\tau_h(x)$ is \eqref{eq:tau-sigma}, and
$|\mathsf S|<N/2$ reads $|S_u(x)|<\sigma_h(x)$: this is the real case.
\end{proof}

\begin{lemma}[Separation by an excluded boundary]\label{lem:separation-principle}
Let $K=\mathcal C_h^{\mathrm{ext}}$. If an open set $U$ meets $K$,
$\partial U\cap K=\varnothing$, and
$\overline U\subset\{\operatorname{Re}c>0\}$, then $K$ has at least
three components.
\end{lemma}
\begin{proof}
The sets $K\cap U$ and $K\cap(-U)$ are disjoint nonempty relatively
clopen sets, since $K\cap U=K\cap\overline U$. A point $iy$ of the
central annulus, $1<y<\sqrt{h+1}$, lies in neither. Every component lies
in one of these three disjoint relatively clopen sets.
\end{proof}

\subsection{An elementary separator at height seven}
\begin{proposition}[Height seven]\label{low:prop:eight}

The lines $\re c=\pm16/5$ miss $\mathcal C_7^{\mathrm{ext}}$. The exterior
core $\mathcal C_7^{\mathrm{ext}}$ has at least three connected components.

\end{proposition}
\begin{proof}
\textit{Idea.} Only the first digits $5$ and $6$ survive, and the second
digit fails a bound, except off the axis, where $\beta_c$ rejects $5$.

Set $x=16/5$, so $R=256/25+y^2$, $\beta=7/(x-1)^2=175/121$ and
$7+\beta=1022/121$; the bound $\beta_c=7/(|c|-1)^2$ of \Cref{sep:lem:outer}
satisfies $\beta_c\le\beta$. By \eqref{sep:eq:recurrence},
$\mathsf V_1=32/5-d_1$ and $\mathsf S_1=768/25-32d_1/5-y^2$; digits
$d_1\le4$ give $\mathsf V_1\ge12/5>175/121$, and $d_1\ge7$ give
$\mathsf S_1\le-352/25-y^2<-1022/121$.
If $d_1=6$, then $\mathsf V_1=2/5$, $\mathsf S_1=-192/25-y^2$, and as
$d_2\ge-7$,
\[
 \mathsf S_2=-\frac{6656}{125}-\frac{32}5d_2-\frac{34}5y^2
 \le-\frac{1056}{125}-\frac{34}{5}y^2,
 \qquad \frac{1056}{125}-\frac{1022}{121}=\frac{26}{15125}>0.
\]
If $d_1=5$, then $\mathsf V_1=7/5\le\beta_c$ forces $(|c|-1)^2\le5$, which
fails strictly if $y^2\ge161/625$: then $|c|^2\ge6561/625$, so
$|c|-1\ge56/25$ and $(56/25)^2=3136/625>5$. Hence $y^2<161/625$, and as
$\mathsf S_1=-32/25-y^2$, digits $d_2\le-3$ give
$\mathsf V_2=-32/25-y^2-d_2>914/625>175/121$, and $d_2\ge-2$ give
\[
 \mathsf S_2=-\frac{2816}{125}-\frac{32}5d_2-\frac{39}5y^2
 \le-\frac{1216}{125}<-\frac{1022}{121}.
\]
So a strict inequality excludes every two-digit word, and by
\Cref{sep:lem:strict} the line misses $\mathcal C_7^{\mathrm{ext}}$;
negation gives $\re c=-16/5$.

The normalized series
\[
 1-6z+\frac{7z^2}{1-z}=\frac{1-7z+13z^2}{1-z}
\]
has height seven and the nonreal zero $\gamma_3^{-1}$, where
$\gamma_3=(7+i\sqrt3)/2$. Polynomial truncation and inversion give
$\gamma_3\in\mathcal C_7^{\mathrm{ext}}$.
The two excluded vertical lines leave nonempty clopen portions containing
$\gamma_3$, $-\gamma_3$, and the central annulus, respectively.
\end{proof}

\subsection{Uniform and finite near-square domains}
\begin{proof}[Proof of Proposition~\ref{prop:residual-separators}]
For $q\ge22$, Proposition~\ref{cert:lower-block} gives the domain
$\Omega_q=\Phi(D_q)$, where
$\Phi(v,Y)=X_q(v,Y^2)+iY/\sqrt q$. Its boundary misses
$\mathcal C_{q^2}$ and its interior meets $\mathcal C_{q^2-2}$.
Since $\partial X_q/\partial v=(q-2)/(2q^3)>0$, the map is a
homeomorphism of the plane and takes $\partial D_q$ to
$\partial\Omega_q$. The defining bounds give
$\overline\Omega_q\subset\{\operatorname{Re}c>q\}$.
Nesting and the separation lemma prove all three heights.
For $5\le q\le21$, Proposition~\ref{cert:twenty-domains} gives the
corresponding bounded domains with closure in the positive half-plane.
For $q=3,4$, height seven was proved above and
Proposition~\ref{cert:small-heights} supplies heights
$8,9,14,15,16$. The same proposition supplies $1,3,4$.
Every supplied domain has an excluded whole boundary and a certified
core point inside it, so the separation lemma applies in each case.
\end{proof}

\section{Exact finite certificates}\label{app:certified}
Every record below states a finite exact obligation used in the preceding
proofs. The acceptance rules of the checks are proved here. The integer and
rational coefficients, cell covers and finite word lists that the checks
consume are supplied in the computational companion~\cite{Companion26},
described in Section~\ref{cert:sec:supplement}. All estimates used at
infinitely many scales are proved in Section~\ref{sec:complete-traces};
a finite numerical cutoff does not replace them. Every finite acceptance
decision uses exact arithmetic.

\begin{table}[htbp]\centering\small
\caption{Proof obligations and their uses.}\label{cert:tab:records}
\begin{tabular}{@{}ll@{}}\toprule
Record & Used for\\\midrule
Short-word strips & Pronic bulk-to-terminal reduction\\
Compact windows & Two remaining pronic windows\\
Terminal finite bands and seams & Finite complement of all-scale estimates\\
Height-six returns and entry & Complete trace at height six\\
Height-two return and entry & Tribonacci trace completion\\
Uniform near-square domains & $q\ge22$ residual separators\\
Finite near-square domains & $5\le q\le21$ residual separators\\
Small bounded domains & $h=1,3,4,8,9,14,15,16$\\
Square-adjacent inequalities, inner prefixes, witnesses & Square-adjacent heights, Section~\ref{sec:gap-heights}\\
Surviving prefixes & Lower bounds for the anchor sets, Section~\ref{sec:remaining-scope}\\
Minimum, slices, constants, moves, landing & Real traces, Section~\ref{sec:real-traces}\\\bottomrule
\end{tabular}\end{table}
\FloatBarrier
\subsection{Reduction lemmas}\label{cert:sec:conventions}

\begin{lemma}[Positivity on a half-line times a box]\label{sep:lem:bernstein}

Let $P(q,t_1,\ldots,t_r)$ be a rational polynomial. After setting $q=q_0+s$,
write
\[
 P(q_0+s,t)=\sum_{k=0}^{d}s^k
       \sum_{\alpha\le m} b_{k,\alpha}
       \prod_{j=1}^r\binom{m_j}{\alpha_j}
       t_j^{\alpha_j}(1-t_j)^{m_j-\alpha_j}.
\]
If every $b_{k,\alpha}\ge0$ and every $b_{0,\alpha}>0$, then
$P>0$ for all $q\ge q_0$ and $t\in[0,1]^r$.
\end{lemma}
\begin{proof}
The Bernstein basis functions are nonnegative on the box and sum to one, so the
$s^0$ term is at least $\min_\alpha b_{0,\alpha}>0$; the other terms are
nonnegative for $s\ge0$.
\end{proof}
In one variable, $p(t)=\sum_ja_jt^j$ may instead be bounded through its
coefficients:
\begin{equation}\label{cert:eq:coefficient-sum}
 p(t)\ge a_0+\sum_{j\ge1}\min\{a_j,0\}\quad(0\le t\le1),\qquad
 p(t)\ge a_0-\sum_{j\ge1}|a_j|\quad(-1\le t\le1).
\end{equation}
Before either test, margins are multiplied only by explicitly positive factors.

\Needspace{11\baselineskip}
For a word $u$ over $B_h$, let $(\mathsf S_u,\mathsf V_u)$ be the
coordinates \eqref{sep:eq:recurrence} after the digits of $u$, starting from
$(2x,1)$; the extensions of $u$ are the words $ud$, $d\in B_h$.

\begin{lemma}[Digit-tree exclusion]\label{cert:lem:tree}

Let $h\ge1$, let $\Lambda\subset\{\re c>1\}$, and let $\mathcal T$ be a finite
set of words over $B_h$ that contains the empty word and every prefix of its
members. Suppose that at every $c\in\Lambda$, with $(x,R)=(\re c,|c|^2)$ and
$\beta=h/(x-1)^2$, every $u\in\mathcal T$ satisfies:
\begin{enumerate}
\item if no extension of $u$ lies in $\mathcal T$, then
 $|\mathsf V_u|>\beta$ or $|\mathsf S_u|>h+\beta$;
\item otherwise, $|\mathsf S_u-d|>\beta$ for every $d\in B_h$ with
 $ud\notin\mathcal T$.
\end{enumerate}
Then $\Lambda\cap\mathcal C_h^{\mathrm{ext}}=\varnothing$.
\end{lemma}
\begin{proof}
Fix $c\in\Lambda$ and $(d_1,\ldots,d_m)\in B_h^m$, where $m$ is the largest
length of a word in $\mathcal T$, and let $u$ be its longest prefix in
$\mathcal T$. If $u$ has no extension in $\mathcal T$, (i) violates
\eqref{sep:eq:coarse} strictly; otherwise the next digit $d$ has
$ud\notin\mathcal T$, and (ii) gives $|\mathsf V_{ud}|=|\mathsf S_u-d|>\beta$.
So strict inequalities exclude every admissible continuation through $m$
steps, and \Cref{sep:lem:strict} gives $c\notin\mathcal C_h^{\mathrm{ext}}$.
\end{proof}

In the records, $\Lambda$ is a closed cell of an arc $t\mapsto x(t)+iy(t)$,
$t\in I$, on which $x$ and $y^2$, hence $R$ and $(x-1)^2$, are rational
polynomials. After multiplication by $(x-1)^2$, (i) is the positivity on $I$
of one of the \emph{exit polynomials} $\pm\mathsf V_u(x-1)^2-h$ and
$\pm\mathsf S_u(x-1)^2-h((x-1)^2+1)$. For (ii), let the extensions of $u$ in
$\mathcal T$ be the $ud$ with $\mathrm{lo}\le d<\mathrm{up}$. One may check
$(\mathsf S_u-\mathrm{lo}+1)(x-1)^2-h>0$ if $\mathrm{lo}>-h$ and
$(\mathrm{up}-\mathsf S_u)(x-1)^2-h>0$ if $\mathrm{up}\le h$. The checks take
$[\mathrm{lo},\mathrm{up}-1]=[-h,h]\cap[\lceil S_*-h/H_*\rceil,
\lfloor S^*+h/H_*\rfloor]$, where $S_*\le\mathsf S_u\le S^*$ and
$(x-1)^2\ge H_*>0$ on $I$. Positivity on $I$ is certified by
\Cref{sep:lem:bernstein} or \eqref{cert:eq:coefficient-sum}, after an affine
change of the cell variable.

\begin{lemma}[Rouch\'e discs]\label{cert:lem:rouche}

Let $P$ be a polynomial, $z_0\in\C$, $\varepsilon>0$,
$P(z_0+w)=\sum_kA_kw^k$ and $|a|_1=|\re a|+|\im a|$. If
$|A_0|_1+\sum_{k\ge2}|A_k|_1\varepsilon^k
<\varepsilon\max\{|\re A_1|,|\im A_1|\}$, then $P$ has exactly one zero in
$|z-z_0|<\varepsilon$, and it is simple; if moreover $P\in\mathcal P_h$ and
the closed disc lies in $\{\im c>0\}$, this zero belongs to $\mathcal C_h$.
\end{lemma}
\begin{proof}
As $|a|\le|a|_1$ and $|A_1|\ge\max\{|\re A_1|,|\im A_1|\}$, the
hypothesis gives $|P(z_0+w)-A_1w|<|A_1w|$ on $|w|=\varepsilon$;
apply Rouch\'e's theorem. A zero in the upper half-plane lies in
$\mathcal R_h\setminus\R\subset\mathcal C_h$.
\end{proof}

\subsection{The records}\label{cert:sec:table}

\begin{proposition}[Certified: short-word strip]\label{cert:short-word-strip}

Let $q\ge4$ be real, read the words $\mathsf A_m,\mathsf B_m,\mathsf U_{m,\ell}$
of \Cref{sec:residual} with real digits, and put $a=q-1$, $b=(q-1)(q-2)$,
$h=q(q-1)$ and $N=2h+1$. For each word $u$ below, each $\delta$ in its interval
and $x=q-\delta$, the four margins
\begin{equation}\label{cert:eq:four-margins}
 N-2x\pm2x^2V_u(x),\qquad N\pm2xS_u(x)
\end{equation}
are positive, with $V_u,S_u$ as in \eqref{eq:VS-def}:
\[
\renewcommand{\arraystretch}{1.15}
\begin{array}{@{}l@{\ \ }l@{\qquad}l@{\ \ }l@{}}
 \mathsf A_2 & [q^{-2}-q^{-3},\,(2q)^{-1}-(4q^2)^{-1}], &
 \mathsf B_2 & [(2q)^{-1}+(2q^2)^{-1},\,q^{-1}],\\
 \mathsf A_3 & [q^{-3},\,(2q^2)^{-1}-(3q^3)^{-1}], &
 \mathsf B_3 & [(2q^2)^{-1}+(2q^3)^{-1},\,q^{-2}],\\
 \mathsf A_2(-h) & [(2q)^{-1}-(4q^2)^{-1},\,(2q)^{-1}-3(4q^3)^{-1}], &
 \mathsf U_{2,1} & [(2q)^{-1}+3(4q^3)^{-1},\,(2q)^{-1}+(2q^2)^{-1}].
\end{array}
\]
\end{proposition}
\begin{certificate}
Put $\delta=\delta_-+(\delta_+-\delta_-)t$ with $0\le t\le1$. Each margin,
multiplied by a positive constant times a power of $q$ and written in
$q=4+s$, is a polynomial in $(s,t)$, positive by \Cref{sep:lem:bernstein} with
$q_0=4$ and $r=1$. The checker computes its Bernstein coefficients, whose constant
row satisfies $b_{0i}\ge p(0,0)\cdot4231364608/127302260011>p(0,0)/32$. \emph{Spot check:} for
$q=4$, $u=\mathsf B_2=(7,-10)$ and $\delta=1/4$, so $x=15/4$: $P_u(x)=-35/16$,
$V_u(x)=-5/16$, $S_u(x)=-5/2$ and $N=25$, and the four margins are
$1115/128$, $3365/128$, $25/4$ and $175/4$.
Two independent programs compute the Bernstein coefficients, one from the
stored polynomials and one by rebuilding them from the recurrence.
\companioncheck{cert:short-word-strip}
\end{certificate}

\begin{proposition}[Certified: compact windows]\label{cert:compact-windows}

Let $q\ge4$ be an integer. Every point of the windows $J_{q,2}$ and $J_{q,3}$
of \Cref{prop:compact-residual-closure} is an interior point of
$\mathcal C_{q(q-1)}$.
\end{proposition}
\begin{certificate}
Each window is covered by consecutive closed cells, each with a word over
$B_{q(q-1)}$ whose margins \eqref{cert:eq:four-margins} are positive on the
cell, where also $x^2+2x<N$. \Cref{prop:strict-trace} gives an interior point
of $\mathcal Z_{q(q-1)}^{\mathrm{ext}}$, and
$\interior\mathcal Z_h=\interior\mathcal C_h$, since an open disc in
$\mathcal Z_h$ minus $\R$ lies in the closed set $\mathcal C_h$. For
$J_{q,2}$ with $q\le43$ and $J_{q,3}$ with $q\in\{4,\ldots,11,13,15,\ldots,
23\}$ the cells are rational intervals, positive by
exact Bernstein coefficients. For the remaining $q$, rational cells handle $44\le q\le1000$
for $J_{q,2}$, followed by symbolic cells on the residue classes
$q=q_0+4k$, $q_0=1001,\ldots,1004$. For $J_{q,3}$, symbolic cells start
at $q_0=12,25,14,27$. Their endpoints and admissible digits depend
polynomially on $q$. Substituting $k=y/(1-y)$ and clearing positive
denominators gives nonnegative tensor Bernstein coefficients with positive
edge coefficients, hence strict positivity for every finite $k\ge0$.
\emph{Spot check:} for
$q=4$, $J_{4,2}=[989/256,995/256]$, and its last cell $[7939/2048,995/256]$
carries the word $(7,-10,-5)$. At $x=97/25$,
$(P_u,V_u)=(-49527/15625,1592/15625)$, $S_u=-9587/3125$ and
$x^2+2x=14259/625<25=N$; the margins are $198317631/9765625$,
$138401119/9765625$, $93247/78125$ and $3813003/78125$.
The cells, words and Bernstein coefficients of every rational and symbolic
case are stored and rechecked.
\companioncheck{cert:compact-windows}
\end{certificate}

\begin{proposition}[Certified: finite terminal checks]\label{cert:terminal-finite}

Let $q\ge4$ be an integer, $h=q(q-1)$, $a=q-1$, $b=(q-1)(q-2)$, and let
$\mathcal B_{q,m}$ be the bands of \eqref{eq:terminal-centre}.
\begin{enumerate}
\item For every pair $(q,m)$ of the finite range
 \eqref{eq:terminal-finite-range}, every point of $\mathcal B_{q,m}$ is an
 interior point of $\mathcal C_h$.
\item Every $x$ with $q-q^{-3}\le x\le q-(1/2+3/(4q^4))q^{-3}$ is an
 interior point of $\mathcal C_h$.
\item With the movements \eqref{eq:terminal-endpoint-movements} and the further
 allowances below, the overlap and containment inequalities in the proof of
 \Cref{prop:terminal-central-hierarchy} hold: for $q\ge8$ (at $m=11$ and the
 small-depth thresholds) those of the residue chains, fixed seams and annulus
 targets, the targets of $\mathsf W_{m,m-2}$ and $\mathsf U_{m+1,m}$ shrunk
 further by $D_m$ and $qD_m$, with $\mathsf R_{m,1,0}$ and $\mathsf R_{m,0,0}$
 strictly capturing $T=-1/2$ after their movements; for $4\le q\le7$ (at
 $m=11$) those of \Cref{tab:low-q-chains} and of the seams
 $\mathsf W_{m,m-3}\leadsto\mathsf W_{m,m-2}\leadsto\mathsf W_{m,m-1}$ and
 $\mathsf P_{m+1}\leadsto\mathsf U_{m+1,m}\leadsto\mathsf U_{m+1,m-1}$, whose
 $\mathsf W$ and $\mathsf U$ cells take the single allowance $1000D_{11}$ and
 whose $\mathsf P$ cell takes $\eta_{P,\le2}$ of
 \eqref{eq:terminal-endpoint-movements}. In both cases the $\mathsf W_{m,m-1}/\mathsf R_{m,1,0}$ overlap, so moved, meets
 $T>-2/3$.
\end{enumerate}
\end{proposition}
\begin{certificate}
(i) In the finite range $m\le10$, so $(q-x)x^{m-1}$ and $(q-x)x^m$ decrease in
$x$ on $[q-q^{-3},q)$, with derivatives $x^{m-2}((m-1)(q-x)-x)<0$ and
$x^{m-1}(m(q-x)-x)<0$; exact endpoint tests enclose $\mathcal B_{q,m}$ in a rational interval, covered by
consecutive closed cells, each with a word over $B_h$ whose margins
\eqref{cert:eq:four-margins} are positive, and \Cref{prop:strict-trace}
applies as for \Cref{cert:compact-windows}. (ii) The same for $4\le q\le11$,
with the extra word $\mathsf B_4(-(a^2+1))$ at $q=4$ only; for $q\ge12$, the
words $\mathsf U_{4,j}$, $0\le j\le4$, on five closed cells of the phase
$q^3(q-x)$, whose cleared margins are positive for all $q\ge12$ by
\Cref{sep:lem:bernstein}. (iii) Each inequality compares endpoints of limiting
intervals, rational in $q$, moved by the allowances; after clearing positive
denominators it has nonnegative coefficients in $q-q_0$ and a positive
constant term on each residue class of $q$ modulo $4$ from the thresholds of
\Cref{prop:terminal-central-hierarchy}; for $4\le q\le7$ the comparisons are
rational, at $m=11$. \emph{Spot check} (iii): at $q=4$ the limiting intervals of
$\mathsf W_{m,m-1}$ and $\mathsf R_{m,1,0}$ are $(-313/384,-55/96)$ and
$(-259/384,-751/1536)$, so their overlap $(-259/384,-220/384)$ meets $T>-2/3=-256/384$
and lies below $-1/2$.
The checks of (i) run over all pairs of \eqref{eq:terminal-finite-range},
those of (ii) over rational cells for $q\le11$ and symbolic cells for
$q\ge12$, and those of (iii) inequality by inequality on each residue class.
For $4\le q\le7$ the checker derives all $42$ printed intervals and tests
the four $\mathsf W_{m,m-1}/\mathsf R_{m,1,0}$ domain seams, $60$ formal and
$58$ domain-clipped consecutive overlaps, and the digit, endpoint, allowance
and monotonicity conditions: $389$ exact rational comparisons in all.
\companioncheck{cert:terminal-finite}
\end{certificate}

\begin{proposition}[Certified: finite checks at height six]\label{cert:height-six}

Let the digits $d\in B_6$ act on states by the maps $F_{x,d}$ of
\eqref{eq:PV-map}, write $F_{x,u}$ for the action of a word $u$, and put
$r=x-1$.
\begin{enumerate}
\item There is a finite cover of $[137/50,2999/1000]$ by consecutive
closed rational parameter bands $I$.  For each band there are a nonempty
finite union $Y_I\subset[-4,4]\times[-2,2]$ of closed rational grid cells,
a positive rational number $\eta_I$, and a digit $d_C\in B_6$ for each
selected cell $C$, with
\[
 F_{x,d_C}(C)+[-\eta_I,\eta_I]^2\subset Y_I\qquad(x\in I).
\]
Moreover, each $I$ has a finite partition into closed rational subintervals,
each supplied with a word $v$ of length at most $11$ satisfying
\[
 F_{x,v}(1,1)+[-\eta_I,\eta_I]^2\subset Y_I
\]
throughout that subinterval.
\item For every $(u,w)$ in the set $\mathcal Q$ of
 \Cref{seven:lem:paired-cover} and $2999/1000\le x\le3$, one of the two states
 there has a word of length at most three whose image satisfies the two
 bounds displayed there.
\item Let $x=3$. The eighteen cells of \Cref{seven:tab:central-cells} cover
 $[-4/5,0]$. Let a cell have the word $\mathsf P_M$, $\mathsf R_{M,\ell}$ or
 $\mathsf C_M$ followed by $s$, and let $Y(\Theta)$ be the limiting state
 $\mathsf P$, $\mathsf R_\ell$ or $\mathsf C$ of
 \eqref{seven:eq:central-limits}. For $\Theta$ in the cell, the images of
 $Y(\Theta)$ under all prefixes of $s$, the empty one included, lie in
 $[-3,3]\times[-3/2,3/2]$, and $F_{3,s}(Y(\Theta))$ in
 $[-11/8,11/8]\times[-39/80,39/80]$.
\item For $2999/1000\le x\le344373767/114791256$, some word maps $(1,1)$ to a
 state with $|P|<3/r$ and $|V|<2/r^2$.
\end{enumerate}
\end{proposition}
\begin{certificate}
(i) For the lower bands, a cell indexed by $(i,j)$ is
$[(i-4N)/N,(i+1-4N)/N]\times[(j-2N)/N,(j+1-2N)/N]$;
the stored integer array selects the cells and their digits.  For the upper
bands the same grid is encoded by consecutive vertical strips, with signed
indices.  The two coordinates of
$F_{x,d}(P,V)=(xP-d,x(P+V)-d)$ have extrema at the endpoints of the parameter
band and the source rectangle.  The verifier rounds these extrema outwards,
adds one grid unit in each direction, and checks that every resulting target
cell is selected.  Thus the displayed inclusion holds with $\eta_I=1/N$.
For the entry words, direct rational Bernstein bounds enclose $P_v(x)$ and
$V_v(x)$ on each closed entry interval; the same grid-unit enlargement is
contained in the selected union.  The entry intervals cover their parameter
bands, whose exact endpoints join from $137/50$ to $2999/1000$.
All integer operations are checked against their storage bounds before
fixed-width arithmetic is used.
(ii) The paired-cone certificate consists of sixty branch--word regions, each
defined by four affine halfspaces in $(u,w)$.  The checker reconstructs each
halfspace from its branch and its word of length at most three.  After
multiplication by $x^3(x-1)^2$, these are polynomial identities of degree at
most six; agreement at seven distinct rational parameters verifies them
identically.  At the rational reference parameter $x=5999/2000$, it reconstructs the polygons
by intersecting the halfspaces and checks their union by a horizontal sweep.
The exact signs of all parallel-line comparisons and triple-line determinants
remain unchanged for $0<3-x<1/1000$: the substitution
$3-x=t/(1000(1+t))$, followed by clearing the positive denominator, gives
coefficients all of one sign, with at least one nonzero coefficient.  The face incidences and halfspace signs of the line arrangement are
therefore constant: a change would force one of these determinants or
parallel-line comparisons to vanish. Coverage by the same fixed
combinations of halfspaces thus persists throughout the open interval.  The two endpoints follow by taking limits in the
finite collection of closed regions.
(iii) States
and images are affine in $\Theta$, and each inequality is checked at both
ends of its cell. (iv) On each closed rational entry interval the checker reconstructs $P_v$ from
the digit word and $V_v=P_v+xP_v'$.  Rational Bernstein coefficients establish
the strict inequalities
\[
 3(x-1)-|(x-1)^2P_v(x)|>0,
 \qquad 2-|(x-1)^2V_v(x)|>0.
\]
Since $x>1$, these are precisely $|P_v(x)|<3/(x-1)$ and
$|V_v(x)|<2/(x-1)^2$.  The entry intervals have exact consecutive endpoints
and cover $[2999/1000,344373767/114791256]$.
\emph{Spot check} (iii): on the first cell,
$\mathsf C(\Theta)=(-9/8,-51/16-9\Theta/2)$ has $V=33/80$ at $\Theta=-4/5$ and
$V=-39/80$ at $\Theta=-3/5$. Being affine in $\Theta$, it lies in
$[-11/8,11/8]\times[-39/80,39/80]$ on the whole cell, with margins $3/8$ and
$1/80$ in $R_0=[-3/2,3/2]\times[-1/2,1/2]$.
The band covers of (i), the sixty word regions of (ii), the eighteen cells
of (iii) and the entry cells of (iv) with their Bernstein coefficients are
stored and rechecked.
\companioncheck{cert:height-six}
\end{certificate}

\begin{proposition}[Certified: return and entry at height two]\label{cert:height-two-return}

Let $Q$, $\|\cdot\|_Q$, $\Theta_{x,w}$, $Y_m$, $s_1,\ldots,s_8$ and $\mu$ be as
in \Cref{app:height-two}. Put $\rho_-=1839286755214161\cdot10^{-15}$; in (ii)
and (iii) let $r=\rho$ and $9/80\le\tau\le9r^3/80$.
\begin{enumerate}
\item For $x\in[433/250,4598217/2500000]$ and $Y\in Q$, some word $w$ over
 $B_2$ with $1\le|w|\le10$ has $\Theta_{x,w}Y\in Q$.
\item If $0\le m\le8$ and $x=r-\tau r^{-3m-3}$, some word $s$ over $B_2$ has
 $\Theta_{x,s}Y_m(x)\in\interior Q$.
\item If $(S-V,V)=\bigl(2-r^2,\,5+2r-2r^2-(16+52r-32r^2)\tau\bigr)$ for a
 state $Y$, some $j$ has $\|\Theta_{r,s_j}Y\|_Q\le1-\mu$.
\item The polynomial $x^3-x^2-x-1$ has a zero in $[\rho_-,\rho_-+10^{-15}]$.
\end{enumerate}
\end{proposition}
\begin{certificate}
The reductions are printed in \Cref{app:height-two}: rational enclosures of
word-image margins on closed parameter bands, with triangulations of the
return diamond, for (i); words $ab^ms$ on exact algebraic phase cells for
(ii); margins affine in $\tau$, checked at both ends of eight phase cells,
for (iii). Positivity is certified by
exact Bernstein bounds over $\mathbb Q(\rho)$ in (ii) and (iii); (iv) is the
sign of $x^3-x^2-x-1$ at the two ends.
\emph{Spot check}, in the coordinates $(P,V)=(S-V,V)$ at $x=9/5$: (i) the
vertex $(3/2,0)$ of $Q$ is mapped by $(2,2)$ to $(-37/50,13/25)$, where
$|u_1|+|u_2|=13/35+94/375=1633/2625<1$; (ii) here $m=0$,
$\tau=(\rho-9/5)\rho^3\in[3/20,3/10]$, and the word $(2,1,-1,-1,-2)$ of the
cell containing $9/5$ maps $(1,1)$ to $(-2786/3125,1219/3125)$, where
$|u_1|+|u_2|=1219/4375+19327/46875=226714/328125<1$.
The return cover of (i) consists of $135$ parameter bands with their
triangulations; (ii) and (iv) are checked over $\mathbb Q(\rho)$ on closed
phase cells, and (iii) by exact capture intervals whose smallest overlap
margin, $16347138486163573/(3255\cdot10^{15})$, exceeds $\mu$.
\companioncheck{cert:height-two-return}
\end{certificate}

\begin{proposition}[Certified: the lower block for $q\ge22$]\label{cert:lower-block}

For an integer $q\ge22$ put $r_q=\lfloor3q/4+5/8\rfloor$,
$C_q=3q/4+5/8-r_q\in\{1/8,3/8,5/8,7/8\}$,
\[
 X_q(v,w)=q+\frac12+\frac1{8q}+\frac{v+w-3/8}{2q^2}+\frac{9/64-v}{q^3},\qquad
 \Omega_q=\Bigl\{X_q(v,Y^2)+\frac{iY}{\sqrt q}:(v,Y)\in D_q\Bigr\},
\]
where $D_q=\{(v,Y)\in\R^2:Y^2<C_q+1/8,\ Y^2-C_q+1/16<v<17/16\}$. Then
$\partial\Omega_q\cap\mathcal C_{q^2}=\varnothing$ and
$\Omega_q\cap\interior\mathcal C_{q^2-2}\ne\varnothing$.
\end{proposition}
\begin{certificate}
Here guards force the first inverse digits
$d_1,d_2,\ldots$ along the curve, and a final margin excludes every
continuation (a tree of depth at most four). With the states $(\mathsf S_j,\mathsf V_j)$ along the forced
digits and the height $h$ of the exclusion, the guards are
\begin{equation}\label{sep:eq:Lguard}
\begin{aligned}
 L_j&=(\mathsf S_{j-1}-d_j+1)(x-1)^2-h, &
 L_j^+&=(\mathsf S_{j-1}-d_j)(x-1)^2-h,\\
 U_j&=(2x-h-\mathsf S_j)(x-1)^2-h, &
 Z_j&=(-\mathsf S_j-h)(x-1)^2-h,\qquad W_j=\mathsf V_j(x-1)^2-h.
\end{aligned}
\end{equation}
Given the earlier digits, $L_j>0$ excludes $j$th digits below $d_j$ and
$L_j^+>0$ those below $d_j+1$, by the $\mathsf V$ bound. Raising the $j$th
digit by $k$ lowers $\mathsf V_j$ by $k$ and $\mathsf S_j$ by $2xk$; so, by
the $\mathsf S$ bound, $U_j>0$ excludes the digits above $d_j$ and $Z_j>0$
excludes $d_j$ and all larger ones, and by the $\mathsf V$ bound $W_j>0$
excludes $d_j$ and all smaller ones. On the boundary, at $h=q^2$, with
$x=X_q(v,w)$, $y^2=w/q$, $w=Y^2$ and
$(d_1,d_2,d_3)=(2q,-q^2+q,-q^2+r_q)$: $L_1,U_1,L_2,U_2>0$ on $\partial D_q$
force $d_1,d_2$; then $W_2>0$ on the right side $v=17/16$, and $L_3,Z_3>0$ on
the left side and the caps. By \Cref{cert:lem:tree},
$\partial\Omega_q\cap\mathcal C_{q^2}^{\mathrm{ext}}=\varnothing$, and
$\partial\Omega_q\subset X$. Each guard, cleared of positive factors, is
positive by \Cref{sep:lem:bernstein} for $q\ge24$ and each value of $C_q$,
with $s=q-24$ and the boundary parameter in $[0,1]$ (the guards of the first
two digits on the whole rectangle $-13/16\le v\le17/16$, $0\le w\le1$), and
directly at $q=22,23$. Interior point: for $q\ge40$,
$c_*=X_q(1-3C_q/2,0)\in\Omega_q$ and the word $(d_1,d_2,d_3)$, whose digits
lie in $B_{q^2-3}$, satisfy \Cref{sep:lem:inner} at $h=q^2-3$, by five margins
positive by \Cref{sep:lem:bernstein} with $s=q-40$; for $22\le q\le39$ an
explicit word does so at $c_*=X_q(v_*,0)\in\Omega_q$ and $h=q^2-2$. Nesting
and $\interior\mathcal Z_h=\interior\mathcal C_h$ give the statement. \emph{Spot check:} for
$q=23$, $r_{23}=17$, $C_{23}=7/8$, and $v_*=-1051/1600$ lies in
$(-13/16,17/16)$. At $x=X_{23}(v_*,0)=915134579/38934400$ the word
$(46,-506,-512,-526)$ has digits in $B_{527}$, and the lens, $\mathsf V$ and
$\mathsf S$ margins of \Cref{sep:lem:inner} at $h=527$ are, rounded, $455.5$,
$0.392$ and $381.3$.
The guard polynomials are generated from the recurrence and their Bernstein
coefficients tested; a second program rechecks the exported coefficients,
and the witness words are stored.
\companioncheck{cert:lower-block}
\end{certificate}

\begin{proposition}[Certified: twenty domains, $5\le q\le21$]\label{cert:twenty-domains}

For all integers $q,h$ with $5\le q\le21$ and $q^2-2\le h\le q^2$, there is an
open set $U\subset\{\re c>1\}$ such that $\partial U\cap\mathcal C_h=\varnothing$
and $U\cap\mathcal C_h\ne\varnothing$.
\end{proposition}
\begin{certificate}
Here $U$ is one of twenty bounded Jordan domains with closure in
$\{\re c>2\}$, each with a range
$h_{\min}\le h\le h_{\max}$ within $q^2-2\le h\le q^2$. Its upper boundary
consists of three or four arcs that are straight in $(x,y^2)$ or in the chart
$(v,w)$ of \Cref{cert:lower-block}, so $x$ and $y^2$ are affine on each arc;
each arc is covered by closed cells excluded at $h_{\max}$ by
\Cref{cert:lem:tree}, and conjugation gives the lower half. A real point
$x_*$ of the domain and a word satisfy \Cref{sep:lem:inner} at $h_{\min}$, so
$x_*\in\interior\mathcal Z_{h_{\min}}^{\mathrm{ext}}\subset
\mathcal C_{h_{\min}}$; nesting gives every $h$ of the range. \emph{Spot
check:} for $q=6$ the domain with range $h=34,35$ has the witness $x_*=162/25$ with the
word $(12,-30)$. At $h=34$ it gives
$(\mathsf S_2,\mathsf V_2)=(-538488/15625,\,282/625)$, and the lens,
$\mathsf V$ and $\mathsf S$ margins are $8781/625$, $2363003/10935000$ and
$1149/31250$.
Every exclusion tree is rebuilt from the stored cells by the exit and
branching rules of \Cref{cert:sec:conventions} with
\Cref{sep:lem:bernstein}, and a second program reconstructs the Bernstein
coefficients by Taylor--Bernstein conversion.
\companioncheck{cert:twenty-domains}
\end{certificate}

\begin{proposition}[Certified: small heights]\label{cert:small-heights}

For every $h\in\{1,3,4,8,9,14,15,16\}$, there is an open set
$U\subset\{\re c>1\}$ such that $\partial U\cap\mathcal C_h=\varnothing$ and
$U\cap\mathcal C_h\ne\varnothing$.
\end{proposition}
\begin{certificate}
Here $U$ is a bounded Jordan domain with closure in the positive half-plane. The boundary arcs are covered by closed cells excluded by
\Cref{cert:lem:tree}; at $h=1$ the cell variable runs over $[-1,1]$ and the
check uses the second bound in \eqref{cert:eq:coefficient-sum} with
the constant $(x_{\min}-1)^{-2}\ge\beta$. 

Core points: for $h=1,3,4,8,9$ a nonreal root of an explicit polynomial in
$\mathcal P_h$, isolated by \Cref{cert:lem:rouche} in a disc inside $U$; for
$h=14,15,16$ the real points $221/50$, $112/25$ and $45327/10000$, where the
words $(8,-13,-9)$, $(8,-12)$ and $(8,-12,-13,-13,-16)$ satisfy
\Cref{sep:lem:inner}. \emph{Spot
check:} at $h=15$, $x=112/25$ and $(8,-12)$ give
$(\mathsf S_2,\mathsf V_2)=(-226688/15625,\,332/625)$, with lens,
$\mathsf V$ and $\mathsf S$ margins $1231/625$, $280159/15680000$ and
$30999/31250$.
For $h=3,4,8,9$ the exclusion trees have depths $7,5,4,6$ and the witness
polynomials degrees $7,6,4,6$, each isolated with Rouch\'e radius $10^{-7}$;
for $h=14,15,16$ the cells and words are rechecked by the two programs of
\Cref{cert:twenty-domains}; for $h=1$ the exclusion propagates
$(\re z,\im z/y)$ instead of $(\mathsf S,\mathsf V)$ and the witness root
is isolated exactly.
\companioncheck{cert:small-heights}
\end{certificate}

\begin{proposition}[Certified: inequalities of the square-adjacent family]\label{cert:sa-inequalities}
For $4\le q\le60$ and $k\in\{0,1,2\}$, with $h=q^2-2+k$ and $m=q^2-2q+k$,
the following hold, where $\varphi(\beta)=2(\beta-q)(\beta-1)^2$ and
$t_-(m')$ denotes the smaller zero of $1-(2q+1)z+(2q+m')z^2$:
\begin{enumerate}
\item $\varphi(q+\delta^*)>h$ with $\delta^*=(\sqrt{q^2+4q-4}-q)/2$, and
$\varphi\bigl(1+\sqrt{h-m/2}\bigr)<m$;
\item $t_-(m)\ge1/(2q-1)$;
\item $t_-(q^2-2q+2)\ge1/(2q-3)$ for $q\ge6$.
\end{enumerate}
\end{proposition}
\begin{certificate}
Items (i)--(iii) are finite instances of the inequalities proved by
hand, for every $q\ge4$, in Lemmas~\ref{lem:sa-anchors},
\ref{lem:sa-cones} and~\ref{lem:sa-existence}; the record
checks the instances $4\le q\le60$ in exact arithmetic as a control on
those proofs, which do not depend on it. Each inequality compares a rational number with a number of
the form $a+b\sqrt D$, $a,b\in\mathbb Q$, $D\in\mathbb Q_{>0}$; its sign is
decided exactly by comparing $a^2$ with $b^2D$ when $a$ and $b$ have
opposite signs.
\companioncheck{cert:sa-inequalities}
\end{certificate}

\begin{proposition}[Certified: inner prefixes]\label{cert:inner-prefix}
Let $(q,m',H',J,\pi)$ be one of
\begin{gather*}
 \bigl(4,10,14,[\tfrac15,\tfrac13],(1,-8,12,13,14,12,13,10,13,13)\bigr),\\
 \bigl(5,17,23,[\tfrac17,\tfrac7{25}],(1,-10,20,22,18,21,22,21,17)\bigr).
\end{gather*}
Then $(H'-m')a=1-a$ at the left endpoint $a$ of $J$; $f^-_\pi>0$ at both
endpoints of $J$; $f^-_\pi(z_0)<0$ for $z_0=138/625$, respectively
$z_0=907/5000$; and there are rationals $z_1<z_2$ in $J$ with
$(f^+_\pi)'(z_1)\le0\le(f^+_\pi)'(z_2)$, $f^+_\pi(z_1)>0$,
$f^+_\pi(z_2)>0$ and $f^+_\pi(z_1)+(f^+_\pi)'(z_1)(z_2-z_1)>0$.
\end{proposition}
\begin{certificate}
By convexity of $f^+_\pi$ on $(0,1)$, the last three inequalities give
$f^+_\pi>0$ on $J$: on $[z_1,z_2]$ the function lies above its tangent
at $z_1$, and outside this interval it is monotone away from the
minimum. All quantities are rational and the inequalities are checked
exactly from the stored prefixes and points.
\companioncheck{cert:inner-prefix}
\end{certificate}

\begin{proposition}[Certified: explicit witnesses]\label{cert:explicit-witnesses}
For each row of Table~\ref{tab:explicit-witnesses}, let $P$ be the
polynomial $1-2qz+\sum_{j=2}^{N-1}a_jz^j$, $T=\sum_{i<r}c_iz^i$ and
$R=(1-z^r)P+z^NT$. Then $g=\gcd(R,R')$ is irreducible of degree three,
four or five, $g^2$ divides $R$, the reversed polynomial of $g$ is the
stated minimal polynomial, $g$ has exactly one zero in $(1/(q+1),1/q)$,
all digits lie in $[m,h]$ for each listed height, and the period is not
constant. Moreover the stated polynomial is monic and all its zeros
other than $\beta$ lie strictly inside the unit disk, so that $\beta$ is a Pisot number.
\end{proposition}
\begin{certificate}
Exact polynomial arithmetic over $\mathbb Z$; the zero count in
$(1/(q+1),1/q)$ is a Sturm count with rational endpoints. The Pisot
property, which the proof of Proposition~\ref{prop:explicit-witnesses}
does not use, is certified by exact isolating intervals and rectangles
for all zeros of the stated polynomial. For a rectangle $I\times J$,
the minimum of $|z|^2$ is $\operatorname{dist}(0,I)^2+\operatorname{dist}(0,J)^2$;
the maximum is attained at a corner. The verifier requires maximum $<1$
or minimum $>1$, refining any undecided box. Real intervals use the same
test with $J=\{0\}$.
\companioncheck{cert:explicit-witnesses}
\end{certificate}

\begin{proposition}[Certified: surviving prefixes]\label{cert:surviving-prefixes}
For $h\in\{25,49,62,63,64,79,81\}$, with $q$, $k$ and $m$ as
in~\eqref{eq:gap-parameters}, the companion lists $n_h$ distinct prefixes
$\pi=(1,-2q,a_2,\ldots,a_{N-1})$ with digits in $[m,h]$, of a common
length $N$, where
\begin{gather*}
 (h,N,n_h)\in\{(25,15,384),\ (49,12,171),\ (62,12,314),\ (63,12,332),\\
 (64,12,299),\ (79,12,463),\ (81,12,369)\},
\end{gather*}
each with rationals $z_0,z_1,z_2\in(0,1)$ such that
$1-(2q+1)z_0+(q^2+k)z_0^2<0$, $f^-_\pi(z_0)<0$, $z_1<z_2$,
$(f^+_\pi)'(z_1)\le0\le(f^+_\pi)'(z_2)$, $f^+_\pi(z_1)>0$,
$f^+_\pi(z_2)>0$ and $f^+_\pi(z_1)+(f^+_\pi)'(z_1)(z_2-z_1)>0$.
Consequently the anchor set $\mathcal T_h$ has at least $n_h$ points.
\end{proposition}
\begin{certificate}
The first inequality places $z_0$ in $I_0$, where $f^-_\pi$ is then
negative; the remaining ones give $f^+_\pi>0$ on $(0,1)$ by convexity, as
in Proposition~\ref{cert:inner-prefix}. Since $f^-_\pi(t_\pm)\ge M_m(t_\pm)=0$ and the contiguity inequality $t_-\ge1/(2q-1)$ holds
(Proposition~\ref{cert:sa-inequalities}), Lemma~\ref{lem:ivt-completion}
extends each prefix to a square-adjacent witness, and distinct prefixes of equal
length extend to distinct witnesses. All quantities are rational.
\companioncheck{cert:surviving-prefixes}
\end{certificate}
The last four records concern Section~\ref{sec:real-traces} and the heights
$h\in\{1,3,4,7\}$. Their programs represent every real quantity by an
interval $[\ell,u]2^{-P}$ with integers $\ell\le u$; sums are exact, and
every product and quotient is rounded outward by exact floor division, so
the computed interval contains the exact one. The parameter $t$ always runs
over an interval, and each inequality is decided by comparing integers for
all $t$ in it at once.

\begin{proposition}[Certified: the minimum of $\mathcal D_h$]\label{cert:rt-minimum}
For $h\in\{1,3,4,7\}$ the companion lists rationals $a_h<b_h$, with
$b_h-a_h=2\cdot10^{-150}$ for $h=1$ and $2\cdot10^{-90}$ otherwise, a tiling of
$[1/(h+1),a_h]$ by closed intervals, a polynomial $Q_h$, and words of lengths
$N_h=850,257,216,162$ and $n_0=49,24,19,15$, such that: no series of $\Fh$ has a
double zero in any tile; $Q_h$ satisfies the hypotheses of
Lemma~\ref{rt:lem:good-prefix} on an interval inside $(a_h,b_h)$; every series
with a double zero in $[a_h,b_h]$ begins with the word of length $N_h$, and
every series with a double zero in $[a_h-10^{-12},b_h+10^{-9}]$ begins with the
word of length $n_0$.
\end{proposition}
\begin{certificate}
On a tile $[x,y]$ put $\beta\in[1/y,1/x]$, enclosed by
$[\beta_0-\delta,\beta_0+\delta]$ with $\beta_0$ and $\delta$ dyadic. A double zero
of $1+\sum a_jz^j$ at $t=1/\beta$ forces, for every $N$,
\[
 |P_N(\beta)|\le\frac h{\beta-1},\qquad |V_N(\beta)|\le\frac h{(\beta-1)^2},
\]
where $P_N=\beta^N+\sum_{j\le N}a_j\beta^{N-j}$ and
$V_N=(N+1)\beta^N+\sum_{j\le N}(N+1-j)a_j\beta^{N-j}$, because
$P_N(\beta)=-\beta^N\sum_{j>N}a_jt^j$, and similarly for the derivative.
The check runs the tree of all digit words, carrying $P_N(\beta_0)$, $V_N(\beta_0)$
and their $\beta$-derivatives as integers at scale $2^K$, with Taylor remainders
$4hN^2\bar\beta^N\delta^2$ and $4hN^3\bar\beta^N\delta^2$ ($\bar\beta=\beta_0+\delta$;
valid since $|P_N''|\le1.37hN^2\beta^N$ for $\beta\ge1.49$), and deletes a word
when its enclosure misses the box. A tile is accepted when no word survives;
the tiles abut exactly. At $h=1$, $469$ tiles cover $[1/2,a_1]$ (depth at most
$850$); at $h=3,4,7$, $291$, $295$ and $308$ tiles. The good-prefix conditions
are integer inequalities over $10^{Dn}$, positivity of $Q_h$ on the interval
following from a lower bound for $Q_h''$ on a grid and the Lipschitz constant
$6h/(1-b)^4$. The same tree on $[a_h,b_h]$ and on the window gives the
common prefixes of the surviving words.
\companioncheck{cert:rt-minimum}
\end{certificate}

\begin{proposition}[Certified: slices and local constants]\label{cert:rt-constants}
For $h\in\{1,3,4,7\}$, with $t$ within $10^{-15}$ of the centre $T_0$ of
$[a_h,b_h]$, the quantities $\Delta_{\mathrm{low}}$, $\Delta_1$ and $q_1^*$ of
Section~\ref{sec:real-traces} satisfy the bounds of Table~\ref{rt:tab:data};
the conditions (C1)--(C7), (D), (E1)--(E7) hold for all $i>n_0$ and all relative
tips within $10^{-80}$ of $t_0$; and the checks (GR) and (EX) hold, together
with $L_{\min}+L_2+20\le L^*(159)$, $1.001v^3\le1-v$ ($h=1$) and $1.001v\le1-v$
($h\ge3$). Here (C3) is $v-(n_0+1)\zeta_{n_0+1}-2hv^{r_0}(r_0/(1-v)+1/(1-v)^2)>0$
with $r_0$ the least $r$ with $v^r<\delta_{n_0+1}+\zeta_{n_0+1}$; (C6), the product of
the correction factors in the proof of Lemma~\ref{rt:lem:taylor}, is
$\bigl(1+(4C_3W(i)+\frac32G_2(i))/\kappa_v\bigr)\bigl(1+ie_i/\bar M\bigr)^2
\bigl(1-(3C_3W(i)+G_2(i))/(2\kappa_v)\bigr)^{-2}\le1+10^{-3}$ with
$e_i=(\kappa_v+4C_3W(i)+\frac32G_2(i))W(i)^2/v^i$; (C7)
is $i\bar\eta_i\le10^{-5}$; (E3) is $(N+1)(\bar Mv^{i-1}/\kappa_v+c_\rho v^{(N+1)/2})/v\le10^{-3}$
for $i+5\le N\le2i+2$, a range that contains the values $N=i+L+5$ of
Lemma~\ref{rt:lem:realization}; (E5) is $(6C_3R+G_2(n_0+1))/(2\kappa_v)\le10^{-4}$; (E6) is
$2h(v+\cdots+v^5)>1.7$; and (E7) is $2hv>1.01$. Conditions (E1), (E2$'$) and (E4) are
$\kappa_vc_\rho^2>1.01h/(1-\tau)$, $0.8hv/(1-v)>1.01$ and $v^4(1.9hv/(1-v)+1)\le0.965$.
\end{proposition}
\begin{certificate}
The slices are taken from an outer table: for $x$ in cells of width $10^{-4}$,
a finite union of intervals containing $S_t(x)$ for every $t$ within
$10^{-12}$ of $T_0$, obtained from the trivial table $[0,2ht/(1-t)^2]$ by
iterating \eqref{rt:eq:slice-recursion} with outward rounding (each iterate is
again an outer table; $45$ iterations, $43$ at $h=7$, where the table is
stationary). A lookup also includes the neighbouring cells, which absorbs
weight shifts below $10^{-4}$. The greedy moment $q^{\mathrm g}_t$ is enclosed
over an interval of weights by following every branch of the greedy
expansion, and the moment of the unknown remainder after $n$ digits by
$t^n(nr+[0,\bar dt/(1-t)^2])$ with $r\in[0,1+10^{-7}]$. $\Delta_{\mathrm{low}}$ and
$\Delta_1$ are minima over $4000$ ($h=1$) or $2000$ cells of $p$, and $q_1^*$ is
a maximum over all branches. Each condition indexed by $i>n_0$ is reduced to
$i=n_0+1$ by monotone majorants, for instance $W(i)\le Ai\tau^{i-1}$ and
$G_2(i)\le Bi^2\tau^{i-2}$, and $\kappa_v$ is enclosed from the $N_h$ certified
digits with the unknown tail bounded by $G_2(N_h+1)/2$. (GR) runs over the
$118$ ($h=1$) and $7$ ($h=3$) admissible pairs $(c,F)$, with worst margins
$0.0708$ and $0.685$; it is empty at $h=4,7$, where $L_G=1$.
\companioncheck{cert:rt-constants}
\end{certificate}

\begin{proposition}[Certified: Move Lemma]\label{cert:rt-move}
Proposition~\ref{rt:prop:move} holds, with the constants of
Table~\ref{rt:tab:data}, and the clearances in (a) and (b) of
Corollary~\ref{rt:cor:move}, that is the distances of the groups from
$v^sL_-(1-\mu)$ in (a) and from $v^s[L_-(1-\mu),L_+(1+\mu)]$ in (b), are at least
$8\cdot10^{-4}v^sL_-$.
\end{proposition}
\begin{certificate}
Fix a start word $d\,u$ of $C_h$. The digits after the pattern are
enumerated by cylinders $uw$ of weakly admissible words over
$\{0,\ldots,\bar d\}$ ($|w|=8,2,2,2$ initially, refined on failure up to
$|w|=40,24,24,14$), the remainder after them in $[0,1+10^{-7}]$. The groups of
the levels inside $uw$ are computed from the known digits; beyond them each
group lies below $t^{|uw|+1}D^*$, where $D^*=\sup\{\max S_t(c+y)-q^{\mathrm g}_t(y)\}$
over $1\le c\le\bar d$ and $y\in[0,1+10^{-7}]$ is at most $7.71,5.56,4.75,4.89$;
this bound is valid for tails with restarts by
Lemma~\ref{rt:lem:restarts}~(ii). The $J_h$ levels before the pattern run over
all weakly admissible context words. For each cylinder, each branch of the
greedy expansion of $1+y_s+[-10^{-7},10^{-7}]$ gives an enclosure of the
first-group position, and the check verifies (a) and (b) against it with the
stated clearance. In all, $71$, $65$, $66$ and $512$ cylinder checks pass for
$h=1,3,4,7$; the smallest clearances, in units of $v^sL_-$, are $1.5\cdot10^{-3}$ and $4.2\cdot10^{-3}$ for
$h=1$, $0.063$ and $0.41$ for $h=3$, $0.014$ and $0.23$ for $h=4$, and
$8.9\cdot10^{-4}$ and $0.54$ for $h=7$.
\companioncheck{cert:rt-move}
\end{certificate}

\begin{proposition}[Certified: Landing Lemma and assembly]\label{cert:rt-landing}\label{cert:rt-assembly}
Proposition~\ref{rt:prop:landing} holds. Moreover, with the constants of
Table~\ref{rt:tab:data} and $F_{\min}=(\mu-10^{-4})L^\circ_{\min}$:
$v^{P_0^*}L^\circ_{\max}(1+\mu)<F_{\min}\le10^{-9}/u_{P_0}$ and $P_0^*\ge J_h+2$;
\begin{itemize}[leftmargin=*]
\item[(K1)] $(\Delta_{\mathrm{low}}-v^{J_h+2}(q_1^*+10^{-5}))\beta^{J_h+1}>L^\circ_{\max}(1+\mu)(1+3\cdot10^{-5})$;
\item[(K2)] the same with $\Delta_1$ in place of $\Delta_{\mathrm{low}}$;
\item[(K4)] $\min(\Delta_{\mathrm{low}},\Delta_1)-v^{L_G}(q_1^*+10^{-5})>0.002$;
\item[(K5)] $1.0001v^{R_h+1}(q_1^*+10^{-5})\le7\cdot10^{-5}L^\circ_{\min}$;
\item[(EX$'$)] the $M_E\ge144$ common digits of (EX) contain a $C_h$-pattern at
 a level at least $\max(M_h,P_0^*)$;
\end{itemize}
the remainders of $t_0$ satisfy $p_{n_0}<2$ and $0\le p_k<1$ for
$n_0<k\le\min(200,N_h-2)$; its certified tail has $C_h$-patterns at levels
$\ge P_0^*$, the deepest at absolute position $i_1\ge158$; and
$X_h=v^{i_1-1}L^\circ_{\max}(1+10^{-4})/(2\kappa_v(1-v^{P_0^*}))$ and
$\bar\eta_{i_1}+X_h(2h+\bar d)/(1-\tau)^2\le\epsilon_h$ are inside the windows of
Propositions~\ref{rt:prop:move} and~\ref{rt:prop:landing}.
\end{proposition}
\begin{certificate}
The landing set is covered by a branch and bound over the core language:
a node fixes a finite core word, encloses $y_s$ over all its continuations
(the constraints of the core are dropped for the unknown part, which only
enlarges the enclosure), and is closed when the digits of the greedy
expansion that are common to the whole node determine a $C_h$-pattern at a
level $\ge P_0^*$ and its first $R_h$ remainders exceed $10^{-20}$; otherwise
it is split. The zero runs $0^k$ are treated one by one for $k\le40,30,40,30$,
and the longer ones with the terminating tails as one node. All
$6263$, $1211$, $40$ and $8678$ nodes close, at depth at most $119$, $25$, $0$
and $23$; the patterns lie at levels $13$--$421$, $7$--$26$, $6$--$10$ and
$4$--$18$. The remainders of $t_0$ are enclosed by a centred form at the
midpoint of $[a_h,b_h]$, using $f(t_0)=f'(t_0)=0$ for the derivative bound.
The remaining items are rational inequalities between the outputs of the
previous records.
\companioncheck{cert:rt-landing}
\end{certificate}

\subsection{Supplementary material}\label{cert:sec:supplement}
The computational companion~\cite{Companion26} contains, for each of
Propositions~\ref{cert:short-word-strip}--\ref{cert:rt-landing}, the frozen
inputs (integer and rational coefficients, cell covers, word lists and
witness data) and a program that accepts them only by exact rational,
symbolic or bounded-integer arithmetic; decimal output is diagnostic and
never decides acceptance. Its documentation lists, record by record, the
inputs and checking programs, and a complete run of all eighty-four checks
on the supplied inputs is recorded there together with the software
versions used; the complete runs of the forty-eight checks of the earlier
records made on 5 and 3 November, the second before the unit-disk test of
record~\ref{cert:explicit-witnesses} was tightened, are preserved alongside
it. The documentation numbers the records D.4--D.19, as in an earlier
version of this paper in which this section was Appendix~D:
Proposition~S4.$k$ here is record D.$k$ there. It also keeps the working name
``finger'' for the square-adjacent family of Section~\ref{sec:gap-heights}
(directory \texttt{certificates/finger}), and its figure \texttt{finger\_region.pdf} is
Figure~\ref{fig:anchor-interval} with ``finger interval'' for the anchor
interval. The companion also contains
the search programs and logs behind Proposition~\ref{prop:explicit-witnesses},
the generators that proposed the certificates of
Proposition~\ref{cert:rt-minimum}, and the data and scripts of the figures
drawn from computed data.

\end{document}